\documentclass[smallextended]{svjour3}       % onecolumn (second format)
\smartqed  % flush right qed marks, e.g. at end of proof
\usepackage[active]{srcltx}
\usepackage{amssymb,amsmath,amssymb,amsopn,amsxtra,amstext,amsbsy}
\usepackage[cp1251]{inputenc}
\usepackage[english]{babel}
\usepackage{setspace}

\allowdisplaybreaks[4]

\makeatletter

\renewcommand\section{\@startsection{section}{1}{\z@}%
    {-21dd plus-8pt minus-4pt}{10.5dd}
     {\large\bfseries\boldmath}}

\renewcommand\subsection{\@startsection{subsection}{2}{\z@}%
    {-21dd plus-8pt minus-4pt}{10.5dd}
     {\normalsize\bfseries\boldmath}}
     
\renewcommand\subsubsection{\@startsection{subsubsection}{3}{\z@}%
    {-13pt plus-8pt minus-4pt}{\z@}{\bfseries\boldmath}}

\makeatother
\DeclareSymbolFont{AMSb}{U}{msb}{m}{n}
\DeclareMathAlphabet{\Bb}{U}{msb}{m}{n}
\DeclareSymbolFont{AMSb}{U}{eur}{m}{n}
\DeclareMathAlphabet{\eurm}{U}{eur}{m}{n}
\numberwithin{equation}{section}
\numberwithin{theorem}{section}
\numberwithin{corollary}{section}
\numberwithin{lemma}{section}
\numberwithin{definition}{section}
\newtheorem{sectlemma}{Lemma}
\numberwithin{sectlemma}{section}
\newtheorem{thmx}{Theorem}

\newtheorem{cormx}{Corollary}

\newtheorem{clash}{}

\newcommand{\fo}[1]{{\mbox{\footnotesize{$#1$}}}}
\newcommand{\tn}[1]{{\mbox{\tiny{$#1$}}}}
\newcommand{\nor}[1]{{\mbox{\normalsize{$#1$}}}}
\newcommand{\im}{{\rm{Im}}}
\newcommand{\re}{{\rm{Re}}}

\newcommand{\ie}{\lambda_{{\tn{\triangle}}}}
\newcommand{\he}{F_{\!{\tn{\triangle}}}}
\newcommand{\het}{F_{\!{\tn{\triangle}}}}
\newcommand{\fet}{\mathcal{F}_{{\tn{\square}}}}

\newcommand{\ihet}{ F_{{{\fo{1/2,1/2;1}}}} }
\newcommand{\hett}{ F_{{{{1/2,1/2;1}}}} }

\newcommand{\jet}{F_{{\fo{3/2,3/2;2}}}}

\newcommand{\Log}{{\rm{Log}}\,}
\newcommand{\Arg}{{\rm{Arg}}\,}
\DeclareMathOperator*{\limsp}{\overline{\lim}}

\newcommand{\diff}{\mathrm{d}}
\newcommand{\imag}{ \hspace{0,02cm}\mathrm{i} \hspace{0,015cm}}

\newcommand\Picklogplus{\mathcal{P}_{\log}(-\infty, 1)}

\newcommand{\Oh}{\mathrm{O}}
\newcommand{\oh}{\mathrm{o}}
\usepackage[linktocpage=true,bookmarksopen=true,bookmarksopenlevel=1]{hyperref}
\journalname{PREPRINT}
\begin{document}

\title{Exponential integral representations  of  theta functions%\thanks{Grants or other notes
}


\author{Andrew Bakan \and H{\aa}kan Hedenmalm %\and Alfonso Montes-Rodr\'{\i}guez???%etc.
}


\institute{Andrew Bakan \at
               Institute of Mathematics, National Academy of Sciences of Ukraine,
01601 Kyiv, Ukraine\\
%              Fax: +123-45-678910\\
              \email{andrew@bakan.kiev.ua}           %  \\
           \and
          H{\aa}kan Hedenmalm \at
              KTH Royal Institute of Technology, SE–10044 Stockholm,
Sweden\\  \email{haakanh@math.kth.se}
}


\maketitle

\begin{abstract}    Let $\Theta_{3} (z):= \sum\nolimits_{n\in\Bb{Z}} \exp (\imag  \pi n^2 z)$ be the
standard Jacobi theta function, which is holomorphic and zero-free in the
upper half-plane $\Bb{H}:=\{z\in\Bb{C}\,|\,\,\im\, z>0\}$, and takes positive
values along $ \imag   \hspace{0,01cm}\Bb{R}_{>0}$, the positive imaginary axis, where $\Bb{R}_{>0}:= (0, +\infty)$.
%\textcolor{red}{We define}
We define
its logarithm $\log\Theta_3(z)$ which is uniquely determined by
the requirements that it should be holomorphic in $\Bb{H}$ and real-valued on
$ \imag  \hspace{0,01cm}\Bb{R}_{>0}$.
We derive an integral representation of $\log\Theta_{3} (z)$  when  $z$ belongs to  the  hyperbolic quadrilateral
\[
\mathcal{F}^{\,{\tn{||}}}_{{\tn{\square}}}:=
\big\{z\in\Bb{C}\,\,\big|\,\,{\im}\, z > 0, \,\,-1\leq{\re}\, z \leq 1,
\,\,|2 z\! -\!  1| > 1,\,\,
| 2 z \!+\!  1|\! > 1\big\}.
\]
Since every point of $\Bb{H}$ is equivalent to at least one point in $\mathcal{F}^{\,{\tn{||}}}_{{\tn{\square}}}$ under the theta subgroup
 of the modular group  on the upper half-plane, this representation carries over
in modified form to all of $\Bb{H}$ via the identity recorded by Berndt.
The logarithms of the related Jacobi theta functions $\Theta_{4}$ and
$\Theta_{2}$ may be conveniently expressed in terms of $\log\Theta_{3}$ via
functional equations, and hence get controlled as well.
Our approach is based on a study of the logarithm of the Gauss hypergeometric
function for a specific choice of the parameters. This connects with the
study of the universally starlike mappings introduced by Ruscheweyh, Salinas,
and Sugawa.

%Insert your abstract here. Include keywords, PACS and mathematical
%subject classification numbers as needed.

\vspace{0.5cm}

\vspace{0.5cm}
% 30C45 Special classes of univalent and multivalent
%%  functions (starlike, convex, bounded rotation, etc.)
%33C05 Classical hypergeometric functions, 2F1
%33E05 Elliptic functions and integrals
%11F27 Theta series; Weil representation; theta correspondences
\end{abstract}

\newpage
\begin{center}
{\tableofcontents}
\end{center}

\vspace{-0,3cm}
\section[\hspace{-0,30cm}. \hspace{0,11cm}Introduction and Main results]{\hspace{-0,095cm}{\bf{.}}  Introduction and Main results}
\label{intro}

Let $\Bb{D}:=\{z\in\Bb{C}\ |\  |z| <1\}$
denote the open unit disk in the complex plane $\Bb{C}$ and
 $\Bb{R}_{>0}:= (0, +\infty)$. We also write ${\Bb{H}}:=
 \left\{ z\in\Bb{C} \mid {{\rm{Im}}} \, z > 0 \right\}$
for the upper half-plane. Given a
domain $D\subset\Bb{C}$, we let ${\rm{Hol}}(D)$ denote the  set
of all holomorphic functions in $D\subset\Bb{C}$.
 We will also need the cone $\mathcal{M}^{+}(\Bb{R})$
  of all nonnegative Borel measures on $\Bb{R}$ such that $ 0\!\leq\!\mu \big((-a , a)\big)\!<\! +\infty$, for each  $a\in \Bb{R}_{>0}$,
  and for a given $\mu\in\mathcal{M}^{+}(\Bb{R})$ we let $
{\mathrm{supp}}\, \mu:= \{ x\in\Bb{R} \mid \mu\big( (x-\varepsilon,\,x+\varepsilon) \big)
 >0\  \ \mbox{for all} \ \varepsilon>0 \}$  denote the support of $\mu$.
  For $0<p<+\infty$, we denote by (see \cite[p.112]{koo})
  \begin{align}\label{inthp}
    &   {H}^p := \left\{ f\in {\rm{Hol}} (\Bb{H}) \ \Big| \ \sup\limits_{y>0}
    \int_{\Bb{R}} |f(x+\imag y)|^p\, {\diff} x < +\infty \right\} \, ,
\end{align}

\noindent the Hardy space  of the upper half-plane,   $m$  the
Lebesgue measure on the real line and $L^{\infty}(\Bb{R})$ the
 space  of all Borel
measurable real-valued functions $f$ on the real line that are
 essentially bounded, equipped with the essential supremum norm
 $\|f\|_{L^{\infty}(\Bb{R})} := \inf\{\, a>0 \,\mid\,
 m( \{ x\in\Bb{R} \, \mid \, |f(x)|> a \} ) = 0 \}$.

 Following the definitions of \cite[pp. 6, 40]{sar}, we denote by
 $\ln : \Bb{R}_{>0}\mapsto \Bb{R}$ the real-valued logarithm defined
 on $\Bb{R}_{>0}$, and let $\Log  ( z )\!=
 \!\ln |z| \!+\! \imag  \Arg  (z)$
  be the principal branch of the logarithm defined for
  $z\!\in \!\Bb{C}\setminus (-\infty, 0]$ with $\Arg  (z)\! \in \!(-\pi, \pi)$.
  Furthermore,  for a simply connected domain $D\!\subset\!\Bb{C}$,
  a point $a\!\in\! D$, and a  function $f\!\in\! {\rm{Hol}}(D)$ which
  is zero-free  in $D$ with $f (a)\!>\!0$, we write  $\log f (z)$  for
  the holomorphic function in $D$ such that  $\exp (\log f (z))\! =\! f (z)$,
  $z\!\in \!D$, and $\log f (a)\! = \!\ln f (a)$ (see \cite[p. 94]{con}).
  Then ${\rm{Re}} \, \log f (z) \!= \!\ln |f (z)|$ and
  $\arg f (z) \!:=\! {\rm{Im}} \, \log f (z)$ for each  $z\!\in \!D$.

As for topology, we denote by $'{\rm{clos}}(A)$ (or $\overline{A}$),
$\mathrm{int}(A)$, and
$\partial A$ the closure, interior, and boundary of a subset $A\subset\Bb{C}$,
respectively.  Moreover, we let $C(A)$ denote the set of all continuous
functions $f : A \to\Bb{C}$.

\medskip
\noindent{\bf Nevanlinna-Pick functions.}
We let $\mathcal{P}$ denote the class of {\emph{Nevanlinna-Pick functions}},
which are holomorphic functions $\Phi$ in $\Bb{C}\setminus\Bb{R}$ with
\begin{align}\label{eq:symm0}
{\rm{Im}}\, \Phi(z)\ge0,\qquad z\in\Bb{H},
\end{align}

\noindent
and the symmetry property
\begin{align}
\label{eq:symm1}
\Phi(z)=\overline{\Phi(\overline{z})},\qquad z\in\Bb{C}\setminus\Bb{R}\,.
\end{align}

\noindent
It is well-known
(see \cite[p.31]{berg})
that unless $\Phi\in\mathcal{P}$ is a real constant, the strict inequality
${\rm{Im}}\, \Phi(z)> 0$ holds for all $z\in \Bb{H}$.
Moreover, each function $\Phi\in\mathcal{P}$ has a unique canonical
representation of the form (see \cite[p. 20, Thm.~1]{D},
\cite[p. 23, Lem.~2.1]{sho})
\begin{align}
\label{f2int}
\Phi(z)
= \alpha z + \beta +
\int_{\Bb{R}} \left( \frac{1}{t-z}-\frac{t}{1+t^2} \right)\, {\diff} \sigma(t)\,,
\quad z\in\Bb{C}\setminus\Bb{R}\,,
\end{align}

\noindent
where $\alpha\geq 0$, $\beta\in\Bb{R}$, and $\sigma\in\mathcal{M}^{+}(\Bb{R})$
is such that
\begin{align*}
\int_{\Bb{R}}\frac{{\diff}\sigma(t)}{1+t^2}<+\infty.
\end{align*}

\noindent
In the converse direction, any function of the given form \eqref{f2int} is in
$\mathcal{P}$.
In the representation \eqref{f2int}, the measure $\sigma$ has the
interpretation of the jump in the imaginary part between the upper and
lower half-planes. If we write, for an open interval $I\subset\Bb{R}$,
\begin{align*}%\label{eq:PickI}
\mathcal{P}(I):=\mathcal{P}\cap{\mathrm{Hol}}\big( I \cup (\Bb{C}\setminus\Bb{R}) \big)
\end{align*} \noindent
it is then a consequence  \cite[p. 26]{D} of the Schwarz reflection principle that
% has been proved in \cite[Lemma~2, p.26]{D} that
\begin{align}
\label{f3int}
\Phi\in\mathcal{P}(I)
\quad\Longleftrightarrow\quad
{\mathrm{supp}}\, \,\sigma\subset\Bb{R}\setminus I.
\end{align} \noindent
%For a background on the Nevanlinna-Pick class, see, e.g., \cite{berg},
%\cite{D}, and \cite{sho}.

\noindent{\bf Logarithms of Nevanlinna-Pick functions.}
For an arbitrary $\Phi\in\mathcal{P}$ with $\Phi(z)\not\equiv a$, $a \leq 0$,
 we may take its
logarithm $\Log  \Phi$ and obtain a holomorphic function which maps
the upper half-plane $\Bb{H}$ into the strip
$\{w\in\Bb{C}\,|\,\,0\leq {\rm{Im}}\, w<\pi\}\subset\Bb{H}\cup\Bb{R}$, and
inherits the
symmetry property \eqref{eq:symm1} from $\Phi$, so that, in particular,
$\Log  \Phi\in\mathcal{P}$. We denote by $\log\mathcal{P}$ the collection of all such
functions $\Log  \Phi$, where $\Phi\in\mathcal{P}$ and $\Phi(z)\not\equiv a$,
$a \leq 0$, and
the observation just made amounts to the inclusion $\log\mathcal{P}\subset\mathcal{P}$.
Such functions $f\!\in\!\log\mathcal{P}\!\subset\!\mathcal{P}$ are characterized
in terms of a
corresponding integral representation (see \cite[p. 27]{D})\vspace{-0,3cm}
\begin{align}\label{f4int}
f(z) = b + \int_{-\infty}^{+\infty}
\left( \frac{1}{t-z}-\frac{t}{1+t^2} \right)\,a(t)\,{\diff} t\,,\quad
z\in\Bb{C}\setminus\Bb{R}\,,
\end{align}

\vspace{-0,1cm}
\noindent
where ${{\diff}}  t := {\rm{d }}  m (t)$, $b\in\Bb{R}$ and
$a\in L^{\infty}(\Bb{R})$ with $0\leq a(x)\leq 1$ on $\Bb{R}$
(almost everywhere with respect to $m$).
On the other hand, any function of the form \eqref{f4int} is in
$\log\mathcal{P}$.

\medskip

\noindent{\bf Universally starlike functions.}
We begin with some notation.
A domain $D$ in the complex plane $\Bb{C}$ is referred to as
\emph{circular} when it is either an open disk or an open half plane.
Moreover, a domain $\Omega$ is said to be \emph{starlike with respect to
the origin} if $0\in\Omega$ and if for each $z_0\in\Omega\setminus\{0\}$
the straight line segment from $0$ to $z_0$ is contained in $\Omega$.

%%%%%%%%%%%%%%%%%%%%%%%%%%%%%%%%%%%%%%%%%%%%%%%%%%%%%%%%%%%%%%%%%%%%%%%%%%%%

%%%%%%%%%%%%%%%%%%%%%%%%%%%%%%%%%%%%%%%%%%%%%%%%%%%%%%%%%%%%%%%%%%%%%%%%%%%%
Associated with the starlike domains we have the notion of starlike univalent
mappings \cite[p. 40]{Dur}. Building on this, Ruscheweyh, Salinas, and Sugawa
\cite[p. 290]{RS} introduced
%\cite[Definition 1.3, p.290]{RS}  introduced
the notion of universal starlikeness in the context of holomorphic
functions in the set $\Bb{C}\setminus[1,+\infty)$.

\begin{definition}\hspace{-0,18cm}{\bf{.}}
A function $\Psi$ is said to be \emph{universally starlike} if
$\Psi$ is holomorphic in $\Bb{C}\setminus[1,+\infty)$, with the normalization
$\Psi(0)=0$, $\Psi^{\,\prime}(0)=1$, and if $\Psi$ maps every circular domain
$D\subset\Bb{C}\setminus[1,+\infty)$ with $0\in D$ one-to-one onto a domain
which is starlike with respect to the origin.
\end{definition}

In \cite[p. 289, Cor. 1.1]{RS},  Ruscheweyh,  Salinas and  Sugawa
characterized the universally starlike functions $\Psi$
as functions of the form
\begin{align}
\label{f5int}
\Psi(z)
=z\,\exp\bigg( \int_{[0,1]} \Log  \frac{1}{1-tz}\
 {\diff} \sigma(t)\bigg),\quad
  z\in\Bb{C}\setminus[1,+\infty),
\end{align}

\noindent
where $\sigma\in\mathcal{M}^{+}(\Bb{R})$ is uniquely determined by the
 requirements
\begin{align}
%\label{eq:univstarlike}
\label{f5aint}
{\mathrm{supp}}\, \,\sigma\subset[0,1],\quad 0\leq\sigma(\Bb{R})\leq 1,\quad \sigma(\{0\})=0.
\end{align}

\noindent
The formulation in \cite[p. 289, Cor. 1.1]{RS} was slightly different, but
 it is easy to see
that it is equivalent to the given one (see \cite[p. 719]{rbs}).

\bigskip

\noindent{\bf Universal starlikeness associated with the hypergeometric
function}.
Let $F_{a,b;c}(z)\!:=\!F (a,b;c;z)$ be the Gauss hypergeometric function
given by
\begin{align}
\label{f6int}
F_{{\fo{a,b;c}}} (z) =\dfrac{\Gamma(c)}{\Gamma(a)\Gamma(b)}
\sum\limits_{n = 0}^{\infty}\
\dfrac{\Gamma(a\!+\!n)\Gamma(b\!+\!n)}{ n!\,\Gamma(c+n)}z^{n}, \qquad
 z \in \Bb{D}\,,
\end{align}

\vspace{-0,1cm}
\noindent
where we restrict the parameters $a,b,c$ to be positive. It is well-known that the
hypergeometric function $F_{{\fo{a,b;c}}}$ extends holomorphically to the set
 $\Bb{C}\!\setminus\![1,+\infty)$. If $0\!<\!b\!<\!c$, we may see
% \blfootnote{\hypertarget{4}{${}^{{\fo{4)}}}$}\cite[pp. 6, 40]{Dur}.
% \hypertarget{5}{${}^{{\fo{5)}}}$}\cite[p. 384, 15.2.1]{olv}.
%\hypertarget{6}{${}^{{\fo{6)}}}$}\cite[p. 59, (10)]{erd1}.}
this from Euler's
integral representation  \cite[p. 59]{erd1}:\vspace{-0,1cm}
\begin{align}
\label{f7int}
F_{{\fo{a,b;c}}} (z)  = \dfrac{\Gamma(c)}{\Gamma(b)\Gamma(c-b)}
\int_{0}^{1}\dfrac{t^{b-1}(1-t)^{c-b-1} {\diff} t }{(1-tz)^{a}}\ .
\end{align}

%%%%%%%%%%%%%%%%%%%%%%%%%%%%%%%%%%%%%%%%%%%%%%%%%%%%%%%%%%%%%%%%%%
%%%%%%%%%%%%%%%%%%%%%%%%%%%%%%%%%%%%%%%%%%%%%%%

%%%%%%%%%%%%%%%%%%%%%%%%%%%%%%%%%%%%%%%%%%%%%%%%%%%%%%%%%%%%%%%%%%%
%%%%%%%%%%%%%%%%%%%%%%%%%%%%%%%%%%%%%%%%%%%%%%

\vspace{-0,1cm}
\noindent
 It was shown in
\cite[p. 292, Thm. 1.8]{RS} that the function $\Psi(z)=
zF_{a,b;c}(z)$ is  universally starlike provided that
$0< b \leq c$ and $0< a \leq \min \{1, c\}$.
%For background material on the hypergeometric function, we refer to e.g.
%\cite{abr}, \cite{erd1}, and \cite{olv}.
%(see \cite[(10), p.59]{erd1})
In particular, it follows from
this theorem of Ruscheweyh, Salinas, and Sugawa  that
for arbitrary  real triples $(a,b,c)$ with $0<b<c$ and
$0< a\leq\min\{1, c\}$, there exists a unique measure
$\sigma_{a,b;c} \in\mathcal{M}^{+}(\Bb{R})$ with
\begin{align}
\label{eq:sigmaprop1}
{\mathrm{supp}}\, \,\sigma_{a,b;c}\subset[0,1],\quad 0<\sigma_{a,b;c}(\Bb{R} )\le1,\quad
\sigma_{a,b;c} \left(\{0\}\right) = 0,
\end{align} \noindent
such that (compare with \eqref{f5int})\vspace{-0,1cm}
\begin{align}
\label{f8int}
F_{{\fo{a,b;c}}} (z)  = \exp\bigg(
\int_{[0,1]} \Log  \frac{1}{1-tz}\ {\diff}\sigma_{a,b;c}(t)\bigg),\quad
z\in\Bb{C}\setminus[1,+\infty).
\end{align} \noindent
This is an existence result and does not tell us how the measure
$\sigma_{a,b;c}$ looks like.

We let $\mathcal{P}_{\log}$ denote the collection of all $f\in\mathcal{P}$ with
$f(z)\not\equiv0$ such that the logarithmic derivative $f'/f\in \mathcal{P}$ as
well, i.e.,
\begin{align}\label{f9int}
    &  \mathcal{P}_{\log}
:= \left\{ \  f\in\mathcal{P}\setminus\{0\}  \quad  \Big|  \quad
  f^{\,\prime}/f\in\mathcal{P} \  \right\}\, .
\end{align}

\noindent
 Moreover, we write $\mathcal{P}_{\log}(-\infty, 1)$ for the subset
  of $\mathcal{P}_{\log}$
consisting of those functions that extend holomorphically to
$\Bb{C}\setminus[1,+\infty)$.
The following result explains the connection with the universally
starlike functions (see  \cite[ Cor. 2.3]{brs}).

\begin{thmx}\hspace{-0,18cm}{\bf{.}}
\label{inttheor1}
The function $\Psi (z) = z\psi (z)$ is universally starlike if and only if
$\psi(0)=1$ and $\psi\in\mathcal{P}_{\log}(-\infty, 1)$.
\end{thmx}

One of our main results is the following theorem.

\begin{theorem}\hspace{-0,18cm}{\bf{.}}\label{th3} For $a=b=1/2$ and $c=1$  the measure
$\sigma_{{{{1/2,1/2;1}}}}$ in the ex\-po\-nen\-tial
integral representation \eqref{f8int} under the conditions
\eqref{eq:sigmaprop1}  has
the following explicit expression:
   \begin{align}\label{f1th0}
   \sigma_{{{{1/2,1/2;1}}}}\big([\,0,\,x]\big)=
    \frac{1}{\pi}\arctan
\dfrac{ \hett   \left(x\right)}{ \hett  \left(1-x\right)},
\qquad  0< x<1.
\end{align}

\noindent
In particular, $\sigma_{{{{1/2,1/2;1}}}}$ is absolutely continuous
 with respect to the Lebesgue measure and has total variation $1/2$.
\end{theorem}

Clearly,
\begin{align*}
0<{\rm{Im}}\,\Log  \frac{1}{1-tz}=\Arg \frac{1}{1-tz}<\pi\,, \qquad
0<t<1\,,\,\,\,
z\in\Bb{H}\,,
\end{align*}

\noindent
so that in \eqref{f8int} for arbitrary $z\in\Bb{H}$ we have
\begin{align}\label{f1cth0}
0<\Arg   \hett  (z)={\rm{Im}}\,
\int_{[0,1]} \Log  \frac{1}{1-tz}\ {\diff}\sigma_{{{{1/2,1/2;1}}}}(t)<
\frac{\pi}{2}\ ,
\end{align}

\noindent
and correspondingly $-\pi/2\! <\!\Arg   \hett  (z)\!<\!0$ for $z$
lying in the lower open half-plane while $\Arg   \hett  (x)=0$ when $x< 1$.
Hence,
\begin{align}\label{f1ath0}\hspace{-0,2cm} \Arg   \hett  (z) \in
 \!\left(-\dfrac{\pi}{2},  \dfrac{\pi}{2}\right) \, , \ \
    {{\rm{Re}}} \, \hett   (z) \!>\! 0 \, , \quad z\!\in\!\Bb{C}\setminus
    [1, +\infty) \,.
\end{align}

\noindent This means that   for arbitrary $z\!\in\!\Bb{C}\!\setminus\!
 [1, +\!\infty)$  the equality \eqref{f8int}
can be written in the forms
%\vspace{-0,1cm}
 \begin{align}\label{f8eir}\hspace{-0,0cm}
    &  \Log   \hett   (z) = \frac{1}{\pi^{2 }} \!\int\limits_{0}^{1}
\frac{\dfrac{1}{t(1-t)}\Log  \left(\dfrac{1}{1-tz}\right)
 }{\vphantom{\dfrac{a}{a}} \hett  (t)^{2} +
 \hett  (1-t)^{2}} \, {\diff} t\ ,\hspace{-0,15cm} \end{align}

   \vspace{-0,2cm}
\noindent or, alternatively,
\vspace{-0,2cm}
\begin{multline}
    \Log    \hett   (z) \!=\!\dfrac{1}{\pi} \int\limits_{0}^{1}
     \dfrac{t}{1+t^2}  \arctan \dfrac{ \hett
       \left(1-t\right)}{ \hett
      \left(t\right)}  {\diff}  t
    \\
    +\! \dfrac{1}{\pi}\! \int\limits_{1}^{+\infty}\!\!
   \left( \frac{1}{t\!-\!z}\! - \!\frac{t}{1\!+\!t^2} \right)\,
   \arctan \dfrac{ \hett
   \left(1\!-\!{1}/{t}\right)}{ \hett
   \left({1}/{t}\right)}\, {\diff} t\, .
\label{f9eir}\end{multline}

\noindent
Next, since  for a positive real $\alpha$,
$\alpha\,\sigma_{{{{1/2,1/2;1}}}}$ is a nonnegative absolutely
 continuous  measure
with total mass $\alpha/2$,  the conditions of \eqref{f5aint}
hold with $\sigma=\alpha\,\sigma_{{{{1/2,1/2;1}}}}$, provided
that $0 < \alpha \leq 2$.
Given  the characterization of the universally starlike functions
\eqref{f5int}, we obtain  in Section~\ref{eir} the following assertion which for the
case when $\alpha=1$ is the special case of \cite[p. 292, Thm. 1.8]{RS}
 when $a=b=1/2$ and $c=1$.{\hyperlink{d1}{${}^{{{\ref*{case1}}}}$}}\hypertarget{bd1}{}
%\vspace{-0,2cm}
\begin{corollary}\hspace{-0,18cm}{\bf{.}}\label{th0}
 We have that $\hett^{\,\alpha}\in \mathcal{P}_{\log}(-\infty, 1)$, so
 that the function
$z\hett (z)^{\alpha}$ is universally starlike, provided
that $0 < \alpha \leq 2$.$\vphantom{A^{A^{A}}}$
\end{corollary}

\bigskip

\noindent{\bf The Schwarz triangle function.} \ We write $\he $ for
the function $F_{1/2,1/2;1}$.

\vspace{0.15cm}
In 1873, Schwarz \cite{sch} established the following fact
(see \cite[p. 97]{erd1}).

\begin{thmx}\hspace{-0,18cm}{\bf{.}}
\label{inttheor2}
The  Schwarz triangle function
\begin{align}
\label{f10int}
       \ie (z) :=   \imag \cdot \frac{\he (1-z)}{\he (z)},
\quad  z \in  (0,1)\cup \left(\Bb{C}\setminus\Bb{R}\right),
\end{align} \noindent
maps the set $(0,1)\cup \left(\Bb{C}\setminus\Bb{R}\right) $ one-to-one
onto the fundamental quadrilateral
\begin{align}\label{f11int}
\fet := \big\{\, z \in\Bb{H}\ \big| \,\  -1 < {\rm{Re}}\,  z < 1 \, , \,  \
|2 z -  1| > 1 , \ | 2 z +  1| > 1 \, \big\}  \, .
\end{align}
\end{thmx}

\noindent
 The function
$\lambda : \fet \to (0,1)\cup \left(\Bb{C}\setminus\Bb{R}\right)$
which is the inverse to $ \ie$
is called the  {\emph{elliptic modular function}}
and is the subject of a large literature (see  \cite[p. 99]{erd1} and
  \cite[p. 579]{olv}).
In the Poincar\'{e} half-plane model of the hyperbolic plane $\Bb{H}$,
$\fet$ is an ideal hyperbolic quadrilateral, and it is  the
set of all interior points of the fundamental domain for the subgroup
$\Gamma(2)$ of the modular group $\Gamma$ on the upper half-plane
$\Bb{H}$ (see \eqref{f6xbeirf}, \cite[p. 20]{ran} and  \cite[p. 115]{cha}).

\vspace{0.15cm} \noindent {\emph{Remark.}} The related function
  $\mu (r) = (\pi/2) \, K(\sqrt{1-r^{2}})
      /K (r) =
    (\pi/2) \, \he  (1-r^{2})/\he  (r^{2}) =
    (\pi/2\,{\rm{i}})\, \ie (r^{2})$
    %\mu (r) = \dfrac{\pi}{2} \, \dfrac{K\left(\sqrt{1-r^{2}}\right)
%       }{K (r)} =
%    \dfrac{\pi}{2} \, \dfrac{\he  (1-r^{2})}{\he  (r^{2})} =
%    \dfrac{\pi}{2 i}\, \ie \left(r^{2}\right)
for $0 < r < 1$ is known in the theory of quasiconformal mapping
and is called
the modulus of the  Gr\"{o}tszch ring \cite[Ch. 5]{vuo}. Here, $K$
 stands for the standard
elliptic integral of the first kind  \cite[p. 569]{abr}.

\vspace{0,2cm}
The intention of this paper is to derive an integral representation for
the logarithm of the theta function $\Theta_3$.
Apart from that, we are also motivated by the desire to write down an
elementary exposition of the basic properties of the theta functions
$\Theta_2,\Theta_3,\Theta_4$.{\hyperlink{d2}{${}^{\ref*{case2}}$}}\hypertarget{bd2}{}

\vspace{0.25cm}\noindent{\bf Outline of the paper.}
 We first describe the elementary properties of the functions $\he $
and $\ie$.
In Section~\ref{bf}, we obtain the basic formulas for the hypergeometric
function $\he $, and show that it is in the Hardy space $H^p$ of
the upper
half-plane for any $p$ with $2 < p < \infty$, and has the properties
\begin{align}
\label{f25preresinvlamhyp}
     \hspace{-0,2cm}{\rm{(a)}} \ \ \he  \in\mathcal{P}(-\infty, 1)\,,
      \quad {\rm{(b)}} \ \
  \he  (z) \neq 0 \,\,\ \text{for all}  \ \  z \in \Bb{C}\setminus
  [1, +\infty)\,.
\end{align}

\noindent
We also obtain that
\begin{align}
\label{f26preresinvlamhyp}
       \imag  \, \ie^{\,\prime} (z)\, \he  (z)^{2}
       =\frac{1}{\pi z(1-z)}
      \ , \quad z \in (0,1)\cup \left(\Bb{C}\setminus\Bb{R}\right).
\end{align}

 In Section~\ref{eir} we prove   Theorem~\ref{th3} and obtain an
 exponential integral representation of  $\ie/\, \imag  \,$.
% (see \eqref{f12eir} and \eqref{f13eir}).

The Schwarz triangle function $\ie$ enjoys the functional relation
\begin{align}
\label{f27preresinvlamhyp}
\ie (z)\ie (1-z) =-1, \quad z \in
\Lambda:= (0,1)\cup \left(\Bb{C}\setminus\Bb{R}\right),
\end{align}

\noindent
and in Section~\ref{sp} we obtain the relationship between the values of
$\ie$ on two sides of the cut along $(-\infty, 0)$:
\begin{align}
\label{f28preresinvlamhyp}
\ie (-x +  \imag  0 ) = 2 +  \ie (-x -
 \imag  0 ) \ , \quad  x>0\,.
\end{align}
\noindent Correspondingly, along the remaining cut $(1, +\infty)$, we obtain
that
\begin{align}
\label{f29preresinvlamhyp}
\ie (1 + x - \imag  0 ) =
\dfrac{\ie (  1+ x +  \imag  0 )}{ 1- 2 \ie (  1+ x +
  \imag  0 ) }\  ,
\quad x>0\,.
\end{align}

\noindent We also obtain the equality of sets
%\vspace{-0,1cm}
\begin{align}\label{f33preresinvlamhyp}
    &  \ie \left((0,1)\cup \big(\Bb{C}\setminus\Bb{R}\big)\right) = \fet \ ,
\end{align}

\noindent which constitutes part of the assertion of  Theorem~\ref{inttheor2}.
 %\vspace{-0,1cm}

In Section~\ref{mth}  we show that
 Lemma~\ref{th2} easily implies the following result which may
 be considered as Liouville-type theorem{\hyperlink{d3}{${}^{\ref*{case3}}$}}\hypertarget{bd3}{} for the
fundamental quadrilateral $\fet$, where $\fet = -1/\fet$
and \\[0,1cm] $\phantom{a}$ \hspace{0,5cm} $\Bb{H}\cap\partial
\fet = (1+  \imag   \Bb{R}_{>0}) \cup (-1+  \imag   \Bb{R}_{>0})\cup
 (1-  \imag   \Bb{R}_{>0})^{-1}\cup
 (- 1-  \imag   \Bb{R}_{>0})^{-1}$.

\begin{sectlemma}\hspace{-0,18cm}{\bf{.}}\label{th1}\hypertarget{hth1}{}
Suppose $f$ is holomorphic on $\fet$ and extends continuously to
its hyperbolic closure $\Bb{H}\, \cap \,{\rm{clos}}\,\fet$. Suppose
in addition
that the boundary values satisfy \
{\rm(a)}\  $f(z)=f(z+2)$ \  and  \ {\rm(b)}\
$ f (-1/z) = f \big(-1/(z+2) \big) $\ \  for each  $z \in -1 +  \imag
  \Bb{R}_{>0}$.
Finally, suppose there exist nonnegative integers $n_{\infty}$,
 $n_{0}$, and $n_{ 1}$ such that
\begin{align*}
     &   (1) \ \  &  \left|f (z)\right|
     =\oh\big(\exp\left(\pi\left( n_{\infty}
+1\right) |z|\right)\big)    \ ,   &     &  & \hspace{-0cm}
  \fet \ni z\to \infty \, ,  \\      &   (2) \ \  &
       \left|f (z)\right|=\oh\big( \exp \left(\pi
\left(n_{0}+1\right)|z|^{-1} \right)\big)  \ , &     &  &
\hspace{-0cm}  \fet \ni z \to 0 \, ,  \\
       &       (3) \ \  &    \left|f (z)\right|=
\oh\big(\exp \left(\pi\left( n_{ 1} +1\right)|z-\sigma|^{-1}\right)\big)
 \ , &     &  &   \hspace{-0cm} \fet \ni z \to\sigma \, ,
\end{align*}
where in $(3)$ we consider both $\sigma\in\{1,-1\}$.
Then there exists an algebraic polynomial $P$ of
degree $\le n_{\infty} + n_{0} + n_{1}$ such that
\begin{align*}
    &  f \big(\ie (z)\big) = \frac{\vphantom{\frac{a}{b}}
P(z)}{\vphantom{\dfrac{a}{b}}
     z^{{\fo{n_{\infty}}}} (1-z)^{{\fo{n_{0}}}} }  \ ,
\quad  z \in  (0,1)\cup \left(\Bb{C}\setminus\Bb{R}\right) \,.
\end{align*}
\end{sectlemma}

In Section~\ref{theta}  we introduce some notation for the standard
theta functions $\Theta_2$, $\Theta_3$, $\Theta_4$.
In Section~\ref{wid} we explain how to obtain Wirtinger's identity
(see \cite[p. 99]{erd1}  with $a\!=\!b\!=\!1/2, c\!=\!1$)
%\vspace{-0,2cm}
    \begin{align}\label{f1theor4}
        &  \Theta_{3}\big(\ie (z)\big)^{2}= \he  (z)   \ ,
         \quad z \in (0,1)\cup \left(\Bb{C}\setminus\Bb{R}\right) \,  ,
    \end{align}

 %\vspace{-0,2cm}
\noindent Together with nonvanishing property \eqref{f25preresinvlamhyp}(b),
 the equality of sets \eqref{f33preresinvlamhyp} and the relationships
 established in Section~\ref{theta},  this gives that
 %\vspace{-0,1cm}
 \begin{align}\label{f0theor3}
    & \Theta_{3}(z)\Theta_{4}(z)\Theta_{2}(z)\neq 0 \ , \quad  z \in \Bb{H}\,.
 \end{align}

%\vspace{-0,1cm}
  In Section~\ref{emf} we show  that $f (z) = ({\Theta_{2}(z)^{4} +
  \Theta_{4}(z)^{4}})/{\Theta_{3}(z)^{4}}$ satisfies the conditions
   of Lemma~{\hyperlink{hth1}{\ref*{th1}}} with $n_{\infty} = n_{0}=
   n_{ 1} = 0$ from which we obtain the Jacobi identity
  %\vspace{-0,2cm}
\begin{align}\label{f8theor3}
      &     \Theta_{2}(z)^{4} + \Theta_{4}(z)^{4} = \Theta_{3}(z)^{4}
       \ , \quad  z \in  \Bb{H} \ .
\end{align}

%\vspace{-0,2cm}
\noindent Then we show that
$f (z) = \Theta_{2}(z)^{4}/{\Theta_{3}(z)^{4}}$ enjoys
the conditions of Lemma~{\hyperlink{hth1}{\ref*{th1}}}
with $n_{1} = 1$
and $ n_{0}= n_{\infty} = 0 $,  from which we deduce the
property (see  \cite[p. 23]{erd1})
\vspace{-0,3cm}
\begin{align}\label{f1theor3}
      &    {\Theta_{2}\big(\ie (z)\big)^{4}}\big/{\Theta_{3}\big(
      \ie (z)\big)^{4}} = z  \ , \quad  z \in
      (0,1)\cup \left(\Bb{C}\setminus\Bb{R}\right) \ .
\end{align}

 \noindent
 As a side remark, this implies that $\ie$ is
univalent in the region $(0,1)\cup(\Bb{C}\setminus\Bb{R})$, which together with
the mapping property \eqref{f33preresinvlamhyp} proved in Section~\ref{sp} below,
completes the proof of Theorem~\ref{inttheor2}.

\vspace{0,1cm}
Since the elliptic modular function $\lambda$ is the inverse of
$\ie$, \eqref{f1theor3}
is the same as the identity $\lambda(z)=
\Theta_{2}(z)^{4}/{\Theta_{3}(z)^{4}}$
for $z\in\fet$, whence it is immediate that
the modular function $\lambda$ extends to a zero-free
holomorphic function in $\Bb{H}$ with period $2$:\vspace{-0,3cm}
\begin{align}\label{f1xtheor3}
    &  \lambda (z) =  {\Theta_{2}(z)^{4}}\big/{\Theta_{3}(z)^{4}}
     \ , \quad  z \in  \Bb{H} \,.
\end{align}

%\vspace{-0,2cm}

By combining \eqref{f1theor3} with \eqref{f8theor3} and
\eqref{f26preresinvlamhyp}, we find that%\vspace{-0,1cm}
 \begin{align}\label{f9theor3}
      &     \lambda^{\,\prime} (z)\!= \!\imag  \pi \lambda (z)
      \big(1-\lambda (z)\big)  \Theta_{3}\left(z\right)^{4}\!=\!
       \imag  \pi {\Theta_{2}\left(z\right)^{4}\Theta_{4}
       \left(z\right)^{4}}\big/{\Theta_{3}\left(z\right)^{4}}
        \, , \  z \!\in\!  \Bb{H} \, ,
\end{align}

%\vspace{-0,1cm}
\noindent which in its turn leads to the following two identities
(see Exercise 16, p. 22 of \cite{law}):
\begin{align}\label{f10theor3}
      &      \frac{\pi}{4\imag }\Theta_{2}(z)^{4}\!=\!
      \frac{\Theta_{4}^{\, \prime} (z)}{\Theta_{4}(z)}\! -\!
      \frac{\Theta_{3}^{\, \prime} (z)}{\Theta_{3}(z)} \, , \ \
      \frac{\pi}{4 \imag  } \Theta_{4}(z)^{4}\!=\!
      \frac{\Theta_{3}^{\, \prime} (z)}{\Theta_{3}(z)}\! - \!
      \frac{\Theta_{2}^{\, \prime} (z)}{\Theta_{2}(z)} \, , \ \
      z \in  \Bb{H} \, .
\end{align}

In addition, we see from \eqref{f1theor3} and \eqref{f1theor4}
that%\vspace{-0,2cm}
  \begin{align}\label{f5theor4}
        &  \Theta_{2}\big(\ie (z)\big)^{4}= z \he  (z)^{2}   \ ,
        \quad z \in (0,1)\cup \left(\Bb{C}\setminus\Bb{R}\right) \,  ,
    \end{align}

%\vspace{-0,2cm}
\noindent This is the function which is universally starlike by
Corollary~\ref{th0}.
\begin{corollary}\hspace{-0,18cm}{\bf{.}}\label{th5}
 The function $ \Theta_{2}\left(\ie \right)^{4}$ is universally
  starlike while $\Theta_{3}(\ie )^{4}$ belongs to the class
  $\mathcal{P}_{\log}(-\infty, 1) $,  and, for every $z\! \in \!\fet
   \!\setminus\!\{ \imag   \Bb{R}_{>0}\} $ we have\vspace{-0,1cm}
 \begin{align}\label{fxth5}
    &    {\rm{(a)}}\  \ ({\rm{Re}} \,  z) \cdot  {\rm{Im}} \,
    \Theta_{3}\left(z\right)^{4} \!>\! 0 \ ,    \quad
  {\rm{(b)}}\  \ ({\rm{Re}} \, z) \cdot  {\rm{Im}}\, \dfrac{\Theta_{3}^{\,\prime}
  \left(z\right)}{\lambda^{\prime}\,(z)\, \Theta_{3}\left(z\right)}\! >\! 0 \, .
 \end{align}
 \end{corollary}

In Section~\ref{log}, we introduce the logarithms  of the theta functions.
Moreover, in Section~\ref{eirth}, we apply Theorem~\ref{th3} in combination
with the Wirtinger identity \eqref{f1theor4} to obtain an integral
representation of $\log\Theta_{3}$ on the set
\begin{align}
\label{f3eirf}
\mathcal{F}^{\,{\tn{||}}}_{{\tn{\square}}}:= \fet\sqcup
(- 1 +  \imag   \Bb{R}_{>0})\sqcup( 1 +  \imag   \Bb{R}_{>0})\,.
\end{align}
\begin{corollary}\hspace{-0,18cm}{\bf{.}}\label{cor1}
 For arbitrary $z \in (0,1)\cup \left(\Bb{C}\setminus\Bb{R}\right)$
  we have $ \ie (z)  \in \fet$,

 \vspace{-0,15cm}
\begin{align}\label{f0yeirth}
    &   \arg \Theta_{3} \big(\ie (z)\big) \in \left(- \frac{\pi}{4}
    \, , \, \frac{\pi}{4}\right)  \, ,
\end{align}

\vspace{-0,5cm}
\begin{align} \label{f1eirf}
    & \hspace{-0,3cm}  \log\Theta_{3} \big(\ie (z)\big) =
    \dfrac{1}{2 \pi^{2 }} \int\limits_{0}^{1}
   \frac{1}{\he (t)^{2} +
\he (1-t)^{2}}\, \Log  \left(\frac{1}{1-tz}\right)\,
\frac{{\diff} t}{t(1-t)} \ , \hspace{-0,2cm}
\end{align}

\vspace{0.25cm}
\noindent where \ $\fet = \ie \big( (0,1)\cup \left(\Bb{C}
\setminus\Bb{R}\right)\big)$, in accordance with \eqref{f33preresinvlamhyp}.
Moreover,
\begin{align}\label{f1aeirf}
    &  \hspace{-0,3cm}  \log\Theta_{3} \left( \pm 1\! +\!  \imag
     y\right)\! =\! \dfrac{1}{2 \pi^{2 }} \int\limits_{0}^{1}
    \frac{1}{\he (t)^{2} +
\he (1-t)^{2}}\, \ln\left(\frac{1}{1+t x}\right)\, \frac{{\diff} t}{t(1-t)}
 \ ,\hspace{-0,2cm}
\end{align}

\vspace{0.15cm}
\noindent where $\Theta_{3} \left( 1\! +\!  \imag   y\right) =
\Theta_{3} \left( - 1\! +\!  \imag   y\right)\in (0,1)$,  \
$\log\Theta_{3} \left( \pm 1\! +\!  \imag   y\right) = \ln\Theta_{3}
\left( \pm 1\! +\!  \imag   y\right)$,

\vspace{-0,15cm}
\begin{align*}
    & y = y (x)= \dfrac{\he   \big( 1/(1 + x)\big)}{\he
    \big( x/(1 + x)\big)}   \ ,
    \quad
    \left\{\begin{array}{l}
    y (0)= +\infty \ ,\\[0,2cm]
    y (+\infty) = 0 \ ,
    \end{array}\right.
    \ \  \ \dfrac{{{\diff}} y (x)}{{{\diff}} x} < 0 \ , \ x > 0 \ .
\end{align*}
\end{corollary}

\vspace{0.25cm}
 Let us consider the periodized set
  \begin{align}
 \label{f2eirf}
\mathcal{F}^{\infty}_{{\tn{\square}}}:=
\bigcup\nolimits_{\,{\fo{m \in \Bb{Z}}}} \left(2m +
\mathcal{F}^{\,{\tn{||}}}_{{\tn{\square}}}\right) \,.
\end{align}
\noindent
 From Corollary~\ref{cor1} we drive an integral
formula for $\log\Theta_{3}$,
\begin{align}
\label{f7eirf}
\log\Theta_{3}\! \left(z\right) =
\frac{1}{2 \pi} \int\limits_{0}^{+\infty} \ \
{\Log  \left(\dfrac{1}{1-\lambda ( \imag
 \tau)\lambda(z)}\right)}
\,\frac{{\diff} \tau}{1+ \tau^{2}},
\qquad z \in \mathcal{F}^{\infty}_{{\tn{\square}}} \ ,
\end{align}

\noindent
or, equivalently (cf. \eqref{f1xtheor3}),
\begin{equation*}
\log\Theta_{3}\! \left(z\right) =  \dfrac{1}{2 \pi}
\int\limits_{0}^{+\infty}\!\!
{\Log  \left(\dfrac{\Theta_{3}( \imag   \tau)^{4}
\Theta_{3}(z)^{4} }
{\Theta_{3}( \imag  \tau)^{4} \Theta_{3}(z)^{4}-\Theta_{2}(
 \imag   \tau)^{4}
\Theta_{2}(z)^{4}}\right)}
\,\frac{{\diff} \tau}{1+ \tau^{2}},
\ \  z \!\in\! \mathcal{F}^{\infty}_{{\tn{\square}}} .
\end{equation*}

\noindent
We should contrast this integral formula with the classical
 series representation (see  \eqref{f3elemtheta},  \eqref{f11elemtheta}
 and compare, e.g., with \cite[p. 338, (4.2)]{ber})
\begin{align*}
    &   \log \Theta_{3}  (z) = \sum\limits_{ n \geq 1}  \dfrac{2}{2n-1}
    \,\dfrac{ {\rm{e}}^{\, \imag
     \pi (2n-1) z}  }{1+{\rm{e}}^{\, \imag   \pi (2n-1) z}  }
    \ , \quad   z \in \Bb{H} \,.
\end{align*}

\medskip

\noindent{\bf The Berndt formula.}
Let ${\rm{SL}}_{2}(\Bb{Z})$ be the multiplicative group of all
$2\times2$ matrices
\begin{align*}
M=\left(\begin{matrix} a & b \\ c & d \end{matrix}\right)
\quad \text{with}\quad a, b, c, d \in\Bb{Z}\,,\,\,\,ad - bc=1\,.
\end{align*}

\noindent
To such a matrix we associate a M\"obius transformation
\begin{align*}
\phi_{M
%{\nor{(\begin{smallmatrix} a & b \\ c & d \end{smallmatrix})}}
}(z) :=
\frac{a  z + b}{c  z + d},  \qquad z\in \Bb{H},
\end{align*}

\noindent
and note that the M\"obius transformation retains all the
 information about the
matrix except that the matrices $M$ and $-M$ give rise to
 the same M\"obius
transformation.
We consider the following subsets of the group
${\rm{SL}}_{2} (\Bb{Z})$:
\begin{align*}  &
{\rm{SL}}_{2} (2, \Bb{Z}) := \left\{  {\fo{\begin{pmatrix}
a & b \\ c & d \\ \end{pmatrix}}}\! \in\!
{\rm{SL}}_{2}(\Bb{Z}) \, \left| \, {\fo{\begin{pmatrix} a & b \\ c & d \\
\end{pmatrix}}}\! \equiv\! {\fo{\begin{pmatrix} 1 & 0 \\ 0 & 1 \\
\end{pmatrix}}} ({\rm{mod}}\, 2)
\, \right\}\right.  , \
 \\
  &  {\rm{SL}}_{2} (\vartheta, \Bb{Z}) :=
{\rm{SL}}_{2} (2, \Bb{Z}) \cup \left\{  {\fo{
\begin{pmatrix} a & b \\ c & d \\ \end{pmatrix}}} \in
{\rm{SL}}_{2}(\Bb{Z}) \, \left| \,   {\fo{
\begin{pmatrix} a & b \\ c & d\\ \end{pmatrix}}} \equiv {\fo{
\begin{pmatrix} 0 & 1 \\ 1 & 0 \\ \end{pmatrix}}}  ({\rm{mod}}\, 2)
\, \right\}\right. .\
\end{align*}

\noindent
We associate with these subsets the induced collections of M\"obius
transformations:\vspace{-0,1cm}
\begin{align}\label{f6xbeirf}
\begin{array}{l}
 \Gamma \,:= \ \{\, \phi_{M} \  | \  M
 \in {\rm{SL}}_{2} (\Bb{Z})\, \} \,, \\[0,15cm]
\Gamma_{\vartheta}:= \left.
\left\{\, \phi_{M} \  \right| \
   M \in
{\rm{SL}}_{2} (\vartheta, \Bb{Z}) \right\}, \quad
    \Gamma(2) := \left.\left\{\, \phi_{M} \  \right| \
   M\in {\rm{SL}}_{2} (2, \Bb{Z}) \right\}.
\end{array}
\end{align}

\noindent
It follows from    \cite[p. 15, Def. 3.3; p. 16, Thm. 3.1]{bkn} that
 \begin{align}
\label{f6beirf}
\Bb{H} =
\bigcup\nolimits_{  {\fo{M=(
\begin{smallmatrix} a & b \\ c & d \end{smallmatrix})
\in  {\rm{SL}}_{2} (\vartheta, \Bb{Z})  }} } \ \
   \phi_{M}\left({\rm{clos}}_{{\tn{\Bb{H}}}}\left(
   \mathcal{F}_{{\fo{\Gamma_{\!\vartheta}}}}\right)\right) , \
  \end{align}

   \noindent
where
\begin{align}\label{f6ybeirf}
{\rm{clos}}_{{\tn{\Bb{H}}}}
\left(\mathcal{F}_{{{\Gamma_{\!\vartheta}}}}\right) :=
\{\, z \in \Bb{H} \, | \,  -1 \!\leq  {\rm{Re}}\,  z \leq  1 ,
 \  | z | \geq 1
  \, \} \! \subset
 \mathcal{F}^{\,{\tn{||}}}_{{\tn{\square}}}.
\end{align}

This means that for arbitrary $z\in\Bb{H}$,
there exists at least one matrix
$M\in {\rm{SL}}_{2} (\vartheta, \Bb{Z})$ such that $\phi_M(z)
\in \mathcal{F}^{\,{\tn{||}}}_{{\tn{\square}}}$.
Now, according to the formula of Berndt (see \cite[p. 339,
Thm.  4.1]{ber}), we have that for any
$M=(\begin{smallmatrix} a & b \\ c & d \end{smallmatrix})\in
{\rm{SL}}_{2} (\vartheta, \Bb{Z}) $ with $c>0$ and $z \in\Bb{H}$,
\begin{align}\label{f6eirf}
    & \log\Theta_{3}\! \left({\nor{\dfrac{a  z\! + \!b}{c
     z\! +\! d}}}\right) \!=\! \log\Theta_{3} (z)\! +\!
     \dfrac{1}{2} \Log
    \dfrac{c z\! +\! d}{ \imag  } \!+\!
    \dfrac{\pi  \imag  }{4} \sum\limits_{k=1}^{c}
    (-1)^{{\fo{k\!+\!1\! +\! \left\lfloor{{\dfrac{k d}{c}}}\right\rfloor  }}} \, ,
\end{align}

\noindent
where, in view of \eqref{f11elemtheta}, the branch of the
logarithm $\log\Theta_3$ is selected
which is real-valued on the positive imaginary axis,
and $\lfloor x\rfloor$ denotes the integer part of $x\in\Bb{R}$.

Hence, \eqref{f1eirf} and \eqref{f1aeirf} combined with \eqref{f6eirf}
% \eqref{f1eirf}, \eqref{f1aeirf}
supply an integral representation for
$\log\Theta_{3}(z)$ when $z\in\Bb{H}$, since we may find
a matrix $M=(\begin{smallmatrix} a & b \\ c & d \end{smallmatrix})\in
{\rm{SL}}_{2} (\vartheta, \Bb{Z}) $ with $c>0$ such that
$\phi_{M}(z) \in\mathcal{F}^{\,{\tn{||}}}_{{\tn{\square}}}$.
The corresponding formulas for $\log\Theta_{2}$ and $\log\Theta_{4}$
may be found from the formula for $\log\Theta_3$, as the following
relations hold for any $z\in\Bb{H}$
 (see \eqref{f12elemtheta} below):
%\vspace{-0,2cm}
\begin{align*}%\label{f6ceirf}
    &   \log \Theta_{4}\left(z\right)\! =\! \log
     \Theta_{3}\left(z\!+\!1\right) \, , \   \log
     \Theta_{2}\left(z\right)\!  =\! \log \Theta_{3}
    \left(1\! -\!1 / z\right)\!  - \!\dfrac{1}{2}
    \Log   \dfrac{z}{ \imag  } \ .
\end{align*}

%%%%%%%%%%%%%%%%%%%%%%%%%%%%%%%%%%%%%%%%%%%%%%%%%%%%%%%%%%%%%%%%%%%%

%\vspace{-0,4cm}
\section[\hspace{-0,30cm}. \hspace{0,11cm}Basic facts about $F$.]{\hspace{-0,095cm}{\bf{.}}  Basic facts about  $\het:= F_{1/2\, ,  1/2\, ; \, 1}$}
\label{bf}

%\vspace{-0,2cm}
To simplify the notation, we denote by $\overline{\Bb{D}}$
 and $\overline{\Bb{H}}$ the closures of $\Bb{D}$ and $\Bb{H}$
  in the complex plane $\Bb{C}$, respectively.
The series \eqref{f6int} and the Euler formula \eqref{f7int}
 for the hypergeometric function in \eqref{f10int} have the
 following form
\begin{align}
     \label{f1preresinvlamhyp}
    &  \he  (z) =\dfrac{1}{\pi} \sum\nolimits_{n = 0}^{\infty}
\dfrac{\Gamma(n+1/2)^{2}}{(n !)^{2}} z^{n}  \ , \quad  z \in
\Bb{D} \, , \\[0cm]
\label{f2preresinvlamhyp}
    &  \he  (z) = \dfrac{1}{\pi }  \int_{0}^{1} \dfrac{ {\diff}
      t}{\sqrt{t(1-t)(1-tz)  }}
    \ , \quad z \in \Bb{C}\setminus [1, +\infty) \ .
\end{align}

\noindent In addition, the function $\he $ satisfies
 the Pfaff formula (see \cite[p. 79]{and})
\begin{align}\label{f5preresinvlamhyp}
    &  \he  (z) = \dfrac{1}{\sqrt{1-z}} \ \he
\left(\dfrac{z}{z-1}\right)   \ , \quad
 z \in \Bb{C}\setminus [1, +\infty)  \ , \
\end{align}

\noindent and for $z \in (1 + \Bb{D})\setminus [1, +\infty) =
  (1 + \Bb{D})\setminus [1, 2]$  has the following expansion (see
\cite[p. 559, 15.3.10]{abr}){\hyperlink{d4}{${}^{\ref*{case4}}$}}\hypertarget{bd4}{} %\vspace{-0,3cm}
\begin{multline}
        \he  (z)    = \dfrac{1 }{\pi}\he  (1-z)\, \Log   \dfrac{1}{1-z} \\
             + \dfrac{2}{\pi^{2}} \sum\limits_{n=0}^{\infty}
              \dfrac{\Gamma (n+1/2)^{2}}{(n !)^{2}}
              \left[\sum\limits_{k\geq n}
               \dfrac{1}{(k +1)(2k +1)}\right](1-z)^{n} \ ,
\label{fa6preresinvlamhyp}\end{multline}

\noindent where \cite[p. 658,  5.1.8.2]{prud1}
$\sum_{k\geq 0} (k +1)^{-1}(2k +1)^{-1}=  2\ln 2$ and the
 summand corresponding to $n=0$ in the series above is equal to
$\pi^{-1}\ln 16$. Therefore{\hyperlink{d5}{${}^{\ref*{case5}}$}}\hypertarget{bd5}{}
\begin{align}
\nonumber %\label{f8preresinvlamhyp}
    & {\rm{(a)}} \   \he  (z) =\dfrac{1}{\pi}\,
\Log  \dfrac{16}{1-z}\ +
{\Oh} \left(|1-z| \ln \dfrac{1}{|1-z|}\right)
 \ , \ \      z \to 1  , \
  z \not\in  [1, 2]  \, ,  \\[0,1cm]
  \label{f7preresinvlamhyp}  &
{\rm{(b)}} \  \he  (z) = 1 + {\Oh} (|z|)  \ ,  \ \
     z \to 0  , \ z \in \Bb{D} \ ,
\\[0,2cm]
    &  \nonumber  %\label{f9preresinvlamhyp}
{\rm{(c)}} \   \he  (z) =\dfrac{\Log   16 (1-z)}{\pi \sqrt{1-z}}
 + {\Oh} \left(\dfrac{\ln |z|}{|z|^{3/2}}\right) \ ,  \ \
    |z| \to +\infty \   , \
z \not\in  [2, +\infty)  \ .
\end{align}

\noindent Here, the property \eqref{f7preresinvlamhyp}(b) is
immediate  from \eqref{f1preresinvlamhyp},  while \eqref{f7preresinvlamhyp}(a)
results from \eqref{fa6preresinvlamhyp} by the asymptotics
of \eqref{f7preresinvlamhyp}(b) applied to the term $\he  (1-z)$
and from the observation that for $z\to 1$ only the term with $n=0$
survives in the series in  \eqref{fa6preresinvlamhyp}.  Finally,
 \eqref{f7preresinvlamhyp}(c) is obtained by substitution of
 \eqref{f7preresinvlamhyp}(a), with $z/(z-1)$ in place of $z$,
  in \eqref{f5preresinvlamhyp} under the condition that $z/(z-1)
   \in (1 + \Bb{D})\setminus [1, 2]$ which is equivalent to
 $z \not\in (1 + \overline{\Bb{D}}) \cup [2, +\infty)$.

It is well-known that the functions $\he  (z)$ and $\he  (1-z)$ are
 two independent solutions of the Euler hypergeometric differential equation
(see \cite[p. 75]{and})
\begin{align}\label{f20preresinvlamhyp}
    &   z (z-1) y^{\,\prime\,\prime}(z) + (2 z-1) y^{\,\prime}(z) +
     y(z)/4 = 0   \ , \quad  z \in \Bb{D} \, ,
\end{align}

\noindent whose linear independence can be easily deduced from the
 formula \eqref{f1preresinvlamhyp}, invariance of
 \eqref{f20preresinvlamhyp} with respect to the change of
 the variable $z\mapsto 1-z$ and the expansion \eqref{fa6preresinvlamhyp}.
  The constant $A$ in the formula  for  the Wronskian in
  \cite[p. 136, Lem. 3.2.6]{and}
\begin{align*}
    &   W_{\!{\tn{\triangle}}} (z) \! :=\! \he^{\,\prime} (z)
    \he  (1-z)\! +\! \he  (z)\he^{\,\prime} (1-z) \!= \!\dfrac{A}{z(1-z)}
\end{align*}

\noindent where (see \cite[p. 557, 15.2.1]{abr})%\vspace{-0,2cm}
\begin{align*}%\label{f22preresinvlamhyp}
    &   \he^{\,\prime} (z) = \dfrac{1}{4}\, \jet (z)  \ , \quad
     z \in \Bb{C}\setminus [1, +\infty) \ ,
\end{align*}

%\vspace{-0,1cm}
\noindent can be calculated by letting  $z \to 1$ in this formula
 and by using the relations $(1-z)\he  (z)\to 0$,
  $(1-z)\he^{\,\prime} (z)\to 1/\pi$
as  $z \to 1$, which are immediate from the expansion
\eqref{fa6preresinvlamhyp} and its differentiated form{\hyperlink{d6}{${}^{\ref*{case6}}$}}\hypertarget{bd6}{}.
This gives (cf.\cite[p. 6, (2.5)]{aga})
\begin{align}\label{f21preresinvlamhyp}
    &     W_{\!{\tn{\triangle}}} (z)  \!=\!   \dfrac{1}{\pi z(1-z)}
      \ , \ \  z \in (0,1)\cup \left(\Bb{C}\setminus\Bb{R}\right) \, .
\end{align}

 For arbitrary $x > 1$ and $t\in [0,1]$ we obviously have
 $|1-t(x \pm \imag   \varepsilon)|\geq |1-tx|$, $\varepsilon > 0$,
 ($\varepsilon \downarrow 0$ means that $\varepsilon \to 0$
 and $\varepsilon > 0$)%\vspace{-0,3cm}
\begin{align}
    &  \label{f10preresinvlamhyp}
 \lim\limits_{\varepsilon \downarrow 0} \sqrt{1-t(x \pm \imag
  \varepsilon)} =
\left\{\begin{array}{ll} \sqrt{1-t x } \ , \   & \hbox{if} \ \
 0 \leq t < 1/x \ ;  \\[0.1cm]
  {\rm{e}}^{{\fo{\mp i \pi/2}}} \sqrt{t x-1}  \ , \   & \hbox{if}
   \ \  1/x < t \leq 1 \ ; \end{array}\right.
\end{align}

\noindent and therefore, by \eqref{f2preresinvlamhyp} and the
Lebesgue dominated convergence theorem \cite[p. 26, 1.34]{rud},
there exist the finite limits ($\he  (x\pm  \imag   0):=
 \lim_{\varepsilon \downarrow 0} \  \he  (x\pm \imag  \varepsilon)$
 are called "radial" limits at $x$)%\vspace{-0,2cm}
\begin{align}
    &   \label{f11preresinvlamhyp}
    \hspace{-0.25cm}
\he  (x\pm  \imag   0)= \dfrac{1}{\pi }  \int\limits_{0}^{1/x}
\dfrac{{{\diff}} t}{\sqrt{t(1-t)(1-tx)  }} \pm  \dfrac{ \imag  }{\pi }
 \int\limits_{1/x}^{1} \dfrac{ {\diff}  t}{\sqrt{t(1-t)(tx -1)  }}  \ , \
\end{align}

%\vspace{-0,3cm}
\noindent where{\hyperlink{d7}{${}^{\ref*{case7}}$}}\hypertarget{bd7}{}%\vspace{-0,1cm}
\begin{align*}
     \int\limits_{0}^{1/x}\dfrac{ {\diff}  t}{\sqrt{t(1-t)(1-tx)  }} &
      = \dfrac{1}{\sqrt{x}}
 \int\limits_{0}^{1/x}\dfrac{ {\diff}  t}{
 \sqrt{t\left(1-t\right)\left(\vphantom{A^{A}}(1/x) - t\right)}}   =
 \dfrac{\pi}{\sqrt{x}} \ \he  \left(\dfrac{1}{x}\right) , \
\end{align*}

%\noindent and
\begin{align*}
      \int\limits_{1/x}^{1} \dfrac{ {\diff}  t}{\sqrt{t(1\!-\!t)(t x\! -\! 1)  }}&
\!= \! \dfrac{1}{\sqrt{x}} \int\limits_{0}^{1\!-\!1/x} \dfrac{ {\diff}
 t}{\sqrt{t \left(1-t\right)\left(1\! -\! \dfrac{1}{x}\!-\!t\right)  }}
 \! = \! \dfrac{\pi}{\sqrt{x}} \ \he \left( 1\!-\! \dfrac{1}{x}\right)   .
\end{align*}

\noindent  Thus
(see \cite[p. 491, 19.7.3; p. 490, 19.5.1]{olv}),%\vspace{-0,1cm}
\begin{align}
    &   \label{f12apreresinvlamhyp}
\he  (x\!\pm\!  \imag   0)\!=\!  \dfrac{1}{\sqrt{x}} \, \he
\left(\dfrac{1}{x}\right)\! \pm \!\dfrac{ \imag  }{\sqrt{x}} \,
\he \left( 1\!- \!\dfrac{1}{x}\right)  ,
\quad x > 1 ,
\end{align}

%\vspace{-0,1cm}
\noindent which can also be written as%\vspace{-0,1cm}
\begin{align}
    &  \label{f12bpreresinvlamhyp}
\he  (1 +x\pm  {\rm{i}}  0)=  \dfrac{1}{\sqrt{x}} \ \he  \left(-
\dfrac{1}{x}\right) \pm  \imag
 \ \he \left( - x \right)  \ , \quad x > 0 \ ,
\end{align}

%\vspace{-0,1cm}
\noindent by virtue of the following equivalent forms of
\eqref{f5preresinvlamhyp}%\vspace{-0,2cm}
\begin{align} \label{f13apreresinvlamhyp}
     \he  (-z) & =\dfrac{1}{\sqrt{1+z}} \
\he  \left(\dfrac{z}{1 + z}\right)
 \, \ ,   &     &    z \in \Bb{C}\setminus (-\infty, -1] \, ,  \\
\dfrac{1}{\sqrt{z}}  \ \he  \left(- \dfrac{1}{z}\right)&   =
\dfrac{1}{\sqrt{1+z}} \  \he  \left(\dfrac{1}{1 + z}\right)
 \,  \ ,   &     &    z \in \Bb{C}\setminus [ -1, 0] \, .
\label{f13bpreresinvlamhyp}
\end{align}

\noindent We  observe that the relation \eqref{f12apreresinvlamhyp}
can also be obtained from  one of Kummer's transformation
rules,{\hyperlink{d26}{${}^{\ref*{case26}}$}}\hypertarget{bd26}{}
\begin{align}\label{f2rem1}
    &  (-z)^{-1/2}\he  (1/z) - \imag
    \he  (z)\, {\rm{sign}} ({\rm{Im}}\, z) =z^{-1/2} \he  (1\!-\!1/z) \ ,
    \ \   z\!\in \! \Bb{C}\setminus\Bb{R} \,,
\end{align}

\noindent   where the principal branch of the square root is used and  ${\rm{sign}} (x)$
 is equal to $-1$ if $x < 0$, $0$ if $x=0$ and $1$ if $x>0$ (cf. \cite[p. 106, (27)]{erd1}).
%$ \sigma (z) = 1$ if $z\in \Bb{H}$ and $ \sigma (z) =-1$ otherwise.

\vspace{0,05cm}
 It follows from \eqref{f7preresinvlamhyp}(a)--(c)
that $\he $ belongs to the Hardy space $ H^{p}  $  for
arbitrary $2 < p < \infty${\hyperlink{d8}{${}^{\ref*{case8}}$}}\hypertarget{bd8}{}.
According to the Schwarz  integral formula
applied to $\imag  \he $ (see \cite[p. 227]{saf}, \cite[p. 128]{koo})
 we get from  \eqref{f12apreresinvlamhyp}  and from the obvious
 consequence  of \eqref{f2preresinvlamhyp},  ${{\rm{Im}}}\, \he  (x)\! =\! 0$,
  $-\infty\! < \!x < \!1 $,
   that%\vspace{-0.1cm}
\begin{align}\label{f24preresinvlamhyp}
    &  \he  (z)=
\dfrac{1}{\pi} \int\nolimits_{1}^{\infty} \he \left( 1-\dfrac{1}{t} \right)
 \dfrac{ {\diff}  t}{(t-z)\sqrt{t}}  \ \ , \quad z \in \Bb{C}\setminus [1, +\infty) \ .
\end{align}

%\vspace{-0.1cm}

%\vspace{-0.1cm}
\noindent
 By using the Cauchy theorem (see \cite[p. 89, 6.6]{con}) and
 \eqref{f7preresinvlamhyp}(a), it can be easily derived  that
 the contour of integration in \eqref{f24preresinvlamhyp} can
 be changed to $1\! + \! \imag \Bb{R}_{\geqslant 0}$ if
 ${\rm{Im}}\, z \! <\!  0$ and to  $1\!  -\!  \imag \Bb{R}_{\geqslant 0}$
  if ${\rm{Im}} \,z \! >\!  0$, where  $\Bb{R}_{\geqslant 0}\!:=\! [\,0,+\infty)$.
  So that{\hyperlink{d9}{${}^{\ref*{case9}}$}}\hypertarget{bd9}{} \vspace{-0.1cm}
 \begin{align}\label{f24xpreresinvlamhyp}
    & \he  (z)=  {\rm{e}}^{{\fo{- \dfrac{\imag\pi\sigma}{2} }}}
    \int\limits_{0}^{\infty}
    \dfrac{\he \left(\dfrac{t}{t+\imag\sigma} \right)
    {{\diff}}  t}{(1 - \imag
     t \sigma -z)\sqrt{1- \imag t \sigma   }}
    \ , \quad
    z \in \sigma\cdot\Bb{H} \ , \  \sigma\in\{1, -1\} \,,
 \end{align}

  \vspace{-0.1cm}
 \noindent
 and hence the function $ \he  $  allows a holomorphic extension
  from the upper half-plane
  $ \Bb{H}$ to
 $ \Bb{C}\setminus (1 - \imag\Bb{R}_{\geqslant 0}) $
  and from  the lower  half-plane $- \Bb{H}$  to
 $ \Bb{C}\setminus (1 + \imag\Bb{R}_{\geqslant 0}) $.
 This means that
$\he $ can be continuously extended from  $ \Bb{H}$
 to $\Bb{H}\!\cup\!(\Bb{R} \!\setminus \!\{1\})$ and from $- \Bb{H}$
 to $(- \Bb{H})\!\cup\!(\Bb{R}\!\setminus \!\{1\})$,
 and that for every point  $x \! \in  \!\Bb{R} \setminus \!\{1\}$ there
 exist two finite "radial" limits satisfying $  \he \! (x\pm\!  \mathrm{i}
  \hspace{0,015cm}   0)
 =\lim_{\, \Bb{H} \ni z \to 0} \  \he  (x\pm z)$,    which can be
 written in the following form, by virtue of \eqref{f12apreresinvlamhyp}
  and \eqref{f13apreresinvlamhyp},\vspace{-0.1cm}
\begin{align}\label{f18preresinvlamhyp}
    & \lim_{{\fo{\, \Bb{H} \ni z\! \to 0}}} \  \he  (x\pm z)\!=\!
    \left\{
      \begin{array}{ll}
\begin{displaystyle}
   \dfrac{1}{\sqrt{x}} \, \he  \left(\!\dfrac{1}{x}\!\right)\!
   \pm \!\dfrac{ \imag  }{\sqrt{x}} \,
   \he \left(\! 1\!- \!\dfrac{1}{x}\!\right) ,
\end{displaystyle}
 & \hbox{if} \ \ x > 1 \ ;  \\[0,35cm]
   \he  (x) \, ,    & \hbox{if} \, \ 0\! \leq\! x\! <\! 1  ; \\[0,0cm]
\begin{displaystyle}
\dfrac{1}{\sqrt{1+|x|}} \  \he  \left(\dfrac{|x|}{1 + |x|}\right) \, ,
\end{displaystyle} & \hbox{if} \ \  x < 0 \ . \
      \end{array}
    \right.
\end{align}

%\vspace{-0.1cm}
\noindent
In view of the obvious consequence  of \eqref{f1preresinvlamhyp}%\vspace{-0.3cm}
\begin{align}\label{f23preresinvlamhyp}
    &  \he  (x) > 0 \ , \quad  0 \leq x < 1 \ , \
\end{align}

\noindent the expressions  \eqref{f18preresinvlamhyp}
yield %%\vspace{-0.2cm}
\begin{align}\label{f19preresinvlamhyp}
    &  \begin{array}{l}
              {{\rm{Re}}}\, \he  (x\!\pm\!  \mathrm{i}
               \hspace{0,015cm}   0) \!>\! 0, \quad  x
              \in \Bb{R}\setminus\{1\}; \quad \ \he  (x)\! >\! 0 \ ,
                \ \ -\infty \!< \!x \!<\! 1;
             \\[0.25cm]
            {{\rm{Im}}}\, \he  (x\!\pm\!  \imag
             0)\! =\! 0  , \quad -\infty\!
             < \!x < \!1 ; \quad
 {{\rm{Im}}}\, \he  (x\!+ \! \imag   0)\! >\! 0
  , \quad x\! >\! 1 \,.
       \end{array}
\end{align}

 The validity of \eqref{f25preresinvlamhyp} follows from
 \eqref{f24preresinvlamhyp}, \eqref{f2int},  \eqref{f3int} and
 \eqref{f19preresinvlamhyp}
(see  \cite[p. 604, Rem. 2.1]{kus}),   which in turn proves the
correctness of the definition \eqref{f10int}.
 This allows to obtain \eqref{f26preresinvlamhyp} from \eqref{f21preresinvlamhyp}
and the identity $\imag  \, \ie^{\,\prime} (z)\, \he  (z)^{2} =
 W_{\!{\tn{\triangle}}} (z)$, which follows from the definition of
  the Wronskian $W_{\!{\tn{\triangle}}}$ and that of the
Schwarz triangle function $\ie$.

 By setting $z=x > 0$ in \eqref{f13apreresinvlamhyp} and
 \eqref{f13bpreresinvlamhyp} we derive from \eqref{f12bpreresinvlamhyp}
  that %\vspace{-0,4cm}
\begin{align}\label{f14preresinvlamhyp}
    &  \dfrac{\he  (1 +x\pm  \imag   0)}{ \he \left( - x \right)}=\pm
    \imag   + \dfrac{\he  \left({1}\big/{(1 + x)}\right)}{\he
     \left( 1- {1}\big/{(1 + x)}\right)}
 \ , \quad x > 0 \ .
\end{align}

%%\vspace{-0.2cm}

%%%%%%%%%%%%%%%%%%%%%%%%%%%%%%%%%%%%%%%%%%%%%%%%%%%%%%%%%%%%%%%%%%%%
%\vspace{-0,4cm}
\section[\hspace{-0,30cm}. \hspace{0,11cm}Exponential integral representation of $F$.]{\hspace{-0,095cm}{\bf{.}}  Exponential integral representation of $F_{1/2\, ,  1/2\, ; \, 1}$}
\label{eir}

%\vspace{-0,2cm}
In the sequel, for $\mu\in\mathcal{M}^{+}(\Bb{R})$  we look at
 the  spaces of Borel measurable real-valued functions
$L^p (\Bb{R}, {\diff} \mu )$, $1\!\leq\! p\! <\! \infty$,  and for
arbitrary function $v\!\in\! L^p(\Bb{R},{\diff} x)$ with $1\!<\! p\! <\! \infty$
we  use  the notation for the (sign changed) Hilbert transform%\vspace{-0,2cm}
\begin{align*}%\label{2.5}
\widetilde{v}(x):=
 \frac{1}{\pi}\lim_{\varepsilon\downarrow 0}\int_{|t-x|>\varepsilon}
\frac{v(t)\, {\diff}  t}{t-x} \ , \   \ x \in \Bb{R} \, ,
\end{align*}

%\vspace{-0,2cm}
\noindent where it is known that $\widetilde{v}\!\in\! L^p(\Bb{R},{\diff} x)$
 by the M.\! Riesz theorem\! (see\! \cite[p. 128]{koo}).
Instead of applying the canonical factorization theorem
 (see \cite[p. 74, Thm. 5.5]{gar}, \cite[p. 119]{koo}) to the
 function $\he  $ in $ H^{p}  $,  we use the property
 \eqref{f25preresinvlamhyp} of  $\he  $ being in
$\mathcal{P}( \,(-\infty , 1)\, )$
and the corresponding representation \eqref{f4int} in the following
 improved form established in \cite[Theorems~2.6, 2.7, 2.8]{brs}.

\begin{thmx}\hspace{-0,18cm}{\bf{.}}\label{theor2}
   Let $v\!\in\!\bigcup_{p>1} L^p(\Bb{R},{\diff} x)$ be nonzero and
   satisfy $ v(x)=0$, if
  $x<1$, $v(x)\geq 0$,  if  $x\geq 1$,%\vspace{-0,2cm}
  \begin{align}\label{f1theor2}
    & {\rm{(a)}} \ \ \int_{1}^{+\infty} \frac{v(t)}{t}\, {\diff} t
    = \pi \ , \quad {\rm{(b)}} \ \ \begin{vmatrix}
\widetilde{v}(x_1) & \widetilde{v}(x_2) \\ v(x_1) & v(x_2)
\end{vmatrix}\geq 0
  \end{align}

  %\vspace{-0,1cm}
\noindent for almost all $1<x_1<x_2<+\infty$.  Then the function%\vspace{-0,1cm}
\begin{align}\label{f2theor2}
    &  \Psi (z) := \frac{1}{\pi} \int_{1}^{+\infty} \frac{v(t)\,
    {\diff} t}{t-z}\,,\quad z\in\Bb{H}\,,
\end{align}

%\vspace{-0,1cm}
\noindent belongs to $\mathcal{P}_{\log}$,  and there exists
  a non-decreasing right-continuous function $\nu$ on $\Bb{R}$ satisfying \
$ 0= \nu(y)\leq\nu(x)\leq 1$, $-\infty \leq y < 1 <x<+\infty$,
such that   for arbitrary $z\in\Bb{C}\setminus [1, +\infty)$
the following equalities hold%\vspace{-0,2cm}
\begin{align}\label{f3theor2}
    &\Log   \Psi (z)\! =\! \beta\! +\! \int\limits_{1}^{+\infty}
   \left( \frac{1}{t-z}\! -\! \frac{t}{1+t^2} \right)\,
   \nu(t)\, {\diff} t\!= \! \int_{[0,1]} \Log  \frac{1}{1-tz}\,
    {\diff} \mu(t) \, ,
\end{align}

%\vspace{-0,1cm}
\noindent where the measure $\mu\in\mathcal{M}^{+}(\Bb{R})$ and the real
 constant $\beta$ are defined by
\begin{align}\label{f4theor2}
    &  \beta := \int_{0}^{1} \dfrac{t \nu(1/t) \,   {\diff}  t }{1+t^2}  \ \, , \quad
    \left\{\begin{array}{rll}\mu \big((-\infty, 0)\big)  &  :=0 \ ,   &   \\[0,1cm]
    \mu\big(\,[0,\,x)\big)  &  := \nu(+\infty) - \nu(1/x) \ ,   &   x > 0 \ .
         \end{array}\right.
\end{align}
\end{thmx}

Denote by the same letter $\nu$ the Lebesgue-Stieltjes measure
 induced on $\Bb{R}$ by a non-decreasing  function
$\nu : \Bb{R} \mapsto [0,1]$ in Theorem~\ref{theor2}
(see \cite[p. 147, Def.~3.9]{mp}).
In addition to that theorem, we will need the following
relationships proved in \cite[Theorems~2.2, 2.6]{brs} and \cite[(2.37)]{brs}.

\begin{cormx}\hspace{-0,18cm}{\bf{.}}\label{corol1}
Under the conditions of Theorem~\ref{theor2},
\begin{align*}%\label{f0corol1}
    & {\mathrm{supp}}\,  \mu  \,\subseteq[0,1] \ , \quad \lim_{{{x \downarrow 0}}}
     \mu \big([0,x)\big) = 0 \ , \quad  \mu([0,1])=\nu(+\infty)\in (0,1],
\end{align*}

\noindent  and  for almost all $x\in\Bb{R}$ we have%\vspace{-0,2cm}
\begin{align}\label{f1corol1}
    & {\rm{(a)}} \ \  \Psi(x+\imag 0)
= \widetilde{v}(x) + \imag  v(x)  \ , \quad {\rm{(b)}} \ \
 \ln\dfrac{\sqrt{1+t^2}}{|t-x|}\in L^1\big(\Bb{R}, {{\diff}}\nu(t)\big) \ , \
   \\    &   \label{f2corol1}
\widetilde{v}(x)=\Big[\cos\pi\nu(x)\Big] \exp\left(  \int_{[1,+\infty) }
   \ln\frac{t}{|x-t|}\, {\diff} \nu(t) \right) \,, \\    &   \label{f3corol1}
    v(x) = \Big[\sin\pi\nu(x)\Big]\exp\left(  \int_{[1,+\infty) }
   \ln\frac{t}{|x-t|}\, {\diff} \nu(t)\right) \,.
\end{align}
  \end{cormx}

We observe that property \eqref{f1corol1}(a) follows directly
 from \eqref{f2theor2} for an arbitrary function
 $v\!\in\! L^p(\Bb{R}, {\diff} x)$,  that vanishes on the interval
  $(-\infty, 1)$ (for $1\!<\! p\! <\! \infty$), in view of
    known consequences of the M. Riesz theorem (see \cite[p. 128]{koo}).
     This fact allows us to deduce from
\eqref{f24preresinvlamhyp} that \eqref{f18preresinvlamhyp} yields%\vspace{-0,3cm}
\begin{align}\label{f1eir}\hspace{-0,16cm}
    &  v (x):= \dfrac{\he  (1-1/x) \chi_{[1, +\infty)} (x)}{\sqrt{x}}
    \,, \ x\in \Bb{R}   \ \Rightarrow \    \widetilde{v}(x)=
    \dfrac{\he  (1/x)}{\sqrt{x}}  \, , \  x > 1 \, .\hspace{-0,05cm}
\end{align}

\noindent
For such $v$ an $\widetilde{v}$ the equality \eqref{f2theor2} for
 $\Psi = \he  $  coincides with \eqref{f24preresinvlamhyp}
and the condition \eqref{f1theor2}(b) holds because it is equivalent to
the non-increasing property of the function $\he  (1-x)/\he  (x)$
 on the interval $(0,1)$ which follows readily
from the following consequence of \eqref{f21preresinvlamhyp},
\begin{align}\label{f2eir}
    & \hspace{-0,2cm} \dfrac{{{\diff}}}{{{\diff}} x} \dfrac{\he
    (1-x)}{\he  (x)} = -  \dfrac{1}{\pi x (1-x) \he  (x)^{2}} < 0 \ ,
     \quad  x \in (0,1) .
\end{align}

\noindent Furthermore,  the condition \eqref{f1theor2}(a) also holds
 in view  of the known
integral relationship (see \cite[p. 399, (4)]{erd2})%\vspace{-0,1cm}
\begin{align*}
    &  \int\nolimits_{1}^{+\infty} \frac{v(t)}{t}\, {\diff} t
    =\!\!\int\nolimits_{1}^{+\infty}t^{-3/2} \he  \left(1-{1}/{t}\right)\,
    {\diff} t =\!\!
    \int\nolimits_{0}^{1}(1-t)^{-1/2}\he  \left(t\right)\, {\diff} t = \pi .
\end{align*}

%\vspace{-0,1cm}
 \noindent Thus, for $\Psi = \he  $ and $v$, $\widetilde{v}$
 defined as in \eqref{f1eir} the conditions of Theorem~\ref{theor2}
 are true and we can   apply the results of Corollary~\ref{corol1}
 to calculate the function~$\nu$. Dividing \eqref{f3corol1} by
  \eqref{f2corol1} for $x > 1$, we obtain,
 by virtue of \eqref{f23preresinvlamhyp},
 \begin{align*}
    &  \tan \pi\nu(x) = \dfrac{v (x)}{\widetilde{v}(x)} =
     \dfrac{\he  \left(1-{1}/{x}\right)}{\he  \left({1}/{x}\right)}>
      0 \ , \quad  x > 1 \ .
 \end{align*}

\noindent Since $\nu(x) \in [0,1]$ we conclude that
\begin{align}\label{f3eir}
    &  \nu(x) =\dfrac{1}{\pi}\arctan \dfrac{\he
    \left(1-{1}/{x}\right)}{\he  \left({1}/{x}\right)} \ ,
    \quad  x > 1  \ , \quad
 \left\{ \begin{array}{l}
     \nu(1+0) = 0 \ , \\
       \nu(+\infty) = {1}/{2} \ ,
  \end{array}    \right.
\end{align}

\noindent and the differentiation of this equality, by taking
account of  \eqref{f2eir}, gives
\begin{align*}%\label{f4eir}
    & \nu^{\,\prime}(x) =
     \dfrac{1}{\pi^{2} (x-1)} \dfrac{1}{\he  \left({1}/{x}\right)^{2}
     + \he  \left(1-{1}/{x}\right)^{2}} > 0 \ , \quad  x > 1 \ ,
\end{align*}

\noindent and $\nu^{\,\prime} \in L^{1}([0,1],  {\diff}  x)$, in view
 of \eqref{f7preresinvlamhyp}(a). This means  that the formulas
 \eqref{f4theor2} (see also \eqref{f3eir}) for the measure $\mu$
 can be written as follows
\begin{align}\label{f5eir}
    &  \mu\big([0,\,x)\big)\!=\!\nu(+\infty)\! -\! \nu(1/x)\! =\!
    \dfrac{1}{2}\! - \!\dfrac{1}{\pi}\arctan \dfrac{\he
     \left(1\!-\!x\right)}{\he  \left(x\right)}  \, , \  x \in (0, 1);
      \\[-0,2cm]    &
    \mu\big(\,[0,\,x)\big)= \dfrac{1}{2}  \ , \quad   x \geq 1  \ ;
     \quad  \quad \mu\big(\,\{0\}\big) = 0  \ .
\label{f6eir}
\end{align}

\noindent This proves the validity of \eqref{f1th0} and shows that
$\mu$ is absolutely continuous with respect to the Lebesgue measure
$m$ on $[0,1]$ and it follows  from  \eqref{f2eir} and
 \eqref{f7preresinvlamhyp}(a) that the Radon-Nikodym derivative
  ${ {\diff}  \mu }/{ {\diff}  x} $ of $\mu$ with respect to $m$
  (see \cite[p. 214]{mp}) has the following form %\vspace{-0,1cm}
\begin{align}\label{f7eir}
    &\ \dfrac{ {\diff}  \mu(x)}{ {\diff}  x} =
    \dfrac{1}{\pi^{2 }  x (1-x)} \dfrac{1}{\he  \left(x\right)^{2}
    + \he  \left(1-x\right)^{2}}
         \ , \ \   x \in (0,1) ,
\end{align}

%\vspace{-0,1cm}
\noindent and ${ {\diff}  \mu }/{ {\diff}  x} \in  L^{1}([0,1], {{\diff}} x)$.
 Therefore for arbitrary $z\in\Bb{C}\setminus [1, +\infty)$ the
  exponential integral representations \eqref{f3theor2} we can
  write in the  forms \eqref{f8eir} and \eqref{f9eir} where
$ \beta :=\int_{0}^{1} ( {t\, \nu(1/t)}/{(1+t^2)}) {\diff}  t $.
It follows from \eqref{f6eir}, \eqref{f7eir} and \eqref{f8eir} that
\begin{align*}%\label{f11eir}
    &   \Arg  \he  (z) \in \left(- \dfrac{\pi}{2}, \dfrac{\pi}{2}\right)
    \ , \quad z\in\Bb{C}\setminus [1, +\infty) \ , \
\end{align*}

\noindent which shows that \eqref{f1ath0} holds, and  we have%\vspace{-0,2cm}
\begin{align}\label{f11aeir}
    & \Log   \he  (1\!-\!z)\! - \!  \Log   \he  (z) \!
    = \! \Log   \dfrac{
    \he  (1\!-\! z)}{\he  (z)} \ , \quad
     z\! \in \! (0,1)\cup \left(\Bb{C}\setminus\Bb{R}\right),
     \\[0,2cm]    &
    \Log   \he  (z)^{\alpha} =  \alpha \, \Log   \he  (z)  \ ,
    \quad 0 < \alpha \leq 2\,, \  z\in\Bb{C}\setminus [1, +\infty) .
\label{f11beir}
\end{align}

%\vspace{-0,1cm}
\noindent Then the representation \eqref{f8eir} written for
$ z\in \Bb{C}\setminus[1,+\infty)$ in the form%\vspace{-0,2cm}
\begin{align*}
    & \frac{z  \he  (z)^{\alpha}}{z}
= \exp\left( \int_{[0,1]} \Log  \frac{1}{1-tz}\  {\diff} \big(\alpha
\mu(t)\big)\right)
   \ , \quad  \alpha \mu([0,1]) = \alpha/2 \in (0, 1]  \ , \
\end{align*}

%\vspace{-0,2cm}
\noindent gives Corollary~\ref{th0} because, by virtue of
 \eqref{f6eir}, \eqref{f5int} and \eqref{f5aint} hold for
 $\Psi =z  \he  (z)^{\alpha} $ and $\alpha \mu$ in place of
 $\sigma$. In addition, the representations \eqref{f8eir}
 and \eqref{f9eir} for arbitrary $ z\! \in \! (0,1)\cup
 \left(\Bb{C}\setminus\Bb{R}\right)$ yield{\hyperlink{d15}{${}^{\ref*{case15}}$}}\hypertarget{bd15}{} %\vspace{-0,2cm}
\begin{align}\label{f12eir}
    &  \Log   \dfrac{ \he  (1\!-\! z)}{\he  (z)}\!=\!\dfrac{1}{
    \pi^{2 }}\! \int\limits_{0}^{1}\!\! \dfrac{
    \Log  \dfrac{1-tz}{1-t+tz}}{t (1-t)\left(\he  \left(t \right)^{2}
    + \he  \left(1-t \right)^{2}\right)} \,  {\diff}  t \ , \\[0,0cm]
      &\label{f13eir}
    \Log   \dfrac{ \he  (1\!-\! z)}{\he  (z)}=\dfrac{1-2z}{\pi}
    \int\limits_{0}^{1} \ \,  \dfrac{\arctan \dfrac{\he
    \left(1-t\right)}{\he  \left(t\right)}}{(1-tz)(1-t +t z)}
    \,  {\diff}  t \ .
\end{align}

%\vspace{-0,1cm}
\noindent To obtain  \eqref{f12eir}, it is sufficient to  use the
fact that for  arbitrary points   $ z\in(0,1)\cup \left(\Bb{C}
\setminus\Bb{R}\right)$
and  $t \in (0,1)$, the  two numbers $1-tz $ and $1-t + t z$ lie
in the open half-plane
     $\left\{ \   a + s \, z \ | \ s \in \Bb{R} \ , \ a > 0
     \ \right\}$
 which implies that%\vspace{-0,1cm}
\begin{align}\label{f14eir}
    &   \big| \Arg  (1\! -\! tz)\! -\!  \Arg  (1\! -\! t\!  +\!
      t z)\big| \! <\!  \pi  \, , \ \ \  t \in (0,1) \ , \  z\!
       \in \! (0,1)\cup \left(\Bb{C}\setminus\Bb{R}\right) \, ,
\end{align}

%\vspace{-0,1cm}
\noindent and therefore for arbitrary $t \in (0,1)$ and $z\! \in \! (0,1)\cup \left(\Bb{C}\setminus\Bb{R}\right) $ we have%\vspace{-0,1cm}
\begin{align}\label{f14aeir}
    &  \Log  \frac{1}{1\!-\!tz}\!-\!\Log  \frac{1}{1\!-\!t \!+\! tz}
     \!=\! \Log  \dfrac{1\!-\!tz}{1\!-\!t\!+\!tz}\, .
\end{align}

%\vspace{-0,1cm}
\noindent Together with \eqref{f6eir} and \eqref{f7eir}, the
inequality \eqref{f14eir} allows us to deduce from \eqref{f12eir}
 that %\vspace{-0,2cm}
\begin{align}\label{f15eir}
    &   \Arg  \dfrac{\he  (1\!-\! z)}{\he  (z)} \in \left(-
    \dfrac{\pi}{2} , \dfrac{\pi}{2}\right)  \ , \quad  z\!
    \in \! (0,1)\cup \left(\Bb{C}\setminus\Bb{R}\right) \ .
\end{align}

%\vspace{-0,1cm}
\noindent By  \eqref{f10int}, the formulas \eqref{f12eir} and
\eqref{f13eir} give the exponential representation of
$\ie/ \imag  $ for all
 $ z\! \in \! (0,1)\cup \left(\Bb{C}\setminus\Bb{R}\right)$.
 Besides that, for arbitrary $t \in (0,1)$ and
 $ z\! \in \! \Bb{C}\setminus\Bb{R}$
 the sign of $\Arg  (1\! -\! tz)\! -\!  \Arg  (1\! -\! t\!  +\!  t z)$
 is obviously equal to $-{\rm{sign}} ({\rm{Im}}\, z)$ and therefore it
 follows from \eqref{f12eir}, \eqref{f6eir} and \eqref{f10int} that
 \begin{align}\label{f16eir}
    &   \Arg  \ie (z) \in \dfrac{\pi}{2} - \left(0 ,
    \dfrac{\pi}{2}\right)\cdot {\rm{sign}} ({\rm{Im}}\, z)  \ , \  \
     z\! \in \! \Bb{C}\setminus\Bb{R}  \ .
\end{align}

%%%%%%%%%%%%%%%%%%%%%%%%%%%%%%%%%%%%%%%%%%%%%%%%%%%%%%%%%%%%%%%%%%%%

%\vspace{-0,4cm}
\section[\hspace{-0,30cm}. \hspace{0,11cm}Special properties of the triangle function]{\hspace{-0,095cm}{\bf{.}} Special properties of   $\ie$}\label{sp}

%\vspace{-0,2cm}
\noindent By  \eqref{f14preresinvlamhyp}, we have %\vspace{-0,3cm}
\begin{align}\label{f14zpreresinvlamhyp}
 \ie (-x \pm  \imag   0 ) = \pm 1 + \ie  \left(
 \dfrac{x}{1 + x}\right)=  \pm 1 +  \imag   \,
 \dfrac{\he   \left( {1}\big/{(1 + x)}\right)}{\he
 \left(  {x}\big/{(1 + x)}\right)}   \, , \ \  x > 0 \, .
\end{align}

\noindent
{\emph{Proof of \eqref{f27preresinvlamhyp},  \eqref{f28preresinvlamhyp}
 and \eqref{f29preresinvlamhyp}.}}
The properties \eqref{f25preresinvlamhyp} and definition \eqref{f10int}
 imply the validity of \eqref{f27preresinvlamhyp}.
 It follows from \eqref{f14zpreresinvlamhyp} that
 \eqref{f28preresinvlamhyp} holds.
Setting $z= -x \pm  \imag   0 $ and  $z= 1+ x \pm
\imag   0 $
in \eqref{f27preresinvlamhyp}  and combining these with
\eqref{f28preresinvlamhyp},  we obtain that \eqref{f29preresinvlamhyp}
 holds for arbitrary $ x>0$.{\hyperlink{d10}{${}^{\ref*{case10}}$}}\hypertarget{bd10}{} $\square$

\vspace{0,1cm}
 For the   function $\ie$ we introduce its  {\emph{remainders from
  singularities}}
 %\vspace{-0,1cm}
\begin{align}\nonumber
    \ie(\infty; z) &  :=\dfrac{1}{  \ie (z) -{\rm{sign}} ({\rm{Im}}\, z)}
    + \dfrac{ \imag  }{ \pi} \ln  |z|  \, ;   &  &
    \ie(0; z)  := - \ie (z)  + \dfrac{\imag  }{\pi}
     \ln\dfrac{1}{|z|} \, ; \\
    \ie(1; z) &:=     \dfrac{1}{\ie (z)}  +
     \dfrac{\imag }{\pi}\ln\dfrac{1}{|1-z|} \ ,   &    &
        z \in  (0,1)\cup \left(\Bb{C}\setminus\Bb{R}\right) \,.
\label{f29ainvlamhyp}\end{align}

%\vspace{-0,1cm}

  A direct calculation of ratios of asymptotic expansions
  \eqref{f7preresinvlamhyp} according to formula \eqref{f10int}
    implies the validity of the following assertion.{\hyperlink{d11}{${}^{\ref*{case11}}$}}\hypertarget{bd11}{}
 \begin{lemma}\hspace{-0,18cm}{\bf{.}}\label{lem3}
 Let $ \Lambda\!:=\! (0,1)\!\cup\! \left(\Bb{C}\!\setminus\!\Bb{R}\right)$.
 Then the following  asymptotic formulas hold:\hypertarget{tempbd9}{}
 \vspace{-0,1cm}
 \begin{align}
   \label{f30preresinvlamhyp}
      \vphantom{a}\hspace{-0,25cm}\ie(\infty; z) & \!=\!
      \dfrac{\ln 16 }{ \imag   \pi} +  \dfrac{\Arg  (1-z)}{\pi}
       + {\Oh} \left(\dfrac{\ln^{2} |z|}{|z|}\right)     \, ,   &
       \hspace{-0,5cm}          \Lambda \ni  z\to \infty  \, ,
      \\    \label{f31preresinvlamhyp}
  \vphantom{a}\hspace{-0,25cm} \ie(1; z)& \!=\!
  \dfrac{\ln 16 }{ \imag   \pi}\! +\!
   \dfrac{\Arg  \big(1/(1\!-\!z)\big)}{\pi}
  \!+ \! {\Oh} \left(|1\!-\!z| \ln \dfrac{1}{|1\!-\!z|}\right)
     \, ,    &   \hspace{-0,5cm}        \Lambda\! \ni \!   z\!\to \! 1  \, ,
                \\
    \label{f32preresinvlamhyp}
    \vphantom{a}\hspace{-0,25cm}\ie(0; z)& \!=\!   \dfrac{ \ln 16 }{
     \imag   \pi} + \dfrac{\Arg  (1/z)}{\pi} +
       {\Oh} \left(|z| \ln
     \dfrac{1}{|z|}\right) \ ,      &    \hspace{-0,5cm}
      \Lambda \ni z\to 0 \, .
   \end{align}

%\vspace{-0,1cm}

\noindent Here, $1-z, 1/(1\!-\!z), 1/z \in \Bb{C}\setminus(-\infty, 0]$
for any $z\in \Lambda$.
\end{lemma}

\vspace{0.15cm}
 The following crucial properties of the remainders from
 singularities hold. As for notation, let ${\rm{sign}} (x)$
 be equal to $-1$ if $x < 0$, $0$ if $x=0$ and $1$ if $x>0$.

\begin{lemma}\hspace{-0,18cm}{\bf{.}}\label{th2}
Let $\sigma (z):= {\rm{sign}} ({\rm{Im}}\, z)$, $z \in \Bb{C}\setminus \Bb{R} $.
Then
%\vspace{-0,1cm}
 \begin{align}\label{f3th2}
    &  \limsp_{{\fo{\Lambda\! \ni\! z\! \to 0}}}
     \big|{\rm{Re}} \,\ie (z)\big| \leq 1 \, , \
    \lim_{{\fo{\Lambda\! \ni\! z\! \to \!1}}} \ie (z)
    =\!\!\! \lim_{{\fo{\Lambda\! \ni\! z\! \to
    \!\infty}}}\big|\ie (z)\! -\!\sigma (z)\big|=0 \, ,
\end{align}

%\vspace{-0,2cm}
 \noindent and there exists a finite positive number
 $\varepsilon_{\!{\tn{\triangle}}}$ such that
 %\vspace{-0,1cm}
 \begin{align}\nonumber
    & {\rm{(a)}} \ \ \big|\ie(\infty; z)\big| \leq 2 \, ,
      \quad     |z| \geq 1/\varepsilon_{\!{\tn{\triangle}}} \ ;
  \quad   {\rm{(b)}} \ \ \big|\ie(0; z)\big|  \leq 2 \, ,   \quad
     |z| \leq \varepsilon_{\!{\tn{\triangle}}} \ ;  \\    &
    {\rm{(c)}} \ \ \big|\ie(1; z)\big|  \leq 2 \, ,  \quad \ \,
      |z-1| \leq \varepsilon_{\!{\tn{\triangle}}}\ ; \quad \ \ \ \
       z \in (0,1)\cup \left(\Bb{C}\setminus\Bb{R}\right)\,.
\label{f1th2} \end{align}
\end{lemma}

\vspace{0.05cm}
 \noindent {\emph{Remark.}} \, The known property $\lambda ( z\pm 1)=
 \lambda(z)/ (\lambda(z) -1)$, $z\in \Bb{H}$, of the elliptic modular
 function $\lambda$ (see \cite[p. 111]{cha}) which is immediate from
 \eqref{f1xtheor3} and \eqref{f9elemtheta} established below,
 together with the identity \eqref{f1theor3}  imply that{\hyperlink{d27}{${}^{\ref*{case27}}$}}\hypertarget{bd27}{}\vspace{-0,1cm}
 \begin{align}\label{f1rem1}
    &  \ie( z) =\ie \big(z/(z-1)\big) + \sigma (z) \ , \quad
    z \in \Bb{C}\setminus \Bb{R} \,.
 \end{align}

\vspace{-0,1cm}
 \noindent
 By letting here $\Lambda \ni z\to \infty$ and by using the asymptotic
 formula  \eqref{f31preresinvlamhyp}, we can improve the remainder
 in the asymptotic formula \eqref{f30preresinvlamhyp} to
 ${\Oh}(|z|^{-1}\ln |z|)$.

\vspace{0.1cm}
\noindent
{\emph{Proof of Lemma~~\ref{th2}.}}
The formulas  \eqref{f30preresinvlamhyp}, \eqref{f31preresinvlamhyp},
\eqref{f32preresinvlamhyp}, the estimate $(1/\pi)\ln 16  <    0,883$
and the fact that $\Arg  y \in (-\pi, \pi)$,
$y \in \Bb{C}\setminus(-\infty, 0]$,  together  imply that
 \eqref{f1th2} and  two limit identities in \eqref{f3th2} hold.
 Regarding the inequality \eqref{f3th2}, it is obtained by
 separation of the real and imaginary  parts in
  \eqref{f32preresinvlamhyp}. This completes the proof
  of Lemma~\ref{th2}.
$\square$

\vspace{0,1cm}
\noindent
{\emph{Proof of \eqref{f33preresinvlamhyp}.}}
By  \eqref{f14zpreresinvlamhyp}, we have for $x>0$
that $|{\rm{Re}} \ie (-x \pm  \imag   0 )|\! = \!1$ ,  while
after substitution of \eqref{f14zpreresinvlamhyp} in
\eqref{f27preresinvlamhyp} we obtain
$|{\rm{Re}} \ie (1\!+\!x\pm \imag   0 )|\! < \!1$. Together with
  the following obvious consequence of \eqref{f18preresinvlamhyp}
   and \eqref{f23preresinvlamhyp},\vspace{-0,1cm}
\begin{align}\label{f4th2}
    & \ie (x\!\pm\!  \imag   0)=\lim\nolimits_{{\fo{\, \Bb{H}
    \ni z \to 0}}} \  \ie (x\pm z)\,, \quad x\in \Bb{R}\setminus
     [0,1]\,,
\end{align}

\vspace{-0,1cm}
 \noindent
and \eqref{f3th2} this entails that \vspace{-0,1cm}
\begin{align*}
    & \limsp\nolimits_{{\fo{\Lambda\! \ni\! z\! \to a}}}
    \pm {\rm{Re}} \,\ie (z) \leq
     \limsp\nolimits_{{\fo{\Lambda\! \ni\! z\! \to a}}}
      \ \left|{\rm{Re}} \,\ie (z)\right| \leq 1   \ , \quad
       a \in  \left\{\infty\right\}\cup \partial\Lambda\ .
\end{align*}

\vspace{-0,2cm}
\noindent  By the maximum principle applied to the two harmonic
 functions\! $-\!1\pm {\rm{Re}} \ie (z)$  in  $\Lambda$,
we obtain that (see \cite[pp. 254, 129, 40]{con}, \cite[p. 47]{hay}){\hyperlink{d12}{${}^{\ref*{case12}}$}}\hypertarget{bd12}{}\vspace{-0,1cm}
\begin{align}\label{add}
    &  \big|{\rm{Re}}\, \ie (z)\big|< 1 \ , \quad  z\in \Lambda\,.
\end{align}

\vspace{-0,2cm}
\noindent
Next, if we assume that there exists a point $\zeta \in
\ie (\Lambda)\setminus \fet $
with ${\rm{Re}} \, \zeta \in (-1,1)$, we appeal to \eqref{f27preresinvlamhyp}
 and the reflection property $1 - \Lambda = \Lambda$
to obtain that $ -1/\zeta \in \ie (\Lambda)$. On the other hand,
 as $ \fet$ is invariant under inversion, the facts
that $\zeta \not\in \fet$ and $-1<{\rm{Re}}\, \zeta <1$ entail that
  $|{\rm{Re}} (-1/\zeta)|> 1$. This,
  however, is impossible because as we just observed
  $ -1/\zeta \in \ie (\Lambda)$, which would yield by \eqref{add}
  the opposite inequality.
So that the inclusion  $  \ie \left(\Lambda\right)\subset \fet$ is immediate.

 To obtain the reverse inclusion as well, we observe that the
set $\ie (\Lambda)$ is open, by  the open mapping
theorem (see \cite[p.99]{con}).
Suppose to the contrary that
$\fet \setminus  \ie (\Lambda)\neq \emptyset$.
Then one may find at
  least one point $\zeta$ in the intersection $\fet \cap \partial
 ( \ie\! \left(\Lambda\right))$, as otherwise
 for each
 point $\xi \in \fet \setminus  \ie (\Lambda)$ there would exist
 $\varepsilon=\varepsilon (\xi) > 0$ such that $(\xi +
 \varepsilon \Bb{D})\cap   \ie (\Lambda)=\emptyset $.  It would then follow that
 $\fet \setminus  \ie (\Lambda)$ is open and  that $\fet$
 may  be represented as the  union of the two non-empty disjoint open
  subsets $\ie (\Lambda)$ and $\fet \setminus  \ie\! (\Lambda)$,
  in contradiction with the connectivity of $\fet$
  (see \cite[p. 92, Thm. 1.6]{lan}).  For this point $\zeta \in \fet \cap \,\partial
 ( \ie\! \left(\Lambda\right))$  there must exist a sequence
  $\{y_{n}\}_{n\geq 1} \subset \Lambda$ such that
  $\zeta = \lim_{n\to \infty}  \ie (y_{n})$ and, by
   replacing this sequence by a suitable subsequence, we may assume
   that $y_{n}$ converges either to
   $\infty$ or to some   boundary  point $y_{\infty} \in \Bb{C} \setminus
   \Lambda = (-\infty, 0] \cup [1, +\infty)$ as $n\to\infty$.
   If the limit of $\{y_{n}\}_{n\geq 1}$ is infinite then
   by \eqref{f3th2}, we have that $\zeta \in \{1,-1\}\subset \partial \fet$.
  If the limit is instead finite and  is equal to $y_{\infty}\in (-\infty, 0)\cup (1, +\infty)$
    then we deduce from \eqref{f4th2} together with \eqref{f14zpreresinvlamhyp}
     and \eqref{f27preresinvlamhyp} that\vspace{-0,25cm}
     \begin{align*}
        & \zeta \in \{\ie (y_{\infty}\pm \imag  0)\}  \subset (\pm1+
    \imag   \Bb{R}_{>0})\cup
     (\pm1-  \imag   \Bb{R}_{>0})^{-1}
    \subset \partial \fet \ .
     \end{align*}

\vspace{-0,2cm}
  \noindent
   Finally, if $y_{\infty}\in
     \{0,1\}$ then according to \eqref{f29ainvlamhyp},
     \eqref{f31preresinvlamhyp} and \eqref{f32preresinvlamhyp}
      we have either $\zeta = \infty$ or $\zeta = 0 \in \partial \fet$.
      So that for all cases we obtain $\zeta \not\in \fet$
      which contradicts the assumption $\zeta\in\fet \cap
      \partial \ie (\Lambda )$. This proves $\fet \subset
      \ie \left(\Lambda\right)$ and completes the proof of
      \eqref{f33preresinvlamhyp}.
$\square$

\vspace{-0,15cm}
%\vspace{-0,4cm}

%\vspace{-0,4cm}

\section[\hspace{-0,30cm}. \hspace{0,11cm}Proof of  Lemma~\ref{th1}.]{\hspace{-0,095cm}{\bf{.}} Proof of  Lemma~\ref{th1}}
\label{mth}

\vspace{-0,15cm}
Introduce the function{\hyperlink{d13}{${}^{\ref*{case13}}$}}\hypertarget{bd13}{}      %\vspace{-0,1cm}
\begin{align*}
    &  \Phi (z) := f \big(\ie (z)\big)  \ , \quad
    z \in \Lambda:= (0,1)\cup \left(\Bb{C}\setminus\Bb{R}\right) \, .
\end{align*}

%\vspace{-0,1cm}
\noindent It follows from \eqref{f14zpreresinvlamhyp} and
\eqref{f27preresinvlamhyp} that  $\ie (1+ x +  \imag   0 )
\in 1 /(1 - \imag  \Bb{R}_{>0})$
for arbitrary $x>0$, and therefore, by  \eqref{f29preresinvlamhyp}
 and the second invariance property Lemma~{\hyperlink{hth1}{\ref*{th1}}}(b)
 of $f$, which can be written as
  $ f (z) = f ( {z}/{(1-2z)}) $ for all $ z \in 1 /(1 - \imag  \Bb{R}_{>0})$,
  we derive that
  $\Phi (1+ x -  \imag   0 ) =
  \Phi (1+ x +  \imag   0 )$, $x>0$.
   While the first invariance property Lemma~{\hyperlink{hth1}{\ref*{th1}}}(a)
    of $f$ together with \eqref{f14zpreresinvlamhyp} and
    \eqref{f28preresinvlamhyp} yields
$\Phi (-x- \imag   0 ) = \Phi (-x+ \mathrm{i}
 \hspace{0,015cm}   0 )$, $x>0$.
 By the Morera theorem  (see \cite[p. 96]{lav}) we get
$\Phi \in  {\rm{Hol}}( \Bb{C}\setminus\{0, 1\})$.
In the notation for the remainders from singularities,
 the properties \eqref{f1th2}  entail that for arbitrary
 $z\in \Lambda$  we have
\vspace{-0,1cm}
 \begin{align}\nonumber
    & {\rm{(a)}} \,  \dfrac{1}{ \left| \ie (z)\! -\!\sigma (z)\right|}\!
     \leq \!2\! + \! \dfrac{ \ln  |z| }{ \pi} \, ,
     \      |z|\! \geq\! 1/\varepsilon_{\!{\tn{\triangle}}} \, ;
  \    {\rm{(b)}} \,  \left|\ie (z)\right| \!
   \leq\! 2 \!+\! \dfrac{1 }{\pi} \ln\dfrac{1}{|z|}\, ,   \
    |z|\! \leq\! \varepsilon_{\!{\tn{\triangle}}} \, ;  \\[0,1cm]    &
    {\rm{(c)}} \ \ {1}/{\left|\ie (z)\right|}     \leq 2 +
    ({1}/{\pi})\ln\left({1}/{|1-z|}\right) \, ,  \quad \ \,
    |z-1| \leq \varepsilon_{\!{\tn{\triangle}}}\, .
\label{f1th1pr} \end{align}

%\vspace{-0,2cm}
\noindent Here $\sigma (z) := {\rm{sign}} ({\rm{Im}}\, z)$, and,
 in view of \eqref{f33preresinvlamhyp} and \eqref{f3th2},
 we obtain %\vspace{-0,2cm}
\begin{align}
  &   \begin{array}{lrll}
       {\rm{(a)}} \,  \ie (z)\! -\!\sigma (z)\to 0  \ ,   &
           \Lambda\ni z \to \infty\,,  &
   \ \  {\rm{(b)}} \ \       \ie (z) \to \infty  \ ,   &
       \Lambda\ni z \to 0\,,\
    \\[0,1cm]
    {\rm{(c)}} \ \ \ie (z) \to 0  \ ,   &
       \Lambda \ni z \to 1  \ ,  &
     \ \    {\rm{(d)}} \ \   \ie (z)\in \fet  \, ,     &
        z \in \Lambda \ .
      \end{array}
\label{f2th1pr}\end{align}

%\vspace{-0,2cm}\noindent
By substituting $\ie (z)$ in place of $z$ in
 Lemma~{\hyperlink{hth1}{\ref*{th1}}}, (1)-(3),
 and letting $z\to 0$ in (1), $z\to 1$ in (2), and
 $z\to \infty$ in (3), and in addition by
using \eqref{f2th1pr} and the inequalities \eqref{f1th1pr},
we find that %\vspace{-0,1cm}
\begin{align}\nonumber
    &  {\rm{(1)}}\, |z|^{{\fo{n_{\infty} +1}}}
    \left|\Phi (z)\right|\! \to 0\!\, ,   &
     \Lambda\! \ni\! z \!\to\! 0 \, ,    &     &
    {\rm{(2)}}\, |1\!-\!z|^{{\fo{n_{0}\!+\!1}}}\left|\Phi (z)\right|\!\to\! 0\, ,
     \  \Lambda \!\ni \!z \!\to\! 1\, ,  \\   &
    {\rm{(3)}}\, |z|^{ {\fo{- n_{ 1} -1}}}\left|\Phi (z)\right| \to 0\ ,
      &     \Lambda \ni z \to  \infty \,  .  &     &
\label{f3th1pr}\end{align}

%\vspace{-0,2cm}
\noindent From the Riemann theorem about removable singularities
(see \cite[p. 103]{con}) it follows that the function
$\Phi_{1} (z):= z^{n_{\infty}} (1-z)^{n_{0}} \Phi (z)$ is
holomorphic at the points $0$ and $1$ and is hence an entire
 function
which, by the property  \eqref{f3th1pr}(3) above, has the
asymptotics
%\vspace{-0,4cm}
\begin{align}\label{f4th1pr}
    &\Phi_{1} (z) = {\oh} \left(|z|^{n_{0} +n_{\infty} +
     n_{ 1} +1}\right) \quad \mbox{as} \quad
     \Lambda \ni z \to \infty\,.
\end{align}

%\vspace{-0,2cm}
\noindent
The continuity of $\Phi_{1}$ ensures  the existence of
$C \in \Bb{R}_{>0}$ such that $|\Phi_{1} (z)| \leq C (1+|z|)^{n_{0}
 +n_{\infty} +  n_{ 1} +1}$, $z\in \Bb{C}$,  which by the
  extended version of the Liouville theorem
   (see \cite[p. 2, Thm. 1]{lev2}) yields that $\Phi_{1} (z)$
    is an algebraic polynomial of degree at most
    $n_{\infty} + n_{0} + n_{1}+1$. But the relationship
    \eqref{f4th1pr} proves that actually its degree cannot
    exceed $n_{\infty} + n_{0} + n_{1}$.
    Lemma~{\hyperlink{hth1}{\ref*{th1}}} follows.

%%%%%%%%%%%%%%%%%%%%%%%%%%%%%%%%%%%%%%%%%%%%%%%%%%%%%%%%%%%%%%%%%%%%
%\vspace{-0,2cm}
\section[\hspace{-0,30cm}. \hspace{0,11cm}Definitions of the theta functions]{\hspace{-0,095cm}{\bf{.}} Definitions of the theta functions}
\label{theta}

%\vspace{-0,2cm}
Introduce the functions
\begin{align}\label{f1elemtheta}
    &
    \begin{array}{ll}
       \theta_{3} (u) \! := \!
1  \!+ \! 2\sum\nolimits_{n\geq 1} u^{ n^2  }      \, ,
    &  \quad   \theta_{2} (u) \!:=
\!1 \!+\! \sum\nolimits_{n\geq 1} u^{n^2 +n} \, ,    \\[0,2cm]
\theta_{4} (u) \! :=  \!1  \!+ \! 2\sum\nolimits_{n\geq 1}
(-1)^{n}u^{ n^2  } \, ,     & \quad u \in \Bb{D} ,
    \end{array}
\end{align}

\noindent which  are obviously holomorphic in the unit disk
$\Bb{D}$ and satisfy
\begin{align}\label{f2elemtheta}
    & \theta_{4} (u) =   \theta_{3} (-u)  \ , \    &     &
     \theta_{2} (-u) = \theta_{2} (u) \ , \    &     &
      u \in \Bb{D} \ .
\end{align}

In order to  apply the Poisson summation formula to the series
in \eqref{f1elemtheta} we assume there $u \in (0,1)$ and replace
$u$ by $\exp (- \pi x)$ with $x > 0$. Then the well-known
integrals \cite[p. 146, (27)]{erd3},
\begin{align*}
    & \sqrt{\dfrac{\pi}{x}} {\rm{e}}^{{\fo{-2a\sqrt{ x}}}}\! =\!
     \int\limits_{0}^{\infty} {\rm{e}}^{{\fo{-x t}}}
     \dfrac{{\rm{e}}^{{\fo{- a^{2}/t}}}}{\sqrt{t}}  {\diff}  t \ , \ \
 \sqrt{\dfrac{\pi}{x}} {\rm{e}}^{{\fo{-a\sqrt{ x}}}}
 \! =\!  \int\limits_{0}^{\infty} {\rm{e}}^{{\fo{-x t}}}
  \dfrac{{\rm{e}}^{{\fo{- a^{2}/(4 t)}}}}{\sqrt{t}}
   {\diff}  t\ , \ \  a, x > 0,
\end{align*}

\noindent  allow us to derive from the known formulas
\begin{align*}
    &  \coth \sqrt{\pi x} =1 + 2 \sum\nolimits_{n \geq 1}
    {\rm{e}}^{{\fo{-2 n \sqrt{\pi x}}}}  \ ,   &     &
\dfrac{1}{\sinh \sqrt{\pi x}} =   2 \sum\nolimits_{n \geq 0}
  {\rm{e}}^{{\fo{-(2 n +1)\sqrt{\pi x}}}} \ , \
\end{align*}

\noindent that for arbitrary $x> 0$ we have
\begin{align}\label{f4elemtheta}
    &  \dfrac{\sqrt{\pi}\coth \sqrt{\pi x}}{\sqrt{x}} =
    \int\nolimits_{0}^{\infty}{\rm{e}}^{{\fo{-x t}}}\
    \dfrac{ 1}{\sqrt{t}}\left(1 + 2
    \sum\nolimits_{n \geq 1}{\rm{e}}^{{\fo{- \pi n^{2}/t}}}
     \right)  {\diff}  t \ , \    \\    &
     \dfrac{\sqrt{\pi}}{\sqrt{x}\sinh \sqrt{ \pi x}} =
      \int\nolimits_{0}^{\infty}{\rm{e}}^{{\fo{- x t}}}\
      \dfrac{2 {\rm{e}}^{{\fo{-\pi/(4 t)}}} }{\sqrt{t}}
      \left(\sum\nolimits_{n \geq 0}{\rm{e}}^{{\fo{-\pi(
      n^{2} +n)/ t}}}\right)  {\diff}  t \ .
\label{f5elemtheta}\end{align}

\noindent The known expansions into the series of the simple
fractions \cite[p. 126]{olv}
\begin{align*}
    &   \sqrt{\pi}\, \dfrac{\coth  \sqrt{\pi x}}{\sqrt{x}} =
      \dfrac{1}{x} + 2 \sum\limits_{n \geq 1} \dfrac{1}{x
      + \pi n^{2}}  \ ,     &     &
\dfrac{\sqrt{\pi}}{\sqrt{ x}\sinh \sqrt{\pi x}} =  \dfrac{1}{x}
+ 2 \sum\limits_{n \geq 1} \dfrac{(-1)^{n}}{x +\pi n^{2}} \ , \
\end{align*}

\noindent for arbitrary $x> 0$ yield readily that
\begin{align*}%\label{f4aelemtheta}
    &  \dfrac{\sqrt{\pi}\coth \sqrt{\pi x}}{\sqrt{x}} =
    \int\nolimits_{0}^{\infty} {\rm{e}}^{{\fo{-x t}}}
    \left(1 + 2 \sum\nolimits_{n \geq 1}{\rm{e}}^{{\fo{-n^{2} \pi  t}}}
     \right)  {\diff}  t \ , \     \\    &
     \dfrac{\sqrt{\pi}}{\sqrt{x}\sinh \sqrt{ \pi x}} =
    \int\nolimits_{0}^{\infty} {\rm{e}}^{{\fo{-x t}}}  \left(1
     + 2 \sum\nolimits_{n \geq 1}(-1)^{n}{\rm{e}}^{{\fo{-n^{2}\pi  t}}}
      \right)  {\diff}  t \ .
%\label{f5aelemtheta}
\end{align*}

\noindent Comparing these equalities with \eqref{f4elemtheta}
 and \eqref{f5elemtheta}  we conclude, by the uniqueness
  theorem for the Laplace transform (see \cite[p. 63, Thm. 6.3]{wid}),
  that for any $t > 0$ we have
 \begin{align}\label{f6elemtheta}
    & \theta_{3} \left({\rm{e}}^{{\fo{-\pi t}}}\right) \!=\!
    t^{-1/2}  \theta_{3} \left({\rm{e}}^{{\fo{-\pi/ t}}}\right)
     \, , \
\theta_{4} \left({\rm{e}}^{{\fo{-\pi t}}}\right) \!= \!
2 \, t^{-1/2} {\rm{e}}^{{\fo{-\pi/(4 t)}}}  \theta_{2}
\left({\rm{e}}^{{\fo{-\pi/ t}}}\right)  \, ,
\end{align}

\noindent and the change of $t$ by $1/t$ in the latter equality
 gives (cf. Exercise 20 in \cite[p. 23]{law})%\vspace{-0,3cm}
\begin{align}\label{f7elemtheta}
    &  2 {\rm{e}}^{{\fo{-\pi t/4 }}}  \theta_{2}
    \left({\rm{e}}^{{\fo{-\pi t}}}\right)  = t^{-1/2} \theta_{4}
\left({\rm{e}}^{{\fo{-\pi / t}}}\right)  \ , \   t > 0 \ .
\end{align}

%\vspace{-0,1cm}
 To have a more simple form of writing the relationships
 \eqref{f2elemtheta}, \eqref{f6elemtheta} and \eqref{f7elemtheta}
  between $\theta_{k} (u)$,
 $2 \leq k \leq 4$,  the following  analytic  functions in
 $\Bb{H}$ are introduced
\begin{align}\label{f8elemtheta}
    &  \Theta_{3}  (z):=   \theta_{3}\left({\rm{e}}^{\imag
     \pi z}\right)\, , &   & \Theta_{4}  (z):=
     \theta_{4}\left({\rm{e}}^{\imag  \pi z}\right)\, , & &
 \Theta_{2}  (z) := 2 {\rm{e}}^{\imag  \pi z /4}
 \theta_{2}\left({\rm{e}}^{\imag  \pi z}\right)\, ,
\end{align}

\noindent where $z \in \Bb{H}$. Regarding these functions,
 the main  relationships can be written for arbitrary
 $z \in \Bb{H}$ as follows, by using the principal branch
  of the square root   (see \cite[p. 531, 20.7.27-29,31-33]{olv}),
\begin{align}\label{f9elemtheta}
    &  \hspace{-0.45cm} \begin{array}{llll}
	{\rm{(a)}}\	\Theta_2(-1/z) &= (z/\imag )^{1/2}\Theta_4 (z)\, ,
 \quad   &\quad {\rm{(b)}}\
	\Theta_3(-1/z) &= (z/\imag )^{1/2}\Theta_3(z)\, ,\ \\[0,1cm]
  {\rm{(c)}} \
	\Theta_4(-1/z) &= (z/\imag )^{1/2}\Theta_2(z) \, , &\quad
{\rm{(d)}} \
	\Theta_2(z+1) &= {\rm{e}}^{\imag \pi/4}\Theta_2(z) \, ,
\\[0,1cm]   {\rm{(e)}} \
	\Theta_3(z+1) &= \Theta_4(z) \, , &\quad {\rm{(g)}} \
	\Theta_4(z+1)&= \Theta_3(z) \, ,
\end{array}\hspace{-0.2cm}
\end{align}

\noindent where \eqref{f9elemtheta}(d),(e),(g) follow readily
from \eqref{f2elemtheta}, while \eqref{f9elemtheta}(a),(b),(c)
follow from \eqref{f6elemtheta} and  \eqref{f7elemtheta}
 because according to these relations the three functions %\vspace{-0,2cm}
\begin{align*}
    & \Theta_2\left(\!-\dfrac{1}{z}\right)\! -\!
    \Theta_4 (z)\left(\dfrac{z}{\imag }\right)^{{\fo{\dfrac{1}{2}}}}, \
     \Theta_3\left(\!-\dfrac{1}{z}\right)\! -\!
     \Theta_3(z)\left(\dfrac{z}{\imag }\right)^{{\fo{\dfrac{1}{2}}}}  , \
\Theta_4\left(\!-\dfrac{1}{z}\right) \! -\!
\Theta_2(z)\left(\dfrac{z}{\imag }\right)^{{\fo{\dfrac{1}{2}}}} ,
\end{align*}

\noindent are all holomorphic on $\Bb{H}$ and
vanish  for all  $ z \in \imag \, \Bb{R}_{>0}$, so
that by the uniqueness theorem  for analytic functions (see
 \cite[p. 78, Thm. 3.7(c)]{con}),    they all
 vanish identically on $\Bb{H}$.{\hyperlink{d16}{${}^{\ref*{case16}}$}}\hypertarget{bd16}{}

\vspace{0.25cm}
 \noindent {\emph{Remark.}} \,
In the notations  $\theta_{k} (z | \tau ) =
\theta_{k} (z, q)$, $2 \leq k \leq 4$,  $z \in \Bb{C}$,
 $ \tau \in \Bb{H}$, $q = {\rm{e}}^{\imag \pi \tau} \in \Bb{D}$
 of \cite[p. 524, 20.2.2-4]{olv}, we have $\theta_{3} (0, u)\!
 =\!\theta_{3} (u) $, $\theta_{4} (0, u)\! =\!\theta_{4} (u)$,
$u \in \Bb{D}$, and%\vspace{-0,2cm}
\begin{align*}%\label{f10elemtheta}
    &  \Theta_{k} (\tau)= \theta_{k} \left(0,  {\rm{e}}^{\imag
     \pi \tau}\right)\ ,  \quad 2 \leq k \leq 4 \ , \quad
       \tau \in \Bb{H}  \ .
\end{align*}

%%%%%%%%%%%%%%%%%%%%%%%%%%%%%%%%%%%%%%%%%%%%%%%%%%%%%%%%%%%%%%%%%%%%

%\vspace{-0,5cm}
\section[\hspace{-0,30cm}. \hspace{0,11cm}Wirtinger's identity]{\hspace{-0,095cm}{\bf{.}} Wirtinger's identity}
\label{wid}

%\vspace{-0,1cm}
{\emph{Proof of the identity \eqref{f1theor4}.}}{\hyperlink{d17}{${}^{\ref*{case17}}$}}\hypertarget{bd17}{}
By \eqref{f25preresinvlamhyp}, the function $\he  (z)$
does not vanish on  $\Bb{C}\setminus [1, +\infty)$ and
therefore we can introduce the function %\vspace{-0,1cm}
\begin{align}\label{f1wid}
    &   \Phi (z) = \dfrac{\Theta_{3}\left(\ie (z)\right)^{2}}{\he  (z)}\ ,
    \quad z \in \Lambda:= (0,1)\cup \left(\Bb{C}\setminus\Bb{R}\right) \,,
   \ \ \Phi \in  {\rm{Hol}}( \Lambda)\,  .
\end{align}

%\vspace{-0,1cm}
\noindent The formulas \eqref{f9elemtheta}(b) and
\eqref{f27preresinvlamhyp} for any $z \in \Lambda$ yield
 that %\vspace{-0,1cm}
\begin{align}\label{f2wid}
    &  \Phi (z) = \dfrac{\imag }{\ie (z)} \dfrac{\Theta_{3}
    \left(-1/\ie (z)\right)^{2}}{\he  (z)}
    =      \dfrac{\Theta_{3}\left(\ie (1-z)\right)^{2}}{\he  (1-z)}
     =  \Phi (1-z)  \ .
\end{align}

%\vspace{-0,1cm}
\noindent By using   \eqref{f9elemtheta}(e),(g) and
\eqref{f28preresinvlamhyp}, for arbitrary  $x > 0$ we deduce
that%\vspace{-0,1cm}
\begin{align}\label{f3wid}
    &  \Phi (-x\! + \! \imag   0 )\! = \!
     \dfrac{\Theta_{3}\left(2\! +\! \ie (-x\! -\!  \imag
      0 )\right)^{2}}{\he  (-x)}\! =  \! \dfrac{\Theta_{3}
      \left( \ie (-x\! - \! \imag   0 )\right)^{2}}{\he  (-x)}
       \!  =\!
      \Phi (-x\! - \!  \imag   0 ) \ .
\end{align}

%\vspace{-0,1cm}
\noindent For arbitrary $z \in \Bb{H}$ and $x > 0$ it follows
from  \eqref{f9elemtheta}(b),(e),(g) and \eqref{f29preresinvlamhyp}
 that%\vspace{-0,1cm}
\begin{align*}
    &  \Theta_{3}\left(\dfrac{z}{1-2z}\right)^{2}\!\! =\!
     (1\!-\!2z) \Theta_{3} (z)^{2}\!, \quad
    1\!-\!2 \ie (1 \!+\! x \!+\!  \imag   0 )\! =\!
    \dfrac{\he  (1\! + \!x \!- \!  \imag   0 )}{\he  (1\!
    +\! x\! +\!  \imag   0 )}  \, ,
\end{align*}

%\vspace{-0,1cm}
\noindent from which and \eqref{f29preresinvlamhyp}
for arbitrary  $x > 0$  we derive
\begin{align}  &
     \Phi (1\! +\! x\! -\!  \imag   0 )\!=\!
    \dfrac{\big(1\!-\!2 \ie (1\! +\! x\! +\!  \imag   0 )\big)
    \Theta_{3} \big( \ie (1\! +\! x\! + \! \imag   0 )\big)^{2}}{\he  (1
    \!+\! x\! -\!  \imag   0 )}=
     \Phi (1\! +\! x \!+ \! \imag
     0 ) \, .
\label{f4wid}\end{align}

%\vspace{-0,1cm}
\noindent By applying the Morera theorem to the properties
\eqref{f3wid} and \eqref{f4wid}  we obtain that
$\Phi \in  {\rm{Hol}}( \Bb{C}\setminus\{0, 1\})$
(see \cite[p. 96]{lav}).

When $\Lambda\ni z\to 0$, by  \eqref{f2th1pr}(b), \eqref{f1th2}(b)
and \eqref{f8elemtheta}, \eqref{f1elemtheta}, \eqref{f11int},
we have  $\he  (z) \to 1$,
$\ie (z)= (\imag /\pi) \ln(1/|z|) - \ie(0; z)$, $|\ie(0; z)| \leq 2,$
 and $\Theta_{3}\left(\ie (z)\right)^{2}= \theta_{3}\left(\exp\left(
 - \ln(1/|z|) \!- \!  \imag \pi \ie(0; z)\right)\right)^{2}$ tends to $1$,
 correspondingly.  Next, by the symmetry property  \eqref{f2wid},
 we obtain the existence of the two limits%\vspace{-0,1cm}
\begin{align}\label{f5wid}
    &  \lim\nolimits_{\Lambda\ni z\to 0} \,  \Phi (z) =
    \lim\nolimits_{\Lambda\ni z\to 1}  \Phi (z) = 1 \,.
\end{align}

%\vspace{-0,2cm}
\noindent By the Riemann theorem about removable singularities
 (see \cite[p. 103]{con}) we deduce that $\Phi$ is an entire
  function satisfying
$ \Phi (0) =  \Phi (1) = 1$.

 But if $z\! \to\! \infty$ lying in the one of the half-planes
  $\sigma\! :=\! {\rm{sign}} ({\rm{Im}}\, z)\! \in \!\{1, -1\}$  then
   by \eqref{f9elemtheta}(c),(e) and \eqref{f2th1pr}(a) we have
   $\Theta_{3}(z)^{2}  = \imag ( z \! - \!\sigma )^{-1}
    \Theta_{2} \left(  \!- \! 1/(z \! - \!\sigma)\right)^{2}$
    and $\ie (z)\! \to\! \sigma$, respectively. In the
    notation of \eqref{f29ainvlamhyp} and in view of
    \eqref{f1th2}(a), by denoting $\delta (z) := \imag
    \pi \ie(\infty; z)$,   we deduce from  \eqref{f8elemtheta}
     and \eqref{f7preresinvlamhyp}(c) that%\vspace{-0,1cm}
\begin{align*}
    &   \Phi (z) =
    \dfrac{\imag  \, \theta_{2}
     \Big(\exp\big({{\nor{ - {\imag  \pi}/{(\ie (z) -
     \sigma)\,}}}}\big)\Big)^{2}\exp\big({{\nor{- {\imag  (\pi /2)}/{
    (\ie (z) -\sigma)}   }}}\big) }{(\ie (z) -\sigma)
     \he  (z)}  \\ & =
    \dfrac{ \delta (z) + \ln  |z|}{\dfrac{\Log   16 (1-z)}{4
    \sqrt{1-z}} + {\Oh} \left(\dfrac{\ln |z|}{|z|^{3/2}}\right)}\
    \dfrac{{\rm{e}}^{{\fo{ - \delta (z)/2  }}}}{\sqrt{ |z| }}\
     \theta_{2}
  \left( {\rm{e }}^{{\fo{ - \delta (z) - \ln  |z|  }}}\right)^{2}
  = {\Oh} (1) \,,
\end{align*}

\noindent as $\Lambda\ni z \to \infty$. As a consequence,
the entire function $\Phi$ is bounded and  by the Liouville
 theorem \cite[p. 77]{con} it is a constant, which must
 equal $1$, by  \eqref{f5wid}.
This establishes the Wirtinger identity \eqref{f1theor4}.
$\square$

{\emph{Proof of \eqref{f0theor3}.}} In view of
\eqref{f33preresinvlamhyp} and \eqref{f25preresinvlamhyp},
 \eqref{f1theor4} yields that $\Theta_{3}(z)\neq 0$ for all
  $z\in\fet$, and  $\Theta_{3}(z)\neq 0$ for all $z \in \pm
  1 + \imag  \, \Bb{R}_{>0}$, by  \eqref{f14zpreresinvlamhyp}{\hyperlink{d18}{${}^{\ref*{case18}}$}}\hypertarget{bd18}{}.
The relations \eqref{f6beirf},  \eqref{f6ybeirf} and the
 equality \eqref{f6eirf} (see also \cite[p. 32, Thm.~7.1]{mum}),
 which is the result of successive applications of the
 transformations $z\mapsto z+2$ and $z \mapsto -1/z$ in
  \eqref{f9elemtheta}(b),(e),(g)
 (see \cite[p. 112, Lem. 2]{cha}), prove that
  $\Theta_{3}(z)\neq 0$ for all $z\in\Bb{H}$. Then
 the following consequences of \eqref{f9elemtheta},
 $\Theta_{4}(z)\!=\! \Theta_{3}(z\!+\!1)$, $\Theta_{2}(z)=
  (\imag /z)^{1/2}\Theta_{3}(1\!-\!1/z)$, for $z \! \in \! \Bb{H}$,
  complete the proof of \eqref{f0theor3}. $\square$

%\vspace{-0.3cm}
\section[\hspace{-0,30cm}. \hspace{0,11cm}Identities for the elliptic modular function]{\hspace{-0,095cm}{\bf{.}} Identities for the elliptic modular function}
\label{emf}

%\vspace{-0.1cm}
{\emph{Proof of \eqref{f8theor3} and \eqref{f1theor3}.}}
We prove that two holomorphic  functions in $\Bb{H}$ %\vspace{-0,2cm}
\begin{align*}
    & f_{1}(z) := ({\Theta_{2}(z)^{4} +
    \Theta_{4}(z)^{4}})/{\Theta_{3}(z)^{4}}  \ , \
     f_{2} (z) :=
      {\Theta_{2}\left(z\right)^{4}}/{\Theta_{3}\left(z\right)^{4}}
      \ , \  z\in \Bb{H} \ , \
\end{align*}

 %\vspace{-0,2cm}
\noindent  satisfies the conditions of Lemma~{\hyperlink{hth1}{\ref*{th1}}}
 with $n_{\infty} = n_{0}= n_{ 1} = 0$  and
$n_{0} = n_{\infty}=0$,  $n_{1}= 1$,
respectively{\hyperlink{d19}{${}^{\ref*{case19}}$}}\hypertarget{bd19}{}.   It follows from \eqref{f9elemtheta} that for arbitrary
$z\in \Bb{H}$
and $2\! \leq \!k \!\leq\! 4$ we have:
$\Theta_{k}\left(z\right)^{4} \! =\! \Theta_{k}\left(z+2\right)^{4}$,
$\Theta_{k}\left({z}/{(1\!-\!2z)}\right)^{4} \! =\!(1\!-\!2z)^{2}
\Theta_{k} (z)^{4}$
and $   \Theta_{k}\left(z\right)^{4}\! =\! -\Theta_{6-k}
\left(-1/z\right)^{4}\!/z^{2} $, while  if $\sigma \in \{1, -1\}$
then
\begin{align*}
    & \Theta_{k}\left(z\right)^{4}\!
    =\!\dfrac{(-1)^{m+k}}{(z\!-\!\sigma)^{2}}\Theta_{m}
    \left(-\dfrac{1}{z\!-\!\sigma}
    %-1/(z\!-\!1)
    \right)^{4}
  %  \!/(z\!-\!1)^{2}
    \, ,  \
    \begin{pmatrix}k \\ m\\ \end{pmatrix} \!\in \!
    \left\{ \begin{pmatrix}2 \\4 \\ \end{pmatrix}, \,
     \begin{pmatrix}3 \\ 2\\ \end{pmatrix} , \,
      \begin{pmatrix}4 \\ 3\\ \end{pmatrix}   \right\},
         %   (k,m)\!\in \!\left\{  (2,4)\,, \, (4,3)\,, \, (3,2) \right\},
       \end{align*}

 %\vspace{-0,1cm}
\noindent from which  for any $z\in \Bb{H}$ and
$\sigma \in \{1, -1\}$ we get %\vspace{-0,2cm}
\begin{align}
    & \nonumber
     {\rm{(a)}} \ \  f_{k} (z+2) = f_{k} (z)  \, ,
      &    &        {\rm{(b)}} \ \  f_{k} \left({z}/{(1-2z)}
      \right)= f_{k} (z)\ ,
     \    \  k \in \{1, 2\} \ ; \\   &   \nonumber
     {\rm{(c)}} \ \   f_{1}(z) = f_{1}(-1/z)  \ , \
       &    &
     {\rm{(d)}} \ \    f_{2} (z) =
     {\Theta_4 (-1/z)^{4}}\big/{\Theta_3 (-1/z)^{4}} \ , \
       \\   &
     {\rm{(e)}} \ \   f_{1}(z) = \dfrac{\Theta_{3}(y)^{4}
      - \Theta_{4}(y)^{4}}{ \Theta_{2}(y)^{4}}   \ , \
      &    &  {\rm{(f)}} \ \ f_{2} (z)
      =\!-\! \dfrac{\Theta_{4}(y)^{4}}{ \Theta_{2}(y)^{4}} \ , \
       y := \!-\!\dfrac{1}{z\!-\!\sigma} \ .
\label{f1emf}\end{align}

%\vspace{-0,1cm}
\noindent The two conditions of invariance (a) and (b)
in Lemma~{\hyperlink{hth1}{\ref*{th1}}} hold for $f_{1}$
and $f_{2}$, in view of \eqref{f1emf}(a),(b) with $-1/z$
in place of $z$.
 It follows from  \eqref{f8elemtheta} that
\begin{align}\label{f2emf}
    &
    \begin{array}{ll}
 {\rm{(a)}} \ \    \Theta_{3}  (z)^{4}\!\!=\! 1\!+\! 8\,
 {\rm{e}}^{\imag  \pi  z}\! + \!
  {\Oh} \left({\rm{e}}^{2 \imag  \pi
  z}\right),   &
  {\rm{(b)}} \ \   \Theta_{4}  (z)^{4}\!\!= \!1\!-\! 8\,
  {\rm{e}}^{\imag  \pi  z}\! +\!
   {\Oh} \left({\rm{e}}^{2 \imag
  \pi  z}\right) , \\[0,075cm]
  {\rm{(c)}} \ \     \Theta_{2}  (z)^{4}\!\!=\!  16\,
   {\rm{e}}^{\imag
  \pi  z} \!  + \!  {\Oh} \left({\rm{e}}^{3 \mathrm{i}
   \hspace{0,015cm}  \pi  z}\right),   &
\hphantom{  {\rm{(c)}} \ \ }  \fet \ni z\to \infty \ ,
    \end{array}
\end{align}

%\vspace{-0,1cm}
\noindent which together with \eqref{f1emf}(c),(d) show that %\vspace{-0,1cm}
\begin{align}\label{f4emf}
    &  \lim\nolimits_{{\fo{\fet\! \ni \!z\!\to\! 0}}} \
     f_{k} (z) = 1 , \
     \lim\nolimits_{{\fo{\fet \!\ni\! z\!\to\! \infty}}}\
      f_{k} (z) = 2-k\ , \ k \in \{1,2\} \ ,
\end{align}

%\vspace{-0,1cm}
\noindent
 and hence  the conditions  (1) and (2) in
 Lemma~{\hyperlink{hth1}{\ref*{th1}}} with $n_{0} =
  n_{\infty} = 0$ hold for $f_{1}$ and $f_{2}$.
The condition (3) in Lemma~{\hyperlink{hth1}{\ref*{th1}}}
 with $n_{1} = 0$ also holds for $f_{1}$ as follows from
 \eqref{f1emf}(a),(e) and \eqref{f2emf}.
Applying the result of Lemma~{\hyperlink{hth1}{\ref*{th1}}}
to $f_{1}$
we obtain  that $f_{1} (\ie(z))= a$ holds on $\Lambda$ for
some constant $a \in \Bb{C}$.
By letting $ \Lambda\ni z \to 0$ we obtain from
 \eqref{f2th1pr}(b) and \eqref{f4emf} that  $a=1$ and
 therefore, by virtue of \eqref{f33preresinvlamhyp},
 $f_{1} (z)= 1$ holds for all $z\in\fet$. Since $f_1$
 is holomorphic on $\Bb{H}$ we get
that  $f_{1} (z)= 1$ throughout $\Bb{H}$, which is the
same as \eqref{f8theor3}.

We  now show  that the condition   (3) in
Lemma~{\hyperlink{hth1}{\ref*{th1}}} holds for $f_{2}$
 with $n_{1} = 1$. In view of \eqref{f1emf}(a),  it is
  sufficient
to verify this  when $\fet \ni z \to 1$. For this case
 it follows from \eqref{f1emf}(e)  and
 \eqref{f2emf}(b),(c) that%\vspace{-0,2cm}
\begin{align*}
    & f_{2} (z)\! =\!  - (1/16)\,{{\rm{e}}^{{\fo{{\imag
    \pi}/{(z\!-\!1)} }}}} \! +\! {\Oh} (1)   \quad
     \ \mbox{as} \ \quad \fet \!\ni\! z\to\!1 \,.
\end{align*}

%\vspace{-0,2cm}
\noindent This shows that the condition   (3) of
Lemma~{\hyperlink{hth1}{\ref*{th1}}} holds with
$n_{1} = 1$, since %\vspace{-0,1cm}
\begin{align*}
    &   {\rm{Re}}\left( \imag {\pi}/{(z\!-\!1)} \right)-
    {2\pi}/{|z\!-\!1|} \leq  - {\pi}/{|z\!-\!1|} \to
    - \infty
     \quad  \ \mbox{as} \ \quad    \fet \ni z \to 1 \ .
\end{align*}

%\vspace{-0,1cm}
\noindent
Applying the result of Lemma~{\hyperlink{hth1}{\ref*{th1}}},
we obtain the existence of $a,b \in \Bb{C}$ such that
$f_{2} (\ie(z))= a z + b$.
By letting $ \Lambda\ni z \to 0$ we obtain by \eqref{f2th1pr}(b)
and by \eqref{f4emf} that $f_{2} (\ie(z))\to 0$ which
yields $b=0$ and therefore  $f_{2} (\ie(z))= a z$.
At the same time, by
letting $ \Lambda\ni z \to 1$ we obtain from \eqref{f2th1pr}(c)
and  \eqref{f4emf} that $f_{2} (\ie(z))\to 1$  and therefore
$a=1$. It follows that $ f_{2} (\ie(z))=z$ for all
 $z\in \Lambda$, which gives \eqref{f1theor3}. $\square$

\vspace{0,1cm}
{\emph{Proof of \eqref{f9theor3} and \eqref{f10theor3}.}}{\hyperlink{d20}{${}^{\ref*{case20}}$}}\hypertarget{bd20}{}
It follows from \eqref{f26preresinvlamhyp}, \eqref{f1theor4}
and \eqref{f1theor3} written in the form
$ \ie (\lambda (z)) = z $, $z \in \fet$, that %\vspace{-0,2cm}
\begin{align*}
    &  \ie^{\,\prime} (\lambda (z)) \lambda^{\,\prime} (z)=1 \ ,
    \quad \imag \, \ie^{\,\prime} (\lambda (z))\,
    \Theta_{3}\left(z\right)^{4} =
    \dfrac{1}{\pi \lambda (z)(1-\lambda (z))} \ , \quad
     z \in \fet \ ,
\end{align*}

%\vspace{-0,1cm}
\noindent from which, by using  \eqref{f1xtheor3} and
\eqref{f8theor3}, we obtain \eqref{f9theor3} for $z \in \Lambda$.
 Since all functions in
 \eqref{f9theor3} are holomorphic in $\Bb{H}$ we obtain
 \eqref{f9theor3} for all  $z\in \Bb{H}$. The equalities
 \eqref{f10theor3}
 follow from \eqref{f9theor3} and  from the following forms of
 writing $\lambda^{\,\prime} (z)$, taking account of
 \eqref{f8theor3}, %\vspace{-0,2cm}
 \begin{align*}
    &  \lambda^{\,\prime} (z) = \dfrac{ {\diff} }{ {\diff}  z}
    \dfrac{\Theta_{2}\left(z\right)^{4}}{\Theta_{3}\left(z\right)^{4}} =4 \dfrac{\Theta_{2}\left(z\right)^{4}}{\Theta_{3}\left(z\right)^{4}}  \left(\dfrac{\Theta_{2}^{\,\prime}\left(z\right)}{\Theta_{2}\left(z\right)} - \dfrac{\Theta_{3}^{\,\prime}\left(z\right)}{\Theta_{3}\left(z\right)}\right)  \ , \\    &
    \lambda^{\,\prime} (z)= - \dfrac{ {\diff} }{ {\diff}  z}
     (1-\lambda (z)) = -\dfrac{ {\diff} }{ {\diff}  z}
     \dfrac{\Theta_{4}\left(z\right)^{4}}{\Theta_{3}\left(z\right)^{4}}  = 4 \dfrac{\Theta_{4}\left(z\right)^{4}}{\Theta_{3}\left(z\right)^{4}}  \left(\dfrac{\Theta_{3}^{\,\prime}\left(z\right)}{\Theta_{3}\left(z\right)} - \dfrac{\Theta_{4}^{\,\prime}\left(z\right)}{\Theta_{4}\left(z\right)}\right)\,. \ \square
 \end{align*}

 %\vspace{-0,2cm}
 \noindent
 {\emph{Proof of Corollary~\ref{th5}.}}{\hyperlink{d21}{${}^{\ref*{case21}}$}}\hypertarget{bd21}{}
By combining Corollary~\ref{th0} for $\alpha =2$  with
\eqref{f5theor4} and \eqref{f1theor4}
we find  that $ \Theta_{2}\left(\ie \right)^{4}$ is
universally starlike and $\Theta_{3}(\ie )^{4}\in\mathcal{P}_{\log}(-\infty, 1) $.
 According to the definition \eqref{f9int} of the class
$\mathcal{P}_{\log} $  the latter property means  that
$\Theta_{3}\left(\ie \right)^{4}\in \mathcal{P}$
and
\begin{align}\label{xxx}
    &  \mathcal{P} \ni \dfrac{\dfrac{{{\diff}}}{{{\diff}} z}
    \Theta_{3}\left(\ie(z) \right)^{4}}{4 \Theta_{3}
    \left(\ie (z)\right)^{4}} = \dfrac{ \Theta_{3}^{\,\prime}
    \left(\ie (z)\right)\ie^{\,\prime} (z) }{\Theta_{3}
    \left(\ie (z)\right) }   =    \dfrac{
    \Theta_{3}^{\,\prime}\left(\ie (z)\right) }{
    \lambda^{\,\prime} ( \ie (z) ) \Theta_{3}\left(\ie
    (z)\right)}  \ ,
  \end{align}

 \noindent because $  \lambda^{\,\prime} ( \ie (z) )
 \ie^{\,\prime} (z)=1  $ for each $z \in \Lambda \supset
   \left(\Bb{C}\setminus\Bb{R}\right)$, as follows from
   \eqref{f1theor3}. We need in the following three facts:
 (a) by \eqref{eq:symm0} and  \cite[p.31]{berg}, every
 nonconstant function $f$ in $ \mathcal{P}$   satisfies
  $({\rm{Im}}\, z ) f (z) >0$ for all $z\in\Bb{C}\setminus\Bb{R}$;
   (b) the relation \eqref{f16eir} yields that the numbers
    ${\rm{Im}}\, z$ and ${\rm{Re}}\, \ie (z)$ have the same sign for
     every $z \in \Bb{C}\setminus\Bb{R}$;  (c) we have
     $\ie \left(\Bb{C}\setminus\Bb{R}\right) =$
     $\fet \!\setminus\!\{ \imag
       \Bb{R}_{>0}\} $, in view of
Theorem~\ref{inttheor2} and the equality $\Bb{R}_{>0} =
\{\, \he (1-x)/ \he (x) \, | \, x \in (0, 1) \, \}$ (see Corollary~\ref{cor1}).
By applying these facts to the properties
$\Theta_{3}\left(\ie \right)^{4}\in \mathcal{P}$ and
\eqref{xxx}
we obtain that \eqref{fxth5}(a) and \eqref{fxth5}(b) hold.
   $\square$

%%%%%%%%%%%%%%%%%%%%%%%%%%%%%%%%%%%%%%%%%%%%%%%%%%%%%%%%%%%%%%%%%%%%

%\vspace{-0,4cm}
\section[\hspace{-0,30cm}. \hspace{0,11cm}Definitions for the logarithms of the theta functions]{\hspace{-0,095cm}{\bf{.}} Definitions for the logarithms of the theta functions}
\label{log}

%\vspace{-0,2cm}
By virtue of \eqref{f0theor3} and \eqref{f8elemtheta},
for every $2 \leq k \leq 4$
the function $\theta_{k}$ in \eqref{f1elemtheta} does
not vanish on $\Bb{D}$ and consequently, by
\cite[p. 94, Cor. 6.17]{con}, there exists a holomorphic
 function  $\log \theta_{k}$ in $\Bb{D}$  such that
 $\exp(\log \theta_{k}) =  \theta_{k}$ on $\Bb{D}$, and
$\log \theta_{k}(0) = \ln \theta_{k}(0) = 0$. In addition,
we know from \eqref{f8elemtheta} with $ z \in \imag
 \, \Bb{R}_{>0}$
 and \eqref{f0theor3} that
   $\theta_{k}( x) > 0$ for $x \in [0, 1)$. Since  for arbitrary
 $x \in (0, 1)$  all other solutions $y$  of the equation
 $\exp (y) = \theta_{k}( x)$ differs from $\ln \theta_{k}(x)$
  in $2\pi i n$ with some $n \in \Bb{Z}\setminus\{0\} $,
 and  both functions $\log \theta_{k} (x)$ and $\ln \theta_{k} (x)$
  are continuous on $[0, 1)$ we obtain that
 $\log \theta_{k}(x) = \ln \theta_{k}(x)$, $x\in [0, 1)$.{\hyperlink{d22}{${}^{\ref*{case22}}$}}\hypertarget{bd22}{}
To obtain  the Maclaurin series for $\log \theta_{k}$  we
 use the following classical Jacobi's expansions into
 infinite products
for arbitrary $u \in \Bb{D}$  (see \cite[pp. 469, 470]{whi},
 \cite[p. 529, 20.4.3, 20.4.4]{olv})
\begin{align*}
    &  \theta_{2} (u) \!=\!\prod\limits_{n \geq 1}
    \left(1\!-\!u^{2n}\right)\left(1\! +\! u^{2n}\right)^{2}
    \, ,  \ \theta_{4} (- u) \! = \!
    \theta_{3} (u)\!=\!\prod\limits_{n \geq 1}
    \left(1\!\!-u^{2n}\right)\left(1\! +\! u^{2n-1}\right)^{2} .
\end{align*}

%\vspace{-0,2cm}
\noindent Taking the  real-valued logarithm of these products
for $u \in (0,1)$, we see that
\begin{align*} &
  \ln  \theta_{3} (u) = \sum\nolimits_{n \geq 1} \ln
   \left(1-u^{2n}\right) + 2 \sum\nolimits_{n \geq 1} \ln
    \left(1 + u^{2n-1}\right) \ , \  \\  &
\ln  \theta_{2} (u) = \sum\nolimits_{n \geq 1} \ln
\left(1-u^{2n}\right) +  2\sum\nolimits_{n \geq 1} \ln  \left(1 + u^{2n}\right) \ .
\end{align*}

\noindent
Next, by expanding $\ln (1\pm x)$ in its Maclaurin series
\cite[p. 68, 4.1.24]{abr} we obtain
(compare, e.g., with \cite[p. 338, (4.2)]{ber}) after
several algebraic manipulations that (see
\cite[p. 65]{pre}){\hyperlink{d23}{${}^{\ref*{case23}}$}}\hypertarget{bd23}{}
\begin{align}\label{hak1}
    &  \ln  \theta_{3} (u) =\sum\nolimits_{ n \geq 1}
     \dfrac{2}{2n-1} \dfrac{   u^{2n-1}}{1+u^{2n-1}}  \ , \quad
    \ln  \theta_{2} (u) =\sum\nolimits_{ n \geq 1}
    \dfrac{1}{n} \dfrac{u^{2n}}{1+u^{2n}} \ , \\    &
     \ln  \theta_{4} (u) =  \ln  \theta_{3} (-u) =
      -\sum\nolimits_{ n \geq 1}  \dfrac{2}{2n-1}
      \dfrac{   u^{2n-1}}{1-u^{2n-1}} \ , \quad u \in (0,1) \ .
\label{hak2}\end{align}

\noindent As the three series in the right-hand sides of
these equalities  converge absolutely and uniformly on
compact subsets of the unit disk,  they represent holomorphic
functions in $\Bb{D}$,  which by the standard  uniqueness
 theorem  (see \cite[p. 78, Thm. 3.7(c)]{con})
shows that the identities \eqref{hak1} and \eqref{hak2}
hold throughout~$\Bb{D}$,
\vspace{-0,1cm}
\begin{align}\nonumber
    &  \log  \theta_{3} (u) =  \sum\limits_{ n \geq 1}\,
    \dfrac{2}{ 2n-1}\,\dfrac{u^{2n-1}}{1+ u^{2n-1}}  \ ,    &     &
 \log  \theta_{4} (u)  =   - \sum\limits_{ n \geq 1}\,
 \dfrac{2}{ 2n-1}\,\dfrac{u^{2n-1}}{1- u^{2n-1}} \ ,   \\    &
 \log  \theta_{2} (u) =  \sum\limits_{ n \geq 1}\,\dfrac{1}{n}\,
 \dfrac{u^{2n}}{1+u^{2n}} \,  ,  &     & u \in \Bb{D}  \ .
\label{f3elemtheta}\end{align}

\noindent
As we recall the relationships \eqref{f8elemtheta} connecting
$\theta_{k}$ with $\Theta_{k}$,
we see that this allows to define the logarithms
 $\log \Theta_{k}$ via
%\vspace{-0.1cm}
\begin{align}\nonumber
    & \log \Theta_{k}  (z):=
    \log\theta_{k}\left({\rm{e}}^{\imag  \pi z}\right)
     \, , \  k \in  \left\{3, 4\right\} \, ,  \\    &
\log \Theta_{2}  (z) := \dfrac{\imag \pi z}{4} +  \ln 2+
\log \theta_{2}\left({\rm{e}}^{\imag  \pi z}\right) \, , \
 z \in \Bb{H} \,.
\label{f11elemtheta}\end{align}

%\vspace{-0.1cm}
\noindent The counterpart of the functional relationships
\eqref{f9elemtheta} reads
  %\vspace{-0.1cm}
\begin{align}\nonumber
    &   \log \Theta_{k}  (2m \!+\!z)\!=\! \log \Theta_{k}
      (z)\, , \  3\! \leq\! k\! \leq\! 4\, ,   &     &\hspace{-0,5cm}
\log \Theta_{2}  (2m\! +\!z)\!=\! \dfrac{\imag \pi m}{2}\!
+\! \log \Theta_{2}  (z)\, ,   \\    &
\log \Theta_{k}\left(-1 / z\right) =
  \log \Theta_{6-k}\left(z\right) +  \dfrac{1}{2}\,
  \Log   \dfrac{z}{\imag } \ ,    &     &
      2 \leq k \leq 4 \ ,
  \label{f12elemtheta} \\[0.2cm]    &
\log \Theta_{k}\left(z-1\right) = \log \Theta_{7-k}\left(z\right)
\ ,  \  3 \leq k \leq 4 \ ,     &     &
\hspace{-0,2cm}\log \Theta_{2}\left(z-1\right) =
 \log \Theta_{2}\left(z\right)-\! \dfrac{\imag \pi }{4}\,.\nonumber
\end{align}

\vspace{-0.25cm}
\section[\hspace{-0,15cm}. \hspace{-0,02cm}Exponential integral representation of theta function]{\hspace{-0,095cm}{\bf{.}} Exponential integral representation of  $\Theta_{3}$}
\label{eirth}

\vspace{-0.2cm} {\emph{Proof of Corollary~\ref{cor1}.}}\,{\hyperlink{d24}{${}^{\ref*{case24}}$}}\hypertarget{bd24}{}
The Wirtinger identity \eqref{f1theor4} and the integral
 representation \eqref{f8eir}, by taking into account the
 notation \eqref{f11elemtheta}, allow to write down
 \eqref{f1eirf} and \eqref{f1aeirf} in the set
 \eqref{f3eirf} because according to
   \eqref{f12elemtheta} we have  $ \Theta_{3} (-1+z) =
    \Theta_{3} (1+z)$, $z \in \Bb{H}$, and it follows
     from \eqref{f1theor4}, \eqref{f1ath0} and
      \eqref{f19preresinvlamhyp} that
\begin{align}\label{f1yeirth}
    &  \arg \Theta_{3} (z) \in (- \pi/4 \, , \, \pi/4)
     \ , \quad  z \in \mathcal{F}^{\,{\tn{||}}}_{{\tn{\square}}} \,.
\end{align}

\noindent
 The expression for $ y = y (x)$ follows from
 \eqref{f14zpreresinvlamhyp} and \eqref{f2eir}. $\square$

\vspace{0,1cm}
\noindent
{\emph{Proof of \eqref{f7eirf}.}}
By  \eqref{f12elemtheta},  we have $ \Theta_{3} (2m+z) =
 \Theta_{3} (z)$ for arbitrary $m \in \Bb{Z}$ and $z \in \Bb{H}$,
 from which it follows that the left-hand sides of the equalities
  \eqref{f1eirf} and \eqref{f1aeirf} can be equivalently
   replaced by $ \log\Theta_{3} \left(2m\! +\!\ie (z)\right)$
    and $\log\Theta_{3} \left( 2m\! +\! 1\! +\! \imag
      y\right)$,
    respectively, for arbitrary integer $m$. This gives the
    integral representation
of $ \log \Theta_{3} (z)$  for all $z$ in the set \eqref{f2eirf},
which with the help  of \eqref{f1theor3}  can be written in the
 form \begin{align}\label{f5eirf}
    &  \log\Theta_{3}\! \left(z\right) = \dfrac{1}{2 \pi^{2 }}
    \int\limits_{0}^{1}
    \dfrac{\dfrac{1}{t (1-t)}
     \Log  \dfrac{1}{\vphantom{\frac{A}{B}}1
     -t\Theta_{2}(z)^{4}\Theta_{3}(z)^{-4}} }{\he  \left(t
     \right)^{2} + \he  \left(1-t \right)^{2}} \
      {\diff}  t \, , \quad  z\in \mathcal{F}^{\infty}_{{\tn{\square}}} \,.
 \end{align}

\noindent
 The representation \eqref{f7eirf} is obtained\,{\hyperlink{d25}{${}^{\ref*{case25}}$}}\hypertarget{bd25}{} from
 \eqref{f5eirf} with the help of making there the
 change of variable
 $ \tau = \ie (t)/\imag  $, $t\in (0,1)$,
  which yields $t = \lambda (\imag \tau)$, $\tau\in (0,+\infty)$,
 because in view of \eqref{f7preresinvlamhyp}(a),(b)
 and \eqref{f2eir} we have
 $\ie (t)/\imag  \to +\infty$, as $t \downarrow 0$,
  $\ie (1)/\imag  = 0$ and
  $\ie^{\,\prime} (t)/\imag < 0$,
   $t\in (0,1)$. $\square$

\newpage\renewcommand{\thesection}{A}

\renewcommand*\thecase{[\arabic{clash}]}

\section[\hspace{-0,2cm}. \hspace{0,03cm}Addendum]{\hspace{-0,095cm}{\bf{.}} Addendum}
\label{detpro}

\setcounter{equation}{2}

%\vspace{0.2cm}
\subsection[\hspace{-0,25cm}. \hspace{0,075cm}Notes on Section~\ref{intro}]{\hspace{-0,11cm}{\bf{.}} Notes on Section~\ref{intro}}
\vspace{-0,15cm}

\begin{clash}{\rm{\hypertarget{d1}{}\label{case1}\hspace{0,00cm}{\hyperlink{bd1}{$\uparrow$}}
 \ In view of \eqref{f9int} and Theorem~\ref{inttheor1},
the property that $\het:=\hett$ is universally starlike obtained by Ruscheweyh, Salinas, and Sugawa in
\cite[p.~292, Thm.~1.8]{RS}, yields that $ \het^{\,\prime} (z)/\het (z)\in\mathcal{P}$. But then  $((d/dz)\het (z)^{2})/\het (z)^{2}$ $=$ $2 \het^{\,\prime} (z)/\het (z)\in\mathcal{P}$ and therefore, by virtue of \eqref{f9int} and Theorem~\ref{inttheor1}, the universal starlikeness of $z \het (z)^{2}$ follows from $\het (z)^{2}\in
\mathcal{P}(-\infty,1)$ which is a direct consequence of \eqref{f1cth0}.
}}\end{clash}

\vspace{-0,2cm}
\begin{clash}{\rm{\hypertarget{d2}{}\label{case2}\hspace{0,00cm}{\hyperlink{bd2}{$\uparrow$}} \
The two crucial
properties of the fundamental quadrilateral $\fet$ are the following:  for each
$z\in \fet$ we have $-1/z\in \fet$ and $z - {\rm{sign}}\, (\re \, z) \in \fet$, where
${\rm{sign}} (x)$ is equal to $-1$ if $x < 0$, $0$ if $x=0$ and $1$ if $x>0$.

It should  be noted that Lemma~{\hyperlink{hth1}{\ref*{th1}}} and the basic relationships between the theta functions \eqref{f1theor4},
\begin{align*}\tag{\ref{f8theor3}}
      &     \dfrac{\Theta_{2}^{4}(z) + \Theta_{4}^{4}(z)}{ \Theta_{3}^{4}(z)} =  1\ , \quad  z \in  \Bb{H} \ ,
\end{align*}

\noindent
 and \eqref{f1theor3}
are obtained following the same lines:

\vspace{0.15cm}

\noindent
{\rm{1)}}  we form the ratio of the expressions which should
be shown to be equal (in Lemma~{\hyperlink{hth1}{\ref*{th1}}} this ratio is denoted by $f$ while in
\begin{align*}\tag{\ref{f1theor3}}
      &   \dfrac{ {\Theta_{2}\big(\ie (z)\big)^{4}}}{{\Theta_{3}\big(\ie (z)\big)^{4}}} = z  \ , \quad  z \in
      (0,1)\cup \left(\Bb{C}\setminus\Bb{R}\right) \ ,
\end{align*}

\noindent
the ratio already appears on the left-hand side);

\vspace{0.15cm}

\noindent
{\rm{2)}} for this ratio we replace its argument by the Schwarz triangle function $\ie(z)$
and prove that the composed function assumes the same values on
different sides of the both cuts along $(-\infty, 0)$ and  $(1, +\infty)$ and hence belongs
to ${\rm{Hol}}(\Bb{C}\setminus\{0,1\})$ by the Morera theorem
(such change  in  \eqref{f1theor3}
 has already been made as well as in
  \begin{align*}\tag{\ref{f1theor4}}
        &  \dfrac{\Theta_{3}\big(\ie (z)\big)^{2}}{\ihet (z)}  =1   \ , \quad z \in (0,1)\cup \left(\Bb{C}\setminus\Bb{R}\right) \,  ,
    \end{align*}

\noindent
 because its  original form  is
$ \Theta_{3}(z)^{2}= \ihet (\lambda(z))$\,),

\vspace{0.15cm}

\noindent
{\rm{3)}} by using the asymptotic formulas we prove that  this function  at $0$ and $1$   either has a pole or removable singularity while at $\infty$
it has  polynomial growth,  and hence it is a rational function by  the Riemann theorem about  removable singularities and by the extended version of  Liouville's theorem;

\vspace{0.15cm}

\noindent
{\rm{4)}} the first terms of asymptotic expansions of  this function at the points $0$, $1$ and at $\infty$  give the values of all  coefficients of  this rational function.

}}\end{clash}

 \begin{subequations}
\begin{clash}{\rm{\hypertarget{d3}{}\label{case3}\hspace{0,00cm}{\hyperlink{bd3}{$\uparrow$}} \ It will be shown below{\hyperlink{d14}{${}^{\ref*{case14}}$}}\hypertarget{bd14}{}
that for arbitrary function $f$ satisfying the conditions of Lemma~{\hyperlink{hth1}{\ref*{th1}}} the functions
\begin{align}\label{f1bd3}
    &   \Phi (z) := f \big(\ie (z)\big)\,,   &     &        z \in \Lambda := (0,1)\cup \left(\Bb{C}\setminus\Bb{R}\right)  \, , \  \\    &
      \Psi (z) := f\Bigg(\ie \bigg(\dfrac{1}{1- z^{2}}\bigg)\Bigg)\,,   &     &       z \in \Bb{H}  \, ,
\label{f2bd3}\end{align}

\noindent satisfy correspondingly the conditions of the following statements which are  simple consequences of the   Liouville theorem \cite[p. 2, Thm. 1]{lev2}  and of the Riemann theorem about  removable singularities \cite[p. 103, Thm. 1.2]{con}.}}\end{clash}
\begin{lemma}\hspace{-0,18cm}{\bf{.}}\label{bd3lem1}
 Let $\Phi$ be holomorphic on $ \Lambda := (0,1)\cup \left(\Bb{C}\setminus\Bb{R}\right)$ and
 can be extended continuously from $\Bb{H}$ to  $(\Bb{H}\cup\Bb{R} ) \setminus [\,0, 1]$ and
  from $-\Bb{H}$
  to $(- \Bb{H}\cup\Bb{R} ) \setminus [\,0, 1]$ such that
\begin{align}\label{f1bd3lem1} {\hypertarget{hbd3lem1}{}}
    &  \Phi (-x- i  0 )= \Phi (-x+ i  0 )  \, , \ \  \Phi (1+ x -  i  0 ) = \Phi (1+ x +  i  0 ) \, , \ \  x > 0 \, ,
\end{align}

\noindent where $\Phi (x\pm  \imag   0):= \lim_{{\fo{\,\varepsilon\! \downarrow\! 0}}} \ \Phi (x\pm \imag  \varepsilon)$,\  $x \in (-\infty, 0) \cup (1, +\infty)$.

\vspace{0,1cm}
 Suppose
that  there exist nonnegative integers $n_{\infty}$,
 $n_{0}$, and $n_{ 1}$ satisfying
\begin{align*}  &{\rm{(1)}}&      |z|^{ {\fo{- n_{ 1} -1}}}\left|\Phi (z)\right|\! \to 0\!\, ,   &   &    \Lambda \! \ni\! z \!\to\!  \infty \, ,\\[0,3cm]
      &  {\rm{(2)}}  &        |z|^{{\fo{n_{\infty} +1}}} \left|\Phi (z)\right|\! \to 0\!\, ,   &   &    \Lambda\! \ni\! z \!\to\! 0 \, ,
    \\[0,3cm]      &
    {\rm{(3)}} &       |1\!-\!z|^{{\fo{n_{0}\!+\!1}}}\left|\Phi (z)\right|\!\to\! 0\, ,   &  &    \Lambda \!\ni \!z \!\to\! 1\,     .
\end{align*}

\noindent Then there exists an algebraic polynomial $P$ of
degree $\le n_{\infty} + n_{0} + n_{1}$ such that \vspace{-0,3cm}
 \begin{align*}
    &  \Phi (z) = \frac{\vphantom{\frac{a}{b}}
P(z)}{\vphantom{\dfrac{a}{b}}
     z^{{\fo{n_{\infty}}}} (1-z)^{{\fo{n_{0}}}} }  \ ,
\quad  z \in  (0,1)\cup \left(\Bb{C}\setminus\Bb{R}\right) \,.
\end{align*}

\end{lemma}

\begin{lemma}\hspace{-0,18cm}{\bf{.}}\label{bd3lem2}
  Let $\Psi$ be holomorphic on $\Bb{H}$ and
can be extended continuously to  $(\Bb{H}\cup\Bb{R} ) \setminus \{-1, 0, 1\}$ such that
\begin{align}\label{f1bd3lem2}  {\hypertarget{hbd3lem2}{}}
    &   \Psi ( x)=  \Psi ( -x)  \ , \quad  x \in \Bb{R}\setminus\{-1, 0, 1\}  \, ,
\end{align}

\noindent where $\Psi (x):= \lim_{{\fo{\,\varepsilon\! \downarrow\! 0}}} \ \Psi (x+ i  \varepsilon)$, \ $x \in \Bb{R}\setminus\{-1, 0, 1\}$.

\vspace{0,1cm}
 Suppose
that  there exist nonnegative integers $n_{\infty}$,
 $n_{0}$, and $n_{ 1}$ satisfying
\begin{align*}
&{\rm{(1)}}&      &    & |z|^{ {\fo{- 2n_{ \infty} -2}}}\left|\Psi (z)\right| \to 0\, ,            &   &       \Bb{H} \! \ni\! z \!\to\!   \infty \,, &   & \\[0,3cm]
&{\rm{(2)}}&     &     &    |z|^{{\fo{\,2n_{0}+2}}}  \left|\Psi (z)\right|\! \to 0\, ,        &    &      \Bb{H}\! \ni\! z \!\to\! 0 \, ,  &   & \\[0,3cm]      &{\rm{(3)}} &     &     &  |\sigma\!-\!z|^{{\fo{\,n_{1}+1}}}\left|\Psi (z)\right|\to 0\, ,   &    &     \Bb{H} \!\ni \!z \!\to\! \sigma\, ,  &   &  \sigma\in\{1,-1\}\, .
\end{align*}

\noindent Then there exists an algebraic polynomial $Q$ of
degree $\le n_{\infty} + n_{0} + n_{1}$ such that\vspace{-0,3cm}
 \begin{align}\label{f2bd3lem2}
    &   \Psi (z) = \dfrac{Q \left(z^{2}\right)}{\vphantom{\dfrac{a}{b}}\left(1-z^{2}\right)^{{\fo{n_{1}}}}  z^{{\fo{2 n_{0}}}}   } \ , \quad  z \in \Bb{H}\,.
\end{align}

\end{lemma}

\end{subequations}

%\vspace{0.25cm}
\subsection[\hspace{-0,25cm}. \hspace{0,075cm}Notes on Section~\ref{bf}]{\hspace{-0,11cm}{\bf{.}} Notes on Section~\ref{bf}}
 \begin{subequations}
\begin{clash}{\rm{\hypertarget{d4}{}\label{case4}\hspace{0,00cm}{\hyperlink{bd4}{$\uparrow$}} \ {\emph{Proof of \eqref{fa6preresinvlamhyp}.}} In view of \cite[15.3.10, p.559]{abr}, for arbitrary $z \in (1 + \Bb{D})\setminus [1, +\infty)$ we have
\begin{align*}
    &     \het (z)    =\dfrac{1}{\pi^{2}} \sum\limits_{n=0}^{\infty} \dfrac{\Gamma (n+1/2)^{2}}{(n !)^{2}}
\Big[2 \psi (n+1) - 2 \psi (n+1/2) - \Log (1-z) \Big] (1-z)^{n} \, ,
\end{align*}

\noindent where  $\psi$ denotes the digamma function  (see \cite[6.3.2, 6.3.3,  p.258]{abr}) and
\begin{align}\label{ad0}
    &   \het (z) := F\left(\dfrac{1}{2} \ , \ \dfrac{1}{2} \ ; \  1  \ , \  z \right) = \dfrac{1}{\pi} \sum\limits_{n = 0}^{\infty}
\dfrac{\Gamma(n+1/2)^{2}}{(n !)^{2}} z^{n}  \ , \quad z \in \Bb{D}\,.
\end{align}

\noindent Thus, for any $z \in (1 + \Bb{D})\setminus [1, +\infty)$ the next equality holds
\begin{multline}\label{ad1}
  \het (z)    = \dfrac{ \Log \dfrac{1}{1-z}}{\pi}\het (1-z) \\[0,3cm] +
  \dfrac{2}{\pi^{2}} \sum\limits_{n=0}^{\infty} \dfrac{\Gamma (n+1/2)^{2}}{(n !)^{2}} \left[\psi (n+1) -  \psi (n+1/2)\right](1-z)^{n}  \,.
\end{multline}

\noindent According to the definition of the digamma function (see \cite[6.3.2, 6.3.4, p.258]{abr})
\begin{align*}
     \psi (1) & = - \gamma \, ,   &         \psi (n)&  = - \gamma + \sum\limits_{k=1}^{n-1} \dfrac{1}{k}  \ , \quad  n \geq 2  \, , \\
 \psi (1/2) & = - \gamma - \ln 4 \, ,      &        \psi (n+1/2) &  = - \gamma - \ln 4  + 2 \sum\limits_{k=0}^{n-1} \dfrac{1}{2k +1}   \\
 &   &     &     =
 - \gamma +2 - \ln 4  +  \sum\limits_{k=1}^{n-1} \dfrac{1}{k +1/2}
 \ , \quad  n \geq 1  \ ,
\end{align*}

\noindent and therefore, assuming $\sum_{k=1}^{0} := 0$, we get
\begin{align*}
    & \psi (n+1) = - \gamma + \sum\limits_{k=1}^{n} \dfrac{1}{k}  \ , \ \  \psi (n+1/2) = - \gamma - \ln 4  +  \sum\limits_{k=1}^{n} \dfrac{2}{2k -1}   \ , \  n \geq 0  \, .
\end{align*}

\noindent Since by \cite[p. 658, 5.1.8.2]{prud1}
\begin{align*}
    &  \sum\limits_{k\geq 0} \dfrac{1}{(k +1)(2k +1)} = \ln 4  \ , \
\end{align*}

\noindent we obtain
\begin{align*}
    &  \psi (n+1) - \psi (n+1/2) = \ln 4 + \sum\limits_{k=1}^{n} \dfrac{1}{k} - \sum\limits_{k=1}^{n} \dfrac{2}{2k -1}=
    \ln 4 + \sum\limits_{k=0}^{n-1} \dfrac{1}{k+1}\\    & - \sum\limits_{k=0}^{n-1} \dfrac{2}{2k +1}    =
   \sum\limits_{k\geq 0} \dfrac{1}{(k +1)(2k +1)} +  \sum\limits_{k=0}^{n-1} \left(\dfrac{1}{k+1} -  \dfrac{2}{2k +1}\right) \\    &   =
   \sum\limits_{k\geq 0} \dfrac{1}{(k +1)(2k +1)}- \sum\limits_{k=0}^{n-1}\dfrac{1}{(k +1)(2k +1)} =
    \sum\limits_{k\geq n} \dfrac{1}{(k +1)(2k +1)}\ ,
\end{align*}

\noindent i.e.,
\begin{align*}
    &  \psi (n+1) -  \psi (n+1/2)= \sum\limits_{k\geq n} \dfrac{1}{(k +1)(2k +1)} \ , \quad  n \geq 0 \ , \\    &
    \psi (1) -  \psi (1/2) = \ln 4 \, .
\end{align*}

\noindent Besides that, the summand corresponding to $n=0$ in the series \eqref{ad1} is equal to ($\Gamma (1/2) = \sqrt{\pi}$,
see \cite[p. 258, 6.1.8]{abr})
\begin{align*}
    &   \dfrac{2}{\pi^{2} }  \Gamma (1/2)^{2}\left(\psi (1) -  \psi (1/2)\right) = \dfrac{\ln 16}{\pi} \ .
\end{align*}

\noindent Thus, for any $z \in (1 + \Bb{D})\setminus [1, +\infty)$ we can write \eqref{ad1} as follows
\begin{align}
    & \nonumber   \het (z)    - \dfrac{ \Log \dfrac{1}{1-z}}{\pi}\het (1-z) \\[0,3cm]  & \nonumber    =
  \dfrac{2}{\pi^{2}} \sum\limits_{n=0}^{\infty} \dfrac{\Gamma (n+1/2)^{2}}{(n !)^{2}} \left[\psi (n+1) -  \psi (n+1/2)\right](1-z)^{n} \\[0,3cm]  & \nonumber    =
  \dfrac{2}{\pi^{2}} \sum\limits_{n=0}^{\infty} \dfrac{\Gamma (n+1/2)^{2}}{(n !)^{2}} \left[\sum\limits_{k\geq n} \dfrac{1}{(k +1)(2k +1)}\right](1-z)^{n}\\[0,3cm]  &    = \dfrac{\ln 16}{\pi} +
  \dfrac{2}{\pi^{2}} \sum\limits_{n=1}^{\infty} \dfrac{\Gamma (n+1/2)^{2}}{(n !)^{2}} \left[\sum\limits_{k\geq n} \dfrac{1}{(k +1)(2k +1)}\right](1-z)^{n} \,,
\label{ad2}\end{align}

\noindent which proves \eqref{fa6preresinvlamhyp}.

}}\end{clash}
\end{subequations}

\begin{clash}{\rm{\hypertarget{d5}{}\label{case5}\hspace{0,00cm}{\hyperlink{bd5}{$\uparrow$}} \  We give the following more precise form of the relations
\eqref{f7preresinvlamhyp}{\rm{(a)}}--{\rm{(c)}}.
}}\end{clash}
\begin{subequations}

\vspace{-0,2cm}
\begin{lemma}\hspace{-0,18cm}{\bf{.}}\label{adth1}
 For the Gauss hypergeometric function
 \begin{align*}
    &   \het (z) := F\left(\dfrac{1}{2} \ , \ \dfrac{1}{2} \ ; \  1  \ , \  z \right) = \dfrac{1}{\pi} \sum\limits_{n = 0}^{\infty}
\dfrac{\Gamma(n+1/2)^{2}}{(n !)^{2}} z^{n}  \ , \quad z \in \Bb{D}\,,
 \end{align*}

\noindent and for  arbitrary $\theta \in (0,1)$ the following inequalities hold
 \begin{align}
      \label{f15preresinvlamhyp}
\Big|\he (z) - 1\Big| \leq \he (\theta)\, |z|  \    \qquad   \qquad
\mbox{for all} \quad     \ z \in \theta\!\!\hphantom{1} \Bb{D} \ ,  &   &   &     &     &     &     &     &     &     &     &     &
\end{align}
\begin{multline}
\label{f16preresinvlamhyp}
\left| \het (z)    -\dfrac{1}{\pi}  \Log \dfrac{16}{1-z} \right|
         \leq      \dfrac{\he (\theta)}{\pi}|1-z|\ln \dfrac{51}{|1-z|} \ ,
 \\[0,3cm]    \mbox{for all} \quad
 z \in (1 \!+ \!\theta \!\!\hphantom{1}\Bb{D})\setminus [1, 1+\theta]  \, ,
 \end{multline}
  \begin{multline}       \label{f17preresinvlamhyp}
  \left|  \he (z)    -\dfrac{1}{\pi}\dfrac{\Log 16 (1-z)}{\sqrt{1-z}} \right| \leq
    \dfrac{\he (\theta)}{\pi} \dfrac{ \left(4 + \ln |z-1|\right)}{|z-1|^{3/2} }\ ,   \\[0,3cm]
\mbox{for all} \quad         \
z \not\in \left(1 + (1/\theta)\!\!\hphantom{1}\overline{\Bb{D}}\,\right) \cup [1+ 1/\theta, +\infty)  \ .
\end{multline}

 \end{lemma}
\end{subequations}
\begin{subequations}
\vspace{0.25cm} {\bf{Proof}} of Lemma~\ref{adth1} and \eqref{f7preresinvlamhyp}{\rm{(a)}}--{\rm{(c)}}.  \ Fix an arbitrary $\theta \in (0,1)$.

\vspace{0.25cm}
\noindent
 {\emph{Proof of\eqref{f7preresinvlamhyp}{\rm{(b)}} and \eqref{f15preresinvlamhyp}.}} It follows from
 \eqref{ad0} that
 \begin{align*}
    &   \het (z)  = \dfrac{1}{\pi} \sum\limits_{n = 0}^{\infty} \dfrac{\Gamma(n+1/2)^{2}}{(n !)^{2}} z^{n} = 1 +
\dfrac{z}{\pi} \sum\limits_{n = 1}^{\infty} \dfrac{\Gamma(n+1/2)^{2}}{\Gamma(n+1)^{2}} z^{n-1} \ , \\    &
\sum\limits_{n = 1}^{\infty} \dfrac{\Gamma(n+1/2)^{2}}{\Gamma(n+1)^{2}} z^{n-1}=
\sum\limits_{n = 0}^{\infty} \dfrac{\Gamma(n+3/2)^{2}}{\Gamma(n+2)^{2}} z^{n}=
\sum\limits_{n = 0}^{\infty} \dfrac{\Gamma(n+3/2)\Gamma(n+3/2)}{\Gamma(n+2) \Gamma(n+2)} z^{n} \ ,
 \end{align*}

\noindent i.e.,
\begin{align*}
    &  \he (z)  = 1 + \dfrac{z}{\pi} \sum\limits_{n = 0}^{\infty} \dfrac{\Gamma(n+3/2)\Gamma(n+3/2)}{\Gamma(n+2) \Gamma(n+2)} z^{n} \ , \quad z \in \Bb{D}\,,
\end{align*}

\noindent where
\begin{align*}
    &  \dfrac{|z|}{\pi} \sum\limits_{n = 0}^{\infty} \dfrac{\Gamma(n+3/2)\Gamma(n+3/2)}{\Gamma(n+2) \Gamma(n+2)} |z|^{n}  \\    &   =
\dfrac{|z|}{\pi} \sum\limits_{n = 0}^{\infty} \left(\dfrac{n+1/2}{n+1}\right)^{2} \dfrac{\Gamma(n+1/2)\Gamma(n+1/2)}{\Gamma(n+1) \Gamma(n+1)} |z|^{n}
 \leq \\ & \leq
  \dfrac{|z|}{\pi} \sum\limits_{n = 0}^{\infty}  \dfrac{\Gamma(n+1/2)^{2}}{(n !)^{2} } |z|^{n} = |z| \he (|z|) \ , \
\end{align*}

\noindent and therefore
\begin{align}\label{ad3}
    &  \Big|\he (z) - 1\Big| \leq \he (\theta)\, |z| \ ,   \quad      \ z \in \theta\!\!\hphantom{1} \Bb{D} \ , \ \theta \in (0,1) \,,
\end{align}

\noindent which proves \eqref{f7preresinvlamhyp}{\rm{(b)}} and \eqref{f15preresinvlamhyp}.

\vspace{0.25cm} \noindent {\emph{Proof of \eqref{f7preresinvlamhyp}{\rm{(a)}} and \eqref{f16preresinvlamhyp}.}}
Let $z \in (1 \!+ \!\theta \,\Bb{D})\setminus [1, 1+\theta]$. Then, by virtue of \eqref{ad3},
\begin{align*}
    & \left| \dfrac{1}{\pi}\het (1-z) \Log \dfrac{1}{1-z} - \dfrac{1}{\pi} \Log \dfrac{1}{1-z}\right| \leq \dfrac{\he (\theta)}{\pi}\, |1-z|
  \left|  \Log \dfrac{1}{1-z}\right| \\    &   \leq
\dfrac{\he (\theta)}{\pi}\,  |1-z| \left(\pi + \ln \dfrac{1}{|1-z|}\right) \ , \
\end{align*}

\noindent while in view of
\begin{align*}
    &  \sum\limits_{k\geq n} \dfrac{1}{(k +1)(2k +1)} \leq \sum\limits_{k\geq 0} \dfrac{1}{(k +1)(2k +1)} -1 = \ln 4 -1 \ , \quad n \geq 1 \ ,
\end{align*}

\noindent we also have that
\begin{align*}
    & \left| \dfrac{2}{\pi^{2}} \sum\limits_{n\geq 1} \dfrac{\Gamma (n+1/2)^{2}}{(n !)^{2}} \left[\sum\limits_{k\geq n} \dfrac{1}{(k +1)(2k +1)}\right](1-z)^{n}\right|  \\    &   \leq
     \dfrac{2}{\pi^{2}} \sum\limits_{n\geq 1} \dfrac{\Gamma (n+1/2)^{2}}{(n !)^{2}} \left[\sum\limits_{k\geq n} \dfrac{1}{(k +1)(2k +1)}\right]|1-z|^{n}
     = \\ & =    \dfrac{2|1-z|\ln (4/e)}{\pi^{2}} \sum\limits_{n=1}^{\infty} \dfrac{\Gamma (n+1/2)^{2}}{\Gamma (n+1)^{2}} |1-z|^{n-1} \leq \\ & \leq
  \dfrac{2|1-z|\ln (4/e)}{\pi^{2}}  \sum\limits_{n=0}^{\infty} \dfrac{\Gamma (n+3/2)^{2}}{\Gamma (n+2)^{2}} |1-z|^{n} \leq \\ & \leq
 \dfrac{2|1-z|\ln (4/e)}{\pi^{2}} \sum\limits_{n=0}^{\infty} \dfrac{\Gamma (n+1/2)^{2}}{\Gamma (n+1)^{2}}  |1-z|^{n} = \\ & =
 \dfrac{2\ln (4/e)}{\pi} |1-z|  \he (|1-z|) \leq \dfrac{2\ln (4/e)}{\pi} |1-z|  \he (\theta) \ , \
\end{align*}

\noindent which together with \eqref{ad2} and the inequality $1 < \he (\theta)$ gives
\begin{align*}
    &  \left| \het (z)    -\dfrac{1}{\pi}  \Log \dfrac{1}{1-z} -\dfrac{\ln 16}{\pi}\right| \\    &   \leq
    \dfrac{\he (\theta)}{\pi}\,  |1-z| \left(\pi + \ln \dfrac{1}{|1-z|}\right) + \dfrac{2\ln (4/e)}{\pi} |1-z|  \he (\theta)  \\ & =
    \dfrac{\he (\theta)}{\pi}|1-z|\ln \dfrac{1}{|1-z|} +  \left(\he (\theta) + (2/\pi)\ln (4/e)\right)|1-z| \\    &   \leq
    \dfrac{\he (\theta)}{\pi}|1-z|\left(2 \ln (4/e) +  \pi +\ln \dfrac{1}{|1-z|} \right)  \\ & =
    \dfrac{\he (\theta)}{\pi}|1-z|\left( \pi  - 2 + \ln \dfrac{16}{|1-z|} \right)
    \ , \
\end{align*}

\noindent i.e., taking account of $16 \exp (\pi - 2) < 50,11 < 51$, we obtain
\begin{align}
    \left| \het (z)    -\dfrac{1}{\pi}  \Log \dfrac{16}{1-z} \right|
         \leq      \dfrac{\he (\theta)}{\pi}|1-z|\ln \dfrac{51}{|1-z|}  \ , \
\label{ad4}\end{align}

\vspace{0.25cm}
\noindent for arbitrary $z \in (1 \!+ \!\theta \,\Bb{D})\setminus [1, 1+\theta]$. This inequality proves \eqref{f7preresinvlamhyp}{\rm{(a)}} and \eqref{f16preresinvlamhyp}.

\vspace{0.25cm} \noindent {\emph{Proof of \eqref{f7preresinvlamhyp}{\rm{(c)}} and \eqref{f17preresinvlamhyp}.}} Since for arbitrary
\begin{align}\label{ad5}
    &  z \not\in \left(1 + (1/\theta)\!\!\hphantom{1}\overline{\vphantom{A^{1}}\Bb{D}}\right)
     \cup [1 +1/\theta, +\infty) \ , \
\end{align}

\noindent we have
\begin{align*}
    & z-1 \not\in (1/\theta)\!\!\hphantom{1}\overline{\vphantom{A^{1}}\Bb{D}} \cup  [1/\theta , +\infty) \ , \\    &
    \dfrac{1}{z-1} \in  \theta \, \Bb{D} \setminus [0, \theta] \ , \\    &
\dfrac{z}{z-1} =     1 + \dfrac{1}{z-1} \in\left(1+ \theta \, \Bb{D} \right)\setminus [1, 1+\theta] \ , \
\end{align*}

\noindent then we can apply \eqref{ad4} with $z/(z-1)$ instead of $z$ for   $z$ satisfying \eqref{ad5},
\begin{align*}
    &   \left| \het \left(\dfrac{z}{z-1}\right)    -\dfrac{1}{\pi}  \Log \dfrac{16}{1-\dfrac{z}{z-1}} \right|
         \leq      \dfrac{\he (\theta)}{\pi}\left|1-\dfrac{z}{z-1}\right|\ln \dfrac{51}{\left|1-\dfrac{z}{z-1}\right|} \ , \
\end{align*}

\noindent to obtain, taking account of $\ln 51 < 3,94 < 4$,
\begin{align*}
    &   \left| \het \left(\dfrac{z}{z-1}\right)    -\dfrac{1}{\pi}\Log 16 (1-z) \right|
         \leq   \dfrac{\he (\theta)}{\pi} \dfrac{4 + \ln |z-1| }{|z-1|}  \ , \
\end{align*}

\noindent from which and the   Pfaff formula
\begin{align*}\tag{\ref{f5preresinvlamhyp}}
    &  \he (z) = \dfrac{1}{\sqrt{1-z}} \ \he
\left(\dfrac{z}{z-1}\right)   \ , \quad
 z \in \Bb{C}\setminus [1, +\infty)  \ , \
\end{align*}

\noindent it follows that
\begin{align*}
    &   \left|  \he (z)    -\dfrac{1}{\pi}\dfrac{\Log 16 (1-z)}{\sqrt{1-z}} \right| \leq
    \dfrac{\left| \het \left(\dfrac{z}{z-1}\right)    -\dfrac{1}{\pi}\Log 16 (1-z) \right|}{\left|\sqrt{1-z}\right|} \\    &
         \leq   \dfrac{\he (\theta)}{\pi} \dfrac{ \left(4 + \ln |z-1|\right)}{|z-1| \left|\sqrt{1-z}\right|}
      = \dfrac{\he (\theta)}{\pi} \dfrac{ \left(4 + \ln |z-1|\right)}{|z-1|^{3/2} }     \ , \
\end{align*}

\noindent  and therefore
\begin{align}\label{ad6}
    &   \left|  \he (z)    -\dfrac{1}{\pi}\dfrac{\Log 16 (1-z)}{\sqrt{1-z}} \right| \leq
    \dfrac{\he (\theta)}{\pi} \dfrac{ \left(4 + \ln |z-1|\right)}{|z-1|^{3/2} }
    \end{align}

\noindent for arbitrary   $z$ satisfying \eqref{ad5}. This proves \eqref{f7preresinvlamhyp}{\rm{(c)}} and \eqref{f17preresinvlamhyp}.
This completes the proof   of Lemma~\ref{adth1}.

\vspace{0,25cm}
Since $\he $ is strictly increasing function on the interval $[0,1)$ the inequalities of  Lemma~\ref{adth1} can be written in another form.
\begin{corollary}\hspace{-0,18cm}{\bf{.}} In the notations of Lemma~\ref{adth1},
the following inequalities hold
 \begin{align}
      \label{f1adcor1}
\Big|\he (z) - 1\Big| \leq \he (|z|)\, |z|  \    \qquad   \qquad
\mbox{for all} \quad     \ z \in \Bb{D} \ ,  &   &   &     &     &     &     &     &     &     &     &     &
\end{align}
\begin{multline}
\label{f2adcor1}
\left| \het (z)    -\dfrac{1}{\pi}  \Log \dfrac{16}{1-z} \right|
         \leq      \dfrac{\he (|1-z|)}{\pi}|1-z|\ln \dfrac{51}{|1-z|} \ ,
 \\[0,3cm]    \mbox{for all} \quad
 z \in (1 + \Bb{D})\setminus [1, 2]  \, ,
 \end{multline}
  \begin{multline}       \label{f3adcor1}
  \left|  \he (z)    -\dfrac{1}{\pi}\dfrac{\Log 16 (1-z)}{\sqrt{1-z}} \right| \leq
     \dfrac{4 + \ln |z-1|}{\pi\,  |z-1|^{3/2} } \he \left(\frac{1}{|z-1|}\right)\ ,   \\[0,3cm]
\mbox{for all} \quad         \
z \not\in \left(1 + \overline{\Bb{\vphantom{a^{a}}D}}\,\right) \cup [2, +\infty)  \ .
\end{multline}
\end{corollary}

\vspace{0.25cm} The following consequence of the inequality
\eqref{f2adcor1} is used below.
\begin{corollary}\hspace{-0,18cm}{\bf{.}}\label{adcor2} The following inequality holds
\begin{align}\label{f1adcor2}
    &  \left| \het \left(1-\dfrac{1}{1+ r e^{{\fo{i \varphi}}} }\right)\right| \leq 3 + \dfrac{\ln r}{\pi}  \, , \quad \varphi \in \left[-\dfrac{\pi}{2}\ , \ \dfrac{\pi}{2} \right]  \,, \  r \geq 3 \,.
\end{align}
{\hypertarget{hadcor2}{}}
\end{corollary}

\vspace{-0.2cm} \noindent {\emph{Proof of Corollary~{\hyperlink{hadcor2}{\ref*{adcor2}}}.}} Observe that \eqref{f1adcor2} follows from the next inequality
\begin{align}\label{f2adcor2}
    &  \left| \het \left(1-\dfrac{ e^{{\fo{i \psi}}}}{ R }\right)\right| \leq 2,\!5 + \dfrac{\ln  R}{\pi}   \ , \quad
     \psi \in [-\pi/2, \pi/2]  \,, \  R \geq 3 \,.
\end{align}

\noindent Actually, for every $r \geq 3$ and $\varphi \in [-\pi/2, \pi/2]$ we have
\begin{align*}
    & r e^{{\fo{i \varphi}}} \in \Bb{C}_{\re\geqslant 0}:= \left\{z\in \Bb{C} \ | \ \re z \geq 0 \right\}   &     &    \Rightarrow \
    1+ r e^{{\fo{i \varphi}}} \in \Bb{C}_{\re\geqslant 0} \\    &    &     &      \ \Rightarrow \
    \Arg \left(1+ r e^{{\fo{i \varphi}}}\right)\in [-\pi/2, \pi/2]  \ ; \end{align*}

\vspace{-0,8cm}
\begin{align*}  &
   \left|1+ r e^{{\fo{i \varphi}}}\right|^{2} = 1+ r^{2} + 2 r \cos \varphi \in \left[1+ r^{2}\ , \ (1+ r)^{2} \right] \ , \  &     &   \phantom{  1+ r^{2} + 2 r \cos \varphi aa}
\end{align*}

\noindent and hence there exist $R$ and $\psi$ such that
\begin{align*}
    &   1+ r e^{{\fo{i \varphi}}} = R e^{{\fo{- i \psi}}}  \ , \quad R \in \left[\sqrt{1+ r^{2}}\ , \ 1+ r\right] \,, \ \psi \in [-\pi/2, \pi/2]\,.
\end{align*}

\noindent Since $r\geq 3$ yields $R \geq \sqrt{1+ r^{2}} \geq 3$, $\ln R \leq \ln (r+1)\leq \ln (1+1/r) + \ln r \leq 0,3  + \ln r$ and obviously
$0,3/\pi < 1/2$ we can apply \eqref{f2adcor2} to get
\begin{align*}
    &  \left| \het \left(1-\dfrac{1}{1+ r e^{{\fo{i \varphi}}} }\right)\right|  = \left| \het \left(1-
    \dfrac{ 1 }{ R e^{{\fo{-i \psi}}}}\right)\right| \leq \dfrac{5}{2} + \dfrac{\ln  R}{\pi} \leq 3 +
    \dfrac{\ln  r }{\pi} \,,
\end{align*}

\noindent which proves the validity of \eqref{f1adcor2}.

\vspace{0.15cm}
To complete the proof of Corollary~{\hyperlink{hadcor2}{\ref*{adcor2}}} it remains to prove  \eqref{f2adcor2}.
Since   for every $R \geq 3$ and $\psi \in [-\pi/2, \pi/2]$ we have
\begin{align*}
    & \dfrac{ e^{{\fo{i \psi}}}}{ R } \in \Bb{D}\setminus [-1, 0] \ \Rightarrow \    1-\dfrac{ e^{{\fo{i \psi}}}}{ R } \in
    (1 + \Bb{D})\setminus [1, 2] \ ,
\end{align*}

\noindent we can apply \eqref{f2adcor1} by using \eqref{ad11e},
\begin{align}{\tag{\ref{ad11e}}}
    & \he \left(\dfrac{1}{2}\right) =  F \left(\dfrac{1}{2},\dfrac{1}{2};1;\dfrac{1}{2}\right) <1,18036 < 2 \ .
\end{align}

%\vspace{-0,3cm}
\noindent  to get%\vspace{-0,5cm}
\begin{align*}
      \left| \het \left(1-\dfrac{e^{{\fo{i \psi}}} }{ R  }\right)    -\dfrac{1}{\pi}  \Log 16 R e^{{\fo{-i \psi}}} \right|
       &  \leq      \dfrac{\he \left({1}/{R}\right)}{\pi R} \ln (51 R)  \\    &   \leq
       \dfrac{\he \left({1}/{2}\right)}{\pi}\left( \dfrac{\ln 51}{R} + \dfrac{\ln R}{R}\right)\\    &   \leq
       \dfrac{2 \ln 51 + 2 \ln 3}{3 \pi} \ ,
\end{align*}

\noindent because the function $x^{-1}\ln x$ is decreasing for $x \geq 3$.  So that
\begin{align*}
    &  \left|\het \left(1-\dfrac{e^{{\fo{i \psi}}} }{ R  }\right)\right| = \left|\het \left(1-\dfrac{e^{{\fo{i \psi}}} }{ R  }\right) - \dfrac{1}{\pi}  \Log 16 R e^{{\fo{-i \psi}}}\right| + \left|\dfrac{1}{\pi}  \Log 16 R e^{{\fo{-i \psi}}}\right|  \\    & \leq
       \dfrac{2 \ln 51 + 2 \ln 3}{3 \pi} + \left| \dfrac{\ln 16 R}{\pi}  -\dfrac{i \psi}{\pi}\right|   \leq
          \dfrac{2 \ln 153 }{3 \pi} + \dfrac{\ln 16 }{\pi} + \dfrac{\ln  R}{\pi} + \dfrac{1}{2} \\    &   =
   \dfrac{\ln  R}{\pi} +    \dfrac{1}{2} +  \dfrac{2 (\ln 153 +  \ln 64)}{3 \pi}  =
     \dfrac{\ln  R}{\pi} +    \dfrac{1}{2} +  \dfrac{2\ln 9792}{3 \pi}\leq  \dfrac{\ln  R}{\pi} +    \dfrac{1}{2} +2\,,
\end{align*}

\noindent which proves the validity of   \eqref{f1adcor2} and completes the proof of Corollary~{\hyperlink{hadcor2}{\ref*{adcor2}}}.$\square$
\end{subequations}

\begin{clash}{\rm{\hypertarget{d6}{}\label{case6}\hspace{0,00cm}{\hyperlink{bd6}{$\uparrow$}} \ For every $z\in (1 + \Bb{D})\setminus [1, 2]$ we have $1-z \in \Bb{D}$ and
therefore the formulas
\begin{multline*}\tag{\ref{fa6preresinvlamhyp}}
         \he (z)    = \dfrac{1 }{\pi}\he (1-z) \Log \dfrac{1}{1-z} \\       + \dfrac{2}{\pi^{2}} \sum\limits_{n=0}^{\infty} \dfrac{\Gamma (n+1/2)^{2}}{(n !)^{2}} \left[\sum\limits_{k\geq n} \dfrac{1}{(k +1)(2k +1)}\right](1-z)^{n} \ ,
\end{multline*}
\begin{align*}\tag{\ref{f1preresinvlamhyp}}
    &   \he (1-z) = 1 + \dfrac{1}{\pi} \sum\limits_{n = 1}^{\infty}
\dfrac{\Gamma(n+1/2)^{2}}{(n !)^{2}} (1-z)^{n}  \ , \quad  1-z \in \Bb{D} \, ,
\end{align*}

\noindent yield $(1-z)\het (z)\to 0$ as $z \to 1$, $z \not\in [1, 2]$, while the differentiation of them gives
\begin{multline*}
         \he^{\,\prime} (z)    = \dfrac{1 }{\pi}\dfrac{\he (1-z)}{1-z}  - \dfrac{1 }{\pi}\he^{\,\prime} (1-z) \Log \dfrac{1}{1-z}
       \\       + \dfrac{2}{\pi^{2}} \sum\limits_{n=1}^{\infty} \dfrac{n \Gamma (n+1/2)^{2}}{(n !)^{2}} \left[\sum\limits_{k\geq n} \dfrac{1}{(k +1)(2k +1)}\right](1-z)^{n-1} \ ,
\end{multline*}
\begin{align*}
    &   \he^{\,\prime} (1-z) =- \dfrac{1}{\pi} \sum\limits_{n = 1}^{\infty}
\dfrac{n \Gamma(n+1/2)^{2}}{(n !)^{2}} (1-z)^{n-1}  \ , \quad  1-z \in \Bb{D} \, ,
\end{align*}

\noindent which means that  $(1-z)\het^{\,\prime} (z)\to 1/\pi$, as $z \to 1$, $z \not\in [1, 2]$.
}}\end{clash}

\begin{subequations}

\begin{clash}{\rm{\hypertarget{d7}{}\label{case7}\hspace{0,00cm}{\hyperlink{bd7}{$\uparrow$}} \ More precisely, for arbitrary $1 < x < +\infty$ we have
\begin{align}\label{ad7}
     \int\limits_{0}^{1/x}\dfrac{d t}{\sqrt{t(1-t)(1-tx)  }} &  = \dfrac{1}{\sqrt{x}}
 \int\limits_{0}^{1/x}\dfrac{d t}{\sqrt{t\left(1-t\right)\left(\vphantom{A^{A}}(1/x) - t\right)}} \\ &
 = \dfrac{1}{\sqrt{x}}   \int\limits_{0}^{1}\dfrac{d t}{\sqrt{t\left(1-t\right)\left(\vphantom{A^{A}}1- (1/x)\, t\right)}} =
 \dfrac{\pi}{\sqrt{x}} \ \he \left(\dfrac{1}{x}\right) , \
\nonumber \end{align}

%\noindent and
\begin{align}\label{ad8}
      \int\limits_{1/x}^{1} \dfrac{d t}{\sqrt{t(1-t)(tx - 1)  }}&
=  \dfrac{1}{\sqrt{x}} \int\limits_{1/x}^{1}\dfrac{d t}{\sqrt{t \left(1-t\right)\left(t - \dfrac{1}{x}\right)  }}
\\ & =  \dfrac{1}{\sqrt{x}} \int\limits_{0}^{1-1/x} \dfrac{d t}{\sqrt{t \left(1-t\right)\left(1 - \dfrac{1}{x}-t\right)  }}
\nonumber  \\ &  =  \dfrac{1}{\sqrt{x}}   \int\limits_{0}^{1}
 \dfrac{d t}{\sqrt{t \left(1 - \dfrac{x-1}{x}\, t\right) (1-t) }} =  \dfrac{\pi}{\sqrt{x}} \ \he\left( 1- \dfrac{1}{x}\right)  \ .
\nonumber \end{align}

\noindent Therefore the "radial" limits $\het (x\pm  \imag   0):= \lim_{\varepsilon \downarrow 0} \  \het (x\pm \imag  \varepsilon)$ exist for every $x>1$ and it can be calculated as follows
\begin{align*}
    &   {\tag{\ref{f12apreresinvlamhyp}}}
\he (x\!\pm\!  \imag   0)\!=\!  \dfrac{1}{\sqrt{x}} \, \he \left(\dfrac{1}{x}\right)\! \pm \!\dfrac{ \imag  }{\sqrt{x}} \, \he\left( 1\!- \!\dfrac{1}{x}\right)  \,,
\quad x > 1 \,.
\end{align*}

\noindent But actually, the integral representation \eqref{f2preresinvlamhyp} allows to prove that for every $x>1$
all non-tangential limits exist  (see \cite[p. 11, Remark 1]{koo}), i.e.,
\begin{align}\label{f1ad7}
    &  \lim\limits_{\substack{{\fo{ \Bb{H} \ni z \to 0}} \\[0,05cm]  {\fo{ \im \, z \geq \delta  |\re\, z|}} }}
\he (x \pm z) = \he (x \pm \imag  0)\, ,
\quad \delta > 0 \,, \ x > 1 \,.
\end{align}

 We prove  the validity of \eqref{f1ad7}. First we observe that for arbitrary $x > 1$ and $t\in [0,1]$
 there exists more general limit
\begin{align}\label{f2ad7}
    &   \lim\limits_{{\fo{\Bb{H} \ni z \to 0}}} \sqrt{1-t(x \pm z)} =
\left\{\begin{array}{ll} \sqrt{1-t x } \ , \   & \hbox{if} \ \  0 \leq t < 1/x \ ;  \\[0.1cm]
  e^{{\fo{\mp i \pi/2}}} \sqrt{t x-1}  \ , \   & \hbox{if} \ \  1/x < t \leq 1 \ , \end{array}\right.
\end{align}

\noindent
 than that of \eqref{f10preresinvlamhyp},
\begin{align}
    &  {\tag{\ref{f10preresinvlamhyp}}}
 \lim\limits_{{\fo{\varepsilon \downarrow 0}}} \sqrt{1-t(x \pm \imag   \varepsilon)} =
\left\{\begin{array}{ll} \sqrt{1-t x } \ , \   & \hbox{if} \ \  0 \leq t < 1/x \ ;  \\[0.1cm]
  e^{{\fo{\mp i \pi/2}}} \sqrt{t x-1}  \ , \   & \hbox{if} \ \  1/x < t \leq 1 \ . \end{array}\right.
\end{align}

\noindent But to apply the Lebesgue dominated convergence theorem \cite[p. 26, 1.34]{rud} in order to get
for arbitrary $\delta > 0$ and $x > 1$ the required relation  \eqref{f1ad7},
\begin{align*}\nonumber
    & \lim\limits_{\substack{{\fo{ \Bb{H} \ni z \to 0}} \\[0,05cm]  {\fo{ \im  \, z\geq \delta  |\re  \, z|}} }} \he (x \pm z) = \dfrac{1}{\pi }  \int_{0}^{1}\left(
   \lim\limits_{\substack{{\fo{ \Bb{H} \ni z \to 0}} \\[0,05cm]  {\fo{ \im \, z \geq \delta  |\re  \, z|}} }}\dfrac{1}{\sqrt{1-t(x \pm  z)}}\right) \dfrac{ d t}{\sqrt{t(1-t)  }}   \\    &
      \stackrel{{\fo{\eqref{f2ad7} }}}{\vphantom{A}=}
   \dfrac{1}{\pi }  \int_{0}^{1}\left(
    \lim\limits_{{\fo{\Bb{H} \ni z \to 0}}}\dfrac{1}{\sqrt{1-t(x  \pm  z)}}\right) \dfrac{ d t}{\sqrt{t(1-t)  }}
   \stackrel{{\fo{\eqref{f2ad7}, \eqref{f11preresinvlamhyp}}}}{\vphantom{A}=}   \he (x \pm \imag  0)\, ,
   \end{align*}

\noindent it is sufficient to satisfy the condition of this theorem about the existence for every $\delta >0$ and $x > 1$  of a function  $g_{\delta, x}(t)
\in L_{1} \left([0,1], d t\right)$ such that
\begin{align}\label{f3ad7}
    &  \left|\dfrac{ 1}{\sqrt{t(1-t)\big(1-t(x  \pm  z)\big)  }}\right| \leq  g_{\delta, x}(t)  \ , \quad
    \begin{array}{l}
       t \in (0,1) \,,  \  z \in \Bb{C}\setminus \Bb{R} \,,\\[0,16cm]
     \im \, z\geq \delta  |\re  \, z| \,.
    \end{array}
\end{align}

\noindent We state that the function
\begin{align}\label{f4ad7}
    & g_{\delta, x}(t) :=  \dfrac{\sqrt{ 1 + \dfrac{1}{\delta}}}{\sqrt{t(1-t)|1-t x|  }}\in L_{1} \left([0,1], d t\right)  \ , \quad
    \delta >0 \,, \ x > 1 \,,
\end{align}

\noindent  satisfies \eqref{f3ad7}.  To prove this, one auxiliary lemma is needed.}}\end{clash}

\vspace{-0,2cm}
\begin{lemma}\hspace{-0,18cm}{\bf{.}}\label{ad7lem1}
     The following identity holds
    \begin{align}\label{f1ad7lem1}
        &  |z|^{2}\,  |x-  z |^{2}  =  x^{2}\,\im^{2}  \, z + \left(\, |z|^{2}      -  x\, \re  \, z\, \right)^{2}
     \ , \quad  x\in \Bb{R}\,, \ z \in \Bb{C} \,.
    \end{align}
\end{lemma}

\noindent
{\emph{Proof of Lemma~\ref{ad7lem1}.}} \ \
  For  $a, b, x \in \Bb{R}$ and $z = a + i b$ we have
 \begin{align*}
    |z|^{2}  |x-  z |^{2}   & = \\ & =   \left(a^{2}+ b^{2}\right)  \left((x-a)^{2}+ b^{2}\right)  = \\ & =
      \left(a^{2}+ b^{2}\right)  \left(a^{2}+ b^{2} + x^{2} - 2a x\right)  = \\ & =
        \left(a^{2}+ b^{2}\right) ^{2} - 2 a x  \left(a^{2}+ b^{2}\right) + x^{2}  \left(a^{2}+ b^{2}\right) = \\ & =
  \left(a^{2}+ b^{2}\right) ^{2} - 2 a x  \left(a^{2}+ b^{2}\right) +  a^{2} x^{2} +    b^{2} x^{2}    = \\ & =
  \left(a^{2}+ b^{2}-  a x \right)^{2}+    b^{2} x^{2}  = \\ & =
  x^{2}\im^{2}  \, z +
 \left( |z|^{2}      -  x \re (z)\right)^{2}\,,
\end{align*}

\noindent which proves the required identity \eqref{f1ad7lem1}. $\square$

\vspace{0,5cm} It follows from \eqref{f1ad7lem1} that
\begin{align}\label{f5ad7}
    &  |z|  \, |x-  z |  \geq  |x| \,\left|\im \, z \right|
     \ , \quad  x\in \Bb{R}\,, \ z \in \Bb{C}\setminus \Bb{R}\,,
\end{align}

\noindent where for $t \in \Bb{R}\setminus\{0\} $ we can set $1-xt$ instead of $x$ and $z t$ instead of $z$ to get
\begin{align*}
    &  | t z|  \, |1-xt- t z |  \geq  |1-xt| \,\left|\im (t z) \right| = |t|  \,|1-xt|  \,\left|\im  \, z \right|  \ , \
\end{align*}

\noindent and hence,  dividing both sides by $|t|$, we see that
\begin{align}\label{f6ad7}
    &   \big|1-t\left(x+z\right) \big| \geq  \dfrac{\left|\im \, z \right|}{|z|}  |1-xt|  \ , \quad  x, t \in \Bb{R}\,, \ z \in \Bb{C}\setminus \Bb{R} \,,
 \end{align}

\noindent because for $t=0$ this inequality is obviously true. Then
\begin{align}\label{f7ad7}
    &   \dfrac{1}{ \big|1-t\left(x\pm z\right) \big| }   \leq  \dfrac{|z|}{\left|\im  \, z \right|}\cdot \dfrac{1}{|1-xt|}   \ , \quad  x, t \in \Bb{R}\,, \ z \in \Bb{C}\setminus \Bb{R}  \,,
\end{align}

\noindent and if $\left|\im  \, z \right|\geq \delta  |\re  \, z|$ we get
\begin{align}\label{f8ad7}
    &   \dfrac{|z|}{\left|\im  \, z \right|} \leq \dfrac{\left|\re  \, z \right| + \left|\im  \, z \right|}{\left|\im  \, z \right|}\leq 1 + \dfrac{1}{\delta} \ , \
\end{align}

\noindent from which and \eqref{f7ad7} we deduce that
\begin{multline*}
    \left|\dfrac{ 1}{\sqrt{t(1-t)\big(1-t(x \pm z)\big)  }}\right|  =  \dfrac{ 1}{\sqrt{t(1-t)\big|1-t(x \pm z)\big|  }} \\          \stackrel{{\fo{\eqref{f7ad7}}}}{\vphantom{A}\leq}
   \dfrac{\sqrt{|z|}}{\sqrt{\left|\im  \, z \right|}}  \dfrac{ 1}{\sqrt{t(1-t)|1-t x|  }}
  \stackrel{{\fo{\eqref{f8ad7}}}}{\vphantom{A}\leq}
  \dfrac{\sqrt{ 1 + \dfrac{1}{\delta}}}{\sqrt{t(1-t)|1-t x|  }}     \ , \
 \end{multline*}

\noindent for arbitrary $x > 1$, $t \in (0,1)$, $\delta >0$ and $ z \in \Bb{C}\setminus \Bb{R}$ satisfying $|\im  \, z| \geq \delta  |\re  \, z|$.
This proves that the function $g_{\delta, x}$ defined as in \eqref{f4ad7} satisfies \eqref{f3ad7} and therefore  the proof of \eqref{f1ad7} is completed.
\end{subequations}

\begin{subequations}

\begin{clash}{\rm{\hypertarget{d26}{}\label{case26}\hspace{0,00cm}{\hyperlink{bd26}{$\uparrow$}} \
We obtain \eqref{f2rem1},
\begin{align*}\tag{\ref{f2rem1}}
    &  (-z)^{-1/2}\he  (1/z) - \imag
    \he  (z)\, {\rm{sign}} ({\rm{Im}}\, z) =z^{-1/2} \he  (1\!-\!1/z) \ ,
    \ \   z\!\in \! \Bb{C}\setminus\Bb{R} \,,
\end{align*}

\noindent  based on the functional  relation \eqref{f15case27},
\begin{align}\tag{\ref{f15case27}}
    &  \ie( z) - {\rm{sign}} ({\rm{Im}}\, z) =\ie\left(\dfrac{z}{z-1}\right)  \ , \quad
    z \in \Bb{C}\setminus \Bb{R} \,,
\end{align}

\noindent
 established below.

 \vspace{0.15cm}
 Observe that the changes of variable  $z$ by $1-z$ and  $z$ by $1/z$ in the Pfaff formula \eqref{f5preresinvlamhyp},
 \begin{align}\label{f0case26}
    &  \he  (z) = \dfrac{1}{\sqrt{1-z}} \ \he
\left(\dfrac{z}{z-1}\right)   \ , \quad
 z \in \Bb{C}\setminus [1, +\infty)  \ , \
\end{align}

\noindent imply
\begin{align}\label{f1case26}
     &  \he (1-z)    = \dfrac{1}{\sqrt{z}} \ \he
\left(1- \dfrac{1}{z}\right)   \ , \quad
 z \in \Bb{C}\setminus (-\infty, 0]  \,,
\end{align}

\noindent and
\begin{align}\label{f2case26}
    &   \ \he \left( \dfrac{1}{z}\right)   =
\dfrac{1}{\sqrt{1-\dfrac{1}{z}}} \  \he \left(\dfrac{1}{1 - z}\right)
  \ , \quad      z \in \Bb{C}\setminus [ 0, 1] \, .
\end{align}

\noindent But for arbitrary $ z \in \Bb{C}\setminus \Bb{R} = \Bb{H} \cup (- \Bb{H})$ two numbers $1-z$ and $-z$ belong to the same half-plane either to $\Bb{H}$ or to $-\Bb{H}$ which yields that
\begin{align*}
    &  \sqrt{1-\dfrac{1}{z}} = \sqrt{\dfrac{1-z}{-z}} = \dfrac{\sqrt{1-z}}{\sqrt{-z}}  \ , \quad  z \in \Bb{C}\setminus \Bb{R} \,,
\end{align*}

\noindent and therefore \eqref{f2case26} can be written in the following form
\begin{align}\label{f3case26}
    &   \ \dfrac{\he \left( \dfrac{1}{z}\right)}{\sqrt{-z}}   =
 \dfrac{\he \left(\dfrac{1}{1 - z}\right)}{\sqrt{1-z}}
  \ , \quad      z \in \Bb{C}\setminus \Bb{R} \, .
\end{align}

Applying to the
 functional  relation \eqref{f15case27},
\begin{align}\tag{\ref{f15case27}}
    &  \ie( z) - \sigma(z) =\ie\left(\dfrac{z}{z-1}\right)  \ , \quad  \sigma(z) :={\rm{sign}} ({\rm{Im}}\, z)
     \ , \  z \in \Bb{C}\setminus \Bb{R} \,,
\end{align}

\noindent
 the definition \eqref{f10int},
\begin{align}\tag{\ref{f10int}}
       \ie (z) :=   \imag \cdot \frac{\he (1-z)}{\he (z)},
\quad  z \in  (0,1)\cup \left(\Bb{C}\setminus\Bb{R}\right),
\end{align}

\noindent of the function $\ie $ we  conclude that
\begin{align*}
    &  i \cdot \frac{\he\left(1-z\right)}{\he\left(z\right)}- \sigma (z) =
i \cdot \frac{\he\left(1-\dfrac{z}{z-1}\right)}{\he\left(\dfrac{z}{z-1}\right)}
  &     &    \Longrightarrow \\[0,3cm] &
\frac{i \cdot \he\left(1-z\right) - \sigma (z)\he\left(z\right) }{\he\left(z\right)} =
\frac{i \cdot\he\left(1-\dfrac{z}{z-1}\right)}{\he\left(\dfrac{z}{z-1}\right)}  &     &     \stackrel{{\nor{\eqref{f0case26}}}}{\vphantom{A}  \Longrightarrow  }
\\[0,3cm] &
 i \cdot \he\left(1-z\right) - \sigma (z)\he\left(z\right)=\frac{i \cdot\he\left(\dfrac{1}{1-z}\right)}{\sqrt{1-z}}  &     &    \Longrightarrow \\[0,3cm] &
 \sigma (z)\he\left(z\right)= i \cdot \he\left(1-z\right)-\frac{i \cdot\he\left(\dfrac{1}{1-z}\right)}{\sqrt{1-z}}
  &     &      \stackrel{{\nor{\eqref{f3case26}}}}{\vphantom{A} \Longrightarrow}   \\[0,3cm] &
 \sigma (z)\he\left(z\right)= i \cdot  \he\left(1-z\right)-
i \cdot  \dfrac{\he \left(\dfrac{1}{z}\right)}{\sqrt{-z}} \
  &     &      \stackrel{{\nor{\eqref{f1case26}}}}{\vphantom{A}\Longrightarrow}
   \\[0,3cm] &
  \sigma (z) \he\left(z\right)= i \cdot
 \dfrac{\he\left(1- \dfrac{1}{z}\right)}{\sqrt{z}}
-i \cdot \dfrac{\he \left(\dfrac{1}{z}\right)}{\sqrt{-z}}
   &     &    \Longrightarrow
    \\[0,3cm] &
    \dfrac{\he \left(\dfrac{1}{z}\right)}{\sqrt{-z}}
-i \cdot \sigma (z)  \he\left(z\right)=
 \dfrac{\he\left(1- \dfrac{1}{z}\right)}{\sqrt{z}} \ , \     &     &    z \in \Bb{C}\setminus \Bb{R} \,.
\end{align*}

\vspace{0.2cm}
\noindent This completes the proof of \eqref{f2rem1}.

}}\end{clash}
\end{subequations}

\begin{subequations}

\begin{clash}{\rm{\hypertarget{d8}{}\label{case8}\hspace{0,00cm}{\hyperlink{bd8}{$\uparrow$}} \ For $z = x + i y$, $ x, y \in \Bb{R}$,
the inequality
\begin{align*}
    &  \left|1 - z t\right|^{2}    =  (1- xt)^{2} + y^{2}t^{2}\geq (1- xt)^{2} \ , \   t \in [0,1] \,,
    \end{align*}

\noindent being applied to the integral representation \eqref{f2preresinvlamhyp},
\begin{align*}
    &  \left|\he (z)\right| \leq  \dfrac{1}{\pi }  \int_{0}^{1} \dfrac{ \diff  t}{\sqrt{t(1-t)|1-tz|  }} \leq
    \dfrac{1}{\pi }  \int_{0}^{1} \dfrac{ \diff  t}{\sqrt{t(1-t)|1-t x|  }} \ ,
\end{align*}

\noindent in view of \eqref{f18preresinvlamhyp}, yield that for any $y \in \Bb{R}$ we have
\begin{align}\label{ad9}
    &  \left|\he (x + i y)\right| \leq  \he (x)=
    \left\{\begin{array}{ll} \he (x) \, ,    & \ \ \hbox{if} \, \ 0\! \leq\! x\! <\! 1  ; \\[0,0cm]
\begin{displaystyle}
\dfrac{1}{\sqrt{1+|x|}} \  \he \left(\dfrac{|x|}{1 + |x|}\right) \, ,
\end{displaystyle} & \ \  \hbox{if} \ \  x < 0 \ ,
    \end{array}\right.
\end{align}

\noindent and, by virtue of \eqref{ad7} and  \eqref{ad8},
\begin{align}\label{ad10}
    &  \left|\he (x + i y)\right| \leq    \dfrac{\he \left(\dfrac{1}{x}\right) + \he \left(1-\dfrac{1}{x}\right) }{\sqrt{x}}    \ , \quad  x > 1 \, , \  y \in \Bb{R}\setminus\{0\} \,.
\end{align}

\noindent According to \cite[p. 557, 15.1.26]{abr},
\begin{align*}
    &  F \left(a, 1-a; b; \dfrac{1}{2}\right)=2^{1-b} \sqrt{\pi} \dfrac{\Gamma(b)}{\Gamma\left(\dfrac{a+b}{2}\right)\Gamma\left(\dfrac{b+1-a}{2}\right)} \ , \quad b \not\in \Bb{Z}_{\leq 0} \ ,
\end{align*}

\noindent $\Bb{Z}_{\leq 0}:= \{n \in \Bb{Z} \,|\,  n \leq 0\}$, which means that
\begin{align*}
    &  F \left(\dfrac{1}{2},\dfrac{1}{2};1;\dfrac{1}{2}\right) =\sqrt{\pi} \dfrac{\Gamma(1)}{\Gamma\left(\dfrac{3}{4}\right)\Gamma\left(\dfrac{3}{4}\right)}=
\dfrac{ \sqrt{ \pi} }{\Gamma\left(\dfrac{3}{4}\right)^{2}} \ , \
\end{align*}

\noindent where in view of $\Gamma( 3/4)= 1,22541...$ (see \cite[p. 255, 6.1.14]{abr}), we have
\begin{align}\label{ad11e}
    &  F \left(\dfrac{1}{2},\dfrac{1}{2};1;\dfrac{1}{2}\right) <1,18036 < 2 \ .
\end{align}

\noindent Setting $\theta = 1/2$ we obtain from \eqref{f15preresinvlamhyp} and \eqref{f16preresinvlamhyp}
\begin{align*}
      \he (x) &\leq 1 + 2 x \leq 2 \ , \  x \in [0,1/2] \ , \  \\
     \he (x)& \leq \dfrac{1}{\pi} \ln \dfrac{16}{1-x} + \dfrac{2}{\pi} (1-x)\ln \dfrac{51}{1-x}
     \leq  \dfrac{1}{\pi} \ln \dfrac{16}{1-x} + \dfrac{1}{\pi} \ln \dfrac{51}{1-x}  \\ & =
      \dfrac{2}{\pi}\ln \dfrac{1}{1-x} + \dfrac{\ln 16 + \ln 51}{\pi} < 3 + \ln \dfrac{1}{1-x}
      \ , \  x \in [1/2, 1) \ ,
\end{align*}

\noindent from which (cf. \cite[p. 494, 19.9(i)]{olv})
\begin{align}\label{ad11}
    &   \he (x) \leq 3 + \ln \dfrac{1}{1-x} \ , \quad   x \in [0, 1) \ .
\end{align}

\noindent Combining this with \eqref{ad9} and \eqref{ad10} we get
\begin{align}\label{ad12}
    &  \left|\he (x + i y)\right| \leq
    \left\{\begin{array}{lll}
      \dfrac{6 + \ln \dfrac{x^{2}}{x-1}}{\sqrt{x}} \, ,    & \ \ \hbox{if} \, \ x > 1  \ , & \ y \in \Bb{R}\setminus\{0\}   ; \\[0,3cm]
    3 + \ln \dfrac{1}{1-x} \, ,    & \ \ \hbox{if} \, \ 0\! \leq\! x\! <\! 1 \ , & \ y \in \Bb{R} ; \\[0,4cm]
\begin{displaystyle}
\dfrac{ 3 + \ln \left(1 + |x|\right)}{\sqrt{1+|x|}}  \, ,
\end{displaystyle} & \ \  \hbox{if} \ \  x < 0 \ , &\  y \in \Bb{R}\,.
    \end{array}\right.
\end{align}

\noindent These inequalities yield that $\het$ belongs to the Hardy space $ H^{p}  $  for
arbitrary $2 < p < \infty$ (see \eqref{inthp}).

}}\end{clash}
\end{subequations}

\begin{subequations}

\begin{clash}{\rm{\hypertarget{d9}{}\label{case9}\hspace{0,00cm}{\hyperlink{bd9}{$\uparrow$}} \ We prove the validity of \eqref{f24xpreresinvlamhyp},
 \begin{align*}{\tag{\ref{f24xpreresinvlamhyp}}}
    & \he (z)=  e^{{\fo{- \dfrac{i\pi\sigma}{2} }}}\int\limits_{0}^{\infty}
    \dfrac{\he\left(\dfrac{t}{t+i\sigma} \right) d  t}{(1 - i t \sigma -z)\sqrt{1- i t \sigma   }} \ , \quad
    z \in \sigma\cdot\Bb{H} \ , \  \sigma\in\{1, -1\} \,.
 \end{align*}

\vspace{0,1cm}
For arbitrary $z_{1}, z_{2} \in \Bb{C}$,
$z_{1}\neq z_{2}$ let $[z_{1}, z_{2}]$ be the straight line segment  from $z_{1}$ to
$z_{2}$.

\vspace{0,1cm}
For given $r > 1$ we denote by $\gamma_{r}^{\pm}$  the contour which passes from $1$ to $1+r$  along
 $[1, 1+r]$, from $1+r$ to $1\pm i r$  along the arc
 \begin{align}\label{f0ad9}
    & \left.\left\{ \,  1+  r e^{{\fo{is}}}\, \, \right| \,  s \in \left[ 0 , \pm \dfrac{\pi}{2} \right]\, \right\}  \, ,
 \end{align}

 \noindent
 and from $1\pm i r$ to $1$  along
 $[1\pm  i r, 1]$, respectively.

\vspace{0.25cm}
Choose  any $z \in  \Bb{H} \cup (- \Bb{H})$ and let
\begin{align*}
    & \sigma := {\rm{sign}} (\im\, z) \in \{1, -1\} \ , \quad
     \Bb{C}_{\re > 0} := \{\, z \in \Bb{C} \,|\, \re\,z> 0\,\} \,.
\end{align*}

\noindent
It follows from $\he \in{\rm{Hol}} \left(\Bb{C}\!\setminus\![1,+\infty)\right)$ that
\begin{align*}
    &  \he\left( 1-\dfrac{1}{t} \right) \in{\rm{Hol}} \Big(\Bb{C}\!\setminus\!(-\infty, 0]\Big) \ ,
\end{align*}

\noindent and therefore\vspace{-0,2cm}
\begin{align*}
    &   \dfrac{ \he\left( 1-\dfrac{1}{t} \right)}{(t-z)\sqrt{t}} \in{\rm{Hol}}
    \Big(\Bb{C}_{\re > 0} \cap \big( - \sigma \Bb{H} + \im\, z  \big)\Big) \ ,
\end{align*}

\noindent where obviously\vspace{-0,2cm}
\begin{align*}
    &   \gamma_{{\fo{r}}}^{{\fo{-\sigma}}} \subset
    \Bb{C}_{\re > 0} \cap \big(   - \sigma \Bb{H} + \im\, z \big) \ , \  r > 1  \,.
\end{align*}

\noindent Applying  the Cauchy theorem (see \cite[p. 89, 6.6]{con}) we get
\begin{align*}
  \int\limits_{{\nor{\gamma_{{\fo{r}}}^{{\fo{-\sigma}}}}}}  \dfrac{ \he\left( 1-\dfrac{1}{t} \right) }{(t-z)\sqrt{t}}d t = 0 \ , \qquad z \in \sigma \Bb{H} \,, \  \ r > 1 \ .
\end{align*}

\noindent This equality  can be written as follows
\begin{align}\label{f1ad9}
    &  \int\limits_{{\fo{1}}}^{{\fo{1+r}}} \  \dfrac{\he\left( 1-\dfrac{1}{t} \right) d  t}{(t-z)\sqrt{t}} =
 (- \sigma  i )\int\limits_{{\fo{0}}}^{{\fo{r}}}
    \dfrac{\he\left( 1-\dfrac{1}{1-\sigma i t} \right) d  t    }{(1 -\sigma i t  -z)\sqrt{1 -\sigma i t   }}  - \Delta_{{\fo{\,r}}} (z)\,, \end{align}

\noindent where
    \begin{align*}
    \Delta_{{\fo{\,r}}} (z)  &:=
    \int\limits_{0}^{{\fo{- \dfrac{\sigma\pi}{2}}}}  \he\left( 1-\dfrac{1}{1+ r e^{{\fo{i s}}} } \right)
    \dfrac{ d \left( 1+ r e^{{\fo{i s}}}  \right)    }{\left(1+ r e^{{\fo{i s}}}-z\right)\sqrt{1+ r e^{{\fo{i s}}}   }}
    \\    &    =
     \dfrac{i}{\sqrt{r}} \int\limits_{0}^{{\fo{- \dfrac{\sigma\pi}{2}}}}  \he\left( 1-\dfrac{1}{1+ r e^{{\fo{i s}}} } \right)
    \dfrac{ e^{{\fo{i s}}}  d s    }{\left( e^{{\fo{i s}}} + \dfrac{1-z}{r}\right)\sqrt{ e^{{\fo{i s}}} + \dfrac{1}{r}  }}
    \ .
\end{align*}

\noindent But for $r \geq 3 + 2 |1-z|$ and arbitrary $\alpha \in \Bb{R}$ we have
\begin{align*}
    & \left|  e^{{\fo{i \alpha}}} + \dfrac{1}{r}\right| \geq 1 - \dfrac{1}{r} \geq \dfrac{2}{3} > \dfrac{4}{9} \ , \\    &   \left| e^{{\fo{i \alpha}}} + \dfrac{1-z}{r}\right| \geq 1 - \dfrac{|1-z|}{r}\geq \dfrac{1}{2}  \ ,
\end{align*}

\noindent and using the result of Corollary~{\hyperlink{hadcor2}{\ref*{adcor2}}},
\begin{align*}{\tag{\ref{f1adcor2}}}
    &  \left| \het \left(1-\dfrac{1}{1+ r e^{{\fo{i \varphi}}} }\right)\right| \leq 3 + \dfrac{\ln r}{\pi}  \, , \quad \varphi \in \left[-\dfrac{\pi}{2}\ , \ \dfrac{\pi}{2} \right]  \,, \  r \geq 3 \,,
\end{align*}

\noindent we obtain
\begin{align*}
      \left|\Delta_{r} (z)\right| & = \left|  \dfrac{(- \sigma i)}{\sqrt{r}} \int\limits_{0}^{{\fo{ \dfrac{\pi}{2}}}}  \he\left( 1-\dfrac{1}{1+ r e^{{\fo{- \sigma i s}}} } \right)
    \dfrac{ e^{{\fo{- \sigma i s}}}  d s    }{\left( e^{{\fo{- \sigma i s}}} + \dfrac{1-z}{r}\right)\sqrt{ e^{{\fo{- \sigma i s}}} + \dfrac{1}{r}  }}\right|  \\[0,3cm]    &    \leq
    \dfrac{1}{\sqrt{r}}
     \int\limits_{0}^{{\fo{ \dfrac{\pi}{2}}}} \left| \he\left( 1-\dfrac{1}{1+ r e^{{\fo{- \sigma i s}}} } \right)\right|
    \dfrac{   d s    }{\left| e^{{\fo{- \sigma i s}}} + \dfrac{1-z}{r}\right|\sqrt{ \left|e^{{\fo{- \sigma i s}}} + \dfrac{1}{r}\right|  }} \\[0,3cm]    &    \leq
     \dfrac{3 \pi }{2 \sqrt{r}} \max\limits_{{\fo{s \in [0, \pi/2]}}} \left| \he\left( 1-\dfrac{1}{1+ r e^{{\fo{- \sigma i s}}} } \right)\right| \  \stackrel{{\fo{\eqref{f1adcor2}}}}{\vphantom{A}=} \
     \dfrac{3 \pi (3 \pi + \ln r)}{2 \pi\sqrt{r}}\ ,
\end{align*}

\vspace{0.25cm}
\noindent which proves that $\lim_{r \to +\infty} \Delta_{\,r} (z) = 0$.  Letting $r \to +\infty$ in \eqref{f1ad9} we obtain
\begin{align*}
    &  \int\limits_{1}^{+\infty} \  \dfrac{\he\left( 1-\dfrac{1}{t} \right) d  t}{(t-z)\sqrt{t}} =
 (- \sigma  i )\int\limits_{0}^{+\infty}
    \dfrac{\he\left( 1-\dfrac{1}{1-\sigma i t} \right) d  t    }{(1 -\sigma i t  -z)\sqrt{1 -\sigma i t   }} \ , \quad
    z \in \sigma\cdot\Bb{H} \ ,
\end{align*}

\noindent which coincides with \eqref{f24xpreresinvlamhyp} because
\begin{align*}
    &  - \sigma  i = e^{{\fo{- \dfrac{i\pi\sigma}{2} }}}  \ , \  1-\dfrac{1}{1-\sigma i t} = \dfrac{-\sigma i t}{1-\sigma i t} = \dfrac{t}{t + i \sigma}  \ , \quad  t > 0 \,, \ \sigma\in\{1, -1\} \,.
\end{align*}

\noindent The equality \eqref{f24xpreresinvlamhyp},
 \begin{align*}{\tag{\ref{f24xpreresinvlamhyp}}}
    & \he (z)=  e^{{\fo{- \dfrac{i\pi\sigma}{2} }}}\int\limits_{0}^{\infty}
    \dfrac{\he\left(\dfrac{t}{t+i\sigma} \right) d  t}{(1 - i t \sigma -z)\sqrt{1- i t \sigma   }} \ , \quad
    z \in \sigma\cdot\Bb{H} \ , \  \sigma\in\{1, -1\} \,,
 \end{align*}

\noindent
 is proved. $\square$

\vspace{0.25cm} Observe that for  $ z \in \Bb{H}$   we deduce from  \eqref{f10int} and \eqref{f24xpreresinvlamhyp} that
 \begin{align*}
    &  \ie (z)=i \frac{\he(1-z)}{\he(z)}=i
    \dfrac{
     i \begin{displaystyle}
     \int\limits_{0}^{\infty}
     \end{displaystyle}
    \dfrac{\he\left(\dfrac{t}{t-i} \right) d  t}{(z + i t  )\sqrt{1+ i t   }}
    }{
    \dfrac{1}{i}\begin{displaystyle}
     \int\limits_{0}^{\infty}
     \end{displaystyle}
    \dfrac{\he\left(\dfrac{t}{t+i} \right) d  t}{(1 - i t  -z)\sqrt{1- i t    }}
    }=
 i    \dfrac{
     \begin{displaystyle}
     \int\limits_{0}^{\infty}
     \end{displaystyle}
    \dfrac{\he\left(\dfrac{t}{t-i} \right) d  t}{(z + i t  )\sqrt{1+ i t   }}
    }{
    \begin{displaystyle}
     \int\limits_{0}^{\infty}
     \end{displaystyle}
    \dfrac{\he\left(\dfrac{t}{t+i} \right) d  t}{(z + i t  -1)\sqrt{1- i t    }}
    }    \ ,
 \end{align*}

 \noindent while if   $ z \in -\Bb{H}$ then\vspace{-0,1cm}
\begin{align*}
    &  \ie (z)=i \frac{\he(1-z)}{\he(z)}=i
    \dfrac{
   \dfrac{1}{i}   \begin{displaystyle}
     \int\limits_{0}^{\infty}
     \end{displaystyle}
    \dfrac{\he\left(\dfrac{t}{t+i} \right) d  t}{(z - i t  )\sqrt{1- i t   }}
    }{
    i\begin{displaystyle}
     \int\limits_{0}^{\infty}
     \end{displaystyle}
    \dfrac{\he\left(\dfrac{t}{t-i} \right) d  t}{(1 + i t  -z)\sqrt{1+ i t    }}
    }=
 i   \dfrac{
      \begin{displaystyle}
     \int\limits_{0}^{\infty}
     \end{displaystyle}
    \dfrac{\he\left(\dfrac{t}{t+i} \right) d  t}{( i t  -z)\sqrt{1- i t   }}
    }{
    \begin{displaystyle}
     \int\limits_{0}^{\infty}
     \end{displaystyle}
    \dfrac{\he\left(\dfrac{t}{t-i} \right) d  t}{(1 + i t  -z)\sqrt{1+ i t    }}
    } \ ,
 \end{align*}

\vspace{-0,1cm}
 \noindent i.e., \eqref{f24xpreresinvlamhyp} for the function $ \ie $ can be written in the following form
 \begin{align}\label{f4athd}
    & \ie (z)=i \frac{\he(1-z)}{\he(z)}=i    \dfrac{
     \begin{displaystyle}
     \int\limits_{0}^{\infty}
     \end{displaystyle}
    \dfrac{\he\left(\dfrac{t}{t-i} \right) d  t}{(z + i t  )\sqrt{1+ i t   }}
    }{
    \begin{displaystyle}
     \int\limits_{0}^{\infty}
     \end{displaystyle}
    \dfrac{\he\left(\dfrac{t}{t+i} \right) d  t}{(z + i t  -1)\sqrt{1- i t    }}
    }  \ , \quad  z\in \Bb{H} \ , \\    &
     \ie (z)=i \frac{\he(1-z)}{\he(z)}=
      i   \dfrac{
      \begin{displaystyle}
     \int\limits_{0}^{\infty}
     \end{displaystyle}
    \dfrac{\he\left(\dfrac{t}{t+i} \right) d  t}{( i t  -z)\sqrt{1- i t   }}
    }{
    \begin{displaystyle}
     \int\limits_{0}^{\infty}
     \end{displaystyle}
    \dfrac{\he\left(\dfrac{t}{t-i} \right) d  t}{(1 + i t  -z)\sqrt{1+ i t    }}
    } \ , \quad  z\in -\Bb{H} \,.
\label{f4bthd} \end{align}
}}\end{clash}

 \vspace{0.25cm}

\end{subequations}
%%%%%%%%%%%%%%%%%%%%%%%%%%%%%%%%%%%%%%%%%%%%%%%%%%%%%%%%%%%%%%%%%%%%%%%%%%%%%%%%%%%%%%%%%%%%%%%%%%%%%%%%%

\subsection[\hspace{-0,25cm}. \hspace{0,075cm}Notes on Section~\ref{eir}]{\hspace{-0,11cm}{\bf{.}} Notes on Section~\ref{eir}}

\begin{subequations}
\vspace{-0.15cm}
\begin{clash}{\rm{\hypertarget{d15}{}\label{case15}\hspace{0,00cm}{\hyperlink{bd15}{$\uparrow$}} \ Actually,
for arbitrary  $ z\! \in \! (0,1)\cup \left(\Bb{C}\setminus\Bb{R}\right)$ we have
$z\!\in\!\Bb{C}\!\setminus\! [1, +\!\infty)$ and $1-z\!\in\!\Bb{C}\!\setminus\! [1, +\!\infty)$ and therefore
 the representations \eqref{f8eir} and \eqref{f9eir} can be written in the form ($\het := F_{1/2,1/2;1}$)
 \begin{align}\label{add1}
    &
    \begin{array}{rl}
  \Log\het (z)  &    =  \begin{displaystyle}\frac{1}{\pi^{2 }} \!\int\limits_{0}^{1}
\frac{\dfrac{1}{t(1-t)}\Log\left(\dfrac{1}{1-tz}\right)\, \diff t  }{\vphantom{\dfrac{a}{a}}\het (t)^{2} +
\het (1-t)^{2}} \ ,
 \end{displaystyle}
   \\[0,7cm]
  \Log\het (1-z)  &    =  \begin{displaystyle}\frac{1}{\pi^{2 }} \!\int\limits_{0}^{1}
\frac{\dfrac{1}{t(1-t)}\Log\left(\dfrac{1}{1-t(1-z)}\right)\, \diff t  }{\vphantom{\dfrac{a}{a}}\het (t)^{2} +
\het (1-t)^{2}} \ ,
 \end{displaystyle}
    \end{array}
  \end{align}

\vspace{0.25cm} \noindent and
\begin{align}\label{add2}
  \begin{array}{l}
  \begin{displaystyle}
\Log \het (z) \!=\!\dfrac{1}{\pi} \int\limits_{0}^{1}  \dfrac{t}{1+t^2}  \arctan \dfrac{\het \left(1-t\right)}{\het
     \left(t\right)}  \diff  t
  \end{displaystyle} \\[0,4cm]
  \begin{displaystyle}
 +\! \dfrac{1}{\pi}\! \int\limits_{1}^{+\infty}\!\!
   \left( \frac{1}{t\!-\!z}\! - \!\frac{t}{1\!+\!t^2} \right)\,
   \arctan \dfrac{\het \left(1\!-\!{1}/{t}\right)}{\het \left({1}/{t}\right)}\, \diff t\,,
  \end{displaystyle}\\[0,7cm]
 \begin{displaystyle}
\Log \het (1-z) \!=\!\dfrac{1}{\pi} \int\limits_{0}^{1}  \dfrac{t}{1+t^2}  \arctan \dfrac{\het \left(1-t\right)}{\het
     \left(t\right)}  \diff  t
  \end{displaystyle} \\[0,4cm]
  \begin{displaystyle}
 +\! \dfrac{1}{\pi}\! \int\limits_{1}^{+\infty}\!\!
   \left( \frac{1}{t\!-\!(1-z)}\! - \!\frac{t}{1\!+\!t^2} \right)\,
   \arctan \dfrac{\het \left(1\!-\!{1}/{t}\right)}{\het \left({1}/{t}\right)}\, \diff t \ .
  \end{displaystyle}
    \end{array}
\end{align}

\vspace{0.25cm} \noindent
Subtracting form the second equality the first one in \eqref{add1} and in \eqref{add2} and using
\begin{align*}\tag{\ref{f11aeir}}
    & \Log \he (1\!-\!z)\! - \!  \Log \he (z) \! = \! \Log \dfrac{ \he (1\!-\! z)}{\he (z)} \ , \quad
     z\! \in \! (0,1)\cup \left(\Bb{C}\setminus\Bb{R}\right),
\end{align*}

\noindent and
\begin{align*}\tag{\ref{f14aeir}}  \hypertarget{addf14aeir}{}
    &  \Log\frac{1}{1-tz}-\Log\frac{1}{1-t + tz} = \Log\dfrac{1-tz}{1-t+tz}
    \, , \quad     z\! \in \! (0,1)\cup \left(\Bb{C}\setminus\Bb{R}\right),
\end{align*}

\vspace{0.25cm}
\noindent where $t \in (0,1)$, we obtain
\begin{align*}
     \Log \dfrac{ \he (1\!-\! z)}{\he (z)} & = \! \Log \he (1\!-\!z)\! - \!  \Log \he (z)   \\ & =
   \frac{1}{\pi^{2 }} \!\int\limits_{0}^{1}
\frac{\dfrac{1}{t(1-t)}\left(\Log\left(\dfrac{1}{1-t +tz}\right) - \Log\left(\dfrac{1}{1-tz}\right) \right)\, \diff t  }{\vphantom{\dfrac{a}{a}}\het (t)^{2} +
\het (1-t)^{2}}    \\ & =
\dfrac{1}{\pi^{2 }}\! \int\limits_{0}^{1}\!\! \dfrac{ \Log\dfrac{1-tz}{1-t+tz}}{t (1-t)\left(\het \left(t \right)^{2} + \het \left(1-t \right)^{2}\right)} \,  \diff  t \,,
\end{align*}

\vspace{0.25cm} \noindent and
\begin{align*}
      \Log \dfrac{ \he (1\!-\! z)}{\he (z)}  & = \! \Log \he (1\!-\!z)\! - \!  \Log \he (z)   \\ & =
    \dfrac{1}{\pi}\! \int\limits_{1}^{+\infty}\!\!
   \left( \frac{1}{t\!-\!(1-z)}\! - \!\frac{1}{t\!-\!z} \right)\,
   \arctan \dfrac{\het \left(1\!-\!{1}/{t}\right)}{\het \left({1}/{t}\right)}\, \diff t   \\ & =
    \dfrac{1}{\pi}\! \int\limits_{1}^{+\infty}\!\!
    \frac{t\!-\!z - \left(t\!-\!1\!+\! z\right)}{\left(t\!-\!1\!+\! z\right)\left(t\!-\!z\right)} \,
   \arctan \dfrac{\het \left(1\!-\!{1}/{t}\right)}{\het \left({1}/{t}\right)}\, \diff t
     \\ & =  \dfrac{1-2z}{\pi} \int\limits_{0}^{1} \ \,  \dfrac{\arctan \dfrac{\het \left(1-t\right)}{\het \left(t\right)}}{(1-tz)(1-t +t z)} \, \diff  t\,,
\end{align*}

\vspace{0.25cm} \noindent which prove the validity of
\begin{align*}\tag{\ref{f12eir}} \hypertarget{addf12eir}{}
    &  \Log \dfrac{ \he (1\!-\! z)}{\he (z)}\!=\!\dfrac{1}{\pi^{2 }}\! \int\limits_{0}^{1}\!\! \dfrac{ \Log\dfrac{1-tz}{1-t+tz}}{t (1-t)\left(\he \left(t \right)^{2} + \he \left(1-t \right)^{2}\right)} \,  \diff  t \ , \\[0,0cm]    &\tag{\ref{f13eir}}
    \Log \dfrac{ \he (1\!-\! z)}{\he (z)}=\dfrac{1-2z}{\pi} \int\limits_{0}^{1} \ \,  \dfrac{\arctan \dfrac{\het \left(1-t\right)}{\het \left(t\right)}}{(1-tz)(1-t +t z)} \,  \diff  t \,,
 \hypertarget{addff13eir}{}\end{align*}

\vspace{0.25cm} \noindent for arbitrary $ z\! \in \! (0,1)\cup \left(\Bb{C}\setminus\Bb{R}\right)$.

\vspace{0.25cm}
It follows from
\begin{align*}
    &  \tag{\ref{f4theor2}}\mu \big((-\infty, 0)\big)  =0 \, , \quad
    \mu\big(\,[0,\,x)\big)   := \nu(+\infty) - \nu(1/x) \ ,  \quad   x > 0 \,,
\end{align*}

\noindent Theorem~\ref{theor2}, according to which $ \nu(1/x)= 0$, $x > 1$, and
\begin{align*}\tag{\ref{f5eir}}
    &  \mu\big([0,\,x)\big)\!=\!\nu(+\infty)\! -\! \nu(1/x)\! =\!
    \dfrac{1}{2}\! - \!\dfrac{1}{\pi}\arctan \dfrac{\he \left(1\!-\!x\right)}{\he \left(x\right)}  \, , \  x \in (0, 1);   \\   &
    \mu\big(\,[0,\,x)\big)= \dfrac{1}{2}  \ , \quad   x \geq 1  \ ; \quad  \quad \mu\big(\,\{0\}\big) = 0  \ ,
 \tag{\ref{f6eir}} \\    &
 \tag{\ref{f7eir}}
    \dfrac{ \diff  \mu(x)}{ \diff  x} =
    \dfrac{1}{\pi^{2 }  x (1-x)} \dfrac{1}{\he \left(x\right)^{2} + \he \left(1-x\right)^{2}}
         \ , \ \   x \in (0,1) ,
\end{align*}

\noindent that
\begin{align}\label{add3}
    &  \mu \left(\Bb{R}\right) = \mu \left([0,1]\right)=
    \dfrac{1}{\pi^{2 } } \int\limits_{0}^{1} \dfrac{1}{ x (1-x)} \dfrac{d x}{\he \left(x\right)^{2} + \he \left(1-x\right)^{2}} = \dfrac{1}{2}\ .
\end{align}

\noindent In view of ({\hyperlink{addf14aeir}{\ref*{f14aeir}}}),
\begin{align*}
    &  \im \ \Log\dfrac{1-tz}{1-t+tz}\in \left(-\pi, \pi \right)
    \, , \quad     z\! \in \! (0,1)\cup \left(\Bb{C}\setminus\Bb{R}\right), \ t \in (0,1)\,,
\end{align*}

\noindent and therefore for arbitrary $z\! \in \! (0,1)\cup \left(\Bb{C}\setminus\Bb{R}\right)$ we deduce from  ({\hyperlink{addf12eir}{\ref*{f12eir}}}) and \eqref{add3} that
\begin{align*}
   \left|\im\, \Log \dfrac{ \he (1\!-\! z)}{\he (z)}\right| & \!=\!\dfrac{1}{\pi^{2 }}\!\left| \int\limits_{0}^{1}\!\! \dfrac{\im\, \Log\dfrac{1-tz}{1-t+tz}}{t (1-t)\left(\he \left(t \right)^{2} + \he \left(1-t \right)^{2}\right)} \,  \diff  t\right| \\   &   \leq
    \dfrac{1}{\pi^{2 }}\! \int\limits_{0}^{1}\!\! \dfrac{ \left|\im\, \Log\dfrac{1-tz}{1-t+tz} \right| }{t (1-t)\left(\he \left(t \right)^{2} + \he \left(1-t \right)^{2}\right)} \,  \diff  t   < \dfrac{\pi}{2} \ ,
\end{align*}

\noindent which proves
\begin{align*}\tag{\ref{f15eir}}
    &   \Arg \dfrac{\het (1\!-\! z)}{\het (z)} \in \left(- \dfrac{\pi}{2} , \dfrac{\pi}{2}\right)  \ , \quad  z\! \in \! (0,1)\cup \left(\Bb{C}\setminus\Bb{R}\right) \ .
\end{align*}

\noindent Besides that, by virtue of
\begin{align*}\tag{\ref{f10int}}
    &        \ie (z) :=   i \, \frac{\he(1-z)}{\he(z)},
\quad  z \in  (0,1)\cup \left(\Bb{C}\setminus\Bb{R}\right),
\end{align*}

\noindent and  ({\hyperlink{addf12eir}{\ref*{f12eir}}}), we obtain
\begin{align*}
    & \Arg \ie (z) = \dfrac{\pi}{2} +  \Arg \dfrac{\het (1\!-\! z)}{\het (z)} =
\dfrac{\pi}{2} +     \im\, \Log \dfrac{ \he (1\!-\! z)}{\he (z)}\\  &   =
         \dfrac{\pi}{2} +\dfrac{1}{\pi^{2 }}\! \int\limits_{0}^{1}\!\! \dfrac{\im \,  \Log\dfrac{1-tz}{1-t+tz}}{t (1-t)\left(\he \left(t \right)^{2} + \he \left(1-t \right)^{2}\right)} \,  \diff  t
         \\  &   =
         \dfrac{\pi}{2} +\dfrac{1}{\pi^{2 }}\! \int\limits_{0}^{1}\!\! \dfrac{\Arg \dfrac{1-tz}{1-t+tz}}{t (1-t)\left(\he \left(t \right)^{2} + \he \left(1-t \right)^{2}\right)} \,  \diff  t \\  &   =
         \dfrac{\pi}{2} +\dfrac{1}{\pi^{2 }}\! \int\limits_{0}^{1}\!\! \dfrac{ \Arg (1\! -\! tz)\! -\!  \Arg (1\! -\! t\!  +\!  t z)}{t (1-t)\left(\he \left(t \right)^{2} + \he \left(1-t \right)^{2}\right)} \,  \diff  t \ , \
\end{align*}

\noindent i.e.,
\begin{align}\label{add4}
    &  \Arg \ie (z) = \dfrac{\pi}{2} +\dfrac{1}{\pi^{2 }}\! \int\limits_{0}^{1}\!\! \dfrac{ \Arg (1\! -\! tz)\! -\!  \Arg (1\! -\! t\!  +\!  t z)}{t (1-t)\left(\he \left(t \right)^{2} + \he \left(1-t \right)^{2}\right)} \,  \diff  t \ , \
\end{align}

\noindent for each $ z\! \in \! (0,1)\cup \left(\Bb{C}\setminus\Bb{R}\right)$. The estimate
\begin{align*}\tag{\ref{f14eir}}
    &   \big| \Arg (1\! -\! tz)\! -\!  \Arg (1\! -\! t\!  +\!  t z)\big| \! <\!  \pi  \, , \ \ \  t \in (0,1) \ , \  z\! \in \! (0,1)\cup \left(\Bb{C}\setminus\Bb{R}\right) \ ,
\end{align*}

\noindent can be strengthened as follows
\begin{align}\label{add5}
    & \Arg (1\! -\! tz)\! -\!  \Arg (1\! -\! t\!  +\!  t z) \in  - \left(0 , \pi\right)\cdot {\rm{sign}} (\im\, z) \ , \
\end{align}

\noindent for any $t \in (0,1)$ and  $ z\! \in \! (0,1)\cup \left(\Bb{C}\setminus\Bb{R}\right)$, because the both numbers
 $\Arg (1\! -\! tz)$ and $ (- 1) \cdot \Arg (1\! -\! t\!  +\!  t z)$ are of the same sign which is equal to $ (- 1) \cdot {\rm{sign}} (\im\, z)$. And therefore we deduce from
 \eqref{add5}, \eqref{add4} and \eqref{add3} that
  \begin{align*}\tag{\ref{f16eir}}
    &   \Arg \ie (z) \in \dfrac{\pi}{2} - \left(0 , \dfrac{\pi}{2}\right)\cdot {\rm{sign}} (\im\, z) \subset \left(0, \pi\right)  \ , \  \  z\! \in \! \Bb{C}\setminus\Bb{R}  \ .
\end{align*}

\noindent Together with ({\hyperlink{addf12eir}{\ref*{f12eir}}}) and ({\hyperlink{addff13eir}{\ref*{f13eir}}}), this implies that
\begin{align}\label{add6}
      &  \Log \ie (z)\!=\!\dfrac{i \pi}{2} +\dfrac{1}{\pi^{2 }}\! \int\limits_{0}^{1}\!\! \dfrac{ \Log\dfrac{1-tz}{1-t+tz}}{t (1-t)\left(\he \left(t \right)^{2} + \he \left(1-t \right)^{2}\right)} \,  \diff  t \, , \\[0,4cm]    &\label{add7}
    \Log \ie (z)=\dfrac{i \pi}{2} +\dfrac{1-2z}{\pi} \int\limits_{0}^{1} \ \,  \dfrac{\arctan \dfrac{\het \left(1-t\right)}{\het \left(t\right)}}{(1-tz)(1-t +t z)} \,  \diff  t \,,
\end{align}

\vspace{0.2cm}
\noindent  for every $ z\! \in \! (0,1)\cup \left(\Bb{C}\setminus\Bb{R}\right)$. Besides that, \eqref{add6} yields
\begin{align}\label{add8}
    &  \ln \left|\ie (z)\right|\!=\! \dfrac{1}{\pi^{2 }}\! \int\limits_{0}^{1}\!\! \dfrac{ \ln\dfrac{|1-tz|}{|1-t+tz|}}{t (1-t)\left(\he \left(t \right)^{2} + \he \left(1-t \right)^{2}\right)} \,  \diff  t \, ,
\end{align}

\noindent  for each $ z\! \in \! (0,1)\cup \left(\Bb{C}\setminus\Bb{R}\right)$ and, in particular,
when $z= i s + 1/2$, $t \in (0,1)$ and $ s \in \Bb{R}$ it follows from
\begin{align*}
   \left|\dfrac{1-t(i s + 1/2)}{1-t+t(i s + 1/2)}\right|= \left|\dfrac{1-t/2 - i t s }{1-t/2 +i t s} \right| =1 \,,
\end{align*}

\noindent that
 \begin{align}\label{add9}
    & \left|\ie \left(\dfrac{1}{2} + i s\right)\right|\!=\!1 \ , \quad  s \in \Bb{R} \, .
\end{align}

}}\end{clash}
\end{subequations}

%%%%%%%%%%%%%%%%%%%%%%%%%%%%%%%%%%%%%%%%%%%%%%%%%%%%%%%%%%%%%%%%%%%%%%%%%%%%%%%%%%%%%%%%%%%%%%%%%%%%%%%%%

\vspace{0.05cm}
\subsection[\hspace{-0,25cm}. \hspace{0,075cm}Notes on Section~\ref{sp}]{\hspace{-0,11cm}{\bf{.}} Notes on Section~\ref{sp}}

\vspace{-0.15cm}
\begin{clash}{\rm{\hypertarget{d10}{}\label{case10}\hspace{0,00cm}{\hyperlink{bd10}{$\uparrow$}} \ We prove \eqref{f29preresinvlamhyp} in more detail.

 \vspace{0,25cm}
 By virtue of \eqref{f12apreresinvlamhyp}, for arbitrary $x>0$ there exist the finite limits
$\ie (1+x \pm i 0)$ and $\ie (-x \pm i 0)$ for which the functional equality

\vspace{-0,1cm}
\begin{align*}
\tag{\ref{f27preresinvlamhyp}}
\ie (1-z) =-\dfrac{1}{\ie (z)}, \quad z \in
\Lambda:= (0,1)\cup \left(\Bb{C}\setminus\Bb{R}\right),
\end{align*}

\noindent gives, taking account that $1-\Lambda = \Lambda$,

\vspace{-0,1cm}
\begin{align*}
    &  \ie (1+x - i 0) =-\dfrac{1}{\ie (- x + i 0)} \ , \  \ie (1+x + i 0) =-\dfrac{1}{\ie (- x - i 0)} \ , \
\end{align*}

\noindent and using
\begin{align*}\tag{\ref{f28preresinvlamhyp}}
\ie (-x +  i  0 ) = 2 +  \ie (-x -  i  0 ) \ , \quad  x>0\,.
\end{align*}

\vspace{0.25cm}
\noindent  we get

\vspace{-0,1cm}
\begin{align*}
    &   \ie (1+x - i 0) =-\dfrac{1}{\ie (- x +  i 0)} = -\dfrac{1}{2 + \ie (- x - i 0)} =
    -\dfrac{1}{2 -  \dfrac{1}{ \ie (1+x + i 0)} } \\[0,3cm]     &   =
    - \dfrac{ \ie (1+x + i 0)}{2  \ie (1+x + i 0) -1} = \dfrac{ \ie (1+x + i 0)}{1- 2  \ie (1+x + i 0)} \ ,
\end{align*}

\vspace{0.25cm}
\noindent which proves \eqref{f29preresinvlamhyp} and its obvious consequence

\vspace{-0,1cm}
\begin{align}\label{ad12a}
    &  \ie (1  +  x  +  i 0)  =  \dfrac{\ie (  1 +  x -  i 0)}{ 2 \ie (  1 +  x  -  i 0)  + 1}  \ , \quad  x > 0 \,.
\end{align}
}}\end{clash}

%\newpage

\begin{subequations}\noindent
\begin{clash}{\rm{\hypertarget{d11}{}\label{case11}\hspace{0,00cm}{\hyperlink{bd11}{$\uparrow$}} \ We prove Lemma~{\hyperlink{tempbd9}{3.1}}.

\vspace{0.15cm}
\noindent
{\emph{Proof of \eqref{f30preresinvlamhyp}.}}
Everywhere below inside this proof we assume that $\sigma := {\rm{sign}} (\im  \, z) \in \left\{1,-1\right\}$, $z \in \Lambda$ and $|z| \to +\infty$.
Then according to the definition
\begin{align*}\tag{\ref{f10int}}
       \ie (z) :=   i \, \frac{\he(1-z)}{\he(z)},
\quad  z \in \Lambda:=  (0,1)\cup \left(\Bb{C}\setminus\Bb{R}\right),
\end{align*}

\noindent and the asymptotic formula
\begin{align*}\tag{\ref{f7preresinvlamhyp}}
    & {\rm{(c)}} \   \he (z) =\dfrac{\Log 16 (1-z)}{\pi \sqrt{1-z}} + O \left(\dfrac{\ln |z|}{|z|^{3/2}}\right) \ ,  \ \      |z| \to +\infty \   , \
z \not\in  [2, +\infty)  \ ,
\end{align*}

\noindent we can write
\begin{align}\nonumber
     \ie (z)  &  = i \dfrac{\dfrac{\Log 16 z}{\pi \sqrt{z}} + O \left(\dfrac{\ln |1-z|}{|1-z|^{3/2}}\right) }{
    \dfrac{\Log 16 (1-z)}{\pi \sqrt{1-z}} + O \left(\dfrac{\ln |z|}{|z|^{3/2}}\right)     } \\[0,3cm]  \nonumber    &   =
    i \dfrac{\sqrt{1-z}}{\sqrt{z}} \cdot \dfrac{\ln | z| + i \Arg (z) +  \ln 16+ O \left(\dfrac{\ln |z|}{|z|}\right) }{\ln |1- z| + i \Arg (1-z) +  \ln 16+ O \left(\dfrac{\ln |z|}{|z|}\right)}\\[0,3cm]     &   =
    i \dfrac{\sqrt{1-z}}{\sqrt{z}} \cdot \dfrac{\ln | z| + i \Arg (z) +  \ln 16+ O \left(\dfrac{\ln |z|}{|z|}\right) }{\ln |z| + i \Arg (1-z) +  \ln 16+ O \left(\dfrac{\ln |z|}{|z|}\right)}  \ , \
\label{ad13}\end{align}

\noindent where we use the obvious relation $\ln |1- z| = \ln | z| + O (|z|^{-1})$ and since $\Lambda \cap \Bb{R} = (0, 1)$ and $|z| \to +\infty$
 it can be considered that $z \in \Bb{C}\setminus \Bb{R} $ and therefore
\begin{align}\label{ad14}   &    \Arg (1-z)-  \Arg (-z) =  \Arg \dfrac{1-z}{(-z)} =\Arg \left(1 - \dfrac{1}{z}\right)  = O \left(|z|^{-1}\right) \ , \\[0,3cm]  \label{ad15}   &
\Arg (z)- \Arg (-z) = \pi \sigma   \ , \\[0,3cm]     & \nonumber
      \dfrac{\sqrt{1-z}}{\sqrt{z}} = \exp \left(\dfrac{1}{2} \Log (1-z) - \dfrac{1}{2} \Log z\right)  \\[0,3cm]   \nonumber   &   =
      \exp \left(\dfrac{1}{2} \Log (1-z)- \dfrac{1}{2} \Log (-z) + \dfrac{1}{2} \Log (-z)  - \dfrac{1}{2} \Log z\right)\\[0,3cm]   \nonumber   &   =
   \exp \Big(\dfrac{1}{2}\big(\ln |1-z| - \ln |z|\big) + \dfrac{i}{2} \big(\Arg (1-z)-  \Arg (-z)\big)   \\[0,3cm]   \nonumber   &   +
   \dfrac{i}{2}  \big(\Arg (-z)-  \Arg (z)\big)     \Big) =
   \exp \left(-\dfrac{i \pi }{2} \sigma + O \left(|z|^{-1}\right)\right) = i^{-\sigma}   +  O \left(|z|^{-1}\right) \ , \
\end{align}

\noindent from which and \eqref{ad13} we deduce that
\begin{align*}
    &   \ie (z)  = i^{1-\sigma} \dfrac{\ln | z| + i \Arg (z) +  \ln 16+ O \left(\dfrac{\ln |z|}{|z|}\right) }{\ln |z| + i \Arg (1-z) +  \ln 16+ O \left(\dfrac{\ln |z|}{|z|}\right)} + O \left(\dfrac{1}{|z|}\right) \ ,
\end{align*}

\noindent which  can also be written as follows
\begin{align}\label{ad16}
    &   \ie (z)  = i^{1-\sigma} \dfrac{\ln | z| + i \Arg (z) +  \ln 16+ O \left(\dfrac{\ln |z|}{|z|}\right) }{\ln |z| + i \Arg (1-z) +  \ln 16+ O \left(\dfrac{\ln |z|}{|z|}\right)}  \ .
\end{align}

\noindent For arbitrary $a \in \Bb{C}$ we obviously have
\begin{align*}
    &  i^{1-\sigma} a - \sigma = \left\{\begin{array}{ll} a-1  = \sigma (a-1)\, ,    & \ \ \hbox{if} \, \ \sigma = 1  ; \\[0,0cm]
 \  1-a  = \sigma (a-1)\, , & \ \  \hbox{if} \ \  \sigma = -1 ,
    \end{array}\right.
\end{align*}

\noindent and therefore \eqref{ad16} together with \eqref{ad14} and \eqref{ad15} yields that
\begin{align} \nonumber
    &   \ie (z) -\sigma  = \sigma \dfrac{ i \Arg (z) -i \Arg (1-z)+ O \left(\dfrac{\ln |z|}{|z|}\right) }{\ln |z| + i \Arg (1-z) +  \ln 16+ O \left(\dfrac{\ln |z|}{|z|}\right)}   \\  \nonumber   &   =
    \sigma \dfrac{ i \pi \sigma + O \left(\dfrac{\ln |z|}{|z|}\right) }{\ln |z| + i \Arg (1-z) +  \ln 16+ O \left(\dfrac{\ln |z|}{|z|}\right)} \\    &   =
     \dfrac{ i \pi  + O \left(\dfrac{\ln |z|}{|z|}\right) }{\ln |z| + i \Arg (1-z) +  \ln 16+ O \left(\dfrac{\ln |z|}{|z|}\right)} \ .
\label{ad17}
\end{align}

\noindent Then according to ($\sigma = {\rm{sign}} (\im  \, z)$)
\begin{align*}
    &  \tag{\ref{f29ainvlamhyp}}
     \ie(\infty; z) =\dfrac{1}{  \ie (z) -\sigma } + \dfrac{ i }{ \pi} \ln  |z| \ ,
\end{align*}

\noindent it follows from \eqref{ad17}  that
\begin{align*}
      \ie(\infty; z)& = \dfrac{ i }{ \pi} \ln  |z| + \dfrac{\ln |z| + i \Arg (1-z) +  \ln 16+ O \left(\dfrac{\ln |z|}{|z|}\right)}{ i \pi  + O \left(\dfrac{\ln |z|}{|z|}\right) } \\[0,3cm]      & =
   \dfrac{ i }{ \pi} \ln  |z| +  \dfrac{\ln |z|}{i \pi } + \dfrac{i \Arg (1-z)}{i \pi} +  \dfrac{\ln 16}{i \pi}+ O \left(\dfrac{\ln^{2} |z|}{|z|}\right)\\[0,3cm]      & = \dfrac{\ln 16 }{i \pi} +  \dfrac{\Arg (1-z)}{\pi}  + O \left(\dfrac{\ln^{2} |z|}{|z|}\right) \ ,
\end{align*}

\noindent which coincides with  \eqref{f30preresinvlamhyp} and completes its proof.

\vspace{0.5cm} \noindent
{\emph{Proof of \eqref{f32preresinvlamhyp}.}} Everywhere below inside this proof we assume that $z \in \Lambda$ and $z \to 0$.
Then according to the definition
\begin{align*}\tag{\ref{f10int}}
       \ie (z) :=   i \, \frac{\he(1-z)}{\he(z)},
\quad  z \in \Lambda:=  (0,1)\cup \left(\Bb{C}\setminus\Bb{R}\right),
\end{align*}

\noindent and the asymptotic formulas
\begin{align*}
    &{\rm{(a)}} \   \he (z) =\dfrac{1}{\pi}\,
\Log\dfrac{16}{1-z}\ +  O \left(|1-z| \ln \dfrac{1}{|1-z|}\right) \ , \ \      z \to 1  , \
  z \not\in  [1, 2]  \, , \\[0,3cm]     &   \tag{\ref{f7preresinvlamhyp}}
  {\rm{(b)}} \  \he (z) = 1 + O (|z|)  \ ,  \ \        z \to 0  , \ z \in \Bb{D} \ ,
\end{align*}

\noindent we can write
\begin{align*}
 - \ie (z)   &  =(- i) \dfrac{
    \dfrac{1}{\pi}\,
\Log\dfrac{16}{z}\ +  O \left(|z| \ln \dfrac{1}{|z|}\right)   }{ 1 + O (|z|)  } \\[0,3cm]      &   =
   \dfrac{1}{\pi i }\,
\Log\dfrac{16}{z}\ +  O \left(|z| \ln \dfrac{1}{|z|}\right) \\[0,3cm]     &   = - \dfrac{i}{\pi }\ln\dfrac{1}{|z|} +
  \dfrac{ \ln 16 }{ i  \pi} + \dfrac{\Arg (1/z)}{\pi} +  O \left(|z| \ln \dfrac{1}{|z|}\right)    \ , \
\end{align*}

\noindent which according to
\begin{align*}
    \tag{\ref{f29ainvlamhyp}}
      \ie(0; z)  := - \ie (z)  + \dfrac{i }{\pi} \ln\dfrac{1}{|z|}  \ ,
\end{align*}

\noindent implies the validity of
\begin{align*}\hypertarget{thrf}{}
    &   \tag{\ref{f32preresinvlamhyp}}
   \ie(0; z) \!=\!   \dfrac{ \ln 16 }{ i  \pi} + \dfrac{\Arg (1/z)}{\pi} +  O \left(|z| \ln \dfrac{1}{|z|}\right) \ ,      \quad      \Lambda \ni z\to 0 \, .
\end{align*}

\noindent

\vspace{0.25cm} \noindent
{\emph{Proof of \eqref{f31preresinvlamhyp}.}} In view of
\begin{align*}
    &   \tag{\ref{f27preresinvlamhyp}}
\ie (z)\ie (1-z) =-1, \quad z \in
\Lambda:= (0,1)\cup \left(\Bb{C}\setminus\Bb{R}\right),
\end{align*}

\noindent we can deduce from the definitions of $\ie(0; z)$ and $\ie(1; z)$ in \eqref{f29ainvlamhyp},
\begin{multline*}
       \ie(0; z)       := - \ie (z)  + \dfrac{\imag  }{\pi} \ln\dfrac{1}{|z|} \, ; \ \
    \ie(1; z)  :=     \dfrac{1}{\ie (z)}  +  \dfrac{\imag }{\pi}\ln\dfrac{1}{|1-z|}  \ ,  \\[0,1cm]
        z \in  (0,1)\cup \left(\Bb{C}\setminus\Bb{R}\right) \,,
\tag{\ref{f29ainvlamhyp}}\end{multline*}

\noindent that
\begin{align} \label{f4adlem1}
    &  \ie(1; z)  =  \ie(0;1- z)  \ , \quad  z \in  (0,1)\cup \left(\Bb{C}\setminus\Bb{R}\right) \,.
\end{align}

\vspace{0.25cm}
\noindent If $z \in \Lambda$ and $z \to 1$ then the property $1-\Lambda=\Lambda$ yields $1-z \in \Lambda$ and $1-z \to 0$ and we can substitute such
$z$ in \eqref{f4adlem1} to get with the help of ({\hyperlink{thrf}{\ref*{f32preresinvlamhyp}}}) the following expression
\begin{align*}
    &
    \ie(1; z)  \!=\! \ie(0;1- z) \!=\!   \dfrac{\ln 16 }{ i  \pi}\! +\! \dfrac{\Arg \big(1/(1\!-\!z)\big)}{\pi} \!+ \! O \left(|1\!-\!z| \ln \dfrac{1}{|1\!-\!z|}\right)      \,,
\end{align*}

\noindent which proves the validity of \eqref{f31preresinvlamhyp} and finishes the proof of Lemma~{\hyperlink{tempbd9}{3.1}}. $\square$

%%%%%%%%%%%%%%%%%%%%%%%%%%%%%%%%%

\vspace{0,25cm}
Observe that the equalities of Lemma~{\hyperlink{tempbd9}{3.1}} can also be written as follows.

\begin{multline}  \label{f30invlamhyp}
- \dfrac{1}{  \ie (z) -{\rm{sign}} (\im z)}      \!=\!\dfrac{  i  }{ \pi} \ln  |z| +\dfrac{ i \ln 16 }{ \pi} -  \dfrac{\Arg (1-z)}{\pi}  + O \left(\dfrac{\ln^{2} |z|}{|z|}\right)     \, , \\[0,25cm]
     \Lambda \ni  z\to \infty  \, ,
\end{multline}
\begin{multline} \label{f31invlamhyp}
- \dfrac{1}{\ie (z)}    \!=\!  \dfrac{ i }{\pi}\ln\dfrac{1}{|1-z|}
+   \dfrac{ i \ln 16 }{ \pi} \! -\! \dfrac{\Arg \big(1/(1\!-\!z)\big)}{\pi} \!+ \! O \left(|1\!-\!z| \ln \dfrac{1}{|1\!-\!z|}\right)      \, ,\\[0,25cm]
\Lambda\! \ni \!   z\!\to \! 1  \, ,
\end{multline}
\begin{multline} \label{f32invlamhyp}\ie (z)    \!=\! \dfrac{ i  }{\pi} \ln\dfrac{1}{|z|}
+ \dfrac{ i \ln 16 }{ \pi}  - \dfrac{\Arg (1/z)}{\pi} +  O \left(|z| \ln \dfrac{1}{|z|}\right) \ , \\[0,25cm]
  \Lambda \ni z\to 0 \, ,
\end{multline}

\vspace{0.25cm}
\noindent and therefore it can be formulated in the following form.}}\end{clash}

\vspace{0.25cm}
\begin{lemma}\hspace{-0,18cm}{\bf{.}}\label{adlem1}
 Let
 \begin{align*}
    &  F\left(\dfrac{1}{2} \ , \ \dfrac{1}{2} \ ; \  1  \ , \  z \right) =\dfrac{1}{\pi }  \int_{0}^{1} \dfrac{ \diff  t}{\sqrt{t(1-t)(1-tz)  }}
    \ , \quad z \in \Bb{C}\setminus [1, +\infty) \,,
 \end{align*}

 \noindent
 be the Gauss hypergeometric function and
 \begin{align*}
    &    \ie (z) :=   i \, \frac{ F\left({{\dfrac{1}{2}}} \ , \ {{\dfrac{1}{2}}} \ ; \  1  \ , \  1-z \right)}{ F\left({{\dfrac{1}{2}}} \ , \ {{ \dfrac{1}{2} }}\ ; \  1  \ , \  z \right)},
\quad  z \in  (0,1)\cup \left(\Bb{C}\setminus\Bb{R}\right),
 \end{align*}

\noindent  be the  Schwarz triangle function. Then the  asymptotic formulas

 \begin{multline}  \label{f1adlem1}
- \dfrac{1}{  \ie (z) -{\rm{sign}} (\im z)}       = \dfrac{  i  }{ \pi} \ln( 16 |z|)  -  \dfrac{\Arg (1-z)}{\pi}  + O \left(\dfrac{\ln^{2} |z|}{|z|}\right)     \, , \\[0,25cm]
     (0,1)\cup \left(\Bb{C}\!\setminus\!\Bb{R}\right) \ni  z\to \infty  \, ,
\end{multline}
\begin{multline} \label{f2adlem1}- \dfrac{1}{\ie (z)}    =  \dfrac{ i }{\pi}\ln\dfrac{16}{|1-z|}
 - \dfrac{\Arg \big(1/(1\!-\!z)\big)}{\pi} +  O \left(|1\!-\!z| \ln \dfrac{1}{|1\!-\!z|}\right)      \, ,\\[0,25cm]
(0,1)\cup \left(\Bb{C}\!\setminus\!\Bb{R}\right)\! \ni \!   z\!\to \! 1  \, ,
\end{multline}
\begin{multline} \label{f3adlem1}\ie (z)    =  \dfrac{ i  }{\pi} \ln\dfrac{16}{|z|}
 - \dfrac{\Arg (1/z)}{\pi} +  O \left(|z| \ln \dfrac{1}{|z|}\right) \ , \\[0,25cm]
 (0,1)\cup \left(\Bb{C}\!\setminus\!\Bb{R}\right) \ni z\to 0 \, ,
\end{multline}

 \vspace{0.25cm}
 \noindent hold, where $1-z, 1/(1\!-\!z), 1/z \in \Bb{C}\setminus(-\infty, 0]$ for any $z\in (0,1)\cup \left(\Bb{C}\!\setminus\!\Bb{R}\right)$.
\end{lemma}

\end{subequations}

\begin{subequations}

\begin{clash}{\rm{\hypertarget{d27}{}\label{case27}\hspace{0,00cm}{\hyperlink{bd27}{$\uparrow$}} \ We prove \eqref{f1rem1},
\begin{align}\tag{\ref{f1rem1}}\label{}
    &  \ie( z) =\ie \big(z/(z-1)\big) + \sigma (z) \ , \quad \sigma (z):= {\rm{sign}} ({\rm{Im}}\, z) \ , \
    z \in \Bb{C}\setminus \Bb{R} \,.
 \end{align}

\vspace{0.25cm} According to \eqref{f1xtheor3},
\begin{align*}\tag{{\ref{f1xtheor3}}}
    &  \lambda (z) =  {\Theta_{2}(z)^{4}}\big/{\Theta_{3}(z)^{4}}
     \ , \quad  z \in  \Bb{H} \,,
\end{align*}

\noindent \eqref{f9elemtheta}(d), (e), written in the form
\begin{align*}\tag{{\ref{f9elemtheta}}}
    &  {\rm{(d)}} \ \ \Theta_2 (z+1)^{4} = -\Theta_2(z)^{4} \, ,  \quad
  {\rm{(e)}} \ \ 	\Theta_3(z+1)^{4} = \Theta_4(z)^{4} \, \ , \quad  z \in  \Bb{H} \,,
\end{align*}

\noindent and \eqref{adt27}(c),(d),(e),
\begin{align*}
    &   {\rm{(c)}} \ \ \Theta_{2} (z+2)^{4} =  \Theta_{2} (z)^{4}\, ,   &     &
     {\rm{(d)}} \ \ \Theta_{3} (z+2) = \Theta_{3} (z)\, , \\   &
      {\rm{(e)}} \ \ \Theta_{4} (z+2) = \Theta_{4} (z)
      \,  ,   &     &     z \in  \Bb{H} \,,
\tag{{\ref{adt27}}}\end{align*}

\noindent together with \eqref{f8theor3},
\begin{align*}\tag{{\ref{f8theor3}}}
      &     \Theta_{2}(z)^{4} + \Theta_{4}(z)^{4} = \Theta_{3}(z)^{4}
       \ , \quad  z \in  \Bb{H} \ ,
\end{align*}

\noindent and \eqref{f0theor3},
\begin{align*}\tag{{\ref{f0theor3}}}
    & \Theta_{3}(z)\Theta_{4}(z)\Theta_{2}(z)\neq 0 \ , \quad  z \in \Bb{H}\,,
 \end{align*}

\noindent  we have  that for any $z \in  \Bb{H}$ the following identities hold,
\begin{align*}
      \lambda (z-1) &  =   \dfrac{\Theta_{2}(z-1)^{4}}{\Theta_{3}(z-1)^{4}}  \    \stackrel{{\fo{\eqref{adt27}(c),(d)}}}{\vphantom{A}=} \ \dfrac{\Theta_{2}(z+1)^{4}}{\Theta_{3}(z+1)^{4}} =  \lambda (z+1)     \ , \\[0,4cm]
         \lambda (z+1)   &=   \dfrac{\Theta_{2}(z+1)^{4}}{\Theta_{3}(z+1)^{4}}  \    \stackrel{{\fo{\eqref{f9elemtheta}(d),(e)}}}{\vphantom{A}=} \ -  \dfrac{\Theta_{2}(z)^{4}}{\Theta_{4}(z)^{4}}
           \ \stackrel{{\fo{\eqref{f8theor3}}}}{\vphantom{A}=}\
          -\dfrac{\Theta_{2}(z)^{4}}{\Theta_{3}(z)^{4} - \Theta_{2}(z)^{4}}
        \\ & =  \dfrac{\dfrac{\Theta_{2}(z)^{4}}{\Theta_{3}(z)^{4}}}{\dfrac{\Theta_{2}(z)^{4}}{\Theta_{3}(z)^{4}}-1}
             \  \stackrel{{\fo{\eqref{f1xtheor3}}}}{\vphantom{A}=}\
        \dfrac{\lambda (z)}{\lambda (z)-1} \ ,
\end{align*}

\noindent i.e.,
\begin{align}\label{f1case27}
     & \lambda (z-1)   = \lambda (z+1)   = \dfrac{\lambda (z)}{\lambda (z)-1}\ , \quad  z \in \Bb{H}\,.
\end{align}

\noindent In view of Theorem~\ref{inttheor2}, for every $z\in \fet$ there exists a unique $y \in   (0,1)\cup \left(\Bb{C}\setminus\Bb{R}\right) $ such that $z = \ie (y)$. By substituting $\ie (y)$ in place of  $z\in \fet \subset \Bb{H}$ in
\eqref{f1case27} we obtain, by virtue of \eqref{f1theor3} written in the form
\begin{align}\label{f2case27}
      &   \lambda \big(\ie (y)\big) = y  \ , \quad  y \in
      (0,1)\cup \left(\Bb{C}\setminus\Bb{R}\right) \ ,
\end{align}

\noindent that
\begin{align}\label{f3case27}
    &   \lambda\left( \ie( y)+ \sigma\right)=\dfrac{y}{y-1} \ , \quad  y \in \Lambda:= (0,1)\cup \left(\Bb{C}\setminus\Bb{R}\right)  \ , \  \sigma \in \{1, -1\} \,.
\end{align}

\noindent Writing \eqref{f2case27} in the equivalent form
\begin{align}\label{f4case27}
      &   \ie (\lambda (z)) = z \,,  \quad z \in \fet \,  ,
\end{align}

\noindent we conclude that by applying the function $\ie $ to both of parts of the identity  \eqref{f3case27}
we obtain
\begin{align}\label{f5case27}
    &  \ie( y)+ \sigma = \ie \left(\dfrac{y}{y-1}\right) \ , \quad  y \in  (0,1)\cup \left(\Bb{C}\setminus\Bb{R}\right)  \ , \  \sigma \in \{1, -1\} \,,
\end{align}

\noindent provided that $ \ie( y)+ \sigma \in \fet $. By \eqref{f33preresinvlamhyp},
 \begin{align*}\tag{{\ref{f33preresinvlamhyp}}}
    &  \ie \left(\Lambda\right) = \fet \ ,
\end{align*}

\noindent the inclusion $ \ie( y)+ \sigma \in \fet $ is equivalent to
\begin{align}\label{f6case27}
    &  \ie( y) \in \fet \cap \left(\fet -\sigma \right) \,.
\end{align}

\noindent It is easy to see that
\begin{align}\label{f7case27}
    & \fet \cap \left(\fet -\sigma \right)=
   \left\{ \begin{array}{cc}
     \left\{\, z \in  \fet\ |\ \re \,z <0  \,\right\} &  \ \mbox{if} \  \sigma = 1  \,,\\[0,3cm]
 \left\{\, z \in  \fet\ | \ \re \,z > 0  \,\right\}  & \ \mbox{if} \  \sigma = -1  \,.
    \end{array}\right.
\end{align}

\noindent while \eqref{f16eir},
\begin{align*}\tag{{\ref{f16eir}}}
    &   \Arg  \ie (z) \in \dfrac{\pi}{2} - \left(0 ,
    \dfrac{\pi}{2}\right)\cdot {\rm{sign}} ({\rm{Im}}\, z)  \ , \  \
     z\! \in \! \Bb{C}\setminus\Bb{R}  \ ,
\end{align*}

\noindent yields that
\begin{align}\label{f8case27}
     \left\{ \begin{array}{rcl}
    \ie \big( \Bb{H}\big)   &  \subset  &    \left\{\, z \in  \fet\ | \ \re \,z > 0  \,\right\} \,,\\[0,3cm]
      \ie \big( - \Bb{H}\big) &   \subset  & \left\{\, z \in  \fet\ | \ \re \,z < 0  \,\right\}  \,.
    \end{array}\right.
    \end{align}

\noindent According to the definition \eqref{f10int},
\begin{align}\tag{\ref{f10int}}
       \ie (z) :=   \imag \cdot \frac{\he (1-z)}{\he (z)},
\quad  z \in  (0,1)\cup \left(\Bb{C}\setminus\Bb{R}\right),
\end{align}

\noindent of the function $\ie $ we have that
\begin{align}\label{f10case27}
    &  \ie \big((0,1)\big) = \left\{\,\imag \cdot \dfrac{\he (1-x)}{\he (x)}  \, \Big| \, x \in (0, 1) \, \right\} =
    \imag        \Bb{R}_{>0} \ , \
\end{align}

\noindent because according to the property \eqref{f2eir}
\begin{align}\tag{\ref{f2eir}}
    & \hspace{-0,2cm} \dfrac{{{\diff}}}{{{\diff}} x} \dfrac{\he
    (1-x)}{\he  (x)} = -  \dfrac{1}{\pi x (1-x) \he  (x)^{2}} < 0 \ ,
     \quad  x \in (0,1) ,
\end{align}

\noindent and $ \he  (0) = 1$, $\lim_{x\in (0,1)\,,\, x \to 1}  \he  (x) = +\infty$ (see \eqref{f7preresinvlamhyp}(a) and (b) ) the function $\he (1-x)/ \he (x)$ decreases  from $+\infty$ to $0$ on the interval $(0,1)$.

\vspace{0.1cm}
 Since
\begin{align*}
    &  \fet  = \left\{\, z \in  \fet\ | \ \re \,z < 0  \,\right\} \sqcup  \imag        \Bb{R}_{>0}\sqcup
    \left\{\, z \in  \fet\ | \ \re \,z > 0  \,\right\}   \ , \\[0,3cm]    &
    \Lambda = \left(-\Bb{H}\right) \sqcup (0,1)  \sqcup \Bb{H}  \ , \
\end{align*}

\noindent then by  Theorem~\ref{inttheor2},
\begin{align*}
    &  \left\{\, z \in  \fet\ | \ \re \,z < 0  \,\right\} \sqcup  \imag        \Bb{R}_{>0}\sqcup
    \left\{\, z \in  \fet\ | \ \re \,z > 0  \,\right\}  =  \fet \\[0,2cm]   &   =  \ie \left(\Lambda\right) =
    \ie \big( -\Bb{H}\big) \sqcup \ie \big((0,1)\big)   \sqcup \ie \big( \Bb{H}\big) \ , \
\end{align*}

\noindent and we deduce from  \eqref{f8case27} and   \eqref{f10case27} that
\begin{align}\label{f11case27}
    &
    \left\{ \begin{array}{rcl}
    \ie \big( \Bb{H}\big)   &   =  &    \left\{\, z \in  \fet\ | \ \re \,z > 0  \,\right\} \,,\\[0,3cm]
      \ie \big( - \Bb{H}\big) &   =  & \left\{\, z \in  \fet\ | \ \re \,z < 0  \,\right\}  \,.
    \end{array}\right.
\end{align}

\noindent Combining this and \eqref{f7case27} yields
\begin{align}\label{f12case27}
    & \fet \cap \left(\fet -\sigma \right)=
   \left\{ \begin{array}{cc}
    \ie \big( - \Bb{H}\big)  &  \ \mbox{if} \  \sigma = 1  \,,\\[0,3cm]
\ie \big( \Bb{H}\big)   & \ \mbox{if} \  \sigma = -1  \,,
    \end{array}\right.
\end{align}

\noindent and therefore, by \eqref{f6case27},
\begin{align}\label{f13case27}
    &  \ie( y) \in \fet \cap \left(\fet -\sigma \right)  \ \Leftrightarrow \ y \in (- \sigma)\cdot\Bb{H}  \ , \quad   \sigma \in \{1, -1\} \,.
\end{align}

\noindent So that \eqref{f5case27} can equivalently be written as follows,
\begin{align}\label{f14case27}
    &   \ie( y)+ \sigma = \ie \left(\dfrac{y}{y-1}\right) \ , \quad  y  \in (- \sigma)\cdot\Bb{H}   \ , \  \sigma \in \{1, -1\} \,.
\end{align}

\noindent This is the same as
\begin{align}\label{f15case27}
    &  \ie( z) - {\rm{sign}} ({\rm{Im}}\, z) =\ie\left(\dfrac{z}{z-1}\right)  \ , \quad
    z \in \Bb{C}\setminus \Bb{R} \,,
\end{align}

\noindent which coincides with \eqref{f1rem1} and completes its proof.

}}\end{clash}
\end{subequations}

 \begin{subequations}
\begin{clash}{\rm{\hypertarget{d12}{}\label{case12}\hspace{0,00cm}{\hyperlink{bd12}{$\uparrow$}} \ We prove that  $|\re\, \ie (z)|\!< \!1$ for all $z\!\in\! \Lambda$.

\vspace{0.25cm}
Since any harmonic function $u : D \mapsto \Bb{R}$ in an open and connected set $D \subset \Bb{R}^{2}$ (see \cite[p. 14, Def. 2.1]{con}) is also subharmonic (see \cite[p. 41]{hay}),  the maximum principle for subharmonic functions in \cite[Theorem 2.3, p.47]{hay} can be formulated for real parts of holomorphic functions as follows.}}\end{clash}

\vspace{-0,15cm}
\begin{thmx}\hspace{-0,18cm}{\bf{.}}
\label{ainttheor1}
Let $f : D \mapsto \Bb{C}$ be a holomorphic function in the open, unbounded and connected set $D \subset \Bb{C}$. Denote by ${\rm{clos}}(D)$ the closure of $D$ in $\Bb{C}$ and  $\partial D :=\left( {\rm{clos}}(D)\right) \setminus D$ the set of all boundary points of $D$. Assume that for arbitrary point $\zeta \in \partial D$ and any $\varepsilon  \in \Bb{R}_{>0}$ there exists $\delta (\zeta, \varepsilon) \in \Bb{R}_{>0}$
such that
\begin{align}\label{f1thd}
    &  \re \, f (z) < \varepsilon  \ , \quad  |z-\zeta | < \delta (\zeta, \varepsilon)\,, \ z \in D \,,
\end{align}

\noindent and that for arbitrary  $\varepsilon  \in \Bb{R}_{>0}$ there exists $\delta (\infty, \varepsilon) \in \Bb{R}_{>0}$
such that
\begin{align}\label{f2thd}
    &  \re \, f (z) < \varepsilon  \ , \quad  |z | > 1/  \delta (\infty, \varepsilon)\,, \ z \in D \,.
\end{align}

\noindent Then either $ \re f \equiv 0$ in $D$ or
\begin{align}\label{f3thd}
    &   \re \, f (z) < 0        \ , \quad  z \in D \,.
\end{align}

\end{thmx}

\vspace{0.25cm} In Theorem~\ref{ainttheor1} we put $D := \Lambda\!=\! (0,1)\cup \left(\Bb{C}\!\setminus\!\Bb{R}\right)$ and
\begin{align}\label{f3athd}
    & f (z) := -1 + \sigma  \ie (z)\,, \quad \sigma \in \{1, -1\}\,.
\end{align}

\noindent  Then
\begin{align*}
    & \partial D =\partial \Lambda = (-\infty, 0] \cup [1,+\infty) \,,
\end{align*}

\noindent and in view of \eqref{f14zpreresinvlamhyp}, written in the form
\begin{align}\tag{\ref{f14zpreresinvlamhyp}}
    &  \ie (-x \pm  \imag   0 ) =  \pm 1 +  \imag   \,
 \dfrac{\he  \left( {1}\big/{(1 + x)}\right)}{\he  \left(  {x}\big/{(1 + x)}\right)}   \, , \ \  x > 0 \, ,
\end{align}

\vspace{-0,25cm}
\noindent we get
\begin{align*}
    &\re f (-x \pm  \imag   0) = -1 + \sigma  \re \ie (-x \pm  \imag   0) =  -1 \pm \sigma \leq 0  \, , \ \  x > 0 \, .
\end{align*}

\noindent Furthermore, by virtue of \eqref{f27preresinvlamhyp}
\begin{align}
\tag{\ref{f27preresinvlamhyp}}
\ie (z) =-\dfrac{1}{\ie (1-z)}, \quad z \in \Lambda:= (0,1)\cup \left(\Bb{C}\setminus\Bb{R}\right),
\end{align}

\noindent and  \eqref{f14zpreresinvlamhyp}, we obtain
\begin{align*}
    &  \re  f (1+x \pm  \imag   0) = -1 + \sigma  \re \ie (1+x \pm  \imag   0) =    -1 -  \sigma \re \dfrac{1}{ \ie (-x \mp  \imag   0)}  \\[0,3cm]    &   =
     -1 -  \sigma \re\dfrac{1}{\mp 1 +  \imag   \,
 \dfrac{\he  \left( {1}\big/{(1 + x)}\right)}{\he  \left(  {x}\big/{(1 + x)}\right)} }
 = -1 +  \sigma \re\dfrac{1}{\pm 1 -  \imag   \,
 \dfrac{\he  \left( {1}\big/{(1 + x)}\right)}{\he  \left(  {x}\big/{(1 + x)}\right)} } \\[0,3cm]    &   =
 -1 \pm  \dfrac{\sigma}{1 +  \dfrac{\he  \left( {1}\big/{(1 + x)}\right)^{2}}{\he  \left(  {x}\big/{(1 + x)}\right)^{2}}}
< 0  \, , \ \  x > 0 \, .
\end{align*}

\noindent So that
\begin{align}\label{f4thd}
    & \re \, f (1+x \pm  \imag   0) < 0  \ , \quad \re \,   f (-x \pm  \imag   0) \leq 0  \ , \quad  x > 0 \,.
\end{align}

\vspace{0,1cm}
\noindent According to \eqref{f18preresinvlamhyp}, the following limits exist
\begin{align}\label{f4cthd}
    &   \lim_{{\fo{\, \Bb{H} \ni z \to 0}}} \  \het (x\pm z) = \he(x\pm i 0)\in \Bb{C}\setminus\{0\}  \ , \quad  x \in \Bb{R}\setminus\{0, 1\} \,,
\end{align}

\noindent and hence, by virtue of \eqref{f10int}, for arbitrary $x \in \Bb{R}\setminus\{0, 1\}$ there exist the finite limits (see also \eqref{f4athd} and \eqref{f4bthd})
\begin{align*}
     \lim_{{\fo{\, \Bb{H} \ni z \to 0}}} \  \ie (x\pm z) &  = i \cdot\lim_{{\fo{\, \Bb{H} \ni z \to 0}}} \
   \frac{\he(1-x\mp z)}{\he(x\pm z)}    = i \cdot\frac{\he(1-x\mp i 0)}{\he(x\pm i 0)}\\[0,2cm]    &= \ie (x\pm i 0)\,.
\end{align*}

\noindent Then obviously   the following finite limits exist
\begin{align}\label{f4ethd}
    &  \lim_{{\fo{\, \Bb{H} \ni z \to 0}}} \  \re \, \ie (x\pm z) = \re\, \ie (x\pm i 0)   \ , \quad
    x \in \Bb{R}\setminus\{0, 1\}\, ,
\end{align}

\noindent and, by virtue of \eqref{f3athd}, for arbitrary $x \in \Bb{R}\setminus\{0, 1\}$ there exist the limits
\begin{align}\label{f4fthd}
    &  \lim_{{\fo{\, \Bb{H} \ni z \to 0}}} \  f (x\pm z) = f (x\pm i 0) \in \Bb{R}  \,.
\end{align}

\noindent This  gives possibility to deduce from \eqref{f4thd} the  existence of $\delta (\zeta, \varepsilon)\in \Bb{R}_{>0}$ satisfying \eqref{f1thd} for arbitrary
$\zeta \in (-\infty, 0) \cup (1,+\infty) = \left(\partial D\right) \setminus\{0, 1\}$ and $\varepsilon \in \Bb{R}_{>0}$.

\vspace{0.25cm} It is easy to see from \eqref{f3adlem1},
\begin{align*}
    &   \tag{\ref{f3adlem1}}\ie (z)    =  \dfrac{ i  }{\pi} \ln\dfrac{16}{|z|}
 - \dfrac{\Arg (1/z)}{\pi} +  O \left(|z| \ln \dfrac{1}{|z|}\right) \, , \  \  \Lambda \ni z\to 0 \, ,
\end{align*}

\noindent that
\begin{align*}
    & \re  \ie (z)    =   - \dfrac{\Arg (1/z)}{\pi} +  O \left(|z| \ln \dfrac{1}{|z|}\right) \, , \  \  \Lambda \ni z\to 0 \, ,
\end{align*}

\noindent and therefore
\begin{align*}
    & \re f (z) = -1 + \sigma \re \ie (z)= -1 - \dfrac{ \sigma \, \Arg (1/z)}{\pi}+  O \left(|z| \ln \dfrac{1}{|z|}\right) \, , \  \  \Lambda \ni z\to 0 \, .
\end{align*}

\noindent This relation together with the inequality
\begin{align*}
    &  -1 - \dfrac{ \sigma \, \Arg (1/z)}{\pi} < 0  \ , \quad  z \in \Lambda \,,
\end{align*}

\noindent leads to the conclusion that for arbitrary $\varepsilon \in \Bb{R}_{>0}$ there exists  $\delta (0, \varepsilon)\in \Bb{R}_{>0}$ satisfying \eqref{f1thd} for $\zeta = 0 \in \partial D$.

\vspace{0.25cm} The two equalities in \eqref{f3th2},
 \begin{align*}\tag{\ref{f3th2}}
    &     \lim_{{\fo{\Lambda\! \ni\! z\! \to \!1}}} \ie (z) =\!\!\! \lim_{{\fo{\Lambda\! \ni\! z\! \to \!\infty}}}\big|\ie (z)\! -\!{\rm{sign}} (\im z)\big|=0 \, ,
\end{align*}

\noindent yield\vspace{-0,3cm}
\begin{align*}
      \lim_{{\fo{\Lambda  \ni  z  \to  1}}}& \  \re \ie (z) = 0  \ ,   &     &
    \lim_{{\fo{\Lambda \cap \Bb{H}  \ni  z  \to  \infty}}} \re \ie (z) = 1   \ , \\[0,3cm]
    \lim_{{\fo{\Lambda \cap \left(-\Bb{H}\right)  \ni  z  \to  \infty}}}  & \ \re \ie (z) = -1   \ , &     &
\end{align*}

\vspace{0.25cm}
\noindent and hence \  $\big(\,\re f (z) := -1 + \sigma \re \ie (z)\,\big)$,
\begin{align*}
      \lim_{{\fo{\Lambda  \ni  z  \to  1}}}&  \ \re f (z) = -1 < 0 \ , &     &
    \lim_{{\fo{\Lambda \cap \Bb{H}  \ni  z  \to  \infty}}} \re f (z) = -1 + \sigma \leq 0   \ , \\[0,3cm]
    \lim_{{\fo{\Lambda \cap \left(-\Bb{H}\right)  \ni  z  \to  \infty}}}  & \ \re f (z) = -1 - \sigma \leq 0    \ .&     &
\end{align*}

\vspace{0.25cm}
\noindent These relations for arbitrary $\varepsilon \in \Bb{R}_{>0}$  prove the existence of $\delta (1, \varepsilon)\in \Bb{R}_{>0}$ satisfying \eqref{f1thd} for $\zeta = 1 \in \partial D$ and the existence of $\delta (\infty, \varepsilon)\in \Bb{R}_{>0}$ satisfying \eqref{f2thd}.

Thus, the conditions of Theorem~\ref{ainttheor1} are satisfied for  $D := \Lambda\!=\! (0,1)\cup \left(\Bb{C}\!\setminus\!\Bb{R}\right)$ and
$f (z) := -1 + \sigma  \ie (z)$ with any $\sigma \in \{1, -1\}$. Since $\re\,  \ie (z)$ is not a constant function in $\Lambda$ we can apply \eqref{f3thd} to get
\begin{align*}
    &  -1 +  \re\,  \ie (z) < 0   \ , \quad  -1 -  \re\,  \ie (z) < 0    \ , \quad  z \in  \Lambda\,,
\end{align*}

\noindent which means the validity of   $|\re\, \ie (z)|\!< \!1$ for each $z\!\in\! \Lambda$, what was to be proved.

\end{subequations}

\vspace{-0.1cm}
\subsection[\hspace{-0,25cm}. \hspace{0,075cm}Notes on Section~\ref{mth}]{\hspace{-0,11cm}{\bf{.}} Notes on Section~\ref{mth}}

\begin{subequations}
\vspace{-0.15cm}
\begin{clash}{\rm{\hypertarget{d13}{}\label{case13}\hspace{0,00cm}{\hyperlink{bd13}{$\uparrow$}} \
We prove Lemma~{\hyperlink{hth1}{\ref*{th1}}} in more detail.

\vspace{0.25cm} Introduce the function
\begin{align}\label{ad18}
    &  \Phi (z) := f \big(\ie (z)\big)  \ , \quad  z \in \Lambda:= (0,1)\cup \left(\Bb{C}\setminus\Bb{R}\right) \ .
\end{align}

%\vspace{-0,1cm}
\noindent It follows from \eqref{f14zpreresinvlamhyp}
\begin{align}\tag{\ref{f14zpreresinvlamhyp}}
 \ie (-x -  i  0 )=  - 1 +  i  \,
 \dfrac{\he  \left( {1}\big/{(1 + x)}\right)}{\he  \left(  {x}\big/{(1 + x)}\right)} \in -1 + i \Bb{R}_{>0}   \, , \ \  x > 0 \, ,
\end{align}

\noindent and \eqref{f27preresinvlamhyp} written in the form
\begin{align*}
    & \tag{\ref{f27preresinvlamhyp}}
\ie (1+x +  i  0) =-\dfrac{1}{\ie (-x -  i  0)}, \quad x > 0 \,,
\end{align*}

\noindent  that  $\ie (1+ x +  i  0 ) \in 1 /(1 - i \Bb{R}_{>0})$
for arbitrary $x>0$. Therefore, by virtue of \eqref{f29preresinvlamhyp},
\begin{align*}\tag{\ref{f29preresinvlamhyp}}
    & \ie (1 + x - i  0 ) = \frac{\ie (  1+ x +  i  0 )}{ 1- 2 \ie (  1+ x +  i  0 ) } , \quad x>0\,,
\end{align*}

\noindent
and by the  property {\rm{(b)}} of $f$ in Lemma~{\hyperlink{hth1}{\ref*{th1}}},  which can be written as
\begin{align*}
    & \hspace{-0,1cm} f \left(-\dfrac{1}{(-1/z)}\right)\! =  \!f (z) \!=\! f \left(\dfrac{z}{1-2z}\right) \!=\! f \left(- \dfrac{1}{2+ (-1/z)}\right)  \, , \ \
      z \!\in \! 1 /(1 -  i  \Bb{R}_{>0}) \, ,
\end{align*}

\noindent and where one can set $ z = \ie (1+ x +  i  0 )  \in 1 /(1 -  i  \Bb{R}_{>0}) $, we derive that
\begin{align}\nonumber
     \Phi (1+ x -  i  0 )& = f \big(\ie (1 + x - i  0 )\big)= f\left(\frac{\ie (  1+ x +  i  0 )}{ 1- 2 \ie (  1+ x +  i  0 ) } \right)
      \\[0,25cm]    &=  f \big(\ie (1 + x + i  0 )\big)  = \Phi (1+ x +  i  0 )     \ , \quad  x>0 \,.
\label{ad19}\end{align}

\noindent
 At the same time the property {\rm{(a)}} of $f$ in Lemma~{\hyperlink{hth1}{\ref*{th1}}}
\begin{align*}
    &  f(z)=f(z+2)  \ , \quad  z=  \ie (-x -  i  0 ) \in -1 +  i  \Bb{R}_{>0} \ , \
\end{align*}

\noindent   together with \eqref{f14zpreresinvlamhyp}
\begin{align}\tag{\ref{f14zpreresinvlamhyp}}
 \ie (-x -  i  0 )=  - 1 +  i  \,
 \dfrac{\he  \left( {1}\big/{(1 + x)}\right)}{\he  \left(  {x}\big/{(1 + x)}\right)} \in -1 + i \Bb{R}_{>0}   \, , \ \  x > 0 \, ,
\end{align}

\noindent
 and \eqref{f28preresinvlamhyp}
\begin{align*}
    & \tag{\ref{f28preresinvlamhyp}}
\ie (-x +  i  0 ) =  \ie (-x -  i  0 ) + 2 \ , \quad  x>0\,,
\end{align*}

\noindent  yields
\begin{align}\nonumber
     \Phi (-x- i  0 )&=f \big(\ie (-x- i  0 )\big)=f \big(\ie (-x- i  0 ) + 2\big) \\[0,25cm]       &
     =f \big(\ie (-x+ i  0 )\big)   = \Phi (-x+ i  0 )     \ , \quad  x > 0 \ .
\label{ad20}\end{align}

\noindent  In view of \eqref{f33preresinvlamhyp}, written as
\begin{align*}
    \tag{\ref{f33preresinvlamhyp}}
    &  \ie \left(\Lambda\right) = \fet \ , \ \ \Lambda := (0,1)\cup \big(\Bb{C}\setminus\Bb{R}\big) \, , \
\end{align*}

\noindent and by the condition $f \in {\rm{Hol}} (\fet)$ of  Lemma~{\hyperlink{hth1}{\ref*{th1}}},  we have

\begin{align*}
    & \Phi (z) = f \big(\ie (z)\big)\,, \     z \in \Lambda \, ; \ \ie \in {\rm{Hol}} (\Lambda) \, ; \  \ie (\Lambda)=\fet \, ; \
    f \in {\rm{Hol}} (\fet) \, , \
\end{align*}

\vspace{0.25cm}
\noindent which yield that $\Phi  \in {\rm{Hol}} (\Lambda)$ (see \cite[p. 34, 2.4]{con}). Thus, by virtue of \eqref{ad19} and \eqref{ad20},  the values of $\Phi$ are the same on two sides of
each cut along $(-\infty, 0]$ and $[1, +\infty)$ of the set $\Lambda$,
\begin{align}\label{ad21}
    &  \Phi (-x- i  0 )= \Phi (-x+ i  0 )  \, , \  \Phi (1+ x -  i  0 ) = \Phi (1+ x +  i  0 ) \, , \ \  x > 0 \ ,
\end{align}

\noindent Applying the Morera theorem  (see \cite[p. 96]{lav}) we get
\begin{align}\label{ad22}
    &  \Phi \in  {\rm{Hol}}( \Bb{C}\setminus\{0, 1\}) \,.
\end{align}

\vspace{0,1cm}
In the notations \eqref{f29ainvlamhyp} for the remainders from singularities written in the form
\begin{align*}
   - \dfrac{1}{  \ie (z) -{\rm{sign}} (\im z)}     & = \dfrac{ i }{ \pi} \ln  |z| - \ie(\infty; z)    \, ;   \\[0,25cm]
  \ie (z)    & = \dfrac{i }{\pi} \ln\dfrac{1}{|z|}  - \ie(0; z)   \, ; \\[0,25cm]
- \dfrac{1}{\ie (z)}    &=  \dfrac{i}{\pi}\ln\dfrac{1}{|1-z|}   -  \ie(1; z)    \ , \quad
        z \in  (0,1)\cup \left(\Bb{C}\setminus\Bb{R}\right) \,,
 \tag{\ref{f29ainvlamhyp}}
\end{align*}

\vspace{0.25cm}
\noindent
the properties \eqref{f1th2}
\begin{align*}
    &      {\rm{(a)}} \ \ \big|\ie(\infty; z)\big| \leq 2 \, ,    \quad     |z| \geq 1/\varepsilon_{\!{\tn{\triangle}}} \ ; \\[0,25cm]
      &    {\rm{(b)}} \ \ \big|\ie(0; z)\big|  \leq 2 \, ,   \quad  \ \,  |z| \leq \varepsilon_{\!{\tn{\triangle}}} \ ;  \\[0,25cm]
      &    {\rm{(c)}} \ \ \big|\ie(1; z)\big|  \leq 2 \, ,  \quad \ \,   |z-1| \leq \varepsilon_{\!{\tn{\triangle}}}\ ; \quad \ \ \ \  z \in (0,1)\cup \left(\Bb{C}\setminus\Bb{R}\right)\,,
 \tag{\ref{f1th2}}
\end{align*}

\noindent mean that for arbitrary $z\in \Lambda$  we have
 \begin{align}\nonumber
    {\rm{(a)}}   &     &      \dfrac{1}{ \left| \ie (z)\! -\!\sigma (z)\right|} &  \leq 2 +  \dfrac{ \ln  |z| }{ \pi} \, ,
    &     &        |z|\! \geq\! 1/\varepsilon_{\!{\tn{\triangle}}} \, ; \\[0,25cm]  \nonumber
      {\rm{(b)}}  &     &      \left|\ie (z)\right| &    \leq 2 + \dfrac{1 }{\pi} \ln\dfrac{1}{|z|}\, ,      &     &     |z|\! \leq\! \varepsilon_{\!{\tn{\triangle}}} \, ;  \\[0,25cm]
    {\rm{(c)}}   &     &     \dfrac{1}{\left|\ie (z)\right|}    & \leq 2 + ({1}/{\pi})\ln\dfrac{1}{|1-z|} \, ,     &     &     |z-1| \leq \varepsilon_{\!{\tn{\triangle}}}\, ,
\label{ad23} \end{align}

%\vspace{-0,2cm}
\noindent where $ \sigma (z) := {\rm{sign}} (\im z)$, and that, in view of  \eqref{f3th2},%\vspace{-0,2cm}
\begin{align}
  &   \begin{array}{lll}
       {\rm{(a)}}   &    \ie (z)\! -\!\sigma (z)\to 0  \ ,   &    \quad      \Lambda\ni z \to \infty\,,    \\[0,3cm]
    {\rm{(b)}} &      \ie (z) \to \infty  \ ,   &    \quad     \Lambda\ni z \to 0\,,    \\[0,3cm]
    {\rm{(c)}} & \ie (z) \to 0  \ ,   &   \quad      \Lambda \ni z \to 1  \ .
      \end{array}
\label{ad24}\end{align}

\noindent
Substituting in  Lemma~{\hyperlink{hth1}{\ref*{th1}}}, 1)-3)
\begin{align*}
    &  \ie (z) \ , \  z\in \Lambda := (0,1)\cup \big(\Bb{C}\setminus\Bb{R}\big) \quad  \ \mbox{ in place of} \ \quad z \ , \  z\in \fet \stackrel{{\fo{\eqref{f33preresinvlamhyp}}}}{\vphantom{A} =} \ie (\Lambda) \ , \
\end{align*}

\vspace{0.25cm}
\noindent and letting  $\Lambda\ni z\to 0$ in 1), $\Lambda\ni z\to 1$ in 2), and $\Lambda\ni z\to \infty$ in 3),
we use \eqref{ad24}    and apply the inequalities \eqref{ad23} to deduce that

\begin{align*}
    &  1) \ \ \Lambda\ni z\to 0 \ \stackrel{{\fo{\eqref{ad24}(b)}}}{\vphantom{A^{A}}\Rightarrow } \ \fet \ni \ie (z) \to \infty
    \ \Rightarrow  \\[0,25cm]    &
      0 = \lim\limits_{{\fo{\fet \ni \ie (z) \to \infty}}}\left|f (\ie (z))\right| \exp \Big(- \pi\left( n_{\infty} +1\right) |\ie (z)|\Big)
    \\    &   \stackrel{{\fo{\eqref{ad23}(b)}}}{\geq }
    \limsp\limits_{{\fo{\Lambda\ni z\to 0}}}\left|\Phi (z)\right| \exp \left(- 2\pi\left( n_{\infty} +1\right)  - \left( n_{\infty} +1\right)\ln\dfrac{1}{|z|}\right) \\[0,25cm]    &
    = e^{{\fo{- 2\pi\left( n_{\infty} +1\right) }}}  \limsp\limits_{{\fo{\Lambda\ni z\to 0}}}\left|\Phi (z)\right| \cdot |z|^{n_{\infty} +1}
     \ \Rightarrow \  \lim\limits_{{\fo{\Lambda\ni z\to 0}}}\left|\Phi (z)\right| \cdot |z|^{n_{\infty} +1}=0 \ ;
\end{align*}

\begin{align*}
    &  2) \ \ \Lambda\ni z\to 1 \ \stackrel{{\fo{\eqref{ad24}(c)}}}{\vphantom{A^{A}}\Rightarrow } \ \fet \ni \ie (z) \to 0  \
    \Rightarrow  \\[0,25cm]    &
      0 = \lim\limits_{{\fo{\fet \ni \ie (z) \to 0}}}\left|f (\ie (z))\right| \exp \Big(- \pi\left( n_{0} +1\right) |\ie (z)|\Big)
    \\    &   \stackrel{{\fo{\eqref{ad23}(c)}}}{\geq }
    \limsp\limits_{{\fo{\Lambda\ni z\to 1}}}\left|\Phi (z)\right| \exp \left(- 2\pi\left( n_{0} +1\right)  - \left( n_{0} +1\right)\ln\dfrac{1}{|1-z|}\right) \\[0,25cm]    &
    = e^{{\fo{- 2\pi\left( n_{0} +1\right) }}}  \limsp\limits_{{\fo{\Lambda\ni z\to 1}}}\left|\Phi (z)\right| \cdot |1-z|^{n_{0} +1}
     \\[0,25cm]    & \Rightarrow \  \lim\limits_{{\fo{\Lambda\ni z\to 1}}}\left|\Phi (z)\right| \cdot |1-z|^{n_{0} +1}=0 \ ;
\end{align*}
\begin{align*}
    &  3) \ \ \Lambda\ni z\to \infty \ , \  {\rm{sign}} (\im z) = \sigma \in \{1, -1\} \ \stackrel{{\fo{\eqref{ad24}(a)}}}{\vphantom{A^{A}}\Rightarrow } \ \fet \ni \ie (z) \to \sigma  \
    \Rightarrow  \\[0,25cm]    &
      0 = \lim\limits_{{\fo{\fet \ni \ie (z) \to \sigma }}}\left|f (\ie (z))\right| \exp \left(- \dfrac{\pi\left( n_{1} +1\right) }{|\ie (z)-\sigma|}\right)
    \\    &   \stackrel{{\fo{\eqref{ad23}(a)}}}{\geq }
    \limsp\limits_{{\fo{\Lambda\ni z\to \infty}}}\left|\Phi (z)\right| \exp \left(- 2\pi\left( n_{1} +1\right)  -
    \left( n_{1} +1\right)\ln|z|\right) \\[0,25cm]    &
    = e^{{\fo{- 2\pi\left( n_{1} +1\right) }}}  \limsp\limits_{{\fo{\Lambda\ni z\to \infty}}}\left|\Phi (z)\right| \cdot |z|^{-n_{1} -1}
    \\[0,25cm]    &\Rightarrow \  \lim\limits_{{\fo{\Lambda\ni z\to \infty}}}\left|\Phi (z)\right| \cdot |z|^{-n_{1} -1}=0 \ ;
\end{align*}

\vspace{0.25cm}
\noindent i.e.,
\begin{align}\nonumber
    &  {\rm{(1)}}  &     &    |z|^{{\fo{n_{\infty} +1}}} \left|\Phi (z)\right|\! \to 0\!\, ,   &     \Lambda\! \ni\! z \!\to\! 0 \, ,
    \\[0,3cm] \nonumber     &
    {\rm{(2)}} &     &  |1\!-\!z|^{{\fo{n_{0}\!+\!1}}}\left|\Phi (z)\right|\!\to\! 0\, ,   &    \Lambda \!\ni \!z \!\to\! 1\, ,  \\[0,3cm]    &
    {\rm{(3)}}&     & |z|^{ {\fo{- n_{ 1} -1}}}\left|\Phi (z)\right| \to 0\ ,   &     \Lambda \ni z \to  \infty \, .
\label{ad25}\end{align}

%\vspace{-0,2cm}
\noindent By the Riemann theorem about removable singularities (see \cite[p. 103]{con}) we obtain that the function
$\Phi_{1} (z):= z^{n_{\infty}} (1-z)^{n_{0}} \Phi (z)$ is holomorphic at $0$ and at $1$ and therefore it
is an entire function satisfying, by virtue of (3) in \eqref{ad25}, the equality
%\vspace{-0,4cm}
\begin{align}\label{ad26}
    &\hspace{-0,3cm} \lim\limits_{{\fo{\Lambda \ni |z| \to \infty}}} \  \dfrac{ \Phi_{1} (z) }{
    |z|^{{\fo{n_{0} +n_{\infty} +  n_{ 1} +1}}}}=  \lim\limits_{{\fo{\Lambda \ni |z| \to \infty}}} \  \dfrac{z^{{\fo{n_{\infty}}}} (1-z)^{{\fo{n_{0}}}} \Phi (z) }{
    |z|^{{\fo{n_{0} +n_{\infty} +  n_{ 1} +1}}}} = 0 \,.\hspace{-0,2cm}
\end{align}

%\vspace{-0,2cm}
\noindent
By the continuity of $\Phi_{1}$, we get the existence of $C \in \Bb{R}_{>0}$ such that $|\Phi_{1} (z)| \leq C (1+|z|)^{n_{0}
 +n_{\infty} +  n_{ 1} +1}$, $z\in \Bb{C}$,  which by the extended version of the Liouville theorem (see \cite[p. 2, Thm. 1]{lev2}) yields that $\Phi_{1} (z)$ is an algebraic polynomial of degree at most $n_{\infty} + n_{0} + n_{1}+1$. But the relationship \eqref{ad26} proves that actually its degree cannot exceed $n_{\infty} + n_{0} + n_{1}$. Lemma~{\hyperlink{hth1}{\ref*{th1}}} is proved.

}}\end{clash}
\end{subequations}

 \begin{subequations}

\begin{clash}{\rm{\hypertarget{d14}{}\label{case14}\hspace{0,00cm}{\hyperlink{bd14}{$\uparrow$}} \ Observe that the relations \eqref{ad25} and \eqref{ad21}
\begin{align*}\tag{\ref{ad21}}
    &  \Phi (-x- i  0 )= \Phi (-x+ i  0 )  \, , \  \Phi (1+ x -  i  0 ) = \Phi (1+ x +  i  0 ) \, , \ \  x > 0 \ ,
\end{align*}

\noindent prove that the function
\begin{align*}\tag{\ref{f1bd3}}
    &   \Phi (z) = f \big(\ie (z)\big)\,,   &     &        z \in \Lambda := (0,1)\cup \left(\Bb{C}\setminus\Bb{R}\right)  \, ,
\end{align*}

\noindent satisfies  the conditions of Lemma~{\hyperlink{hbd3lem1}{\ref*{bd3lem1}}}, while its conclusion  has already been deduced from its conditions
in the paragraph after the formulas \eqref{ad25}.

\vspace{0.35cm} Prove now that the function
\begin{align*}\tag{\ref{f2bd3}}
    &    \Psi (z) = f\Bigg(\ie \bigg(\dfrac{1}{1- z^{2}}\bigg)\Bigg) =  \Phi \left(\dfrac{1}{1- z^{2}}\right) \,,   &     &       z \in \Bb{H}  \, ,
\end{align*}

\noindent satisfies  the conditions of Lemma~{\hyperlink{hbd3lem2}{\ref*{bd3lem2}}}, where the function $\Phi$ satisfies  the conditions of Lemma~{\hyperlink{hbd3lem1}{\ref*{bd3lem1}}}.

\vspace{0.15cm}
The  formula for $"z_{2}"$ in  \cite[item 4, p.155]{lav} with  $a = 1/2$  and $1/2 - z$ in place of $z $ states that $\sqrt{1 - 1/z}$ maps $(0, 1) \cup \left(\Bb{C}\setminus\Bb{R}\right)$ one-to-one onto $\Bb{H} $. Therefore the inverse mapping $1/ (1-z^{2})= 1/ (1+(z/i)^{2})$ maps  $\Bb{H} $ one-to-one onto  $(0, 1) \cup \left(\Bb{C}\setminus\Bb{R}\right)$. It follows easily from
\begin{align*}
   & \dfrac{1}{1 - \left(x + i \varepsilon\right)^{2}} = \dfrac{1\,+ \,\varepsilon^{2}\, -\, x^{2} \,+\, 2 i \hspace{0,05cm} \varepsilon x}{\left(1+ \varepsilon^{2} - x^{2}\right)^{2} + 4 \varepsilon^{2} x^{2}} \ , \qquad
   x \in \Bb{R} \ , \ \varepsilon > 0 \ ,
\end{align*}

\noindent that\vspace{-0,3cm}
\begin{align*}
    \begin{array}{rlcl}
   x \in  (0,1) \ ,    &   \Bb{H} \ni z \to \,\pm \,x    &  \quad \Leftrightarrow \quad  &
    \,\pm \, \Bb{H} \ni  \dfrac{1}{1- z^{2}} \to \dfrac{1}{1- x^{2}} \in
     (1, + \infty)\,,
     \\[0,6cm]
 x \in  (1,+\infty)   \ ,   &    \Bb{H} \ni z \to \,\pm \,x    &  \quad \Leftrightarrow  \quad  &
   \,\pm \, \Bb{H} \ni  \dfrac{1}{1- z^{2}} \to \dfrac{1}{1- x^{2}} \in
      (- \infty, 0) \ ,
     \end{array}
\end{align*}

\noindent which in the notation $\eurm{s} (z) := 1/(1-z^{2})$ can symbolically be written as follows
\begin{align*}
    \eurm{s} \Big( (0,1)+ i 0\Big)  & = (1, + \infty) + i 0  \ ,   &
     \eurm{s} \Big( (-1,0)+ i 0\Big)  & = (1, + \infty) - i 0  \ , \\
       \eurm{s} \Big( (1,+\infty)+ i 0\Big)  & = (- \infty, 0) + i 0  \ ,  &
        \eurm{s} \Big( (- \infty, -1)+ i 0\Big)  & = (- \infty, 0) - i 0  \ .
\end{align*}

\noindent Therefore the condition of Lemma~{\hyperlink{hbd3lem1}{\ref*{bd3lem1}}} about the possibility to extend  $\Phi  \in {\rm{Hol}} \left((0, 1) \cup \left(\Bb{C}\setminus\Bb{R}\right)\right)$ continuously from $\Bb{H}$ to  $(\Bb{H}\cup\Bb{R} ) \setminus [\,0, 1]$ and
  from $-\Bb{H}$
  to $(- \Bb{H}\cup\Bb{R} ) \setminus [\,0, 1]$
such that
\begin{align*}
    &   \Phi (-x- i  0 )= \Phi (-x+ i  0 )  \, , \  \Phi (1+ x -  i  0 ) = \Phi (1+ x +  i  0 ) \, , \ \  x > 0 \, ,
\end{align*}

\noindent for the function
\begin{align*} &
    \Psi (z) = \Phi\left(\dfrac{1}{1- z^{2}}\right) \in {\rm{Hol}} \left(\Bb{H}\right) \, , \quad   \Psi: \Bb{H} \mapsto  \Bb{C}\,,
\end{align*}

\noindent  is equivalent to the possibility of extending   $\Psi \in {\rm{Hol}} \left(\Bb{H}\right)$ continuously to the set  $(\Bb{H}\cup\Bb{R} ) \setminus \{-1, 0, 1\}$  such that
\begin{align*}\tag{\ref{f1bd3lem2}} &
\Psi (x+ i  0) =  \Psi (-x+ i  0) \ , \quad  x \in \Bb{R}\setminus\{-1, 0, 1\}  \, ,
\hypertarget{hf1bd3lem2}{}\end{align*}

\noindent holds. Furthermore, the conditions (1), (2) and (3) of  Lemma~{\hyperlink{hbd3lem1}{\ref*{bd3lem1}}}
\begin{align}& {\rm{(1)}}&      &    & |z|^{ {\fo{- n_{ 1} -1}}}\left|\Phi (z)\right| \to 0\, ,   &     &     \Lambda \ni z \to  \infty \,,     \label{f1bd12}\\[0,3cm]
      &  {\rm{(2)}}  &     &     &    |z|^{{\fo{n_{\infty} +1}}} \left|\Phi (z)\right| \to 0\, ,   &    &      \Lambda\! \ni\! z \!\to\! 0 \, ,
    \label{f2bd12} \\[0,3cm]     &
    {\rm{(3)}} &     &     &  |1\!-\!z|^{{\fo{n_{0}\!+\!1}}}\left|\Phi (z)\right|\to 0\, ,   &    &     \Lambda \!\ni \!z \!\to\! 1\, ,
\label{f3bd12}\end{align}

\noindent  after substituting $1/(1-z^{2})$ instead of $z$ lead to  the conditions (1), (2) and (3) of  Lemma~{\hyperlink{hbd3lem2}{\ref*{bd3lem2}}}
\begin{align}\label{f4bd12}\hypertarget{hf4bd12}{}
&{\rm{(1)}}&      &    & |z|^{ {\fo{- 2n_{ \infty} -2}}}\left|\Psi (z)\right| \to 0\, ,            &   &       \Bb{H} \! \ni\! z \!\to\!   \infty \,, &   & \\[0,3cm]\label{f5bd12}\hypertarget{hf5bd12}{}
&{\rm{(2)}}&     &     &    |z|^{{\fo{\,2n_{0}+2}}}  \left|\Psi (z)\right|\! \to 0\, ,        &    &      \Bb{H}\! \ni\! z \!\to\! 0 \, ,  &   & \\[0,3cm]      &{\rm{(3)}} &     &     &  |\sigma\!-\!z|^{{\fo{\,n_{1}+1}}}\left|\Psi (z)\right|\to 0\, ,   &    &     \Bb{H} \!\ni \!z \!\to\! \sigma\, ,  &   &  \sigma\in\{1,-1\}\, .
\label{f6bd12}\hypertarget{hf6bd12}{}\end{align}

\noindent Actually, since
\begin{align*}
    \lim\limits_{{\fo{\Bb{H}\! \ni \! z \!\to \!\infty}}} \dfrac{1}{1-z^{2}} = 0   \stackrel{{\fo{\eqref{f2bd12}}}}{\vphantom{A}\Rightarrow}
  0 &  =\lim\limits_{{\fo{\Bb{H}\! \ni \! z \!\to \!\infty}}} \left|\dfrac{1}{1-z^{2}}\right|^{{\fo{n_{\infty} +1}}} \left|\Phi \left(\dfrac{1}{1-z^{2}}\right)\right| \\[0,2cm]    &     =
   \lim\limits_{{\fo{\Bb{H}\! \ni \! z \!\to \!\infty}}}  \dfrac{\left|\Phi \left(\dfrac{1}{1-z^{2}}\right)\right|}{|z|^{{\fo{2 n_{\infty} +2}}}} =
     \lim\limits_{{\fo{\Bb{H}\! \ni \! z \!\to \!\infty}}}  \dfrac{\left|\Psi (z)\right|}{|z|^{{\fo{2 n_{\infty} +2}}}}\ ,
\end{align*}

\vspace{0.25cm}
\noindent   we obtain that \eqref{f4bd12} holds, whereas
\begin{align*}
     \lim\limits_{{\fo{\Bb{H}\! \ni \! z \!\to \!0}}}\dfrac{1}{1-z^{2}} = 1  \   \stackrel{{\fo{\eqref{f3bd12}}}}{\vphantom{A}\Rightarrow}  \
  0 & =  \lim\limits_{{\fo{\Bb{H}\! \ni \! z \!\to \!0}}} \left|1\!-\!\dfrac{1}{1-z^{2}}\right|^{{\fo{n_{0}\!+\!1}}}\left|\Phi \left(\dfrac{1}{1-z^{2}}\right)\right| \\[0,2cm]    &     =
   \lim\limits_{{\fo{\Bb{H}\! \ni \! z \!\to \!0}}} \left|\dfrac{z^{2}}{1-z^{2}}\right|^{{\fo{n_{0}\!+\!1}}}\left|\Phi \left(\dfrac{1}{1-z^{2}}\right)\right|\\[0,2cm]    &     =
   \lim\limits_{{\fo{\Bb{H}\! \ni \! z \!\to \!0}}} |z|^{{\fo{2n_{0}\!+\!2}}}  \left|\Psi \left(z\right)\right| \, ,
\end{align*}

\noindent   proves \eqref{f5bd12}. Finally,  for arbitrary $\sigma \in \{1, -1\}$ we have
\begin{align*}
    \lim\limits_{{\fo{\Bb{H}\! \ni \! z \!\to \!\sigma}}} \dfrac{1}{1-z^{2}} = \infty  \  \   \stackrel{{\fo{\eqref{f1bd12}}}}{\vphantom{A}\Rightarrow}  \ \   0 &   =   \lim\limits_{{\fo{\Bb{H}\! \ni \! z \!\to \!\sigma}}}  \dfrac{\left|\Phi \left(\dfrac{1}{1-z^{2}}\right)\right| }{\left| \dfrac{1}{1-z^{2}}\right|^{ {\fo{ n_{ 1} +1}}}} = \\[0,2cm]  & =
  \lim\limits_{{\fo{\Bb{H}\! \ni \! z \!\to \!\sigma}}}  \left|1-z^{2}\right|^{ {\fo{ n_{ 1} +1}}}  \left|\Psi \left(z\right)\right|
  \\[0,2cm]  & =
2^{ {\fo{ n_{ 1} +1}}}  \lim\limits_{{\fo{\Bb{H}\! \ni \! z \!\to \!\sigma}}}  \left|z-\sigma \right|^{ {\fo{ n_{ 1} +1}}}  \left|\Psi \left(z\right)\right|   \, ,
 \end{align*}

\noindent and hence \eqref{f6bd12} is true. Thus, the function $\Psi $  satisfies  all the conditions of Lemma~{\hyperlink{hbd3lem2}{\ref*{bd3lem2}}}.

\vspace{0.25cm} We  give now an independent proof of Lemma~{\hyperlink{hbd3lem2}{\ref*{bd3lem2}}}.

\vspace{0.15cm}{\emph{Proof of Lemma~{\hyperlink{hbd3lem2}{\ref*{bd3lem2}}}.}} \ \
By virtue of the condition ({\hyperlink{hf1bd3lem2}{\ref*{f1bd3lem2}}}), the Morera theorem applied to the function
\begin{align}\label{f0bd12}
     & \widehat{\Psi} (z) := \left\{
                     \begin{array}{ll}
                     \Psi (z)\,, & \hbox{if} \ z \in (\Bb{H}\cup\Bb{R} ) \setminus \{-1, 0, 1\}\,, \\[0.25cm]
                     \Psi (-z)\, , & \hbox{if} \ z \in (-\Bb{H}\cup\Bb{R} ) \setminus \{-1, 0, 1\}\,,
                     \end{array}
                   \right.
\end{align}

\noindent gives $\widehat{\Psi} \in {\rm{Hol}} \left(\Bb{C}\setminus \{-1, 0, 1\}\right)$. In other words,   the function
  $\Psi \in {\rm{Hol}} (\Bb{H})$ permits a holomorphic extension $\Psi := \widehat{\Psi}$   from $\Bb{H}$ to $\Bb{C}\setminus\{-1, 0, 1\}$ satisfying
\begin{align}\label{f7bd12} \Psi \in {\rm{Hol}} \left(\Bb{C}\setminus \{-1, 0, 1\}\right) \, , \quad
  \Psi(z) =  \Psi (-z)  \ , \ z \in \Bb{C}\setminus \{-1, 0, 1\}  \ ,
\end{align}

\noindent
  which means that the obtained analytic extension  of $\Psi $ is even on $\Bb{C}\setminus\{-1, 0, 1\}$.

  \vspace{0,1cm}
  By the Riemann theorem about removable singularities (see \cite[p. 103]{con}) it follows from  the conditions
  (\hyperlink{hf5bd12}{2})   and   (\hyperlink{hf5bd12}{3}) of  Lemma~{\hyperlink{hbd3lem2}{\ref*{bd3lem2}}}
that the function
$\Psi_{1} (z):= (1-z^{2})^{n_{1}} z^{2n_{0}+1}  \Psi (z)$ is holomorphic at $0$, $-1$ and at $1$,  and consequently, in view of \eqref{f7bd12}, it
is an odd entire function satisfying, by virtue of   the condition   (\hyperlink{hf4bd12}{1})  in Lemma~{\hyperlink{hbd3lem2}{\ref*{bd3lem2}}}, the relations
\begin{align}\label{f8bd12} &
\begin{array}{ll}
{\rm{(a)}}   &    \Psi_{1} (-z) = - \Psi_{1} (z) \, , \ \ \Psi_{1} (z):= (1-z^{2})^{{\fo{n_{1}}}} z^{{\fo{2n_{0}+1}}}  \Psi (z) \, , \  \  z \in \Bb{C}\,;
\\[0,4cm]
{\rm{(b)}}   &  |z|^{{\fo{- 2 n_{\infty}-2 - 2n_{1} - 2n_{0}-1}}}\, \left|\Psi_{1} (z)\right| \to 0  \ , \ \   z \to \infty \ .
\end{array}
    \ \
\end{align}

\noindent By the properties \eqref{f8bd12}(b) and $ \Psi_{1} \in {\rm{Hol}} (\Bb{C})$, there exists a constant $C \in \Bb{R}_{>0}$ such that
\begin{align*}
    & |\Psi_{1} (z)| \leq C (1+|z|)^{2 n_{\infty}+2 + 2n_{1} + 2n_{0}+1} \ , \quad  z\in \Bb{C}\,,
\end{align*}
\noindent  which by the extended version of the Liouville theorem (see \cite[p. 2]{lev2}) yields that $\Psi_{1} (z)$ is an algebraic polynomial $P$ of degree at most $2 n_{\infty}+2 + 2n_{1} + 2n_{0}+1$.
But the relationship \eqref{f8bd12}(b) proves that actually  the degree of $P$ cannot exceed $2 n_{\infty}+2 + 2n_{1} + 2n_{0}$. In view of the property \eqref{f8bd12}(a),
$P$ is an odd algebraic polynomial ($P=\Psi_{1}$), and hence its degree cannot be greater than   $2 n_{\infty}+1 + 2n_{1} + 2n_{0}$. Consequently,  there exists an algebraic polynomial of degree $\leq n_{\infty} + n_{1} + n_{0}$ such that
\begin{align*}
    &  \Psi_{1} (z) = P (z) = z \, Q \left(z^{2}\right) \, , \  \  z \in \Bb{C}\,,
\end{align*}

\noindent from which it follows that
\begin{align*}
    &  \Psi (z) =   \dfrac{\Psi_{1} (z)}{(1-z^{2})^{{\fo{n_{1}}}} z^{{\fo{2n_{0}+1}}} } =
    \dfrac{ z \, Q \left(z^{2}\right)}{(1-z^{2})^{{\fo{n_{1}}}} z^{{\fo{2n_{0}+1}}} } =
    \dfrac{  Q \left(z^{2}\right)}{(1-z^{2})^{{\fo{n_{1}}}} z^{{\fo{2n_{0}}}} } \, , \  \  z \in \Bb{C}\, .
\end{align*}

\noindent Thus, \eqref{f2bd3lem2} holds and Lemma~{\hyperlink{hbd3lem2}{\ref*{bd3lem2}}} is proved. $\square$
}}\end{clash}
\end{subequations}

\subsection[\hspace{-0,25cm}. \hspace{0,075cm}Notes on Section~\ref{theta}]{\hspace{-0,11cm}{\bf{.}} Notes on Section~\ref{theta}}

\begin{subequations}
\vspace{-0.15cm}
\begin{clash}{\rm{\hypertarget{d16}{}\label{case16}\hspace{0,00cm}{\hyperlink{bd16}{$\uparrow$}} \ Several additional to \eqref{f9elemtheta} relationships between the theta functions
will be derived in this note.

\vspace{0.25cm}
For arbitrary $z \in \Bb{H}$ we have
\begin{align*}
    \Theta_{2} (z+2)       &   \stackrel{{\fo{\eqref{f9elemtheta}(d)}}}{\vphantom{A}=}       e^{i\pi/4}\Theta_{2}(z+1)       &
    \hspace{-1,95cm}  \stackrel{{\fo{\eqref{f9elemtheta}(d)}}}{\vphantom{A}=}\,  i \Theta_{2} (z) \ , \\
   \Theta_{3} (z+2)       &  \stackrel{{\fo{\eqref{f9elemtheta}(e)}}}{\vphantom{A}=}
     \Theta_{4}(z+1)      &    \hspace{-1,95cm}   \stackrel{{\fo{\eqref{f9elemtheta}(g)}}}{\vphantom{A}=}\,   \Theta_{3} (z) \ , \\
    \Theta_{4} (z+2)       &  \stackrel{{\fo{\eqref{f9elemtheta}(g)}}}{\vphantom{A}=}         \Theta_{3}(z+1)       &
    \hspace{-1,95cm}   \stackrel{{\fo{\eqref{f9elemtheta}(e)}}}{\vphantom{A}=}           \Theta_{4} (z) \ ,
\end{align*}

\noindent and therefore
\begin{align}\label{adt27}
    &     \hspace{-0.45cm} \begin{array}{llll}
{\rm{(a)}}\	 \Theta_{2} (z+2) &= \, i\, \Theta_{2} (z)\, ,
& \quad {\rm{(d)}}\ \Theta_{3} (z+2) &= \Theta_{3} (z)	\, ,\
\\[0,1cm]
{\rm{(b)}} \  \Theta_{2} (z+2)^{2}  &=\! -\Theta_{2} (z)^{2}	\, ,
&\quad {\rm{(e)}} \ 	\Theta_{4} (z+2) &= \Theta_{4} (z) \, ,
\\[0,1cm]
 {\rm{(c)}} \  \Theta_{2} (z+2)^{4}  &=\,\,\, \Theta_{2} (z)^{4}	\, ,
 &
\end{array}\hspace{-0.2cm}
\end{align}

\noindent for any $z \in \Bb{H}$.

Since $(-\pi/2, \pi/2)\pm (-\pi/2, \pi/2)\subset (-\pi, \pi)$ then for the principal branch of the square root
\begin{align*}
    &  \sqrt{z} = \exp \left(\dfrac{1}{2}\ln |z| + \dfrac{i}{2} \Arg (z)\right)  \ , \quad
 \Arg (z) \in (-\pi, \pi)  \ , \     z \in \Bb{C}\setminus (-\infty, 0] \,,
\end{align*}

\noindent we have
\begin{align*}
    &  \sqrt{z_{1}} \sqrt{z_{2}} = \sqrt{z_{1}z_{2}}  \ , \quad  \dfrac{\sqrt{z_{1}}}{\sqrt{z_{2}}} = \sqrt{\dfrac{z_{1}}{z_{1}}}  \ ,
    \qquad  z_{1} \, , \,  z_{2} \in \Bb{H} \ ,
\end{align*}

\noindent
and for arbitrary $z \in \Bb{H}$ we obtain
\begin{align*}
    \Theta_{2} (z) &         \stackrel{{\fo{\eqref{f9elemtheta}(c)}}}{\vphantom{A}=}   (i/z)^{1/2}\Theta_{4}(-1/z)
  \stackrel{{\fo{\eqref{adt27}(e)}}}{\vphantom{A}=}\, (i/z)^{1/2}\Theta_{4}(2-1/z) \\[0,3cm]    &    \stackrel{{\fo{\eqref{f9elemtheta}(a)}}}{\vphantom{A}=}\,
  \left(\dfrac{i}{z}\right)^{1/2}   \left(\dfrac{i}{2-1/z}\right)^{1/2} \Theta_{2} \left(- \dfrac{1}{2-1/z}\right)  \\[0,25cm]    & \ \quad =
\dfrac{i^{1/2} i^{1/2}}{ z^{1/2}  (2-1/z)^{1/2}}  \, \Theta_{2}\left(\dfrac{z}{1-2z}\right)
       =  \dfrac{i}{(2z-1)^{1/2} } \, \Theta_{2}\left(\dfrac{z}{1-2z}\right)      \, ,
\end{align*}

\begin{align*}
    \Theta_{3} (z) &         \stackrel{{\fo{\eqref{f9elemtheta}(b)}}}{\vphantom{A}=}   (i/z)^{1/2}\Theta_{3}(-1/z)
  \stackrel{{\fo{\eqref{adt27}(d)}}}{\vphantom{A}=}\, (i/z)^{1/2}\Theta_{3}(2-1/z) \\[0,3cm]    &    \stackrel{{\fo{\eqref{f9elemtheta}(b)}}}{\vphantom{A}=}\,
  \left(\dfrac{i}{z}\right)^{1/2}   \left(\dfrac{i}{2-1/z}\right)^{1/2} \Theta_{3} \left(- \dfrac{1}{2-1/z}\right)  \\[0,25cm]    & \ \quad =
\dfrac{i^{1/2} i^{1/2}}{ z^{1/2}  (2-1/z)^{1/2}}  \, \Theta_{3}\left(\dfrac{z}{1-2z}\right)
       =  \dfrac{i}{(2z-1)^{1/2}} \, \Theta_{3}\left(\dfrac{z}{1-2z}\right)      \, ,
\end{align*}

\begin{align*}
    \Theta_{4} (z) &         \stackrel{{\fo{\eqref{f9elemtheta}(a)}}}{\vphantom{A}=}   (i/z)^{1/2}\Theta_{2}(-1/z)
  \stackrel{{\fo{\eqref{adt27}(a)}}}{\vphantom{A}=}\,\, - \,  i (i/z)^{1/2}\Theta_{2}(2-1/z) \\[0,3cm]    &    \stackrel{{\fo{\eqref{f9elemtheta}(c)}}}{\vphantom{A}=}\,
\, - \, i  \left(\dfrac{i}{z}\right)^{1/2}   \left(\dfrac{i}{2-1/z}\right)^{1/2} \Theta_{4} \left(- \dfrac{1}{2-1/z}\right)  \\[0,25cm]    & \ \quad =
\, - \, \dfrac{i \, i^{1/2} i^{1/2}}{ z^{1/2}  (2-1/z)^{1/2}}  \, \Theta_{4}\left(\dfrac{z}{1-2z}\right)
       =  \dfrac{1}{(2z-1)^{1/2}} \, \Theta_{4}\left(\dfrac{z}{1-2z}\right)      \, ,
\end{align*}

\noindent from which  it follows that
\begin{align}\label{adt28}
&\hspace{-0.45cm} \begin{array}{llll}
{\rm{(a)}}\	\ \Theta_{2} (z)   &= i\,\dfrac{ \Theta_{2}\left(\dfrac{z}{1-2z}\right)}{(2z-1)^{1/2} } \,   \, ,
 &\quad {\rm{(d)}}\ \ \Theta_{2} (z)^{2}   &= -\,\dfrac{\Theta_{2}\left(\dfrac{z}{1-2z}\right)^{\!2}}{(2z-1) } \,   \, ,
\\[0,4cm]
{\rm{(b)}} \  \ 	\Theta_{3} (z)   &=i\, \dfrac{ \Theta_{3}\left(\dfrac{z}{1-2z}\right)}{(2z-1)^{1/2} } \, \, ,
&\quad {\rm{(e)}} \  \ \Theta_{3} (z)^{2}   &= -\,\dfrac{ \Theta_{3}\left(\dfrac{z}{1-2z}\right)^{\!2}}{(2z-1) } \,   \,  ,
\\[0,4cm]
{\rm{(c)}} \  \ \Theta_{4} (z)   &=\, \,\dfrac{\Theta_{4}\left(\dfrac{z}{1-2z}\right)}{(2z-1)^{1/2} } \, \, ,
& \quad {\rm{(f)}} \ \  \Theta_{4} (z)^{2}   &=\, \, \, \, \, \dfrac{\Theta_{4}\left(\dfrac{z}{1-2z}\right)^{\!2}}{(2z-1) } \,    \, ,
\end{array}\hspace{-0.2cm}
\end{align}

\noindent for any $z \in \Bb{H}$. Thus,
\begin{align}\label{adt28a}
    &  \Theta_{k} (z)^{4}   = \dfrac{\Theta_{k}\left(\dfrac{z}{1-2z}\right)^{\!4}}{(2z-1)^{2} } \ , \quad  2 \leq k \leq 4 \ , \ z \in \Bb{H}\,.
\end{align}

\vspace{0.25cm} In the similar manner we deduce  from \eqref{f9elemtheta} that  for arbitrary $z \in \Bb{H}$
\begin{align*}
    &   \Theta_{2} (z)  \,\stackrel{{\fo{\eqref{f9elemtheta}(d)}}}{\vphantom{A}=}\,
    e^{i\pi/4}\Theta_2(z-1)
      \,\stackrel{{\fo{\eqref{f9elemtheta}(c)}}}{\vphantom{A}=}\,  e^{i\pi/4} (i/(z-1))^{1/2}\Theta_{4}(-1/(z-1))    \\
       &  = i \dfrac{\Theta_{4} \left(\dfrac{1}{1-z}\right)}{(z-1)^{1/2}} \ \Rightarrow \
       i \dfrac{\Theta_{4} \left(-\dfrac{1}{1+z}\right)}{(z+1)^{1/2}} =  \Theta_{2} (z+2)
       \,\stackrel{{\fo{\eqref{adt27}(a)}}}{\vphantom{A}=}\, i \Theta_{2} (z)  \ \Rightarrow \   \\    &
    \Theta_{2} (z)= i \dfrac{\Theta_{4} \left(\dfrac{1}{1-z}\right)}{(z-1)^{1/2}} = \dfrac{\Theta_{4} \left(-\dfrac{1}{1+z}\right)}{(z+1)^{1/2}} \ ,
\end{align*}

\begin{align*}
    &   \Theta_{3} (z)  \,       \stackrel{{\fo{\eqref{f9elemtheta}(e)}}}{\vphantom{A}=} \,\Theta_4(z-1)
      \,\stackrel{{\fo{\eqref{f9elemtheta}(a)}}}{\vphantom{A}=}\,   (i/(z-1))^{1/2}\Theta_{2}(-1/(z-1))    \\
       &  = i^{1/2} \dfrac{\Theta_{2} \left(\dfrac{1}{1-z}\right)}{(z-1)^{1/2}} \ \Rightarrow \
       i^{1/2} \dfrac{\Theta_{2} \left(-\dfrac{1}{1+z}\right)}{(z+1)^{1/2}} =  \Theta_{3} (z+2)
       \,\stackrel{{\fo{\eqref{adt27}(d)}}}{\vphantom{A}=}\,  \Theta_{3} (z)  \ \Rightarrow \   \\    &
    \Theta_{3} (z)= i^{1/2} \dfrac{\Theta_{2} \left(\dfrac{1}{1-z}\right)}{(z-1)^{1/2}} =i^{1/2} \dfrac{\Theta_{2} \left(-\dfrac{1}{1+z}\right)}{(z+1)^{1/2}} \ ,
\end{align*}

\begin{align*}
    &   \Theta_{4} (z)  \,       \stackrel{{\fo{\eqref{f9elemtheta}(g)}}}{\vphantom{A}=} \,\Theta_3(z-1)
      \,\stackrel{{\fo{\eqref{f9elemtheta}(a)}}}{\vphantom{A}=}\,   (i/(z-1))^{1/2}\Theta_{3}(-1/(z-1))    \\
       &  = i^{1/2} \dfrac{\Theta_{3} \left(\dfrac{1}{1-z}\right)}{(z-1)^{1/2}} \ \Rightarrow \
       i^{1/2} \dfrac{\Theta_{3} \left(-\dfrac{1}{1+z}\right)}{(z+1)^{1/2}} =  \Theta_{4} (z+2)
       \,\stackrel{{\fo{\eqref{adt27}(e)}}}{\vphantom{A}=}\,  \Theta_{4} (z)  \ \Rightarrow \   \\    &
    \Theta_{4} (z)= i^{1/2} \dfrac{\Theta_{3} \left(\dfrac{1}{1-z}\right)}{(z-1)^{1/2}} =i^{1/2} \dfrac{\Theta_{3} \left(-\dfrac{1}{1+z}\right)}{(z+1)^{1/2}} \ ,
\end{align*}

\vspace{0.35cm}
\noindent from which  it follows that

\begin{align}\label{adt29}
&\hspace{-0.45cm} \begin{array}{l}
{\rm{(a)}}\	\ \ \Theta_{2} (z)= i \ \dfrac{\Theta_{4} \left(- \dfrac{1}{z-1}\right)}{(z-1)^{1/2}}\ \  = \dfrac{\Theta_{4} \left(-\dfrac{1}{1+z}\right)}{(z+1)^{1/2}} \ ,
\\[0,4cm]
{\rm{(b)}} \  \ 	  \Theta_{3} (z)= i^{1/2} \dfrac{\Theta_{2} \left(- \dfrac{1}{z-1}\right)}{(z-1)^{1/2}} =i^{1/2} \dfrac{\Theta_{2} \left(-\dfrac{1}{1+z}\right)}{(z+1)^{1/2}} \ ,
\\[0,4cm]
{\rm{(c)}} \  \  \Theta_{4} (z)= i^{1/2} \dfrac{\Theta_{3} \left(- \dfrac{1}{z-1}\right)}{(z-1)^{1/2}} =i^{1/2} \dfrac{\Theta_{3} \left(-\dfrac{1}{1+z}\right)}{(z+1)^{1/2}} \ ,
\end{array}\hspace{-0.2cm}
\end{align}

\vspace{0.25cm}
\noindent and

\begin{align}\label{adt30}
&\hspace{-0.45cm} \begin{array}{l}
{\rm{(a)}}\	\ \ \Theta_{2} (z)^{2}=  -\dfrac{\Theta_{4} \left(- \dfrac{1}{z-1}\right)^{2}}{z-1}\ \  = \dfrac{\Theta_{4} \left(-\dfrac{1}{z+1}\right)^{2}}{z+1} \ ,
\\[0,4cm]
{\rm{(b)}} \  \ 	  \Theta_{3} (z)^{2}=  \, i \, \dfrac{\Theta_{2} \left(- \dfrac{1}{z-1}\right)^{2}}{z-1} = \, i \,  \dfrac{\Theta_{2} \left(-\dfrac{1}{z+1}\right)^{2}}{z+1} \ ,
\\[0,4cm]
{\rm{(c)}} \  \  \Theta_{4} (z)^{2}=  \, i  \, \dfrac{\Theta_{3} \left(- \dfrac{1}{z-1}\right)^{2}}{z-1} = \, i \,  \dfrac{\Theta_{3} \left(-\dfrac{1}{z+1}\right)^{2}}{z+1} \ , \quad z \in \Bb{H}\,.
\end{array}\hspace{-0.2cm}
\end{align}

\vspace{0.5cm}
\noindent Therefore for arbitrary  $\sigma \in \{1, -1\}$ we have

\begin{align}\nonumber
&
{\rm{(a)}} \  \ 	  \Theta_{3} (z)^{4}=  \, - \, \dfrac{\Theta_{2} \left(- \dfrac{1}{z-\sigma}\right)^{4}}{(z-\sigma)^{2}} \ ,
  &  &
{\rm{(b)}} \  \  \Theta_{4} (z)^{4}=  \, -  \, \dfrac{\Theta_{3} \left(- \dfrac{1}{z-\sigma}\right)^{4}}{(z-\sigma)^{2}} \ ,
\\[0,4cm]  &
{\rm{(c)}}\	\ \ \Theta_{2} (z)^{4}=  \dfrac{\Theta_{4} \left(- \dfrac{1}{z-\sigma}\right)^{4}}{(z-\sigma)^{2}}\, ,  &    &     \quad z \in \Bb{H}\,.
\label{adt30a}\end{align}

}}\end{clash}
\end{subequations}

\subsection[\hspace{-0,25cm}. \hspace{0,075cm}Notes on Section~\ref{wid}]{\hspace{-0,11cm}{\bf{.}} Notes on Section~\ref{wid}}

\begin{subequations}

\vspace{-0.15cm}
\begin{clash}{\rm{\hypertarget{d17}{}\label{case17}\hspace{0,00cm}{\hyperlink{bd17}{$\uparrow$}} \
 We prove \eqref{f1theor4} in more detail.

\vspace{0.25cm}   By \eqref{f24preresinvlamhyp}, \eqref{f23preresinvlamhyp} and \eqref{f18preresinvlamhyp}, the function $\het (z)$ does not vanish on  $\Bb{C}\setminus [1, +\infty)$ and therefore for $z \in \Lambda:= (0,1)\cup \left(\Bb{C}\setminus\Bb{R}\right)$
 we can introduce the function\vspace{-0,2cm}
\begin{align*}\tag{\ref{f1wid}}
    &   \Phi (z) := \dfrac{\Theta_{3}\left(\ie (z)\right)^{2}}{\het (z)}
    = \dfrac{\Theta_{3}\left(i \, \dfrac{\he(1-z)}{\he(z)}\right)^{2}}{\het (z)}
      \ .
\end{align*}

\vspace{0.25cm}
\noindent Since $\Theta_{3}\in  {\rm{Hol}}( \Bb{H})$, $\het\in {\rm{Hol}}( \Bb{C}\setminus  [1, +\infty)$  and $\ie \in {\rm{Hol}}( \Lambda)$ we get $\Phi \in  {\rm{Hol}}( \Lambda)$. The formulas \eqref{f9elemtheta}(b), \eqref{f10int} and  \eqref{f27preresinvlamhyp} for any $z \in \Lambda$ yield that %\vspace{-0,1cm}
\begin{align*}\tag{\ref{f2wid}}
    &  \Phi (z)=\dfrac{\Theta_{3}\left(\ie (z)\right)^{2}}{\het (z)} \,\stackrel{{\fo{\eqref{f9elemtheta}(b)}}}{\vphantom{A}=}\, \dfrac{i}{\ie (z)} \dfrac{\Theta_{3}\left(-1/\ie (z)\right)^{2}}{\het (z)}\\[0,3cm]    &
    \stackrel{{\fo{\eqref{f10int}}}}{\vphantom{A}=}\,
       \dfrac{i}{i \, \dfrac{\he(1-z)}{\he(z)}} \dfrac{\Theta_{3}\left(-1/\ie (z)\right)^{2}}{\het (z)}     \,\stackrel{{\fo{\eqref{f27preresinvlamhyp}}}}{\vphantom{A}=}\,    \dfrac{\Theta_{3}\left(\ie (1-z)\right)^{2}}{\het (1-z)} =  \Phi (1-z)  \ ,
\end{align*}

\noindent i.e.,
\begin{align}\label{ad27}
    &  \Phi (z)=  \Phi (1-z)   \ , \quad  z \in \Lambda:= (0,1)\cup \left(\Bb{C}\setminus\Bb{R}\right)\,.
\end{align}

\noindent Since $\het\in {\rm{Hol}}\big( \Bb{C}\setminus  [1, +\infty)\big)$ we have
\begin{align}\label{ad27a}
    &  \het (-x +  i  0 ) =  \het (-x -  i  0 )   \ , \quad  x > 0 \,,
\end{align}

\noindent and by using   \eqref{adt27}(c) and  \eqref{f28preresinvlamhyp}, for arbitrary  $x > 0$ we deduce that
\begin{align*}
    &  \Phi (-x\! + \! i  0 )\! = \! \dfrac{\Theta_{3}\left(\ie (-x\! + \! i  0)\right)^{2}}{\het (-x\! + \! i  0)}
   \,\stackrel{{\fo{\eqref{f28preresinvlamhyp}, \eqref{ad27a}}}}{\vphantom{A}=}\,
     \dfrac{\Theta_{3}\left(2\! +\! \ie (-x\! -\!  i  0 )\right)^{2}}{\het (-x)} \\[0,3cm]    &
     \,\stackrel{{\fo{\eqref{adt27}(d)}}}{\vphantom{A}=}\,\dfrac{\Theta_{3}\left( \ie (-x\! - \! i  0 )\right)^{2}}{\het (-x)} \!  =\!
      \Phi (-x\! - \!  i  0 ) \ ,
\tag{\ref{f3wid}}\end{align*}

\noindent i.e.,
\begin{align}\label{ad28}
    &  \Phi (-x  +   i  0 )  =   \Phi (-x  -    i  0 )  \ , \quad  x > 0 \ .
\end{align}

 For arbitrary  $x > 0$ it follows from
\begin{align*}
    &  \tag{\ref{f29preresinvlamhyp}}
1- 2 \ie (  1+ x +  i  0 ) =
\dfrac{\ie (  1+ x +  i  0 )}{\ie (1 + x - i  0 )  } ,
\quad x>0.
\end{align*}

\noindent and
\begin{align*}\tag{\ref{f10int}}
       \ie (z) =   i \, \dfrac{\he(1-z)}{\he(z)},
\quad  z \in \Lambda:=  (0,1)\cup \left(\Bb{C}\setminus\Bb{R}\right),
\end{align*}

\noindent that
\begin{align*}
    &  1- 2 \ie (  1+ x +  i  0 ) = \dfrac{ i \, \dfrac{\he(-x - i 0)}{\he( 1+ x +  i  0)}}{
 i \, \dfrac{\he(-x + i 0)}{\he( 1+ x -  i  0)}     }=
  \dfrac{  \dfrac{\he(-x )}{\he( 1+ x +  i  0)}}{
  \dfrac{\he(-x )}{\he( 1+ x -  i  0)}     }=\dfrac{\he( 1+ x -  i  0)}{\he( 1+ x +  i  0)} \ ,
\end{align*}

\noindent i.e.,
\begin{align}\label{ad29}
    &  1- 2 \ie (  1+ x +  i  0 ) = \dfrac{\he( 1+ x -  i  0)}{\he( 1+ x +  i  0)} \ , \quad  x > 0 \ .
\end{align}

\noindent Besides that, \eqref{adt28}(e) implies
\begin{align}\label{ad30}
    &  \Theta_{3}\left(\dfrac{z}{1-2z}\right)^{2}\!\! =\!  (1\!-\!2z) \Theta_{3} (z)^{2} \ , \quad  z\in \Bb{H}\,,
\end{align}

\noindent and therefore, by \eqref{f29preresinvlamhyp}, we obtain
\begin{align*}
    &   \Phi (1\! +\! x\! -\!  i  0) = \dfrac{\Theta_{3}\big(\ie (1\! +\! x\! -\!  i  0)\big)^{2}}{\het (1\! +\! x\! -\!  i  0)}
     \,\stackrel{{\fo{\eqref{f29preresinvlamhyp}}}}{\vphantom{A}=}\,
     \dfrac{\Theta_{3}\left(\dfrac{\ie (  1+ x +  i  0 )}{ 1- 2 \ie (  1+ x +  i  0 ) }\right)^{2}}{\het (1\! +\! x\! -\!  i  0)}
     \\[0,3cm]    &
 \stackrel{{\fo{\eqref{ad30}}}}{\vphantom{A}=}\,
 \big(1\!-\!2 \ie (1\! +\! x\! +\!  i  0 )\big) \dfrac{\Theta_{3}\big(\ie (1\! +\! x\! +\!  i  0)\big)^{2}}{\het (1\! +\! x\! -\!  i  0)}
  \\[0,3cm]    &
 \stackrel{{\fo{\eqref{ad29}}}}{\vphantom{A}=}\,
  \dfrac{\he( 1+ x -  i  0)}{\he( 1+ x +  i  0)}
 \dfrac{\Theta_{3}\big(\ie (1\! +\! x\! +\!  i  0)\big)^{2}}{\het (1\! +\! x\! -\!  i  0)}=
 \dfrac{\Theta_{3}\big(\ie (1\! +\! x\! +\!  i  0)\big)^{2}}{\het (1\! +\! x\! +\!  i  0)}
 \\[0,4cm]    &    = \Phi (1\! +\! x\! +\!  i  0)  \ , \quad  x > 0  \ , \
\end{align*}

\noindent i.e.,
\begin{align}\label{ad31}
    &  \Phi (1\! +\! x\! -\!  i  0) =  \Phi (1\! +\! x\! +\!  i  0)  \ , \quad  x > 0  \ .
\end{align}

\noindent Applying the Morera theorem to the function $\Phi \in  {\rm{Hol}}( (0,1)\cup \left(\Bb{C}\setminus\Bb{R}\right))$ with the properties \eqref{ad31} and \eqref{ad28} (see \cite[p. 96]{lav}) we obtain that
\begin{align}\label{ad32}
    &  \Phi \in  {\rm{Hol}}( \Bb{C}\setminus\{0, 1\}) \,.
\end{align}

\vspace{0.25cm} Let $\Lambda\ni z\to 0$. Then, by \eqref{f29ainvlamhyp},
\begin{align*}
    & \ie (z)   =   \dfrac{i }{\pi} \ln\dfrac{1}{|z|} - \ie(0; z) \ , \
\end{align*}

\noindent where in view of \eqref{f1th2}(b) there exists a finite positive number $\varepsilon_{\!{\tn{\triangle}}}$ such that
 $\big|\ie(0; z)\big|  \leq 2$  for all   $|z| \leq \varepsilon_{\!{\tn{\triangle}}}$. Therefore
\eqref{f8elemtheta} and \eqref{f1elemtheta} yield that
\begin{align*}
    &  \Theta_{3}\left(\ie (z)\right)^{2}= \theta_{3} \left(e^{{\fo{ i\pi \ie (z)  }}}\right)^{\!2} =
     \theta_{3}\left(e^{{\fo{  - \ln\dfrac{1}{|z|} \!- \!  i\pi \ie(0; z)     }}}\right)^{\!2} \\    &
   \left(  1  \!+ \! 2\sum\limits_{n\geq 1} e^{{\fo{  - n^{2}\ln\dfrac{1}{|z|} \!- \!  i\pi n^{2}\ie(0; z)     }}} \right)^{2}=
   \left(  1  \!+ \! 2\sum\limits_{n\geq 1}|z|^{n^{2}} e^{{\fo{ - \!  i\pi n^{2}\ie(0; z)     }}} \right)^{2}
    \to 1\,,
\end{align*}

\noindent because
\begin{align*}
    &  \left|\,|z|^{n^{2}} e^{{\fo{ - \!  i\pi n^{2}\ie(0; z)     }}}\,\right| \leq\left( e^{2 \pi}|z|\right)^{n^{2}}
     \ , \quad  n \geq 1  \ , \  |z| \leq \varepsilon_{\!{\tn{\triangle}}} \,  .
\end{align*}

\noindent Since in view of \eqref{f7preresinvlamhyp}(b) we have  $\het (z) \to 1$, it follows that there exists the limit
\begin{align*}
    & \lim\limits_{{\fo{\Lambda\ni z\to 0}}} \Phi (z) =
    \lim\limits_{{\fo{\Lambda\ni z\to 0}}}  \dfrac{\Theta_{3}\left(\ie (z)\right)^{2}}{\het (z)} = 1 \,.
\end{align*}

\noindent By using \eqref{ad27}, we obtain the existence of two limits
\begin{align} \tag{\ref{f5wid}}
    &  \lim\limits_{{\fo{\Lambda\ni z\to 0}}} \,  \Phi (z) =  \lim\limits_{{\fo{\Lambda\ni z\to 1}}}  \Phi (z) = 1 \,.
\end{align}

\noindent Applying to $\Phi \in  {\rm{Hol}}( \Bb{C}\setminus\{0, 1\})$ (see \eqref{ad32})  the Riemann theorem about removable singularity (see \cite[p. 103]{con}) under the properties \eqref{f5wid} we obtain that $\Phi$ is entire function satisfying
\begin{align}\label{ad32a}
   \Phi (0) =  \Phi (1) = 1 \,.
\end{align}

\vspace{0.25cm} Now let $\Lambda\ni z\! \to\! \infty$ approaching from  one of the half-planes $\sigma\! :=\! {\rm{sign}} (\im z)\! \in \!\{1, -1\}$. To study  the behaviour of
\begin{align*}
     \Phi (z) = \dfrac{\Theta_{3}\big(\ie (z)\big)^{2}}{\het (z)}
\end{align*}

\noindent  we apply \eqref{adt30}(b)  to get
\begin{align}\nonumber
     & \Theta_{3} \big(\ie (z)\big)^{2}=   \dfrac{i}{\ie (z)-\sigma}\, \Theta_{2}
     \left(- \dfrac{1}{\ie (z)-\sigma}\right)^{2}  \\[0,4cm]  &
     \,\stackrel{{\fo{\eqref{f8elemtheta}}}}{\vphantom{A}=}\,
     \dfrac{4 i e^{{\fo{-\dfrac{i \pi}{2}  \dfrac{1}{\ie (z)-\sigma} }}}  }{\ie (z)-\sigma}
    \, \theta_{2}\left( e^{{\fo{-  \dfrac{i \pi}{\ie (z)-\sigma} }}} \right)^{2} \ , \
\label{ad33}\end{align}

\noindent where  by virtue of \eqref{f29ainvlamhyp},
\begin{align*}
    &- \dfrac{1}{  \ie (z) -\sigma }   =  \dfrac{ i }{ \pi} \ln  |z| - \ie(\infty; z) \ , \
\end{align*}

\noindent and in view of \eqref{f1th2}(b) there exists a finite positive number $\varepsilon_{\!{\tn{\triangle}}}$ such that
 $\big|\ie(\infty; z)\big|  \leq 2$  for all   $|z| \geq 1/\varepsilon_{\!{\tn{\triangle}}}$, so that
 \begin{align}\label{ad33a}
    &  \ie(\infty; z) = O (1) \ , \quad  \Lambda\ni z\! \to\! \infty\,.
 \end{align}

 \noindent
  This means that
\begin{align*}
      e^{{\fo{-\dfrac{i \pi}{2}  \dfrac{1}{\ie (z)-\sigma} }}} &=
     e^{{\fo{ -\dfrac{1}{2}\ln  |z| -  \dfrac{i \pi}{2}\ie(\infty; z) }}}
     = \dfrac{ e^{{\fo{  -  \dfrac{i \pi}{2}\ie(\infty; z) }}}}{\sqrt{|z|}}
       \ , \  \\
      e^{{\fo{-  \dfrac{i \pi}{\ie (z)-\sigma} }}}& =
     e^{{\fo{ -\ln  |z| -  i \pi\ie(\infty; z) }}}=
          \dfrac{  e^{{\fo{  -  i \pi\ie(\infty; z) }}}}{|z| }
          \ , \  \\
     \dfrac{4 i}{  \ie (z) -\sigma }   & =  \dfrac{ 4 }{ \pi} \ln  |z| + 4 i \ie(\infty; z) \ , \
\end{align*}

\noindent and after substituting these expressions in \eqref{ad33} we obtain
\begin{align}\label{ad34}
   & \hspace{-0,3cm}\Theta_{3} \big(\ie (z)\big)^{2}= \dfrac{\dfrac{ 4 }{ \pi} \ln  |z| + 4 i \ie(\infty; z)}{\sqrt{|z|}}
   \,e^{{\fo{  -  \dfrac{i \pi}{2}\ie(\infty; z) }}}
    \theta_{2}\left(\dfrac{  \,e^{{\fo{  -  i \pi\ie(\infty; z) }}}}{|z| }\right)^{\!2}\! \!\!. \hspace{-0,25cm}
    \end{align}

\noindent At the same time  \eqref{f7preresinvlamhyp}(c) yields
  \begin{align*}
        \he (z)& =\dfrac{\Log 16 (1-z)}{\pi \sqrt{1-z}} + O \left(\dfrac{\ln |z|}{|z|^{3/2}}\right)\\[0,3cm]  & =
    \dfrac{\ln |1-z| + i \Arg (1-z) + \ln 16}{\pi \sqrt{1-z}}  + O \left(\dfrac{\ln |z|}{|z|^{3/2}}\right)
    \\  & =
  \dfrac{\ln |z| }{\pi \sqrt{1-z}} + O \left(\dfrac{1}{|z|^{1/2}}\right)\\[0,3cm]  & =
   \dfrac{\ln |z| }{\pi } \,e^{{\fo{- \dfrac{1}{2} \ln |1-z| -  \dfrac{i}{2}\Arg (1-z)     }}} + O \left(\dfrac{1}{|z|^{1/2}}\right)\\[0,3cm]  & =
    \dfrac{\ln |z| }{\pi } \,e^{{\fo{- \dfrac{1}{2} \ln |z| -  \dfrac{i}{2}\Arg (1-z)     }}} + O \left(\dfrac{1}{|z|^{1/2}}\right)\\[0,3cm]  & =
    \dfrac{\ln |z| }{\pi\sqrt{|z|} } \,e^{{\fo{-  \dfrac{i}{2}\Arg (1-z)     }}} + O \left(\dfrac{1}{|z|^{1/2}}\right)\,.
\end{align*}

\noindent Dividing \eqref{ad34} by this expression and taking account of \eqref{ad33a} and\\ $\lim_{u\to 0} \theta_{2} (u) = \lim_{u\to 0} (1 \!+\! \sum\nolimits_{n\geq 1} u^{n^2 +n}) =1$, in view of \eqref{f1elemtheta}, we obtain
\begin{align*}
     \Phi (z)  & = \dfrac{\Theta_{3}\big(\ie (z)\big)^{2}}{\het (z)} \\[0,3cm]  & =
   \dfrac{\dfrac{\dfrac{ 4 }{ \pi} \ln  |z| + 4 i \ie(\infty; z)}{\sqrt{|z|}}
   \,e^{{\fo{  -  \dfrac{i \pi}{2}\ie(\infty; z) }}}
    \theta_{2}\left(\dfrac{  \,e^{{\fo{  -  i \pi\ie(\infty; z) }}}}{|z| }\right)^{\!2}}{  \dfrac{\ln |z| }{\pi\sqrt{ |z|} } \,e^{ {\fo{ - \dfrac{i}{2} \Arg (1-z)  }}} + O \left(\dfrac{1}{|z|^{1/2}}\right)}\\[0,3cm]  & =
    \dfrac{\dfrac{\dfrac{ 4 }{ \pi} \ln  |z| + 4 i \ie(\infty; z)}{\sqrt{|z|}}
      \left(1 + O \left(\dfrac{1}{|z|}\right)\right)}{  \dfrac{\ln |z| }{\pi\sqrt{ |z|} } + O \left(\dfrac{1}{|z|^{1/2}}\right)}     \,e^{{\fo{ \dfrac{i}{2} \Arg (1-z) -  \dfrac{i \pi}{2}\ie(\infty; z) }}}\\[0,3cm]  & =
      4 \dfrac{\ln  |z| + i \pi \ie(\infty; z)}{\ln |z| + O (1)}    \,e^{{\fo{ \dfrac{i}{2} \Arg (1-z) -  \dfrac{i \pi}{2}\ie(\infty; z) }}}\left(1 + O \left(\dfrac{1}{|z|}\right)\right)\\[0,3cm]  & =
      4 \,e^{{\fo{ \dfrac{i}{2} \Arg (1-z) -  \dfrac{i \pi}{2}\ie(\infty; z) }}}
      \dfrac{1 + \dfrac{i \pi \ie(\infty; z)}{\ln |z|}}{ 1 + O \left(\dfrac{1}{\ln |z|}\right)}\left(1 + O \left(\dfrac{1}{|z|}\right)\right)
     \\[0,3cm]  & \,\stackrel{{\fo{\eqref{f30preresinvlamhyp}}}}{\vphantom{A}=}\,
         4 \,e^{{\fo{ \dfrac{i}{2} \Arg (1-z) -  \dfrac{i \pi}{2}\left( \dfrac{\ln 16 }{ i  \pi} +  \dfrac{\Arg (1-z)}{\pi}  + O \left(\dfrac{\ln^{2} |z|}{|z|}\right) \right) }}}
      \!\!\left(1\! +\! O \!\left(\!\dfrac{1}{\ln |z|}\!\right)\!\right)\\[0,3cm]  & =
      4 \,e^{{\fo{ -  \dfrac{1}{2}\ln 16   + O \left(\dfrac{\ln^{2} |z|}{|z|}\right)  }}}
      \!\!\left(1\! +\! O \!\left(\!\dfrac{1}{\ln |z|}\!\right)\!\right) = 1\! +\! O \!\left(\!\dfrac{1}{\ln |z|}\!\right) \,.
\end{align*}

\noindent In particular, this implies that  the entire function $\Phi$ is bounded and  by the Liouville theorem \cite[p. 77]{con} it is a constant, which equals to $1$, by \eqref{ad32a}. The Wirtinger identity \eqref{f1theor4} follows.
}}\end{clash}
\end{subequations}

\begin{subequations}

\begin{clash}{\rm{\hypertarget{d18}{}\label{case18}\hspace{0,00cm}{\hyperlink{bd18}{$\uparrow$}} \ We explain the proof in greater detail.
According to \eqref{f9elemtheta}(b) and \eqref{adt27}(d),
\begin{align}\label{ad35}
    & {\rm{(a)}} \ \  \Theta_{3} (z+2) = \Theta_{3} (z)	\, ,
    \ \   {\rm{(b)}} \ \ \Theta_3(-1/z)= (z/i)^{1/2}\Theta_3(z)  \ , \  z \in \Bb{H} \ ,
\end{align}

\noindent and it has already been proved that
\begin{align}\label{ad36}
    &  \Theta_{3} (z) \neq 0  \ , \quad  z\in \mathcal{F}^{\,{\tn{||}}}_{{\tn{\square}}} \ , \    \\[0,25cm]    &    \nonumber \mathcal{F}^{\,{\tn{||}}}_{{\tn{\square}}}:= \big\{z\in\Bb{C}:\,\,{\im}\, z > 0, \,\,-1\leq{\re}\, z \leq 1,
\,\,|2 z\! -\!  1| > 1,\,\, | 2 z \!+\!  1|\! > 1\big\}\,.
\end{align}

\vspace{0.25cm}
\noindent
In the notations
\begin{align*}
    &  {\rm{SL}}_{2}(\Bb{Z}) := \left.\left\{\ \left(\begin{matrix} a & b \\ c & d \end{matrix}\right) \ \right| \
\quad  a, b, c, d \in\Bb{Z}\,,\,\,\,ad - bc=1\,. \ \right\}  \ , \   \\    &
\phi_{
{\nor{(\begin{smallmatrix} a & b \\ c & d \end{smallmatrix})}}
}(z) :=
\dfrac{a  z + b}{c  z + d}\ ,  \quad z\in \Bb{H}\,, \ \  \left(\begin{matrix} a & b \\ c & d \end{matrix}\right) \in {\rm{SL}}_{2}(\Bb{Z}) \ , \\    &   I := \begin{pmatrix} 1 & 0 \\ 0 & 1 \\ \end{pmatrix}  \, , \
 T := \begin{pmatrix} 1 & 1 \\ 0 & 1 \\ \end{pmatrix}  \, , \  S := \begin{pmatrix} 0 & 1 \\ -1 & 0 \\ \end{pmatrix} \ , \  S^{2} =  - I \, , \  T^{n} = \begin{pmatrix} 1 & n \\ 0 & 1 \\ \end{pmatrix} \, , \  n \in \Bb{Z} \, , \\    &
 \phi_{{\fo{I}}} (z) = z  \, , \ \phi_{{\fo{S}}} (z) =  - \dfrac{1}{z}  \ ,  \  \phi_{{\fo{\,T^{{\tn{n}}}}}} (z) =  z+n  \, , \  n \in \Bb{Z} \, ,
\end{align*}

\noindent we deduce from \eqref{ad35}, that
\begin{align*}
    &  \Theta_{3} \left( \phi_{{\fo{\,T^{{\tn{2 n}}}}}} (z)\right) = \Theta_{3} (z)	\, , \ n  \in \Bb{Z} \, ,
    \ \    \Theta_3\left(\phi_{{\fo{S}}} (z)\right) = (z/i)^{1/2}\Theta_3 (z)  \ , \  z \in \Bb{H} \ .
\end{align*}

\noindent This means that
\begin{align*}
    & z \in \Bb{H} \ , \   \Theta_{3} (z) \neq 0   \ \Rightarrow \  \Theta_3\left(\phi_{{\fo{S}}} (z)\right)\neq 0  \
    \mbox{and    } \   \Theta_{3} \left( \phi_{{\fo{\,T^{{\tn{2 n}}}}}} (z)\right) \neq 0  \ \mbox{for all} \ n  \in \Bb{Z} \,.
\end{align*}

\noindent Taking the repeated applications of $\phi_{{\fo{S}}}$ or $ \phi_{{\fo{\,T^{{\tn{2 n}}}}}}$ with arbitrary $n  \in \Bb{Z}\setminus\{0\} $
to a given $z\in \mathcal{F}^{\,{\tn{||}}}_{{\tn{\square}}}$ and taking account  the basic property of the M{\"o}bius transformations
\begin{multline*}
    \phi_{A_{1}A_{2}\cdot\ldots\cdot A_{N}} (z) = \phi_{A_{1}} \Big(\phi_{A_{2}} \big(\ldots\big( \phi_{A_{N}} (z) \big)\ldots\big)\Big) \ ,  \\       A_{1}, A_{2}, \ldots,  A_{N} \in {\rm{SL}}_{2} (\Bb{Z})  \, , \  z \in  \Bb{H}  \, , \  N \geq 2 \, ,
\end{multline*}

\noindent we obtain that
\begin{align}\label{ad37}
    &  \Theta_3\left(\phi_{{\fo{M}}} (z)\right)\neq 0   \ , \quad   z\in \mathcal{F}^{\,{\tn{||}}}_{{\tn{\square}}} \ , \
\end{align}

\noindent for arbitrary matrix $M$ which can be represented as a product of finite number of matrices $S$, $T^{2}$ and $T^{-2}$.
Such set of matrices $M$ forms the subgroup (see \cite[p. 112, Lem. 2]{cha})
\begin{align*}
   &  {\rm{SL}}_{2} (\vartheta, \Bb{Z})\! :=\!
 \left\{ \! {\fo{
\begin{pmatrix} a & b \\ c & d \\ \end{pmatrix}}}\! \in\!
{\rm{SL}}_{2}(\Bb{Z}) \, \left| \,  {\fo{\begin{pmatrix} a & b \\ c & d \\
\end{pmatrix}}}\! \equiv\! {\fo{\begin{pmatrix} 1 & 0 \\ 0 & 1 \\
\end{pmatrix}}} ({\rm{mod}}\, 2)  \ \mbox{or} \  {\fo{
\begin{pmatrix} a & b \\ c & d\\ \end{pmatrix}}}\! \equiv\! {\fo{
\begin{pmatrix} 0 & 1 \\ 1 & 0 \\ \end{pmatrix}}}  ({\rm{mod}}\, 2)
\! \right\}\right.
\end{align*}

\noindent of the group ${\rm{SL}}_{2}(\Bb{Z})$ and it has been proved in \cite[p. 16, Thm. 3.1]{bkn} (where it is necessary to set $\lambda = 2$
 and apply  \cite[p. 15, Def. 3.3]{bkn}) that
\begin{align*}
\tag{\ref{f6beirf}}
\Bb{H} =
\bigcup\nolimits_{  {\fo{M=(
\begin{smallmatrix} a & b \\ c & d \end{smallmatrix})
\in  {\rm{SL}}_{2} (\vartheta, \Bb{Z})  }} } \ \
   \phi_{M}\left({\rm{clos}}_{{\tn{\Bb{H}}}}\left( \mathcal{F}_{{\fo{\Gamma_{\!\vartheta}}}}\right)\right) , \
  \end{align*}

\noindent where
\begin{align*}
    &  {\rm{clos}}_{{\tn{\Bb{H}}}}\left(\mathcal{F}_{{{\Gamma_{\!\vartheta}}}}\right) :=
\{\ z \in \Bb{H} \ | \  -1 \!\leq  \re\,  z \leq  1 ,   \  | z | \geq 1
  \ \} \! \subset
 \mathcal{F}^{\,{\tn{||}}}_{{\tn{\square}}}\,.
\end{align*}

\noindent Together with property \eqref{ad37} this yields that $\Theta_{3} (z) \neq 0$ for all $z \in \Bb{H}$.

\vspace{0.25cm}  Then
 the following consequences of \eqref{f9elemtheta}:
\begin{align*}
    &  \Theta_{4}(z) \,\stackrel{{\fo{\eqref{f9elemtheta}(e)}}}{\vphantom{A}=}\,  \Theta_{3}(z\!+\!1) \ , \\[0,3cm]    &
    \Theta_{2}(z)\,\stackrel{{\fo{\eqref{f9elemtheta}(c)}}}{\vphantom{A}=}\, (i/z)^{1/2}\Theta_{3}(-1/z)\,\stackrel{{\fo{\eqref{f9elemtheta}(e)}}}{\vphantom{A}=}\, (i/z)^{1/2}  \Theta_{3}(1\!-\!1/z) \, ,
    \quad z \! \in \! \Bb{H}\,,
\end{align*}

\vspace{0.25cm}
\noindent prove that $\Theta_{2} (z) \neq 0$,  $\Theta_{4} (z) \neq 0$ for all $z \in \Bb{H}$, which completes the proof of \eqref{f0theor3}.
}}\end{clash}
\end{subequations}

\subsection[\hspace{-0,25cm}. \hspace{0,075cm}Notes on Section~\ref{emf}]{\hspace{-0,11cm}{\bf{.}} Notes on Section~\ref{emf}}

\begin{subequations}

\vspace{-0.15cm}
\begin{clash}{\rm{\hypertarget{d19}{}\label{case19}\hspace{0,00cm}{\hyperlink{bd19}{$\uparrow$}} \ We proceed with the proof of \eqref{f8theor3} and \eqref{f1theor3} by supplying more details regarding the functions
\begin{align}\label{ad38}
    & f_{1}(z) := \dfrac{\Theta_{2}(z)^{4} + \Theta_{4}(z)^{4}}{\Theta_{3}(z)^{4}}  \ , \  f_{2} (z) :=   \dfrac{\Theta_{2}\left(z\right)^{4}}{\Theta_{3}\left(z\right)^{4}}  \ , \quad    z\in \Bb{H} \ .
\end{align}

We state that for  any $z\in \Bb{H}$ the following relations hold
\begin{align*}
    &
     {\rm{(a)}} \ \  f_{k} (z+2) = f_{k} (z)  \, ,     &    &        {\rm{(b)}} \ \  f_{k} \big({z}/{(1-2z)}\big)= f_{k} (z)\ ,
     \    \  k \in \{1, 2\} \ , \\[0,25cm]   &
     {\rm{(c)}} \ \   f_{1}(z) = f_{1}(-1/z)  \ , \    &    &
     {\rm{(d)}} \ \    f_{2} (z) =  {\Theta_4 (-1/z)^{4}}\big/{\Theta_3 (-1/z)^{4}} \ , \   \\[0,25cm]   &
     {\rm{(e)}} \ \   f_{1}(z) = \dfrac{\Theta_{3}(y)^{4} - \Theta_{4}(y)^{4}}{ \Theta_{2}(y)^{4}}   \ , \
      &    &  {\rm{(f)}} \ \ f_{2} (z) = -\dfrac{\Theta_{4}(y)^{4}}{ \Theta_{2}(y)^{4}} \, , \  y :=\dfrac{1}{\sigma\!-\!z} \, .
\tag{\ref{f1emf}}\hypertarget{syst}{}\end{align*}

\vspace{0.25cm}
\noindent
{\emph{Proof of \eqref{f1emf}{\rm{(a)}}.}} According to \eqref{adt27}(c),(d),(e) each function $\Theta_{2}^{4}$,  $\Theta_{3}^{4}$ and
$\Theta_{4}^{4}$ is periodic in $\Bb{H}$ with period $2$.

\vspace{0,1cm}
\noindent {\emph{Proof of \eqref{f1emf}{\rm{(b)}}.}} Substituting \eqref{adt28a} in \eqref{ad38} we obtain the equalities \eqref{f1emf}{\rm{(b)}}.

\vspace{0,1cm}
\noindent {\emph{Proof of  \eqref{f1emf}{\rm{(c), (d)}}.}} Substituting  \eqref{f9elemtheta}(a),(b),(c) in \eqref{ad38} we obtain the equalities \eqref{f1emf}{\rm{(c)}}, {\rm{(d)}}.

\vspace{0,1cm}
\noindent {\emph{Proof of  \eqref{f1emf}{\rm{(e), (f)}}.}} Substituting  \eqref{adt30a}(a),(b),(c) in \eqref{ad38} we obtain
\begin{align*}
    &  f_{1}(z) = \dfrac{\Theta_{2}^{4}(z) + \Theta_{4}^{4}(z)}{\Theta_{3}^{4}(z)} =
    \dfrac{\dfrac{\Theta_{4} \left(- \dfrac{1}{z-\sigma}\right)^{4}}{(z-\sigma)^{2}}  \, -  \, \dfrac{\Theta_{3} \left(- \dfrac{1}{z-\sigma}\right)^{4}}{(z-\sigma)^{2}}  }{\, - \, \dfrac{\Theta_{2} \left(- \dfrac{1}{z-\sigma}\right)^{4}}{(z-\sigma)^{2}}}
     \\[0,2cm]    &   = \left|\,y := - \dfrac{1}{z-\sigma}\,\right| = \dfrac{\Theta_{4}^{4}(y) - \Theta_{3}^{4}(y)}{-\Theta_{2}^{4}(y)}=
     \dfrac{\Theta_{3}(y)^{4} - \Theta_{4}(y)^{4}}{ \Theta_{2}(y)^{4}} \ , \\[0,2cm]    &
     f_{2} (z) =   \dfrac{\Theta_{2}\left(z\right)^{4}}{\Theta_{3}\left(z\right)^{4}} =
     \dfrac{ \dfrac{\Theta_{4} \left(- \dfrac{1}{z-\sigma}\right)^{4}}{(z-\sigma)^{2}}}{\, - \, \dfrac{\Theta_{2} \left(- \dfrac{1}{z-\sigma}\right)^{4}}{(z-\sigma)^{2}} }= \left|\,y := - \dfrac{1}{z-\sigma}\,\right| =-\dfrac{\Theta_{4}(y)^{4}}{ \Theta_{2}(y)^{4}} \, ,
\end{align*}

\vspace{0.2cm}
\noindent which prove the equalities \eqref{f1emf}{\rm{(e)}}, {\rm{(f)}}.

\vspace{0.25cm} The equalities \eqref{f1emf}(a) mean that the condition Lemma~{\hyperlink{hth1}{\ref*{th1}}}(a) holds for $f_{1}$ and $f_{2}$, while the equality \eqref{f1emf}(b) for every $k \in \{1, 2\}$ and $z\in \Bb{H}$  yields that
\begin{align*}
    &  f_{k} \left(-\dfrac{1}{z}\right)= f_{k} \big((-1/z)\big) \stackrel{{\fo{({\hyperlink{syst}{\ref*{f1emf}}})(b)}}}{\vphantom{A}=} f_{k} \left(\dfrac{(-1/z)}{1-2(-1/z)}\right) =  f_{k} \left(- \dfrac{1}{2+ z}\right)  \, ,
\end{align*}

\noindent i.e., the condition Lemma~{\hyperlink{hth1}{\ref*{th1}}}(b) also holds for $f_{1}$ and $f_{2}$.

\vspace{0.25cm}  Besides that, it follows from \eqref{f1elemtheta} and \eqref{f8elemtheta} that as $\fet \ni z\to \infty$ we have
\begin{align*}
    \Theta_{3}  (z)^{4} & \!=\! \left(1  \!+ \! 2\sum\limits_{n\geq 1} e^{{\fo{i \pi n^2 z}}}  \right)^{4}\!\!=\!
    \left(1  \!+ \!2  e^{{\fo{i \pi  z}}}  \!+ \! 2 e^{{\fo{4 i \pi  z}}}\sum\limits_{n\geq 2} e^{{\fo{i \pi (n^2 -4) z}}}  \right)^{\!4}
     \\[0,25cm]    &   =  1\!+\! 8 e^{i \pi  z}\! + \!  O \left(e^{2 i \pi  z}\right)
     \ , \\[0,35cm]
      \Theta_{4}  (z)^{4} & \!=\! \left(1  \!+ \! 2\sum\limits_{n\geq 1} (-1)^{n} e^{{\fo{i \pi n^2 z}}} \right)^{\!\!4}\!\!=\!
      \left(\!1 \! -\! 2  e^{{\fo{i \pi  z}}} \! +\! 2 e^{{\fo{4 i \pi  z}}}\sum\limits_{n\geq 2} (-1)^{n} e^{{\fo{i \pi (n^2 -4) z}}}\! \right)^{\!\!4}
        \\[0,25cm]     &   =1\!-\! 8 e^{i \pi  z}\! +\!   O \left(e^{2 i \pi  z}\right)   \ , \\[0,35cm]
  \Theta_{2}  (z)^{4} & \!=\! 16 e^{{\fo{i \pi z }}}\!\!   \left(1 \!+\! \sum\limits_{n\geq 1} e^{{\fo{i \pi (n^2 +n) z}}} \right)^{\!\!4}\! \!=\!
      16 e^{{\fo{i \pi z }}}\!\!   \left(\!1 \!+\! e^{{\fo{2 i \pi z}}} \sum\limits_{n\geq 1} e^{{\fo{i \pi (n^2 +n -2) z}}}\! \right)^{\!\!4}
       \\[0,25cm]     &   = 16 e^{i \pi  z} \!  + \!  O \left(e^{3 i \pi  z}\right) \, ,
\end{align*}

\noindent and therefore
\begin{align*}\tag{\ref{f2emf}}
    &
    \begin{array}{ll}
 {\rm{(a)}} \ \    \Theta_{3}  (z)^{4}\!\!=\! 1\!+\! 8 e^{i \pi  z}\! + \!  O \left(e^{2 i \pi  z}\right),   &
 \  {\rm{(b)}} \ \   \Theta_{4}  (z)^{4}\!\!= \!1\!-\! 8 e^{i \pi  z}\! +\!   O \left(e^{2 i \pi  z}\right) , \\[0,25cm]
  {\rm{(c)}} \ \     \Theta_{2}  (z)^{4}\!\!=\!  16 e^{i \pi  z} \!  + \!  O \left(e^{3 i \pi  z}\right),   &
\hphantom{  {\rm{(c)}} \ \ \,}  \fet \ni z\to \infty \, .
    \end{array}\hypertarget{sist}{}
\end{align*}

\noindent Using these relations together with ({\hyperlink{syst}{\ref*{f1emf}}})(c), (d) we obtain

\vspace{-0,1cm}
\begin{align}
    &  \lim\limits_{{\fo{\fet \!\ni\! z\!\to\! \infty}}}\  f_{1} (z) =
  \lim\limits_{{\fo{\fet \!\ni\! z\!\to\! \infty}}}\  \dfrac{\Theta_{2}^{4}(z) + \Theta_{4}^{4}(z)}{\Theta_{3}^{4}(z)}
\nonumber   \\[0,25cm] &
 \,\stackrel{{\fo{({\hyperlink{sist}{\ref*{f2emf}}})}}}{\vphantom{A}=}\,
    \lim\limits_{{\fo{\fet \!\ni\! z\!\to\! \infty}}}\  \dfrac{ 16 e^{i \pi  z} \!  + \!  O \left(e^{3 i \pi  z}\right) + \!1\!-\! 8 e^{i \pi  z}\! +\!   O \left(e^{2 i \pi  z}\right)}{1\!+\! 8 e^{i \pi  z}\! + \!  O \left(e^{2 i \pi  z}\right)} = 1 \label{ad39}
    \ ,    \\[0,65cm] &  \nonumber
    \lim\limits_{{\fo{\fet \!\ni\! z\!\to\! \infty}}}\  f_{2} (z) = \! \!  \lim\limits_{{\fo{\fet \!\ni\! z\!\to\! \infty}}}\  \dfrac{\Theta_{2}\left(z\right)^{4}}{\Theta_{3}\left(z\right)^{4}}
 \,\stackrel{{\fo{({\hyperlink{sist}{\ref*{f2emf}}})}}}{\vphantom{A}=}\! \!
     \lim\limits_{{\fo{\fet \!\ni\! z\!\to\! \infty}}}\  \dfrac{  16 e^{i \pi  z} \!  + \!  O \left(e^{3 i \pi  z}\right)}{ 1\!+\! 8 e^{i \pi  z}\! + \!  O \left(e^{2 i \pi  z}\right)} = 0
       \, ,    \\[0,65cm] &  \nonumber
       \lim\limits_{{\fo{\fet \!\ni\! z\!\to\! 0}}}\  f_{1} (z)
       \,\stackrel{{\fo{({\hyperlink{syst}{\ref*{f1emf}}})(c)}}}{\vphantom{A}=}\, \lim\limits_{{\fo{\fet \!\ni\! z\!\to\! 0}}}\  f_{1} (-1/z)
     \,\stackrel{{\fo{ (\fet = -1/\fet)}}}{\vphantom{A}=}\,  \lim\limits_{{\fo{\fet \!\ni\! z\!\to\! \infty}}}\  f_{1} (z)
            \\[0,25cm] &  \nonumber   \,\stackrel{{\fo{\eqref{ad39}}}}{\vphantom{A}=}\, 1
         \ ,    \\[0,65cm] & \nonumber
 \lim\limits_{{\fo{\fet \!\ni\! z\!\to\! 0}}}\  f_{2} (z)   \,\stackrel{{\fo{({\hyperlink{syst}{\ref*{f1emf}}})(d)}}}{\vphantom{A}=}\,
 \lim\limits_{{\fo{\fet \!\ni\! z\!\to\! 0}}}\ \dfrac{\Theta_4 (-1/z)^{4}}{\Theta_3 (-1/z)^{4}}\ \stackrel{{\fo{ (\fet = -1/\fet)}}}{\vphantom{A}=} \!
 \! \! \!  \lim\limits_{{\fo{\fet \!\ni\! z\!\to\! \infty}}}\ \dfrac{\Theta_4 (z)^{4}}{\Theta_3 (z)^{4}}
   \\[0,25cm] & \nonumber
 \,\stackrel{{\fo{({\hyperlink{sist}{\ref*{f2emf}}})}}}{\vphantom{A}=}\,
  \lim\limits_{{\fo{\fet \!\ni\! z\!\to\! \infty}}}\ \dfrac{1\!-\! 8 e^{i \pi  z}\! +\!   O \left(e^{2 i \pi  z}\right)}{ 1\!+\! 8 e^{i \pi  z}\! + \!  O \left(e^{2 i \pi  z}\right)} = 1 \ ,
 \end{align}

\vspace{0.25cm}
\noindent so that,
 \begin{align*}\tag{\ref{f4emf}}
    & {\rm{(a)}}\  \lim\limits_{{\fo{\fet\! \ni \!z\!\to\! 0}}} \ f_{k} (z) = 1 , \ {\rm{(b)}}\
     \lim\limits_{{\fo{\fet \!\ni\! z\!\to\! \infty}}}\  f_{k} (z) = 2-k\ , \ \ \ k \in \{1,2\} \ .
\hypertarget{sst}{}\end{align*}

\vspace{0.25cm}
Moreover, for arbitrary $\sigma \in \{1, -1\}$ if $\fet \!\ni\! z\!\to\! \sigma $ then we can consider
${\rm{sign}} ({\re}\, z )= \sigma$  so that $z - \sigma \in \fet $ and consequently $y:= - 1/(z - \sigma ) \in \fet$. Hence,
\begin{align*}
    & \lim\limits_{{\fo{\fet \!\ni\! z\!\to\! \sigma}}}\  f_{1} (z) =
    \lim\limits_{\substack{{\fo{\fet \!\ni\! z\!\to\! \sigma}}\\ {\fo{{\rm{sign}} ({\re}\, z) = \sigma}}}}\  f_{1} (z) \stackrel{{\fo{({\hyperlink{syst}{\ref*{f1emf}}})(e)}}}{\vphantom{A}=}
    \lim\limits_{{\fo{\fet \!\ni\! y\!\to\! \infty}}}\  \dfrac{\Theta_{3}(y)^{4} - \Theta_{4}(y)^{4}}{ \Theta_{2}(y)^{4}} \\[0,35cm]     &
    \,\stackrel{{\fo{({\hyperlink{sist}{\ref*{f2emf}}})}}}{\vphantom{A}=}\, \lim\limits_{{\fo{\fet \!\ni\! y\!\to\! \infty}}}\
     \dfrac{1\!+\! 8 e^{i \pi  y}\! + \!  O \left(e^{2 i \pi  y}\right) - \big(1\!-\! 8 e^{i \pi  y}\! +\!   O \left(e^{2 i \pi  y}\right)\big)}{ 16 e^{i \pi  y} \!  + \!  O \left(e^{3 i \pi  y}\right)} = 1  \ , \\[0,45cm]      &
     \lim\limits_{{\fo{\fet \!\ni\! z\!\to\! \sigma}}}\ e^{{\fo{- \dfrac{2 \pi}{|z-\sigma|}}}} f_{2} (z) =
    \lim\limits_{\substack{{\fo{\fet \!\ni\! z\!\to\! \sigma}}\\ {\fo{{\rm{sign}} ({\re}\, z) = \sigma}}}}\ e^{{\fo{- \dfrac{2 \pi}{|z-\sigma|}}}} f_{2} (z)= \left|\,y := - \dfrac{1}{z-\sigma}\,\right| \\[0,45cm]     &
    \stackrel{{\fo{({\hyperlink{syst}{\ref*{f1emf}}})(f)}}}{\vphantom{A}=}\!\!\!
    \lim\limits_{{\fo{\fet \!\ni\! y\!\to\! \infty}}}\   -\dfrac{\Theta_{4}(y)^{4}}{ \Theta_{2}(y)^{4}}\, e^{{\fo{-2 \pi |y|}}}
    \,\stackrel{{\fo{({\hyperlink{sist}{\ref*{f2emf}}})}}}{\vphantom{A}=}\!\!\! \lim\limits_{{\fo{\fet \!\ni\! y\!\to\! \infty}}}
   \!\!\!  -\dfrac{1\!-\! 8 e^{i \pi  y}\! +\!   O \left(e^{2 i \pi  y}\right) }{  16 e^{i \pi  y} \!  + \!  O \left(e^{3 i \pi  y}\right)}\,e^{{\fo{-2 \pi |y|}}}\\[0,5cm]     &  =
    \lim\limits_{{\fo{\fet \!\ni\! y\!\to\! \infty}}}
   \!\!\!  \dfrac{-1\!+\! 8 e^{i \pi  y}\! +\!   O \left(e^{2 i \pi  y}\right) }{  16  \!  + \!  O \left(e^{2 i \pi  y}\right)}\,
   e^{{\fo{-2 \pi |y| - i \pi y}}} \\    &   =
    \lim\limits_{{\fo{\fet \!\ni\! y\!\to\! \infty}}}\left(- \dfrac{1}{16} +  O \left(e^{ i \pi  y}\right)\right)
      e^{{\fo{- \pi |y| -  \pi (|y| + i y ) }}}
    = 0 \,,
\end{align*}

\vspace{0.25cm}
\noindent i.e.,
\begin{align}\nonumber
    & {\rm{(a)}}\ \lim\limits_{{\fo{\fet \!\ni\! z\!\to\! \sigma}}}\  f_{1} (z) = 1 \ ,  \\    &   \ {\rm{(b)}}\   \lim\limits_{{\fo{\fet \!\ni\! z\!\to\! \sigma}}}\ e^{{\fo{- \dfrac{2 \pi}{|z-\sigma|}}}} f_{2} (z) =0  \ , \quad    \sigma \in \{1, -1\} \,.
\label{ad40}\end{align}

\vspace{0.25cm} The relations ({\hyperlink{sst}{\ref*{f4emf}}}) for $k=1$ and \eqref{ad40}(a) prove that the conditions  1), 2) and 3) in Lemma~{\hyperlink{hth1}{\ref*{th1}}}
hold for $f_{1}$ with $n_{0} = n_{\infty} = n_{1}= 0$.

\vspace{0.15cm}
Applying the result of Lemma~{\hyperlink{hth1}{\ref*{th1}}} to $f_{1}$ we obtain the existence of $a \in \Bb{C}$ such that $f (\ie(z))= a$, $z\in \Lambda$.

\vspace{0.15cm}
By letting $ \Lambda\ni z \to 0$ we obtain, by \eqref{f2th1pr}(b),(d),  that $\fet \!\ni\!\ie (z) \to \infty$ and, by ({\hyperlink{sst}{\ref*{f4emf}}})(b) with $k=1$,
 that $a=1$.  Thus, by virtue of \eqref{f33preresinvlamhyp}, $f_{1} (z)= 1$ for all $z\in\fet$. Since $f_1 \in {\rm{Hol}} (\Bb{H})$ we get
that $f_{1} (z)= 1$ for all $z\in\Bb{H}$, by virtue of the uniqueness theorem  for analytic functions (see
 \cite[p. 78, Thm. 3.7(c)]{con}).  This completes the proof of \eqref{f8theor3}.

\vspace{0.45cm} Finally, the relations ({\hyperlink{sst}{\ref*{f4emf}}}) for $k=2$ and \eqref{ad40}(b) prove that the conditions  1), 2) and 3) in Lemma~{\hyperlink{hth1}{\ref*{th1}}}
hold for $f_{2}$ with $n_{0} = n_{\infty}=0$,  $n_{1}= 1$.

\vspace{0.15cm}
Applying the result of Lemma~{\hyperlink{hth1}{\ref*{th1}}} we obtain the existence of $a,b \in \Bb{C}$ such that $f_{2} (\ie(z))= a z + b$.

\vspace{0.15cm}
By letting $ \Lambda\ni z \to 0$ we obtain, by \eqref{f2th1pr}(b),(d),  that $\fet \!\ni\!\ie (z) \to \infty$ and, by
({\hyperlink{sst}{\ref*{f4emf}}})(b) with $k=2$, that $f_{2} (\ie(z))\to 0$ which  implies that $b=0$ and hence  $f_{2} (\ie(z))= a z$ for $z\in \Lambda$.

\vspace{0.15cm}
But if we let $ \Lambda\ni z \to 1$ we obtain, by \eqref{f2th1pr}(c),(d),  that $\fet \!\ni\!\ie (z) \to 0$  and, by ({\hyperlink{sst}{\ref*{f4emf}}})(a) with $k=2$, that $f_{2} (\ie(z))\to 1$  and therefore
$a=1$. Thus, $ f_{2} (\ie(z))=z$ for all $z\in \Lambda$, what was to be proved for the validity of \eqref{f1theor3}.
}}\end{clash}
\end{subequations}

\begin{subequations}

\begin{clash}{\rm{\hypertarget{d20}{}\label{case20}\hspace{0,00cm}{\hyperlink{bd20}{$\uparrow$}} \ We prove \eqref{f9theor3} and \eqref{f10theor3} in more detail.

\vspace{0.25cm} \noindent {\emph{Proof of\eqref{f9theor3}.}} \
Recall that
\begin{align*}
    &  \Lambda:=
      (0,1)\cup \left(\Bb{C}\setminus\Bb{R}\right) \ .
\end{align*}

By virtue of
\begin{align*}\tag{\ref{f1xtheor3}}
    &  \lambda (z) =  {\Theta_{2}^{4}(z)}\big/{\Theta_{3}^{4}(z)} \ , \quad  z \in  \Bb{H} \,,
\end{align*}

\noindent we can write \eqref{f1theor3}
\begin{align*}\tag{\ref{f1theor3}}
      &  \lambda (\ie(z)) = {\Theta_{2}\big(\ie (z)\big)^{4}}\big/{\Theta_{3}\big(\ie (z)\big)^{4}} = z  \ , \quad  z \in
     \Lambda \ .
\end{align*}

\noindent in the form
\begin{align*}
    &  \ie \left(\lambda (\ie(z))\right)= \ie(z)  \ , \quad  z \in \Lambda \ , \
\end{align*}

\noindent and since according to \eqref{f33preresinvlamhyp}
\begin{align*}\tag{\ref{f33preresinvlamhyp}}
    &  \ie \left(\Lambda\right) = \fet \ , \
\end{align*}

\noindent we get
\begin{align}\label{ad41}
    &  \ie (\lambda (z)) = z \,,  \quad z \in \fet \,,
\end{align}

\noindent from which, in view of $\lambda \in {\rm{Hol}}(\fet)$ and $\ie \in {\rm{Hol}}(\Lambda)$, we obtain by differentiation that
\begin{align}\label{ad42}
    &  \ie^{\,\prime} (\lambda (z)) \lambda^{\,\prime} (z)=1  \,,  \quad z \in \fet \,.
\end{align}

\noindent Besides that, substituting in \eqref{f1theor4}
\begin{align*}\tag{\ref{f1theor4}}
        &  \Theta_{3}\big(\ie (y)\big)^{2}= \het (y)   \ , \quad y \in \Lambda \,  ,
    \end{align*}

\noindent $y = \lambda (z)$ with $z\in \fet$ we obtain by  \eqref{ad41}
\begin{align}\label{ad43}
    &  \Theta_{3}\big(z\big)^{2}= \het (\lambda (z))   \ , \quad z\in \fet \,  ,
\end{align}

\noindent while,  substituting in \eqref{f26preresinvlamhyp}
\begin{align*}\tag{\ref{f26preresinvlamhyp}}
       i  \, \ie^{\,\prime} (y)\, \he (y)^{2} =\frac{1}{\pi y(1-y)}
      \ , \quad y \in \Lambda \,,
\end{align*}

\noindent $y = \lambda (z)$ with $z\in \fet$  we deduce that for arbitrary  $z\in \fet$ we have
\begin{align*}
    &  \dfrac{1}{\pi \lambda (z)(1-\lambda (z))} =  i  \, \ie^{\,\prime} ( \lambda (z))\, \he ( \lambda (z))^{2}
    \stackrel{{\fo{\eqref{ad42}}}}{\vphantom{A}=} \dfrac{i  \he ( \lambda (z))^{2}}{\lambda^{\,\prime} (z)}
     \stackrel{{\fo{\eqref{ad43}}}}{\vphantom{A}=}\dfrac{i  \Theta_{3}\big(z\big)^{4}}{\lambda^{\,\prime} (z)}  \ , \
\end{align*}

\noindent from which it follows that
\begin{align}\label{ad44}
    &   \lambda^{\,\prime} (z)\!= \!i \pi \lambda (z)\big(1-\lambda (z)\big)  \Theta_{3}\left(z\right)^{4} \ , \quad  z\in \fet\,.
\end{align}

\noindent Since the functions in both sides of this equality are holomorphic on $\Bb{H}$,  by  the uniqueness theorem  for analytic functions (see
 \cite[p. 78, Thm. 3.7(c)]{con}) we get the validity of \eqref{ad44} on the whole $\Bb{H}$. This proves the left-hand side equality in \eqref{f9theor3}
 \begin{align}\label{ad45}
    &   \lambda^{\,\prime} (z)\!= \!i \pi \lambda (z)\big(1-\lambda (z)\big)  \Theta_{3}\left(z\right)^{4} \ , \quad  z\in \Bb{H}\,.
\end{align}

\noindent
Substituting in the right-hand side of  \eqref{ad45} the expression
\begin{align*}\tag{\ref{f1xtheor3}}
    &  \lambda (z) =  \dfrac{\Theta_{2}^{4}(z)}{\Theta_{3}^{4}(z)} \ , \quad  z \in  \Bb{H} \,,
\hypertarget{ssthr}{}\end{align*}

\noindent and using \eqref{f8theor3}
\begin{align*}\tag{\ref{f8theor3}}
      &     \Theta_{2}^{4}(z) + \Theta_{4}^{4}(z) = \Theta_{3}^{4}(z)  \ , \quad  z \in  \Bb{H} \ ,
\end{align*}

\noindent according to which
\begin{align}\label{ad47}
    &  1-\lambda (z) = 1 -  \dfrac{\Theta_{2}^{4}(z)}{\Theta_{3}^{4}(z)} = \dfrac{\Theta_{3}^{4}(z)- \Theta_{2}^{4}(z)}{\Theta_{3}^{4}(z)}=
     \dfrac{\Theta_{4}^{4}(z)}{\Theta_{3}^{4}(z)}  \ , \
\end{align}

\noindent we deduce from \eqref{ad45} that for every $z \in  \Bb{H}$ we have
\begin{align*}
    & \lambda^{\,\prime} (z)\!= \!i \pi \lambda (z)\big(1-\lambda (z)\big)  \Theta_{3}\left(z\right)^{4} =
    i \pi  \dfrac{\Theta_{2}^{4}(z)}{\Theta_{3}(z)^{4}}\dfrac{\Theta_{4}(z)^{4}}{\Theta_{3}(z)^{4}}\Theta_{3}\left(z\right)^{4} =
  i \pi  \dfrac{\Theta_{2}(z)^{4}\Theta_{4}(z)^{4}}{\Theta_{3}(z)^{4}} \,.
\end{align*}

\noindent
This proves the right-hand side equality in  \eqref{f9theor3}
\begin{align}\label{ad46}
    &  \lambda^{\,\prime} (z)\!= \!  i \pi  \dfrac{\Theta_{2}(z)^{4}\Theta_{4}(z)^{4}}{\Theta_{3}(z)^{4}} \ , \quad  z \in  \Bb{H}  \,.
\end{align}

\noindent \hfill $\square$

\vspace{0.25cm} \noindent {\emph{Proof of \eqref{f10theor3}.}} \ Combining \eqref{ad46}  with ({\hyperlink{ssthr}{\ref*{f1xtheor3}}})
 and \eqref{ad47} we obtain
\begin{align*}
    &  i \pi  \dfrac{\Theta_{2}(z)^{4}\Theta_{4}(z)^{4}}{\Theta_{3}(z)^{4}}=  \lambda^{\,\prime} (z) = \dfrac{d}{d z}\dfrac{\Theta_{2}\left(z\right)^{4}}{\Theta_{3}\left(z\right)^{4}} =4 \dfrac{\Theta_{2}\left(z\right)^{4}}{\Theta_{3}\left(z\right)^{4}}  \left(\dfrac{\Theta_{2}^{\,\prime}\left(z\right)}{\Theta_{2}\left(z\right)} - \dfrac{\Theta_{3}^{\,\prime}\left(z\right)}{\Theta_{3}\left(z\right)}\right)  \\[0,3cm]    &  \Rightarrow \
      \frac{i \pi}{4  } \Theta_{4}(z)^{4}\!=\!   \dfrac{\Theta_{2}^{\,\prime}\left(z\right)}{\Theta_{2}\left(z\right)} - \dfrac{\Theta_{3}^{\,\prime}\left(z\right)}{\Theta_{3}\left(z\right)} \ ,
      \\[0,5cm]    &
     i \pi  \dfrac{\Theta_{2}(z)^{4}\Theta_{4}(z)^{4}}{\Theta_{3}(z)^{4}}= - \dfrac{d}{d z} (1-\lambda (z)) = -\dfrac{d}{d z}\dfrac{\Theta_{4}\left(z\right)^{4}}{\Theta_{3}\left(z\right)^{4}}  \\[0,3cm]    &    = 4 \dfrac{\Theta_{4}\left(z\right)^{4}}{\Theta_{3}\left(z\right)^{4}}  \left(\dfrac{\Theta_{3}^{\,\prime}\left(z\right)}{\Theta_{3}\left(z\right)} - \dfrac{\Theta_{4}^{\,\prime}\left(z\right)}{\Theta_{4}\left(z\right)}\right) \ \Rightarrow \
     \frac{i \pi}{4}\Theta_{2}(z)^{4}\!=\!  \dfrac{\Theta_{3}^{\,\prime}\left(z\right)}{\Theta_{3}\left(z\right)} - \dfrac{\Theta_{4}^{\,\prime}\left(z\right)}{\Theta_{4}\left(z\right)} \ ,
     \end{align*}

\vspace{0.25cm}
\noindent which proves the validity  \eqref{f10theor3}.  $\square$
}}\end{clash}
\end{subequations}

\begin{subequations}

\begin{clash}{\rm{\hypertarget{d21}{}\label{case21}\hspace{0,00cm}{\hyperlink{bd21}{$\uparrow$}} \ Observe that Corollary~\ref{th0} and
 %\begin{align*} \tag{\ref{f5theor4}}
%        &  \Theta_{2}\big(\ie (z)\big)^{4}= z \het (z)^{2}   \ , \quad z \in (0,1)\cup \left(\Bb{C}\setminus\Bb{R}\right) \,  ,
%    \end{align*}
\begin{align*} \tag{\ref{f1theor4}}
        &  \Theta_{3}\big(\ie (z)\big)^{4}= \het (z)^{2}   \ , \quad z \in (0,1)\cup \left(\Bb{C}\setminus\Bb{R}\right) \,  ,
 \end{align*}

\noindent yields that $\Theta_{3}\big(\ie (z)\big)^{4} \in \Picklogplus $. According to the definition
\begin{align*}\tag{\ref{f9int}}
    &  \mathcal{P}_{\log}
:= \left\{ \  f\in\mathcal{P}\setminus\{0\}  \quad  \Big|  \quad
  f^{\,\prime}/f\in\mathcal{P} \  \right\}\, ,
\end{align*}

\noindent this implies
\begin{align}\label{ad50}
    & \Theta_{3}\left(\ie \right)^{4}\in \mathcal{P}(-\infty, 1) \ , \
\end{align}

\noindent  and
\begin{align}\label{ad48}
    & \mathcal{P}(-\infty, 1) \ni \dfrac{1}{\Theta_{3}\left(\ie (z)\right)^{4}}\dfrac{d}{d z}\Theta_{3}\left(\ie(z) \right)^{4}=
     \dfrac{4 \Theta_{3}^{\,\prime}\left(\ie (z)\right)\ie^{\,\prime} (z) }{\Theta_{3}\left(\ie (z)\right) }\ .
\end{align}

\noindent Since the class $\mathcal{P}(-\infty, 1)$ is invariant with respect to multiplication of its functions by any positive constant
and in view of the equality
\begin{align*}
    &  \lambda^{\,\prime} ( \ie (z) ) \ie^{\,\prime} (z)=1  \, , \quad    z \in \Lambda := (0,1)\cup \left(\Bb{C}\setminus\Bb{R}\right)\,,
\end{align*}

\noindent which follows from the differentiation of \eqref{f1theor3}
\begin{align*}\tag{\ref{f1theor3}}
      &  \lambda (\ie(z)) =  z  \ , \quad  z \in
     \Lambda \ ,
\end{align*}

\noindent we deduce from \eqref{ad48} that
\begin{align}\label{ad49}
    &  \mathcal{P}(-\infty, 1) \ni
     \dfrac{ \Theta_{3}^{\,\prime}\left(\ie (z)\right)\ie^{\,\prime} (z) }{\Theta_{3}\left(\ie (z)\right) }=
     \dfrac{\Theta_{3}^{\,\prime}\left(\ie (z)\right) }{ \lambda^{\,\prime} ( \ie (z) ) \Theta_{3}\left(\ie (z)\right)}\ .
\end{align}

\noindent Since the both functions in \eqref{ad49} and \eqref{ad50} are nonconstant we deduce from   \eqref{ad49}, \eqref{ad50}
and the definition of the Nevanlinna-Pick functions (see \eqref{eq:symm0} and \eqref{eq:symm1}) that
\begin{align}\label{ad51}
    &\hspace{-0,3cm}  \left(\im\, z\right) \cdot \im\, \Theta_{3}\left(\ie (z)\right)^{4} > 0  \ , \
    \left(\im\, z\right) \cdot \im\, \dfrac{\Theta_{3}^{\,\prime}\left(\ie (z)\right) }{ \lambda^{\,\prime} ( \ie (z) ) \Theta_{3}\left(\ie (z)\right)} > 0  \ , \hspace{-0,3cm}
\end{align}

\noindent for every $z\in\Bb{C}\setminus\Bb{R}$. It readily follows from \eqref{f16eir}
 \begin{align*}\tag{\ref{f16eir}}
    &   \Arg \ie (z) \in \dfrac{\pi}{2} - \left(0 , \dfrac{\pi}{2}\right)\cdot {\rm{sign}} (\im\, z)  \ , \  \  z\! \in \! \Bb{C}\setminus\Bb{R}  \ ,
\end{align*}

\noindent that
\begin{align}\label{ad52}
    &  (\im\, z) \cdot \re \, \ie (z) > 0  \ , \quad  z\! \in \! \Bb{C}\setminus\Bb{R}  \ ,
\end{align}

\noindent and therefore \eqref{ad51} can be written as follows
\begin{align}\label{ad53}
     &\hspace{-0,3cm}  \left(\re\, z\right) \cdot \im\, \Theta_{3}\left(z\right)^{4} > 0   \ , \quad
    \left(\re\, z\right) \cdot \im\, \dfrac{\Theta_{3}^{\,\prime}\left(z\right) }{ \lambda^{\,\prime} ( z ) \Theta_{3}\left(z\right)} > 0  \ , \hspace{-0,3cm}
\end{align}

\noindent for arbitrary $z \in   \mathcal{F}_{{\tn{\square}}} \setminus\{ i  \Bb{R}_{>0}\}$, by virtue of \eqref{f33preresinvlamhyp}
\begin{align*}\tag{\ref{f33preresinvlamhyp}}
    &  \ie \left((0,1)\cup \big(\Bb{C}\setminus\Bb{R}\big)\right) = \fet \ .
\end{align*}

\noindent Corollary~\ref{th5} follows as \eqref{ad50} and \eqref{ad53} coincide with the statements of Corollary~\ref{th5}(a),(b),(c), correspondingly.  $\square$
}}\end{clash}
\end{subequations}

\subsection[\hspace{-0,25cm}. \hspace{0,075cm}Notes on Section~\ref{log}]{\hspace{-0,11cm}{\bf{.}} Notes on Section~\ref{log}}

\begin{subequations}

\vspace{-0.15cm}
\begin{clash}{\rm{\hypertarget{d22}{}\label{case22}\hspace{0,00cm}{\hyperlink{bd22}{$\uparrow$}} \ We prove that $\log \theta_{k}(x) = \ln \theta_{k}(x)$, $x\in [0, 1)$.

 \vspace{0.15cm}
 Let  $2 \leq k \leq 4$ be fixed. Denote
\begin{align*}
    & A_{k} := \left\{ \  x\in [\,0, 1) \ \big| \  \log \theta_{k}(x) = \ln \theta_{k}(x)  \ \right\} \,,
\end{align*}

\noindent where $0 \in A_{k}$ because $\log \theta_{k}(0) = \ln \theta_{k}(0) = 0$.
Since $\exp\big(\log \theta_{k}(x)\big) =  \theta_{k}(x)$ for all $x\in [\,0, 1)$ we have
\begin{align}\label{adad1}
    &   \log \theta_{k}(x) \in \bigcup\nolimits_{{\fo{n \in \Bb{Z}}}}   \ \big\{ 2 \pi i n + \ln \theta_{k}(x)\big\}  \ , \quad  x \in (0,1) \,.
\end{align}

\vspace{-0.15cm} It follows from  $\log \theta_{k}$, $\ln \theta_{k}\! \in \!C \big([\,0, 1)\big)$ that $\log \theta_{k}\!-\!\ln \theta_{k} \!\in\! C \big([\,0, 1)\big)$ and the set $A_{k}$ contains each its limit point lying in the segment $[\,0, 1)$.
Assume that \\[0,15cm] $\phantom{a}$
\hspace{4cm}$B_{k}:= [\,0, 1) \setminus A_{k} \neq \emptyset \,$.

\noindent
Then \vspace{-0,2cm}
\begin{align}\label{adad2}
    &  b_{k}:= \inf B_{k} \in [\,0, 1)  \,.
\end{align}

 \vspace{-0,1cm} \noindent If $b_{k}=0$ then $b_{k} \in A_{k}$, but if $b_{k}>0$ then
   $[\,0, b_{k}) \subset A_{k}$ and we also have $b_{k} \in A_{k}$ because $b_{k}$ is the limit point of $A_{k}$ lying on $(0,1)$. Thus, $[\,0, b_{k}] \subset A_{k}$ and the property $\log \theta_{k} - \ln \theta_{k} \in C \big([\,0, 1)\big)$ gives the existence of $\varepsilon_{k} \!\in\! (0, 1-b_{k})$  such that
   $|\log \theta_{k} (x)\!-\!\ln \theta_{k} (x)| \! <\! 2 \pi$ for all $x \in [b_{k}, b_{k}\!+\!\varepsilon_{k})$. By virtue of  \eqref{adad1}, this
    yields  $[b_{k}, b_{k}\!+\!\varepsilon_{k}) \! \subset\!  A_{k}$ and therefore $\inf B_{k} \!\geq\! b_{k}+\varepsilon_{k}$, which contradicts \eqref{adad2}.
 This contradiction implies $B_{k}\! =\! \emptyset$, i.e.,  $\log \theta_{k}(x)\! =\! \ln \theta_{k}(x)$ for all $x\!\in\! [\,0, 1)$, what was to be proved.
}}\end{clash}
\end{subequations}

 \begin{subequations}
\begin{clash}{\rm{\hypertarget{d23}{}\label{case23}\hspace{0,00cm}{\hyperlink{bd23}{$\uparrow$}} \ Actually,
 for $u \in (0,1)$ it follows from
\begin{align*} &
  \ln  \theta_{3} (u) = \sum\nolimits_{n \geq 1} \ln  \left(1-u^{2n}\right) + 2 \sum\nolimits_{n \geq 1} \ln  \left(1 + u^{2n-1}\right) \ , \  \\  &
\ln  \theta_{2} (u) = \sum\nolimits_{n \geq 1} \ln  \left(1-u^{2n}\right) +  2\sum\nolimits_{n \geq 1} \ln  \left(1 + u^{2n}\right) \ , \
\end{align*}

\noindent that\vspace{-0,2cm}
\begin{align*}
    &  \ln  \theta_{3} (u) = - \sum\limits_{n, m \geq 1} \dfrac{u^{2nm}}{m}  + 2 \sum\limits_{n, m \geq 1} (-1)^{m-1} \dfrac{u^{2nm - m}}{m}
    \\ & =   2 \sum\limits_{ m \geq 1}  \dfrac{(-1)^{m-1}}{m} u^{-m} \dfrac{u^{2m}}{1-u^{2m}} - \sum\limits_{ m \geq 1}  \dfrac{1}{m}\dfrac{u^{2m}}{1-u^{2m}}  \\ & =      2 \sum\limits_{ m \geq 1}  \dfrac{ (-1)^{m-1}}{m} \dfrac{u^{m}}{1-u^{2m}}-   \sum\limits_{ m \geq 1} \dfrac{1}{m}\dfrac{u^{2m}}{1-u^{2m}}  \\ & =
 \sum\limits_{ m \geq 1}  \dfrac{2 }{2m-1} \dfrac{u^{2m-1}}{1-u^{2(2m-1)}}-   \sum\limits_{ m \geq 1}  \dfrac{1 }{m} \dfrac{u^{2m}}{1-u^{4m}}-\sum\limits_{ m \geq 1} \dfrac{1}{m}\dfrac{u^{2m}}{1-u^{2m}}\\ & =
\left|- u^{2m}- u^{2m} (1+u^{2m})   = - 2 u^{2m} (1+u^{2m}) + u^{2m} (1+u^{2m}) - u^{2m}\right|
 \\ & =
\sum\limits_{ m \geq 1}  \dfrac{2 }{2m-1} \dfrac{u^{2m-1}}{1-u^{2(2m-1)}}-  \sum\limits_{ m \geq 1} \dfrac{2}{m}\dfrac{u^{2m}}{1-u^{2m}}
 \\ & +\sum\limits_{ m \geq 1} \dfrac{1}{m}\dfrac{u^{2m}(1+u^{2m})}{1-u^{4m}}-   \sum\limits_{ m \geq 1}  \dfrac{1 }{m} \dfrac{u^{2m}}{1-u^{4m}}
 \\ & =
\sum\limits_{ m \geq 1}  \dfrac{2 }{2m-1} \dfrac{u^{2m-1}}{1-u^{2(2m-1)}}-
\sum\limits_{ m \geq 1} \dfrac{2}{2m-1}\dfrac{u^{2(2m-1)}}{1-u^{2(2m-1)}} \\ & -
 \sum\limits_{ m \geq 1} \dfrac{1}{m}\dfrac{u^{4m}}{1-u^{4m}}
+\sum\limits_{ m \geq 1}  \dfrac{1 }{m} \dfrac{u^{4m}}{1-u^{4m}}
 \\ & =
\sum\limits_{ m \geq 1}  \dfrac{2}{2m-1} \dfrac{   u^{2m-1}(1-u^{2m-1})}{1-u^{2(2m-1)}}  =
\sum\limits_{ m \geq 1}  \dfrac{2}{2m-1} \dfrac{   u^{2m-1}}{1+u^{2m-1}} \ ,
\end{align*}

\vspace{0.05cm}
\noindent and
\begin{align*}
     & \ln  \theta_{2} (u) =
- \sum\limits_{n, m \geq 1} \dfrac{u^{2nm}}{m}  +  2 \sum\limits_{n, m \geq 1} (-1)^{m-1} \dfrac{u^{2nm }}{m}  \\ & =
-\sum\limits_{ m \geq 1}\dfrac{1}{m} \dfrac{u^{2m}}{1-u^{2m}}  +  2 \sum\limits_{ m \geq 1}\dfrac{(-1)^{m-1}}{m}\dfrac{u^{2m}}{1-u^{2m}} \\ & =
-\sum\limits_{ n \geq 1}\dfrac{1}{n} \dfrac{u^{2n}}{1-u^{2n}}  +  2 \sum\limits_{ n \geq 1}\dfrac{1}{2n-1}\dfrac{u^{2(2n-1)}}{1-u^{2(2n-1)}} -
 2 \sum\limits_{ n \geq 1} \dfrac{1}{2n}\dfrac{u^{4n}}{1-u^{4n}}   \\ & =
-\sum\limits_{ n \geq 1}\dfrac{1}{n} \dfrac{u^{2n}}{1-u^{2n}}  +   2 \sum\limits_{ n \geq 1} \dfrac{1}{n}\dfrac{u^{2n}}{1-u^{2n}} - 2
\sum\limits_{ n \geq 1} \dfrac{1}{2n}\dfrac{u^{4n}}{1-u^{4n}} -  2 \sum\limits_{ n \geq 1} \dfrac{1}{2n}\dfrac{u^{4n}}{1-u^{4n}}  \\ & =
\sum\limits_{ n \geq 1}\dfrac{1}{n} \dfrac{ u^{2n}}{1-u^{2n}} -  \sum\limits_{ n \geq 1} \dfrac{1}{n}\dfrac{2u^{4n}}{1-u^{4n}} =
\sum\limits_{ n \geq 1}\dfrac{1}{n} \dfrac{u^{2n} (1+u^{2n}) -  2u^{4n}}{1-u^{4n}}\\ & =
\sum\limits_{ n \geq 1}\dfrac{1}{n} \dfrac{u^{2n} -u^{4n}}{1-u^{4n}}   =
\sum\limits_{ n \geq 1}\dfrac{1}{n} \dfrac{u^{2n}(1 -u^{2n})}{1-u^{4n}} =
 \sum\limits_{ n \geq 1}\dfrac{1}{n} \dfrac{u^{2n}}{1+u^{2n}} \ ,
\end{align*}

%\vspace{0.25cm}
\noindent while\vspace{-0,2cm}
\begin{align}\label{ad54}
    &   \ln  \theta_{4} (u) =  \ln  \theta_{3} (-u) = -\sum\limits_{ m \geq 1}  \dfrac{2}{2m-1} \dfrac{   u^{2m-1}}{1-u^{2m-1}} \ .
\end{align}

\vspace{-0,3cm}
\noindent Therefore
\begin{align}\nonumber
    &  \ln  \theta_{3} (u) =  \sum\limits_{ n \geq 1}\, \dfrac{2}{ 2n-1}\,\dfrac{u^{2n-1}}{1+ u^{2n-1}}  \ ,    &     &
 \ln  \theta_{4} (u)  =   - \sum\limits_{ n \geq 1}\, \dfrac{2}{ 2n-1}\,\dfrac{u^{2n-1}}{1- u^{2n-1}} \ ,   \\    &
 \ln  \theta_{2} (u) =  \sum\limits_{ n \geq 1}\,\dfrac{1}{n}\, \dfrac{u^{2n}}{1+u^{2n}} \,  ,  &     & u \in (0,1)  \ .
\label{ad55}
\end{align}

}}\end{clash}
\end{subequations}

\subsection[\hspace{-0,1cm}. \hspace{-0,06cm}Notes on Section~\ref{eirth}]{\hspace{-0,11cm}{\bf{.}} Notes on Section~\ref{eirth}}

\begin{subequations}

\vspace{-0.15cm}
\begin{clash}{\rm{\hypertarget{d24}{}\label{case24}\hspace{0,00cm}{\hyperlink{bd24}{$\uparrow$}} \
We prove Corollary~\ref{cor1} in more detail. Recall that
\begin{align*}
    &  \Lambda:=
      (0,1)\cup \left(\Bb{C}\setminus\Bb{R}\right)
\end{align*}

\noindent and that holomorphic in the upper half-plane function $\log \Theta_{3} \in {\rm{Hol}} (\Bb{H})$ has been defined in \eqref{f11elemtheta}
as holomorphic extension  from $i\cdot \Bb{R}_{>0}$ to $\Bb{H}$ of the function $\ln \Theta_{3} (z)= \ln\theta_{3}(\exp(i \pi z))$
and, by virtue of \eqref{f3elemtheta},
\begin{align}\label{ad56a}
    &   \log \Theta_{3}  (z) = \sum\limits_{ n \geq 1}  \dfrac{2}{2n-1} \,\dfrac{ e^{ i  \pi (2n-1) z}  }{1+e^{ i  \pi (2n-1) z}  }
    \ , \quad   z \in \Bb{H}  \,.
\end{align}

\noindent

In view of analytic property of the composition of two analytic functions  (see \cite[p. 34, 2.4]{con}), it follows from
\begin{align*}
    \tag{\ref{f33preresinvlamhyp}}
    & \ie \in {\rm{Hol}} (\Lambda) \, , \quad   \ie \left(\Lambda\right) = \fet \subset \Bb{H}\, ,
\end{align*}

\noindent that
\begin{align}\label{ad56}
    &   \log\Theta_{3} \big(\ie (z)\big)\in {\rm{Hol}} (\Lambda) \, .
\end{align}

\noindent  According to  \eqref{f10int}, \eqref{f2preresinvlamhyp} and \eqref{f1elemtheta}, we have
\begin{align*}
       \ie (y)\! = \!  i \, \frac{\he(1-y)}{\he(y)} \in i \Bb{R}_{> 0} \, , \quad \he (y) \!= \!\sum\limits_{n = 0}^{\infty}
\dfrac{\Gamma(n+1/2)^{2}}{\pi\, (n !)^{2}} y^{n} \!>\! 0  \, ,
\quad   y\! \in\! (0,1)\,,
\end{align*}

\noindent $\Theta_{3}(i x)= \theta_{3}(\exp(- \pi x)) > 0$, $x > 0$,  and therefore
\begin{align}\label{ad57}
    &   \log\Theta_{3} \big(\ie ( y)\big) = \ln \Theta_{3} \big(\ie ( y)\big)  \ , \quad   y \in (0,1) \ .
\end{align}

 We recall that, in view of \eqref{f1cth0} and \eqref{f1ath0}, the function  $\Log \het $ can  be considered as  the holomorphic extension  from the interval $(0,1)$ to $\Lambda$ of the function $\ln \het$ with
\begin{align}\label{ad57a}
    &   \Arg \het (z)  \in  \left(-\dfrac{\pi}{2},  \dfrac{\pi}{2}\right)
     \ , \quad  z \in \Lambda \,.
\end{align}

  Applying the logarithm to the Wirtinger identity
  \begin{align}\tag{\ref{f1theor4}}
        &  \Theta_{3}\big(\ie (z)\big)^{2}= \het (z)   \ , \quad z \in (0,1)\cup \left(\Bb{C}\setminus\Bb{R}\right) \,  ,
    \end{align}

\noindent for $z = y \in (0,1)$ we obtain   that
\begin{align*}
    &   \Log \het (y) =   \ln \het (y) = \ln \Theta_{3}\big(\ie (y)\big)^{2} = 2 \ln \Theta_{3}\big(\ie (y)\big) = 2 \log\Theta_{3} \big(\ie ( y)\big)  \, ,
\end{align*}

\noindent i.e., the  two functions   $ \Log \het$ and $ 2\, \log\Theta_{3} \big(\ie \big) $ coincide on the interval $(0,1)$ and at the same time they are analytic in $\Lambda$.  By  the uniqueness theorem  for analytic functions (see
 \cite[p. 78, Thm. 3.7(c)]{con}),  they coincide on $\Lambda $,
\begin{align}\label{ad59}
    &  2\, \log\Theta_{3} \big(\ie (z)\big)= \Log \het(z) \ , \quad  z \in \Lambda \,,
\end{align}

\noindent
 and, by virtue of \eqref{ad57a},
\begin{align*}
    & 2 \, \im  \log\Theta_{3} \big(\ie (z)\big) = \im \, \Log \het (z) =
    \Arg \het (z) \in \left(-\dfrac{\pi}{2},  \dfrac{\pi}{2}\right)  \ , \quad  z \in \Lambda \, ,
\end{align*}

\noindent from which it follows that
\begin{align}\label{ad58}
    &    \arg \Theta_{3} (\ie (z)) \in \left(- \frac{\pi}{4} \, , \, \frac{\pi}{4}\right)   \ , \quad  z \in (0,1)\cup \left(\Bb{C}\setminus\Bb{R}\right)\,,
\end{align}

\noindent which means that \eqref{f0yeirth} holds. Since (\,by \eqref{f33preresinvlamhyp}\,)  $ \ie \left(\Lambda\right) = \fet$,   then \eqref{ad58}  can be written as follows
\begin{align}\label{ad60}
    &  \arg \Theta_{3} \big(z\big) \in \left(-\dfrac{\pi}{4} \, , \, \dfrac{\pi}{4}\right)  \, , \quad  z  \in \fet \, ,
\end{align}

\noindent  while  \eqref{ad59} implies that
  \begin{align*}
        & 2 \, \log \Theta_{3}\big(\ie (- x\pm i \varepsilon)\big)= \Log \het (- x \pm i \varepsilon)   \ , \quad x, \varepsilon > 0 \,  .
    \end{align*}

\noindent  Letting  here $\varepsilon\downarrow 0$ and drawing up the facts that  $\log \Theta_{3} \in {\rm{Hol}} (\Bb{H})$ and the limits $\het (- x\pm i 0) > 0$ and $\ie (- x\pm i 0)  \in \Bb{H}$ exist and finite,  we deduce from the expressions \eqref{f14zpreresinvlamhyp} and \eqref{f18preresinvlamhyp} that
\begin{align}\label{ad61}
    & 2 \, \log \Theta_{3}\left(\pm 1 +  i  \,
 \dfrac{\he  \left( {1}\big/{(1 + x)}\right)}{\he  \left(  {x}\big/{(1 + x)}\right)}\right)= \ln \he (-x) \     \ , \quad x > 0 \,  ,
\end{align}

\noindent where
\begin{align*}
    &  \he (-x) =  \dfrac{\he \left(\dfrac{x}{1 + x}\right)}{\sqrt{1+x}} > 0  \ , \quad   x > 0 \,.
\end{align*}

\noindent If in \eqref{ad61} we introduce the notation
\begin{align*}
    &  y (x):= \dfrac{\het  \big( 1/(1 + x)\big)}{\het  \big( x/(1 + x)\big)}   \ ,
\end{align*}

\noindent then \eqref{f7preresinvlamhyp}(a),(b) imply $ y (0):= \lim_{\,t\downarrow 0} y (t) = +\infty $ and
$y (+\infty):=\lim_{t\to +\infty} y (t) = 0$, while \eqref{f2eir} yields $y^{\,\prime}(x) < 0$ for all $x>0$.
By substituting in \eqref{ad61} the expression \eqref{f8eir} with $z = - x \!\in\!\Bb{C}\!\setminus\! [1, +\!\infty)$, we obtain
\begin{align*}
    &  2\, \log\Theta_{3} \left( \pm 1\! +\!  i  y (x)\right)\! =\! \dfrac{1}{ \pi^{2 }} \int\limits_{0}^{1}
    \frac{\dfrac{1}{t(1-t)} \ln\left(\dfrac{1}{1+t x}\right) }{\he(t)^{2} +
\he(1-t)^{2}}\, \, \diff t\ , \quad   x > 0 \,,
\end{align*}

\noindent which completes the proof of \eqref{f1aeirf} and together with \eqref{adt27}(d) establishes that
\begin{align}\label{ad62}
    & \arg \Theta_{3} \left(\pm 1\! +\!  i  x\right) =0 \ , \
     \Theta_{3} \left( 1\! +\!  i  x\right) = \Theta_{3} \left( - 1\! +\!  i  x\right) \in \left(0,1\right) \, , \quad   x > 0 \, .
\end{align}

\noindent It also follows from \eqref{ad62} and \eqref{ad60} that
\begin{align*}\tag{\ref{f1yeirth}}
    &   \arg \Theta_{3} (z) \in \left(-\dfrac{\pi}{4} \, , \, \dfrac{\pi}{4}\right)  \ ,
    \quad  z  \in  \mathcal{F}^{\,{\tn{||}}}_{{\tn{\square}}}  \, .
\end{align*}

By substituting in \eqref{ad59} the expression \eqref{f8eir} with $z \in \Lambda \subset \Bb{C}\!\setminus\! [1, +\!\infty)$ we obtain
\begin{align*}\tag{\ref{f1eirf}}
    & \hspace{-0,3cm} 2\, \log\Theta_{3} \big(\ie (z)\big) = \dfrac{1}{ \pi^{2 }} \int\limits_{0}^{1}
   \frac{\dfrac{1}{t(1-t)}\Log\left(\dfrac{1}{1-tz}\right)}{\he(t)^{2} +
\he(1-t)^{2}}\, \,\diff t  \ , \ \ z \in \Lambda\,, \hspace{-0,2cm}
\end{align*}

\noindent we obtain the validity of \eqref{f1eirf}. Corollary~\ref{cor1} now follows.

}}\end{clash}
\end{subequations}

 \begin{subequations}
\begin{clash}{\rm{\hypertarget{d25}{}\label{case25}\hspace{0,00cm}{\hyperlink{bd25}{$\uparrow$}} \ More precisely, in view of \eqref{f1xtheor3}, the representation
\eqref{f5eirf} can be written as follows
\begin{align}\label{ad64}
    &  \log\Theta_{3}\! \left(z\right) = \dfrac{1}{2 \pi^{2 }} \int\limits_{0}^{1}
    \dfrac{\dfrac{1}{t (1-t)}  \Log\dfrac{1}{\vphantom{\frac{A}{B}}1-t\lambda (z)} }{\het \left(t \right)^{2} + \het \left(1-t \right)^{2}} \   \diff  t \, , \quad  z\in \mathcal{F}^{\infty}_{{\tn{\square}}} \,.
 \end{align}

\noindent Then the function
\begin{align*}
    &  \tau =  \tau (t) :=  \dfrac{\ie (t)}{i} =  \frac{\he(1-t)}{\he(t)} \in  \Bb{R}_{> 0}  \ , \quad  t \in (0,1)\,,
\end{align*}

\noindent by virtue of  \eqref{f7preresinvlamhyp}(a),(b),  satisfies $ \tau (0):= \lim_{\,t\downarrow 0} \tau (t) = +\infty $ and
$\tau (+\infty):=\lim_{t\to +\infty} \tau (t) = 0$, while by \eqref{f2eir},
\begin{align*}\tag{\ref{f2eir}}
    & \hspace{-0,2cm}\tau^{\,\prime} (t)= \dfrac{ \diff }{ \diff  t} \dfrac{\het (1-t)}{\het (t)} = -  \dfrac{1}{\pi t (1-t) \het (t)^{2}} < 0 \ , \quad  t \in (0,1)\, .
\end{align*}

\noindent Therefore, for any $t \in (0,1)$ we have
\begin{align*}
    & \dfrac{\dfrac{1}{\pi t (1-t)}}{\het \left(t \right)^{2} + \het \left(1-t \right)^{2}} =
      \dfrac{1}{\pi t (1-t) \het (t)^{2}} \cdot \dfrac{1}{1+ \dfrac{\het \left(1-t \right)^{2}}{\het \left(t \right)^{2}}} =
      - \dfrac{\tau^{\,\prime} (t)}{1+ \tau (t)^{2}} \ , \
\end{align*}

\noindent and for arbitrary $z\in \mathcal{F}^{\infty}_{{\tn{\square}}}$ we can write \eqref{ad64} in the form of \eqref{f7eirf}
\begin{align*}
    &   \log\Theta_{3}\! \left(z\right) = - \dfrac{1}{2 \pi} \int\limits_{0}^{1}
    \dfrac{ \Log\dfrac{1}{\vphantom{\frac{A}{B}}1-t\lambda (z)}}{1+ \tau (t)^{2}}  \diff  \tau (t) =
\frac{1}{2 \pi} \int\limits_{0}^{+\infty} \ \
\,\frac{\Log\dfrac{1}{1-\lambda ( i  \tau)\lambda(z)}}{1+ \tau^{2}} \diff \tau\ ,
\end{align*}

\noindent because $  \ie (t)= i \tau$, $t\in (0,1)$, implies $t = \lambda (\ie (t)) = \lambda (i\tau)$, $\tau\in (0,+\infty)$, by
 \eqref{f1theor3} with  $z \in (0,1)$  combined with \eqref{f1xtheor3}.
}}\end{clash}
\end{subequations}

%\vspace{0.5cm}
%\subsection{Notes on Section~\ref{eirth}}
%
%
%\noindent
%%%%%%%%%%%%%%%%%%%%%%%%%%%%%%%%%%%%%%%%%%%%%%%%%%%%%%%%%%%%%%%%%%%%%%%%%%%%%%%%%%%%%%%%%%%%%%%%%%%%%%%%%

\vspace{1cm}

\newpage

\addcontentsline{toc}{section}{References}

%\vspace{-0,5cm}

\end{document}